\documentclass{amsart}
\usepackage{graphicx}
\usepackage[all]{xy}
\usepackage{latexsym}
\usepackage{amssymb}
\usepackage{amsmath}
\usepackage{color}

\usepackage{tikz, tikz-cd}
 \tikzstyle{int}=[circle, draw,fill=black,outer sep=0,minimum size=3pt, inner sep=0]
  \tikzstyle{ext}=[circle, draw=black,outer sep=0,inner sep=1pt]

\xyoption{arc}

 \newfam\cyrfam

  \font\tencyr=wncyr10

  \font\sevencyr=wncyr7

  \font\fivecyr=wncyr5

  \textfont\cyrfam=\tencyr \scriptfont\cyrfam=\sevencyr

    \scriptscriptfont\cyrfam=\fivecyr

  \newfam\cyifam

  \font\tencyi=wncyi10

  \font\sevencyi=wncyi7

  \font\fivecyi=wncyi5

  \textfont\cyifam=\tencyi \scriptfont\cyifam=\sevencyi

    \scriptscriptfont\cyifam=\fivecyi

\def\id{{\mbox{1 \hskip -7pt 1}}}
\newcommand{\sgn}{{\mathit s  \mathit g\mathit  n}}
 \newcommand{\lon}{\longrightarrow}
 \newcommand{\bu}{\bullet}
 
 \newcommand{\rar}{\rightarrow}
 \newcommand{\hook}{\hookrightarrow}

\newcommand{\p}{{\partial}}

\newcommand{\Der}{\mathrm{Der}}

\newcommand{\A}{{\mathbb A}}

 \newcommand{\Z}{{\mathbb Z}}
 \newcommand{\bS}{{\mathbb S}}
 \renewcommand{\P}{{\mathbb P}}
 \newcommand{\C}{{\mathbb C}}
 \newcommand{\R}{{\mathbb R}}
 \newcommand{\N}{{\mathbb N}}
 \newcommand{\K}{{\mathbb K}}

 \newcommand{\bbH}{{\mathbb H}}
\newcommand{\Conf}{{\mathit{Conf}}}
 \newcommand{\ot}{\otimes}

\newcommand{\bfX}{{\mathbf X}}
\newcommand{\bfY}{{\mathbf Y}}

\newcommand{\sC}{{\mathsf C}}

\newcommand{\sG}{{\mathsf G}}
\newcommand{\sK}{{\mathsf K}}

\newcommand{\sP}{{\mathsf P}}

\newcommand{\Def}{\mathsf{Def}}
\newcommand{\fGC}{\mathsf{fGC}}
\newcommand{\Arg}{{\mathrm{Arg}}}

\newcommand{\GCor}{\mathsf{GC}^\mathit{or}}
\newcommand{\GC}{\mathsf{GC}}

\newcommand{\dfGC}{\mathsf{dfGC}}

\newcommand{\Holie}{\mathcal{H}\mathit{olie}}
\newcommand{\Holied}{\mathcal{H}\mathit{olie}_d}

\newcommand{\Assb}{\mathcal{A}\mathit{ssb}}

  \newcommand{\LB}{\mathcal{L}\mathit{ieb}}

\newcommand{\wqLB}{\widehat{\LB}^{\mathrm{quant}}}
\newcommand{\wLB}{\widehat{\LB}}
\newcommand{\wLBm}{\widehat{\LB}^{\mathrm{min}}_\infty}
\newcommand{\LBm}{{\LB}^{\mathrm{min}}_\infty}

 \newcommand{\Beq}{\begin{equation}}
 \newcommand{\Eeq}{\end{equation}}
 \newcommand{\Beqr}{\begin{eqnarray}}
 \newcommand{\Eeqr}{\end{eqnarray}}
 \newcommand{\Beqrn}{\begin{eqnarray*}}
 \newcommand{\Eeqrn}{\end{eqnarray*}}
 \newcommand{\Ba}{\begin{array}}
 \newcommand{\Ea}{\end{array}}
 \newcommand{\Bi}{\begin{itemize}}
 \newcommand{\Ei}{\end{itemize}}
 \newcommand{\Bc}{\begin{center}}
 \newcommand{\Ec}{\end{center}}
 \newcommand{\fg}{{\mathfrak g}}

\newcommand{\fr}{{\mathfrak r}}

\newcommand{\ft}{{\mathfrak t}}

\newcommand{\fA}{{\mathfrak A}}
 
\newcommand{\fC}{{\mathfrak C}}
 \newcommand{\fG}{{\mathfrak G}}

 \newcommand{\f}{{\mathcal O}}
 \newcommand{\cA}{{\mathcal A}}
 \newcommand{\cB}{{\mathcal B}}
 \newcommand{\cC}{{\mathcal C}}
 \newcommand{\caD}{{\mathcal D}}
 \newcommand{\cE}{{\mathcal E}}
 \newcommand{\cF}{{\mathcal F}}
 \newcommand{\cG}{{\mathcal G}}
 \newcommand{\caH}{{\mathcal H}}

 \newcommand{\cK}{{\mathcal K}}
 \newcommand{\caL}{{\mathcal L}}
 \newcommand{\cM}{{\mathcal M}}
 
 \newcommand{\cP}{{\mathcal P}}

 \newcommand{\cT}{{\mathcal T}}

 \newcommand{\cU}{{\mathcal U}}

 \newcommand{\al}{\alpha}
 \newcommand{\be}{\beta}
 \newcommand{\ga}{\gamma}
 
 \newcommand{\Ga}{\Gamma}

 \newcommand{\var}{\varepsilon}
 \newcommand{\la}{\lambda}
 \newcommand{\om}{\omega}

 \newcommand{\bx}{{\mathbf x}}
\newcommand{\by}{{\mathbf y}}
\newcommand{\bt}{{\mathbf t}}

\newcommand{\bz}{{\mathbf z}}

 \newcommand{\Hom}{{\mathrm H\mathrm o\mathrm m}}
 
 \newcommand{\sip}{\smallskip}
 \newcommand{\bip}{\bigskip}
 \newcommand{\mip}{\vspace{2.5mm}}
\theoremstyle{plain}
\swapnumbers

\newtheorem{prop-def}[theorem]{Proposition-definition}

\newtheorem{f-theorem}{Formality Theorem}[section]
\newtheorem{main-theorem}{Main~Theorem}[section]
\newtheorem{section-theorem}{Theorem}[section]

\theoremstyle{definition}

\renewcommand{\thesubsection}{\bf\arabic{section}.\arabic{subsection}}
\renewcommand{\thesubsubsection}{\bf\arabic{section}.\arabic{subsection}.\arabic{subsubsection}}

\begin{document}

 \sloppy

 \newenvironment{proo}{\begin{trivlist} \item{\sc {Proof.}}}
  {\hfill $\square$ \end{trivlist}}

\long\def\symbolfootnote[#1]#2{\begingroup%
\def\thefootnote{\fnsymbol{footnote}}\footnote[#1]{#2}\endgroup}

\title{An explicit two step quantization of Poisson structures\\ and Lie bialgebras}

\author{Sergei~Merkulov}
\address{Sergei~Merkulov:  Mathematics Research Unit, Luxembourg University,  Grand Duchy of Luxembourg }
\email{sergei.merkulov@uni.lu}

\author{Thomas~Willwacher}
\address{Thomas~Willwacher: Institute of Mathematics, University of Zurich, Zurich, Switzerland}
\email{thomas.willwacher@math.uzh.ch}

\begin{abstract} We develop a new approach to deformation quantizations of Lie bialgebras and Poisson structures  which goes in two steps.

\sip

In the first step one associates to any Poisson  (resp.\ Lie bialgebra) structure a so called {\em quantizable}\,
Poisson (resp.\ Lie bialgebra) structure.  We show explicit transcendental formulae for this correspondence.

\sip

In the second step one deformation quantizes a {\em quantizable}\, Poisson (resp.\ Lie bialgebra) structure. We show again explicit  transcendental formulae for this second step correspondence  (as a byproduct we obtain configuration space models for biassociahedron and bipermutohedron).

\sip

In the Poisson case the first step is the most non-trivial one and requires a choice of an associator while the second step quantization is essentially unique, it is independent of a choice of an associator and  can be done by a trivial induction. We conjecture that  similar statements hold true in the case of Lie bialgebras.

\sip

The  main new result is a surprisingly simple explicit universal formula (which uses only smooth differential forms) for universal quantizations of finite-dimensional Lie bialgebras.

%\sip
%\noindent {\sc Mathematics Subject Classifications} (2000). 17B37, 16W30, 51M20.

%\noindent {\sc Key words}. Strongly homotopy bialgebras, Hopf algebras,  Lie bialgebras,
% deformation quantization, configuration spaces.
\end{abstract}
 \maketitle
\markboth{Sergei Merkulov and Thomas Willwacher}{Explicit deformation quantization of Lie bialgebras}

{\small
{\small
\tableofcontents
}
}

{\Large
\section{\bf Introduction}
}

\bip

\subsection{Two classical deformation quantization problems} There are two famous deformation quantization problems, one deals with quantization of Poisson structures on finite dimensional manifolds and another with quantization of Lie bialgebras.

\sip

A lot is known by now about the first deformation quantization problem: we have an explicit formula for a universal deformation quantization \cite{Ko}, we also know  that all homotopy inequivalent universal deformation quantizations are classified by the set of Drinfeld associators and that, therefore, the Grothendieck-Teichm\"uller group acts on such quantizations.

\sip

Also much is known about the second quantization problem.
Thanks to  Etingof and Kazhdan in \cite{EK} it is proven that, for any choice of a Drinfeld associator, there exists a universal quantization of an arbitrary Lie bialgebra. Later Tamarkin gave a second proof of the Etingof-Kazhdan
deformation quantization theorem in \cite{Ta}, and very recently Severa found a third proof \cite{Se}. The theorem follows furthermore from the more general results of \cite{GY}. All these proofs give us existence theorems
for deformation quantization maps, but show no hint on how such a quantization might look like explicitly to any order in $\hbar$.

\sip

In this paper we show a new  transcendental explicit formula for universal quantization
of finite-dimensional Lie bialgebras. This work is based on the study
of compactified configuration spaces in $\R^3$ which was motivated  by (but not identical to) an earlier work of Boris Shoikhet \cite{Sh1}; it
gives in particular a new proof of the Etingof-Kazhdan existence theorem.
The methods used in the construction of that formula work well also in two dimensions,
and give us  new explicit formulae for a universal quantization of Poisson structures. Let us explain main ideas of the paper first in this very popular case.

\subsection{Deformation quantization of Poisson structures}
Let  $C^\infty(\R^d)$ be the commutative algebra of smooth functions in $\R^n$.
A {\em star product} in $C^\infty(\R^n)$ is an associative product,
$$
\Ba{rccc}
*_\hbar: & C^\infty(\R^n) \times C^\infty(\R^n) & \lon & C^\infty(\R^n)\\
   &   (f(x),g(x)) & \lon & f *_\hbar g = fg + \sum_{k\geq 1}^\infty \hbar^k B_k(f,g)
\Ea
$$
where all operators $B_k$ are bi-differential. One can check that the associativity condition
for $*_\hbar$ implies that $\pi(f,g):= B_1(f,g) - B_1(g,f)$ is a Poisson structure in $\R^n$; then
$*_\hbar$ is called a {\em deformation quantization}\, of $\pi\in \cT_{poly}(\R^n)$.

\sip

The {\em deformation quantization problem}\, addresses the question: given a Poisson structure in $\R^n$, does there exist a
star product $*_\hbar$ in $C^\infty(\R^n)$ which is its deformation quantization?

\sip

This problem was solved by Maxim Kontsevich \cite{Ko} by giving an explicit direct map between the two sets

$$
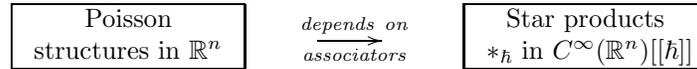
\xymatrix{
\Ba{c}\frame{\mbox{$\Ba{c} \ \mathrm{Poisson} \ \\
 \  \mathrm{structures\ in}\ \R^n\   \Ea$}}\Ea \ \ \ \ \ \
 \ar[r]^{\mathit{depends\ on}}_{\mathit{associators}} &
\ \ \ \ \ \
\
\Ba{c}\frame{\mbox{$\Ba{c} \ \mathrm{Star\ products} \ \\
\ *_\hbar\ \mathrm{in}\ C^\infty(\R^n)[[\hbar]]  \Ea$}}
\Ea
 }
$$
In fact, a stronger correspondence was proven --- the formality theorem.
Later Dmitry Tamarkin proved \cite{Ta2} an existence theorem for deformation quantizations which exhibited a non-trivial role of Drinfeld's associators.

\sip

In this paper we consider an intermediate object --- a {\em quantizable}\, Poisson structure --- so that the quantization process splits in two steps as follows
$$
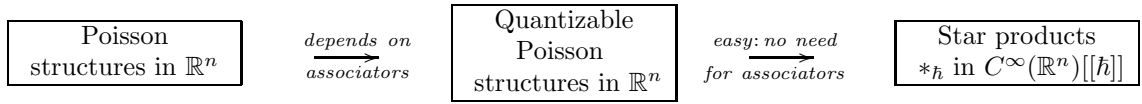
\xymatrix{
\Ba{c}\frame{\mbox{$\Ba{c} \ \mathrm{Poisson} \ \\
 \  \mathrm{structures\ in}\ \R^n\   \Ea$}}\Ea \ \ \ \ \ \
 \ar[r]^{\mathit{depends\ on}}_{\mathit{associators}} &
\ \ \ \ \ \
\Ba{c}
\frame{\mbox{$\Ba{c} \ \mathrm{Quantizable} \ \\  \mathrm{Poisson} \ \\
 \  \mathrm{structures\ in}\ \R^n  \Ea$}}
\Ea \ \ \ \ \ \
\ar[r]^{\mathit{easy:\, no\,\, need}}_{\mathit{for\,\, associators}} &
\ \ \ \ \ \ \
\Ba{c}\frame{\mbox{$\Ba{c} \ \mathrm{Star\ products} \ \\
\ *_\hbar\ \mathrm{in}\ C^\infty(\R^n)[[\hbar]]  \Ea$}}
\Ea
 }
$$

If an ordinary Poisson structure
is a Maurer-Cartan element $\pi \in \cT_{poly}(\R^n)$ of the classical Schouten bracket $[\ ,\ ]_S$,
$$
[\pi,\pi]_S=0,
$$
a quantizable Poisson structure $\pi^{\diamond}$ is a bivector field in $\cT_{poly}(\R^n)[[\hbar]]$ which is Maurer-Cartan element,
\Beq\label{1: eqn for pi^diamond}
\frac{1}{2}[\pi^{\diamond},\pi^{\diamond}]_S + \frac{\hbar}{4!} [\pi^{\diamond},\pi^{\diamond},\pi^{\diamond},\pi^{\diamond}]_4 +
\frac{\hbar^2}{6!} [\pi^{\diamond},\pi^{\diamond},\pi^{\diamond},\pi^{\diamond},\pi^\diamond,\pi^\diamond]_6
+ \ldots =0,
\Eeq
a certain $\caL ie_\infty$ structure in $\cT_{poly}(\R^n)$,
$$
\left\{[\, \ , \ldots ,\ ]_{2k}: \cT_{poly}(\R^n)^{\ot 2k}\rar \cT_{poly}(\R^n)[3-4k] \right\}_{k\geq 1}
$$
which we call a {\em Kontsevich-Shoikhet $\caL ie_\infty$ structure}\, as it was  was  introduced by Boris Shoikhet in \cite{Sh} with a reference to an important contribution by Maxim Kontsevich via an informal communication. As the Schouten bracket, this $\caL ie_{\infty}$ structure makes sense in infinite dimensions.
It was proven in \cite{Wi2} that the Kontsevich-Shoikhet structure is {\em the unique}\, non-trivial deformation of the standard
Schouten bracket in $\cT_{poly}(\R^n)$ in the class of universal $\caL ie_\infty$ structures which makes sense in {\em any}\, (including infinite) dimension (it is a folklore conjecture that in finite dimensions the Schouten bracket $[\ ,\ ]$ is rigid, i.e. admits no universal homotopy non-trivial deformations).

\sip

A map
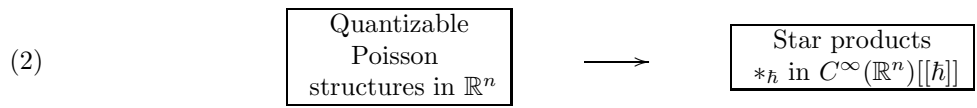
\Beq\label{1: from qua Poisson to star products}
\xymatrix{
\Ba{c}
\frame{\mbox{$\Ba{c} \ \mathrm{Quantizable} \ \\  \mathrm{Poisson} \ \\
 \  \mathrm{structures\ in}\ \R^n  \Ea$}}
\Ea \ \ \ \ \ \
\ar[r]^{} &
\ \ \ \ \ \ \
\Ba{c}\frame{\mbox{$\Ba{c} \ \mathrm{Star\ products} \ \\
\ *_\hbar\ \mathrm{in}\ C^\infty(\R^n)[[\hbar]]  \Ea$}}
\Ea
}
\Eeq
was constructed in \cite{Sh} for any $n$ (including the case $n=+\infty$) with the help of the hyperbolic geometry and transcendental formulae. It was shown in
\cite{Wi2,B} that this universal map (which comes in fact from a $\caL ie_\infty$ morphism)
is essentially unique and can, in fact, be constructed by a trivial (in the sense that no choice
of an associator is needed) induction.

\sip

What is new in our paper is the following Theorem proven in Section 4 below.

\subsubsection{\bf Theorem}\label{1: Theorem on ordinary and qua Poisson} {\em For any finite $n$ and any choice of an associator, there is 1-1
correspondence between the two sets,
$$
\xymatrix{
\Ba{c}\frame{\mbox{$\Ba{c} \ \mathrm{Poisson} \ \\
 \  \mathrm{structures\ in}\ \R^n\   \Ea$}}\Ea \ \ \ \
 \leftrightarrow &
\Ba{c}
\frame{\mbox{$\Ba{c} \ \mathrm{Quantizable} \ \\  \mathrm{Poisson} \ \\
 \  \mathrm{structures\ in}\ \R^n  \Ea$}}
\Ea
}
$$
More precisely, there is a $\caL ie_\infty$ isomorphism,
$$
F: \left(\cT_{poly}(\R^n), [\ ,\ ]_S\right) \lon \left(\cT_{poly}(\R^n)[[\hbar]], \{ [\ , \ , \ldots ,\ ]_{2n}\}_{n\geq 1}\right)
$$
 from the Schouten algebra to the Kontsevich-Shoikhet one.
}

\sip

We show explicit transcendental formulae for the $\caL ie_\infty$ morphism $F$ in (\ref{4: F_k,l from_g}). Composing this morphism
with the essentially unique arrow in (\ref{1: from qua Poisson to star products}) we obtain an acclaimed
new explicit formula for universal quantization of Poisson structures. In fact we obtain a family of such formulae parameterized by smooth functions on $S^1:=\{(x,y)\in \R^2| x^2+y^2=1\}$ with compact support
in the upper $(y>0)$ half-circle; all the associated maps $F$ are homotopy equivalent to each other.

\subsection{Deformation quantization of Lie bialgebras}
 Let $V$ be a $\Z$-graded %finite-dimensional
real vector space, and let
$\f_V:= {\odot^{\bullet}} V= \oplus_{n\geq 0} \odot^n V$
be the space of polynomial functions on $V^*$ equipped with the standard graded commutative and
cocommutative bialgebra structure. If $\cA ss\cB $ stands for the prop of bialgebras, then the standard product and coproduct in $\f_V$ give us a representation,
\Beq\label{1: rho_0}
\rho_0: \cA ss\cB \lon \cE nd_{\f_V}.
\Eeq
A {\em formal deformation}\, of the standard bialgebra structure in $\f_V$ is a continuous
morphisms of props,
$$
\rho_\hbar: \cA ss\cB[[\hbar]] \lon \cE nd_{\f_V}[[\hbar]],
$$
$\hbar$ being a formal parameter, such that $\rho_\hbar|_{\hbar=0}=\rho_0$. It is well-known
\cite{D}
that if $\rho_\hbar$ is a formal deformation of $\rho_0$, then $\frac{d \rho_\hbar}{d\hbar}|_{\hbar=0}$ makes the vector space $V$ into a Lie bialgebra, that is, induces a representation,
$$
\nu: \caL ie \cB \lon \cE nd_V
$$
of the prop of Lie bialgebras $\caL ie \cB$ in $V$. Thus Lie bialgebra structures, $\nu$,
in $V$ control infinitesimal formal deformations of $\rho_0$. Drinfeld formulated a deformation quantization problem:
given $\nu$ in $V$, does $\rho_\hbar$ exist such that $\frac{d \rho_\hbar}{d\hbar}|_{\hbar=0}$ induces $\nu$? This problem was solved affirmatively in \cite{EK,Ta2,Se}. In this paper we give a new proof of the Etingof-Kazhdan theorem  which shows such
an explicit formula in the form $\sum_\Ga c_\Ga \Phi_\Ga$, where the sum runs over a certain family of graphs, $\Phi_\Ga$ is a certain operator uniquely determined by
each graph $\Ga$ and
$c_\Ga$ is an absolutely convergent integral, $\int_{{C}_{\bu,\bu}(\Ga)}\Omega_\Ga$,
of a smooth differential form $\Omega_\Ga$ over a certain configuration space of points in a 3-dimensional subspace, $\caH$, of the Cartesian product, $\overline{\bbH}\times \overline{\bbH}$, of two copies of the closed upper-half plane. Our construction goes in two steps,
$$
\xymatrix{
\Ba{c}\frame{\mbox{$\Ba{c} \ \mathrm{Lie\ bialgebra} \ \\
 \  \mathrm{structures\ in}\ \R^n\   \Ea$}}\Ea \ \ \ \ \ \
 \ar[r]^{} &
\ \ \ \ \ \
\Ba{c}
\frame{\mbox{$\Ba{c} \ \mathrm{Quantizable} \ \\  \mathrm{Lie\ bialgebra} \ \\
 \  \mathrm{structures\ in}\ \R^n  \Ea$}}
\Ea \ \ \ \ \ \
\ar[r]^{} &
\ \ \ \ \ \ \
\Ba{c}\frame{\mbox{$\Ba{c} \ \mathrm{Bialgebra\ structures}\\
\ (*_\hbar, \Delta_\hbar)\ \mathrm{in}\ \odot^\bu(\R^n)[[\hbar]]  \Ea$}}
\Ea
 }
$$
as in the case of quantization of Poisson structures.
We show in \S 5 an explicit universal formula for the first
 arrow (behind which lies a $\caL ie_\infty$ morphism in the full analogy to the Poisson case), and then in \S 6 an explicit universal formula for the second arrow. The composition of the two gives us an explicit formula
for a universal quantization of an arbitrary finite-dimensional Lie bialgebra,  one of the main results of our paper. This result raises, however, open questions on the classification theory of both maps above, and on the graph cohomology description
of a quantizable Lie bialgebra structure; here the situation is much less clear than in the Poisson case discussed above.

\sip

We remark that an explicit configuration space integral formula (based on a propagator which is a generalized function rather than a smooth differential form) for the quantization of finite dimensional Lie bialgebras was claimed in  B.\ Shoikhet's preprint \cite{Sh1}. Furthermore, an odd analog of the properad governing quantizable Lie bialgebras has been investigated in \cite{KMW}.

\subsection{Structure of the paper}
\S 2 is a self-consistent reminder on graph complexes and configuration space models for the 1-coloured operad $\caH olie_d$ of (degree shifted) strongly homotopy Lie algebras, and for the 2-coloured operad $\cM or(\caH olie_d)$ of their morphisms.

\sip

In \S 3 we obtain explicit universal formulae for $\caL ie_\infty$ morphisms
relating Poisson (resp., Lie bialgebra) structures with their {\em quantizable}\, counterparts.

\sip

\S 4 shows a new explicit two-step formula for universal quantization of Poisson structures (depending only on a choice of a smooth function on the circle $S^1$ with support in the upper half of $S^1$), and proves classification claims (made in \S 1.2)
about every step in that construction.

\sip

\S 5 reminds key facts about the minimal resolutions, $\Assb_\infty$ and $\LB_\infty^{\mathrm{min}}$, of the  prop $\Assb$ of associative bialgebras and, respectively, of the prop $\LB$ of Lie bialgebras, and introduces a prop $\wLB_\infty^{\mathrm{quant}}$ of strongly homotopy {\em quantizable}\, Lie bialgebras. We use results of \S 3 to give an explicit transcendental morphism of dg props
$\wLB_\infty^{\mathrm{quant}} \rar \wLB_\infty^{\mathrm{min}, \circlearrowright}$, where $\wLB_\infty^{\mathrm{min}, \circlearrowright}$ is the wheeled closure of the completed version of the dg prop $\wLB_\infty^{\mathrm{min}}$, and hence an explicit morphism   $ \wLB^{\mathrm{quant}} \rar \wLB^\circlearrowright$ from the prop of quantizable Lie bialgebras into the wheeled closure of the completed prop of ordinary Lie bialgebras.

\sip

In \S 6 we show an explicit transcendental formula for a morphism of props
$\Assb \lon \caD(\wLB^{\mathrm{quant}})$, where $\caD$ is the polydifferential endofunctor on props introduced in \cite{MW2}, and show that it lifts by induction to a morphism of dg props $\Assb_\infty \lon \caD(\wLB^{\mathrm{quant}}_\infty)$. This gives us explicit formulae for a universal quantization of quantizable Lie bialgebras. Combining this formula with the explicit formula from \S 5, we obtain finally an explicit transcendental  formula for a universal quantization of ordinary finite-dimensional Lie bialgebras.

\sip

In Appendix A we prove a number of  Lemmas on vanishing of some classes of integrals involved into our formula for quantization of Lie bialgebras.

\sip

In Appendix B we construct surprisingly simple configuration space models for the bipermutahedron and biassociahedron posets introduced by Martin Markl in \cite{Ma} following an earlier work by  Samson Saneblidze and Ron Umble \cite{SU}.

%%%%%%%%%%%%%%%%%%%%%%%%

 \subsection{Some notation}
 The set $\{1,2, \ldots, n\}$ is abbreviated to $[n]$;  its group of automorphisms is
denoted by $\bS_n$;
the trivial one-dimensional representation of
 $\bS_n$ is denoted by $\id_n$, while its one dimensional sign representation is
 denoted by $\sgn_n$. The cardinality of a finite set $A$ is denoted by $\# A$. For a graph $\Ga$ its set of vertices (resp., edges) is denote by $V(\Ga)$ (resp., $E(\Ga)$).

\sip

We work throughout in the category of $\Z$-graded vector spaces over a field $\K$
of characteristic zero.
If $V=\oplus_{i\in \Z} V^i$ is a graded vector space, then
$V[k]$ stands for the graded vector space with $V[k]^i:=V^{i+k}$ and
and $s^k$ for the associated isomorphism $V\rar V[k]$; for $v\in V^i$ we set $|v|:=i$.
For a pair of graded vector spaces $V_1$ and $V_2$, the symbol $\Hom_i(V_1,V_2)$ stands
for the space of homogeneous linear maps of degree $i$, and
$\Hom(V_1,V_2):=\bigoplus_{i\in \Z}\Hom_i(V_1,V_2)$; for example, $s^k\in \Hom_{-k}(V,V[k])$.

\sip

For a
prop(erad) $\cP$ we denote by $\cP\{k\}$ a prop(erad) which is uniquely defined by
 the following property:
for any graded vector space $V$ a representation
of $\cP\{k\}$ in $V$ is identical to a representation of  $\cP$ in $V[k]$.
 The degree shifted operad of Lie algebras $\caL \mathit{ie}\{d\}$  is denoted by $\caL ie_{d+1}$ while its minimal resolution by $\caH olie_{d+1}$; representations of $\caL ie_{d+1}$ are vector spaces equipped with Lie brackets of degree $-d$.

\sip

For a right (resp., left) module $V$ over a group $G$ we denote by $V_G$ (resp.\
$_G\hspace{-0.5mm}V$)
 the $\K$-vector space of coinvariants:
$V/\{g(v) - v\ |\ v\in V, g\in G\}$ and by $V^G$ (resp.\ $^GV$) the subspace
of invariants: $\{\forall g\in  G\ :\  g(v)=v,\ v\in V\}$. If $G$ is finite, then these
spaces are canonically isomorphic as $char(\K)=0$.

%The symbol $\R^+$ stands for the multiplicative group of positive real numbers.
%We work in the category of semialgebraic
%  manifolds $X$, and denote by $\Omega_X^\bu$ the de Rham algebra of  piecewise %semialgebraic differential forms as defined in \cite{KoSo} and further developed %in \cite{HLTV}
% If $\om_1$ and $\om_2$
%are such forms on semialgebraic manifolds $X_1$ and, respectively, $X_2$, then %the form $p_1^*(\om_1)\wedge p_2^*(\om_2)$
%on $X_1\times X_2$,
%where $p_1: X_1\times X_2\rar X_1$ and  $p_2: X_1\times X_2\rar X_2$ are natural %projections,
%is often abbreviated to $\om_1\wedge \om_2$.

%%%%%%%%%%%%%%%%%%%%%%%%%%%%%%%%%%%%%%%%%%%%%%%%%%%%%%%%%%%%%%%%%%%%%
%%%%%%%%%%%%%%%%%%%%%%%%%%%%%%%%%%%%%%%%%%%%%%%%%%%%%%%%%%%%%%%%%%%%%

\bip

{\Large
\section{\bf Graph complexes and configuration spaces}
}

\mip

\subsection{Directed graph complexes}\label{2: subsec on DRGC}
Let $\cG_{k,l}$ be a set of directed graphs $\Ga$ with $k$ vertices and $l$  edges such that
some bijections $V(\Ga)\rar [k]$ and $E(\Ga)\rar [l]$ are fixed, i.e.\ every edge and every vertex of $\Ga$ has a numerical label. There is
a natural right action of the group $\bS_k \times  \bS_l$ on the set $\cG_{k,l}$ with $\bS_k$ acting by relabeling the vertices and  $\bS_l$ by relabeling the
edges. % and the $(\bS_2)^l$ via the reversing the directions of the edges.
For each fixed integer $d$, consider a collection of $\bS_k$-modules $\caD\cG ra_{d}=\{\caD\cG ra_d(k)\}_{k\geq 1}$, where
$$
\caD\cG ra_d(k):= \prod_{l\geq 0} \K \langle \cG_{k,l}\rangle \ot_{ \bS_l}  \sgn_l^{\ot |d-1|} [l(d-1)].
$$

It has an operad structure with the composition rule,
$$
\Ba{rccc}
\circ_i: &  \caD\cG ra_d(p) \times \caD\cG ra_d(q) &\lon & \caD\cG ra_d(p+q-1),  \ \ \forall\ i\in [n]\\
         &       (\Ga_1, \Ga_2) &\lon &      \Ga_1\circ_i \Ga_2,
\Ea
$$
given by substituting the graph $\Ga_2$ into the $i$-labeled vertex $v_i$ of $\Ga_1$ and taking the sum over  re-attachments of dangling edges (attached before to $v_i$) to vertices of $\Ga_2$
in all possible ways.

\sip

For any operad $\cP=\{\cP(k)\}_{n\geq 1}$  in the category of graded vector spaces,
the linear the map
$$
\Ba{rccc}
[\ ,\ ]:&  \sP \ot \sP & \lon & \sP\\
& (a\in \cP(p), b\in \cP(q)) & \lon &
[a, b]:= \sum_{i=1}^p a\circ_i b - (-1)^{|a||b|}\sum_{i=1}^q b\circ_i a\ \in \cP(p+q-1)
\Ea
$$
makes a graded vector space
$
\sP:= \prod_{k\geq 1}\cP(k)$
into a Lie algebra \cite{KM}; moreover, these brackets induce a Lie algebra structure on the subspace
of invariants
$
\sP^\bS:=  \prod_{n\geq 1}\cP(k)^{\bS_k}$. In particular,
the graded vector space
$$
\mathsf{dfGC}_{d}:= \prod_{k\geq 1} \cG ra_{d}(k)^{\bS_k}[d(1-k)]
$$
is a Lie algebra with respect to the above Lie brackets, and as such it can be identified
with the deformation complex $\Def(\caL ie_d\stackrel{0}{\rar} \cG ra_{d})$ of the zero morphism. Hence non-trivial Maurer-Cartan elements of $(\mathsf{dfGC}_{d}, [\ ,\ ])$ give us non-trivial morphisms of operads
\Beq\label{2: morphism from Lie to dGra}
i:\caL ie_d {\lon} \caD\cG ra_{d}.
\Eeq
 One such non-trivial morphism $f$ is given explicitly on the generator of $\caL ie_{d}$ by \cite{Wi}
\Beq\label{2:  map from Lie to dgra}
i \left(\Ba{c}\begin{xy}
 <0mm,0.66mm>*{};<0mm,3mm>*{}**@{-},
 <0.39mm,-0.39mm>*{};<2.2mm,-2.2mm>*{}**@{-},
 <-0.35mm,-0.35mm>*{};<-2.2mm,-2.2mm>*{}**@{-},
 <0mm,0mm>*{\circ};<0mm,0mm>*{}**@{},
   %<0mm,0.66mm>*{};<0mm,3.4mm>*{^1}**@{},
   <0.39mm,-0.39mm>*{};<2.9mm,-4mm>*{^{_2}}**@{},
   <-0.35mm,-0.35mm>*{};<-2.8mm,-4mm>*{^{_1}}**@{},
\end{xy}\Ea\right)=
\Ba{c}\resizebox{6.3mm}{!}{\xy
(0,1)*+{_1}*\cir{}="b",
(8,1)*+{_2}*\cir{}="c",
\ar @{->} "b";"c" <0pt>
\endxy}
\Ea  + (-1)^d
\Ba{c}\resizebox{7mm}{!}{\xy
(0,1)*+{_2}*\cir{}="b",
(8,1)*+{_1}*\cir{}="c",
\ar @{->} "b";"c" <0pt>
\endxy}
\Ea=:\xy
 (0,0)*{\bullet}="a",
(5,0)*{\bu}="b",
\ar @{->} "a";"b" <0pt>
\endxy
\Eeq
Note that elements of $\mathsf{dfGC}_{d}$ can be identified with graphs from $\caD\cG ra_d$ whose vertices' labels are symmetrized (for $d$ even) or skew-symmetrized (for $d$ odd) so that in pictures we can forget about labels of vertices  and denote them by unlabelled black bullets as in the formula above. Note also that graphs from  $\mathsf{dfGC}_{d}$ come equipped with an orientation, $or$, which is a choice of ordering of edges (for $d$ even) or a choice of ordering of vertices (for $d$ odd) up to an even permutation in both cases. Thus every graph $\Ga\in \mathsf{dfGC}_{d}$  has at most two different orientations, $or$ and $or^{opp}$, and one has
the standard relation, $(\Ga, or)=-(\Ga, or^{opp})$; as usual, the data $(\Ga, or)$ is abbreviated by $\Ga$ (with some choice of orientation implicitly assumed).  Note that the homological degree of graph $\Ga$ from $\mathsf{dfGC}_{d}$ is given by
$
|\Ga|=d(\# V(\Ga) -1) + (1-d) \# E(\Ga).
$

\sip

The above morphism (\ref{2:  map from Lie to dgra}) makes
 $(\mathsf{dfGC}_{d}, [\ ,\ ])$ into a {\em differential}\, Lie algebra with the differential
 $$
 \delta:= [\xy
 (0,0)*{\bullet}="a",
(5,0)*{\bu}="b",
\ar @{->} "a";"b" <0pt>
\endxy ,\ ].
 $$
 This dg Lie algebra contains a  dg subalgebra $\mathsf{dGC}_{d}$ spanned by connected graphs
with at least bivalent vertices.
% and with no bivalent vertices of the form  $\xy
%
% (0,0)*{}="a",
%(4,0)*{\bu}="b",
%(8,0)*{}="c",
%
%\ar @{->} "a";"b" <0pt>
%\ar @{->} "b";"c" <0pt>
%\endxy$.
It was proven in \cite{Wi} that
$$
H^\bu(\mathsf{dfGC}_{d})= \odot^{\bu\geq 1}\left(\mathsf{dGC}_{d}[-d]\right)[d]
$$
so that there is no loss of generality of working with $\mathsf{dGC}_{d}$ instead of
$\mathsf{dfGC}_{d}$. Moreover, one has an isomorphism of Lie algebras \cite{Wi},
$$
H^0(\mathsf{dGC}_{d})=\fg\fr\ft_1,
$$
where $\fg\fr\ft_1$ is the Lie algebra of the Grothendieck-Teichm\"u ller group $GRT_1$ introduced by Drinfeld in the context of deformation quantization of Lie bialgebras. Nowadays, this group play an important role in many  areas of mathematics.

\subsection{Oriented graph complexes}\label{2: subsec on oriented graph complexes} A graph $\Ga$ from the set $\cG_{k,l}$
is called {\em oriented}\, if it contains no {\em wheels}, that is, directed paths of edges
forming a closed circle; the subset of $\cG_{k,l}$ spanned by oriented graphs is denoted by $\cG^{or}_{k,l}$. It is clear that the subspace
$\cG ra_d^{or}\subset \caD \cG ra_d$ spanned by oriented graphs is a suboperad.
The morphism (\ref{2:  map from Lie to dgra}) factors through the inclusion $\cG ra_d^{or}\subset
\caD\cG ra_d$ so that one can consider a graph complex
$$
\mathsf{fGC}^{or}_d:=\Def\left(\caL ie_d \stackrel{i}{\rar} \cG ra_d^{or}\right)
$$
and its subcomplex $\GCor_d$ spanned by connected graphs with at least bivalent vertices
and with no bivalent vertices of the form  $\xy
 (0,0)*{}="a",
(4,0)*{\bu}="b",
(8,0)*{}="c",
\ar @{->} "a";"b" <0pt>
\ar @{->} "b";"c" <0pt>
\endxy$. This subcomplex determines the cohomology of the full graph complex,
$H^\bu(\mathsf{fGC}^{or}_d)=\odot^{\bu\geq 1} (H^\bu(\GCor_d)[-d])[d]$.
It was proven in \cite{Wi2} that
$$
H^\bu(\GCor_{d+1})=H^\bu(\mathsf{dGC}_d).
$$
In particular, one has a remarkable isomorphism of Lie algebras,
$
H^0(\GCor_3)=\fg\fr\ft$.
It was also proven in \cite{Wi2} that the cohomology group $H^1(\GCor_2)=H^1(\mathsf{dGC}_1)$ is one-dimensional and is spanned by the following graph
$$
\Upsilon_4:=
\la\Ba{c}\resizebox{11mm}{!}{\xy
(0,0)*{\bullet}="1",
(-7,16)*{\bullet}="2",
(7,16)*{\bullet}="3",
(0,10)*{\bullet}="4",
\ar @{<-} "2";"4" <0pt>
\ar @{<-} "3";"4" <0pt>
\ar @{<-} "4";"1" <0pt>
\ar @{<-} "2";"1" <0pt>
\ar @{<-} "3";"1" <0pt>
\endxy}\Ea
+
%\frac{1}{12}
2\la
\Ba{c}\resizebox{11mm}{!}{\xy
(0,0)*{\bullet}="1",
(-6,6)*{\bullet}="2",
(6,10)*{\bullet}="3",
(0,16)*{\bullet}="4",
\ar @{<-} "4";"3" <0pt>
\ar @{<-} "4";"2" <0pt>
\ar @{<-} "3";"2" <0pt>
\ar @{<-} "2";"1" <0pt>
\ar @{<-} "3";"1" <0pt>
\endxy}\Ea
+
%\frac{1}{24}
\la
 \Ba{c}\resizebox{11mm}{!}{\xy
(0,16)*{\bullet}="1",
(-7,0)*{\bullet}="2",
(7,0)*{\bullet}="3",
(0,6)*{\bullet}="4",
\ar @{->} "2";"4" <0pt>
\ar @{->} "3";"4" <0pt>
\ar @{->} "4";"1" <0pt>
\ar @{->} "2";"1" <0pt>
\ar @{->} "3";"1" <0pt>
\endxy}\Ea
 \ \ \ \ \ \forall \la\in \R\setminus 0.
$$

Moreover $H^2(\GCor_2)=\K$ and is spanned by a graph with four vertices. This means that one can construct by induction a new Maurer-Cartan element in the Lie algebra $\GCor_2$ (the integer subscript in the summand $\Upsilon_n$ stands for the number of vertices of graphs)
\Beq\label{2: KS MC element Upsilon}
\Upsilon_{KS}= \xy
 (0,0)*{\bullet}="a",
(6,0)*{\bu}="b",
\ar @{->} "a";"b" <0pt>
\endxy  + \Upsilon_4
+ \Upsilon_6 + \Upsilon_8 + \ldots
\Eeq
as all obstructions have more than $7$ vertices and hence do not hit the unique cohomology class
in $H^2(\GCor_2)$. Up to gauge equivalence, this new Maurer-Cartan element $\Upsilon$ is the {\em only}\, non-trivial deformation of the standard Maurer-Cartan element $\xy
 (0,0)*{\bullet}="a",
(6,0)*{\bu}="b",
\ar @{->} "a";"b" <0pt>
\endxy$.
We call this element  {\em Kontsevich-Shoikhet}\, one as it was first found by Boris Shoikhet in \cite{Sh} with a reference to an important contribution by Maxim Kontsevich via an informal communication.

\subsection{On a class of representations of graph complexes} Consider a formal power series algebra
$$
\f_n:=\K[[x^1,\ldots, x^n]]
$$
 in $n$ formal homogeneous variables and let
$\mathrm{Der}(\f_n)$ be the Lie algebra of continuous derivations of  $\f_n$.
Then, for any integer $d\geq 2$, the completed vector space
$$
\A_d^{(n)}:= \widehat{\odot^\bu} \left( \mathrm{Der}(\f_n)[d-1]\right)
$$
is canonically a $d$-algebra, that is, a graded commutative algebra equipped with compatible Lie brackets $[\ ,\ ]_S$ of degree $1-d$ (here $\widehat{\odot^\bu}$ stands for the completed graded symmetric tensor algebra functor).
One can identify  $\A_d^{(n)}$ with the ring of formal power series,
$$
\A_{d}^{(n)}:=\K[[x^1,\ldots, x^n, \psi_1,\ldots, \psi_n]]
$$
 generated by formal variables satisfying the condition
$$
|x^i| + |\psi_i|=d-1, \ \ \ d\in \Z,
$$
Then Lie bracket (of degree $1-d$) is given explicitly by
\Beq\label{2: standard Lie bracket in A_d}
[f_1, f_2]_S=\sum_{i=1}^n\frac{f_1\overleftarrow{\p}}{{\p} \psi_i}\frac{\overrightarrow{\p} f_2}{\p x^i} + (-1)^{|f_1||f_2|+ (d-1)(|f_1|+|f_2|)}
\frac{f_2\overleftarrow{\p}}{{\p} \psi_i}\frac{\overrightarrow{\p} f_1}{\p x^i}
\Eeq
 A degree $d$ element $\ga\in {\A}_d^{(n)}$ is called {\em Maurer-Cartan}\, if it satisfies
the condition $[\ga,\ga]_S=0$.

\sip

We are interested in an $n\rar \infty$ version of $\A_d^{(n)}$ which retains the above canonical $d$-algebra structure. Clearly, the sequence of canonical projections of graded vector spaces,
$$
\ldots \lon \A_{d}^{(n+2)} \lon \A_{d}^{(n)} \lon \A_{d}^{(n-1)}
$$
does not respect the above Lie bracket, so that the associated inverse limit $\displaystyle\lim_{\leftarrow} \A_{d}^{(n)}$ can not be a $d$-algebra.
There is a chain of injections of formal power series algebras,
$$
\ldots \lon \f_n \lon \f_{n+1} \lon \f_{n+2}\lon  \ldots
$$
and we denote the associated {\em direct}\, limit by
$$
\f_\infty:= \lim_{n\lon \infty} \f_{n}.
$$
Let $V_\infty$ stand for the infinite-dimensional graded vector space with the infinite basis
$\{x_1,x_2,\dots \}$ and set
$$
 \A_d^{\infty}:= \prod_{m\geq 0} \Hom\left(\odot^m(V_\infty[1-d]), \f_\infty\right)
$$
This is a vector subspace of the inverse limit
$$
\lim_{\leftarrow} \A_{d}^{(n)}=\K[[x^1,x^2,\ldots, \psi_1,\psi_2,\ldots]]
$$
spanned by formal power series in two infinite sets of graded commutative generators $X=\{x^1,x^2,\ldots\}$ and
$\Psi=\{ \psi_1,\psi_2,\ldots\}$ with the property that every monomial in generators from the set
$\Psi$ has as a coefficient a formal power series from the ring $\f_k$ for some finite number $k$.
Clearly, the subspace  $\A_d^{\infty}$ is a well-defined $d$-algebra.

\sip

\bip

The first interesting for application case has $d=2$, $|x^i|=0$ and $|\psi_i|=1$. The associated 2-algebra ${\A}_{2}^{(n)}$ can be identified
with the Gerstenhaber algebra  $\cT_{poly}(\R^n)$ of formal polyvector fields on $\R^n$
%\Beq\label{2: space V span of x^i}
%V:=\mathrm{linear\ span}\langle x^1,x^2,\ldots, x^n\rangle
%\Eeq
so that its Maurer-Cartan are formal power series  Poisson structures on $\R^n$. Its $n\rar \infty$ version  ${\A}_{2}^{\infty}$ gives us the Gerstenhaber algebra of polyvector fields  on the infinite-dimensional space $\R^\infty$.

\sip

The second interesting example has $d=3$ and  $|x^i|=|\psi_i|=1$. In this case Maurer-Cartan elements of $\A_{3}^{(n)}$
satisfying the conditions $\ga|_{x^i=0}=0$ and $\ga|_{\psi_i=0}=0$, $\forall\ i\in [n]$, are cubic polynomials
$$
\ga:=\sum_{i,j,k\in I}\left( C_{ij}^k \psi_kx^i x^j + \Phi_k^{ij}x^k \psi_i\psi_j
\right),
$$
and the equation $[\ga,\ga]=$ implies that the associated to the structure
constants  $\Phi_k^{ij}$  and, respectively, $C_{ij}^k$ linear maps,
$$
\bigtriangleup:\R^n\rar \wedge^2 \R^n,\ \ \  [\ ,\ ]:\wedge^2\R^n \rar \R^n
$$
define a Lie bialgebra structure in $\R^n$.

\sip

The above Lie brackets $[\ ,\ ]_S$ give us a representation
$$
\caL ie_d \lon \cE nd_{\A_d^{(n)}}
$$
for any $n\geq 1$.
 In fact, this representation factors through  morphism
(\ref{2: morphism from Lie to dGra}) and a canonical representation $\Phi$
\Beq\label{2: lie to dGra to End}
\Ba{cccc}
 \Phi: & d\cG ra_d  & {\lon} & \cE nd_{\A_{d}^{(n)}}\\
     &      \Ga     &      \lon             & \Phi_\Ga
          \Ea
\Eeq
of the operad $d\cG ra_d$ in $\A_{d}^{(n)}$ defined, for any $\Ga\in d\cG ra_d(k)$, by a linear map
\Beq\label{2: representation Phi_Gamma}
\Ba{rccc}
\Phi_\Ga: & \ot^k \A_{d}^{(n)} & \lon & \A_{d}^{(n)}\\
          & (f_1,f_2,\ldots,f_k) & \lon & \rho_\Ga(f_1,f_2,\ldots,f_k)
\Ea
\Eeq
where
$$
\Phi_\Ga(f_1,\ldots, f_k) :=m\left(\prod_{e\in E(\Ga)} \Delta_e \left(f_1(x, \psi)\ot f_2(x, \psi)\ot \ldots\ot f_k(x, \psi) \right)\right)
$$
and, for a directed edge $e$ connecting vertices labeled by integers $i$ and $j$,
$$
\Delta_e:= \sum_{a=1}^n \frac{\p}{\p x_{(i)}^a} \ot \frac{\p}{\p \psi_{a(j)}}
$$
with the subscript $(i)$ or $(j)$ indicating that the derivative operator is to be applied to the $i$-th or $j$-th factor in the tensor product.
The symbol $m$  denotes the multiplication map,
$$
\Ba{rccc}
m:&   \ot^k \A_{d}^{(n)} & \lon & \A_{d}^{(n)}\\
   & f_1\ot f_2\ot \ldots \ot f_k &\lon & f_1f_2\cdots f_k.
\Ea
$$
The morphism of dg operads (\ref{2: lie to dGra to End}) induces a morphism of the dg Lie algebras
$$
s: \mathsf{dfGC}_d:=\Def\left(\caL ie_d \stackrel{i}{\rar}  d\cG ra_d\right) \lon
CE^\bu\left(\A_{d}^{(n)}, \A_d^{(n)}\right):= \Def\left(\caL ie_d \stackrel{\Phi \circ i}{\lon}  \cE nd_{\A_d^{(n)}}\right).
$$
Here
$$
CE^\bu\left(\A_d^{(n)}, \A_d^{(n)}\right)=\mbox{Coder}\left(\odot^{\bu\geq 1}(\A_d^{(n)}[d])\right)
$$
is the standard Chevalley-Eilenberg deformation complex of the Lie algebra $\A_d^{(n)}$, that is,
the dg Lie algebra of coderivations of the graded co-commutative coalgebra $\odot^{\bu\geq 1}(\A_d^{(n)}[d])$. Therefore any Maurer-Cartan element $\ga$ in the graph complex $\mathsf{dfGC}_d$
gives a Maurer-Cartan element  $s(\Ga)$ in $\mbox{Coder}(\odot^{\bu\geq 1}(A_{d}^{(n)}[d]))$, that is a $\caH olie_d$ algebra structure in $\A_d^{\infty}$, for any {\em finite}\, number $n$. Moreover, if $\ga$ belongs to the Lie subalgebra $\fGC^{or}_d$, then the associated  $\caH olie_d$ structure remains well-defined
 in the limit  $n\rar +\infty$, i.e.\ it is well-defined in $\A_d^{\infty}$.

\subsubsection{\bf Example} The Maurer-Cartan element $\xy
 (0,0)*{\bullet}="a",
(5,0)*{\bu}="b",
\ar @{->} "a";"b" <0pt>
\endxy  \in \fGC_d^{or}\subset \fGC_d$ (see (\ref{2: lie to dGra to End})) gives rise to the standard Lie brackets (\ref{2: standard Lie bracket in A_d}) in
$\A_d^{(n)}$.

\subsubsection{\bf Example} The Maurer-Cartan element $\Upsilon_{KS}\in \fGC_2^{or}$ from (\ref{2: KS MC element Upsilon})
gives rise to a{\em Kontsevich-Shoikhet}\, $\caL ie_\infty$ structure in $\A_2^{(n)}=\cT_{poly}(\R^n)$,
$$
\left\{[\, \ , \ldots ,\ ]_{2k}: (\A_2^{(n)})^{\ot 2k}\rar \A_2^{(n)}[3-4k] \right\}_{k\geq 1}
$$
where
$$
[\, \ , \ldots ,\ ]_{2k}:=\Phi_{\Upsilon_{2k}}.
$$
 It was introduced by Boris Shoikhet in \cite{Sh} with a reference to an important contribution
by Maxim Kontsevich via a private communication. This  $\caL ie_\infty$ structure is well defined in the limit $n\rar +\infty$.

\mip

We shall consider in the next section some transcendental constructions of Maurer-Cartan elements in
$\fGC_d$ and
 $\fGC_d^{or}$ in which we shall use heavily the following configuration space models of classical operads.

\subsection{Configuration space model for the operad $\Holied$}
Let
$$
\Conf_k(\R^d):=\{p_1, \ldots, p_k\in \R^d\ |\, p_i\neq p_j\ \mbox{for}\ i\neq j\}
$$
be the configuration space of $k$ pairwise distinct points in $\R^d$, $d\geq 2$. The group $\R^+ \ltimes \R^d$ acts freely on each configuration space $\Conf_k(\R^d)$ for $k\geq 2$,
$$
(p_1, \ldots, p_k) \lon (\la p_1 + a, \ldots, \la p_k+a),\ \ \ \ \forall \la\in \R^+, a\in \R^d,
$$
so that the space of orbits,
 $$
C_{k}(\R^d):=\Conf_k(\C)/{\R^+ \ltimes \R^d},
$$
a smooth real $(kd-d-1)$-dimensional manifold.
The space $C_2(\R^d)$ is homeomorphic to the sphere $S^{d-1}$ and hence is compact.

For $\geq 3$ the compactified configuration
space $\overline{C}_k(\R^{d})$ is defined as the closure of an embedding \cite{Ko, Ga}
$$
\Ba{ccccc}
C_k(\R^d) & \lon &  (S^{d-1} )^{k(k-1)} &\times& (\R\P^2)^{k(k-1)(k-2)}\\
(p_1, \ldots, p_k) & \lon &  \prod_{i\neq j} \frac{p_i-p_j}{|p_i-p_j|} &\times&
 \prod_{i\neq j\neq l\neq i}\left[|p_{i}-p_{j}| :
|p_{j}-p_{l}|: |p_{i}-p_{l}|\right]
\Ea
$$
The space $\overline{C}_k(\C)$ is a smooth (naturally oriented) manifold with corners.
Its codimension 1 strata is given by
$$
\p \overline{C}_k(\C) = \bigsqcup_{A\subset [k]\atop \# A\geq 2} C_{k - \# A + 1}(\C)\times
 C_{\# A}(\C)
$$
where the summation runs over all possible  proper subsets of $[k]$ with cardinality $\geq 2$.
Geometrically, each such  stratum corresponds to the $A$-labeled elements of the set $\{p_1, \ldots, p_k\}$ moving
very close
to each other. If we represent $\overline{C}_k(\R^d)$ by the (skew)symmetric $k$-corolla of
degree\footnote{We prefer working with {\em co}chain complexes, and hence adopt gradings accordingly.}
 $k+1-kd$
\Beq\label{2: Lie_inf corolla}
\Ba{c}\resizebox{21mm}{!}{ \xy
(1,-5)*{\ldots},
(-13,-7)*{_1},
(-8,-7)*{_2},
(-3,-7)*{_3},
(7,-7)*{_{k-1}},
(13,-7)*{_k},
 (0,0)*{\circ}="a",
(0,5)*{}="0",
(-12,-5)*{}="b_1",
(-8,-5)*{}="b_2",
(-3,-5)*{}="b_3",
(8,-5)*{}="b_4",
(12,-5)*{}="b_5",
\ar @{-} "a";"0" <0pt>
\ar @{-} "a";"b_2" <0pt>
\ar @{-} "a";"b_3" <0pt>
\ar @{-} "a";"b_1" <0pt>
\ar @{-} "a";"b_4" <0pt>
\ar @{-} "a";"b_5" <0pt>
\endxy}\Ea=(-1)^{d|\sigma|}
\Ba{c}\resizebox{21mm}{!}{\xy
(1,-6)*{\ldots},
(-13,-7)*{_{\sigma(1)}},
(-6.7,-7)*{_{\sigma(2)}},
%(-3,-7)*{_{\sigma(3)}},
%(7,-8)*{_{n-1}},
(13,-7)*{_{\sigma(k)}},
 (0,0)*{\circ}="a",
(0,5)*{}="0",
(-12,-5)*{}="b_1",
(-8,-5)*{}="b_2",
(-3,-5)*{}="b_3",
(8,-5)*{}="b_4",
(12,-5)*{}="b_5",
\ar @{-} "a";"0" <0pt>
\ar @{-} "a";"b_2" <0pt>
\ar @{-} "a";"b_3" <0pt>
\ar @{-} "a";"b_1" <0pt>
\ar @{-} "a";"b_4" <0pt>
\ar @{-} "a";"b_5" <0pt>
\endxy}\Ea,
\ \ \ \forall \sigma\in \bS_k,\ k\geq2
\Eeq
then the boundary operator in the associated face complex of $\overline{C}_\bu(\R^d)$ takes a familiar form
\Beq\label{3: Lie_infty differential}
\p\hspace{-3mm}
\Ba{c}\resizebox{21mm}{!}{ \xy
(1,-5)*{\ldots},
(-13,-7)*{_1},
(-8,-7)*{_2},
(-3,-7)*{_3},
(7,-7)*{_{k-1}},
(13,-7)*{_k},
 (0,0)*{\circ}="a",
(0,5)*{}="0",
(-12,-5)*{}="b_1",
(-8,-5)*{}="b_2",
(-3,-5)*{}="b_3",
(8,-5)*{}="b_4",
(12,-5)*{}="b_5",
\ar @{-} "a";"0" <0pt>
\ar @{-} "a";"b_2" <0pt>
\ar @{-} "a";"b_3" <0pt>
\ar @{-} "a";"b_1" <0pt>
\ar @{-} "a";"b_4" <0pt>
\ar @{-} "a";"b_5" <0pt>
\endxy}\Ea
=
\sum_{A\varsubsetneq [k]\atop
\# A\geq 2} \pm
\Ba{c}\resizebox{25mm}{!}{
\begin{xy}
<10mm,0mm>*{\circ},
<10mm,0.8mm>*{};<10mm,5mm>*{}**@{-},
<0mm,-10mm>*{...},
<14mm,-5mm>*{\ldots},
<13mm,-7mm>*{\underbrace{\ \ \ \ \ \ \ \ \ \ \ \ \  }},
<14mm,-10mm>*{_{[k]\setminus A}};
<10.3mm,0.1mm>*{};<20mm,-5mm>*{}**@{-},
<9.7mm,-0.5mm>*{};<6mm,-5mm>*{}**@{-},
<9.9mm,-0.5mm>*{};<10mm,-5mm>*{}**@{-},
<9.6mm,0.1mm>*{};<0mm,-4.4mm>*{}**@{-},
<0mm,-5mm>*{\circ};
<-5mm,-10mm>*{}**@{-},
<-2.7mm,-10mm>*{}**@{-},
<2.7mm,-10mm>*{}**@{-},
<5mm,-10mm>*{}**@{-},
<0mm,-12mm>*{\underbrace{\ \ \ \ \ \ \ \ \ \ }},
<0mm,-15mm>*{_{A}},
\end{xy}}
\Ea
\Eeq
implying the following

\subsubsection{\bf Proposition \cite{GJ}}  {\em The fundamental chain complex
 of the family of
compactified configurations spaces,
$\{\overline{C}_k(\R^d)\}_{k\geq 2}$, has the structure of a dg free non-unital  operad canonically isomorphic
to the operad $\Holied$ of degree shifted strongly homotopy Lie algebras.}

\subsection{Configuration space model for the operad $\cM or(\Holied)$}
Let $\cM or(\Holied)$ be a two-coloured operad whose representations in a pair of dg vector spaces $V_{in}$ and $V_{out}$ is a triple $(\mu_{in}, \mu_{out}, F)$ consisting of  $\Holied$ structures $\mu_{in}$ on $V_{in}$
and $\mu_{out}$ on $V_{out}$, and  of a $\Holied$ morphism, $F:(V_{in},\mu_{in})\rar (V_{out},\mu_{out})$,
between them. There is a configuration space model \cite{Me1} for this operad which plays one of central roles in this paper.

\sip

The Abelian group
$\R^d$ acts freely,
$$
\Ba{ccc}
\Conf_k(\C) \times \R^d& \lon & \Conf_A(\C)\\
(p=\{p_i\}_{i\in [k]},a) &\lon & p+a:= \{p_i+a\}_{i\in [k]}
\Ea
$$
on the configuration space $\Conf_k(\R^d)$ for any $k\geq 1$
 so that
the quotient
$$
\fC_A(\R^d):=\Conf_A(\R^d)/\R^d
$$
is a $k(d -1)$-dimensional manifold.
There is a diffeomorphism,
$$
\Ba{rccccc}
\Psi_k: & \fC_k(\R^d) & \lon & C_k(\R^d) & \times & (0,1)\\
  &           p     & \lon & \frac{p-p_c}{|p-p_c|} && \frac{|p-p_c|}{1+ |p-p_c|}
\Ea
$$
where
$$
p_c:=\frac{1}{k}(p_1+\ldots + p_k).
$$

Note that the configuration  $\frac{p-p_c}{|p-p_c|}$ is invariant under the larger group $\R^+\ltimes \R^d$
and hence belongs to $C_k(\R^d)$. For any  non-empty subset
$A\subseteq [n]$ there is a natural map
$$
\Ba{rccc}
\pi_A : & \fC_k(\C) & \lon & \fC_A(\C)\\
        & p=\{p_i\}_{i\in [k]} & \lon & p_A:=\{p_i\}_{i\in A}
\Ea
$$
which forgets all the points labeled by elements of the complement $[k]\setminus A$.

\sip

The space $\fC_{1}(\R^d)$ is a point and hence is compact. For $k\geq 2$
a {\em semialgebraic compactification} $\overline{\fC}_k(\R^d)$ of $\fC_{k}(\R^d)$
can be defined as the closure of a composition \cite{Me1},
\Beq\label{2: first compactifn of fC(C)}
\fC_{k}(\R^d)\stackrel{\prod \pi_A}{\lon} \prod_{A\subseteq [k]\atop \# A\geq 2} \fC_{\# A}(\R^d)
\stackrel{\prod \Psi_A}{\lon} \prod_{A\subseteq [k]\atop  \# A\geq 2}
 C_{\# A}(\R^d)\times (0, 1) \hook \prod_{A\subseteq [k]\atop  \# A\geq 2}
 \overline{C}_{\# A}(\R^d)\times [0,1].
\Eeq
Thus all the limiting points in this compactification
come  from configurations in which a group or groups of points move too {\em close}\,  to each other
 within each group and/or a group or
groups of points which are moving too {\em far}\, (with respect to the standard Euclidean distance) away
from each other.
The codimension one boundary strata in $\widehat{\fC}_{n}(\R^d)$ correspond to the limit values $0$ or $+\infty$
of the parameters $|p-p_c|$, and are given by \cite{Me1}
\Beq\label{2:codimension 1 boundary in widehat{C}_n(C)}
\displaystyle
\p \overline{\fC}_{k}(\R^d) = \bigsqcup_{A\subseteq [n]\atop \# A\geq 2} \left(\overline{\fC}_{n - \# A + 1}(\R^d)\times
 \overline{C}_{\# A}(\R^d)\right)\
 \bigsqcup_{[k]=B_1\ \sqcup \ldots \sqcup B_i\atop{
 2\leq l\leq k \atop \#B_1,\ldots, \#B_l\geq 1}}\left( \overline{C}_{k}(\R^d)\times \overline{\fC}_{\# B_1}(\R^d)\times \ldots\times
 \overline{\fC}_{\# B_l}(\R^d)
 \right)
\Eeq
where
\Bi
\item
the first summation runs over all possible  subsets $A$ of $[k]$ and each summand corresponds to  $A$-labeled elements of the set $\{p_1, \ldots, p_k\}$ moving {\em close}\,
to each other,
\item  the
second summation runs over all possible decompositions of $[k]$ into $l\geq 2$  disjoint
non-empty subsets $B_1, \ldots, B_l$, and each summand  corresponds to $l$ groups of points (labeled, respectively,
 by disjoint ordered subsets $B_1, \ldots B_l$ of $[k]$) moving {\em far}\,  from each other while keeping relative
  distances within each group finite.
  \Ei

\sip

Note that the faces of the type $\overline{C}_\bu(\C)$ appear in (\ref{2:codimension 1 boundary in widehat{C}_n(C)})
 in two different ways --- as the strata describing collapsing points and as the strata controlling groups of
points going infinitely away from each other --- and they do {\em not}\, intersect in $\widehat{\fC}_{\bu}(\C)$. For that reason one has to assign to these two groups
of faces  different colours and represent collapsing $\overline{C}_k(\C)$-stratum  by, say,  white corolla
with straight legs as in (\ref{2: Lie_inf corolla}), the $\overline{C}_k(\R)$-stratum at ``infinity" by, say,
   a version of (\ref{2: Lie_inf corolla}) with ``broken" legs,
$
\hspace{-3mm}
\Ba{c}\resizebox{21mm}{!}
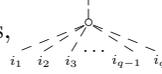
{\xy
(1,-5)*{\ldots},
(-13,-7)*{_{i_1}},
(-8,-7)*{_{i_2}},
(-3,-7)*{_{i_3}},
(7,-7)*{_{i_{q-1}}},
(13,-7)*{_{i_q}},
 (0,0)*{\circ}="a",
(0,5)*{}="0",
(-12,-5)*{}="b_1",
(-8,-5)*{}="b_2",
(-3,-5)*{}="b_3",
(8,-5)*{}="b_4",
(12,-5)*{}="b_5",
\ar @{--} "a";"0" <0pt>
\ar @{--} "a";"b_2" <0pt>
\ar @{--} "a";"b_3" <0pt>
\ar @{--} "a";"b_1" <0pt>
\ar @{--} "a";"b_4" <0pt>
\ar @{--} "a";"b_5" <0pt>
\endxy}\Ea
$. The face $\overline{\fC}_{n}(\C)$ can be represented pictorially by a 2-coloured (skew)symmetric corolla with black vertex,
$\hspace{-3mm}
\Ba{c}\resizebox{21mm}{!}
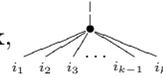
{\xy
(1,-5)*{\ldots},
(-13,-7)*{_{i_1}},
(-8,-7)*{_{i_2}},
(-3,-7)*{_{i_3}},
(7,-7)*{_{i_{k-1}}},
(13,-7)*{_{i_k}},
 (0,0)*{\bullet}="a",
(0,5)*{}="0",
(-12,-5)*{}="b_1",
(-8,-5)*{}="b_2",
(-3,-5)*{}="b_3",
(8,-5)*{}="b_4",
(12,-5)*{}="b_5",
\ar @{--} "a";"0" <0pt>
\ar @{-} "a";"b_2" <0pt>
\ar @{-} "a";"b_3" <0pt>
\ar @{-} "a";"b_1" <0pt>
\ar @{-} "a";"b_4" <0pt>
\ar @{-} "a";"b_5" <0pt>
\endxy}\Ea
$ of degree $d(1-k)$.

Each space $\overline{\fC}_{k}(\R^d)$ has a natural structure of a smooth manifold with corners.

\mip

\subsubsection{\bf Proposition \cite{Me1}}\label{2: Propos on the face complex of Mor(Lie_infty)}
{\em The disjoint union
\Beq\label{3: Lie_infty config topol operad}
\underline{\fC}(\R^d):=\overline{C}_\bu(\R^d)\sqcup \overline{\fC}_{\bu}(\R^d)\sqcup
\overline{C}_\bu(\R^d)
\Eeq
is a 2-coloured operad in the category of semialgebraic manifolds (or smooth manifolds with corners). Its complex of fundamental chains
can be identified with the operad $\cM or(\Holied)$ which is a dg free  non-unital 2-coloured operad
generated by the  corollas,
$$
\cM or(\Holied):= \cF ree
\left\langle
\Ba{c}\resizebox{21mm}{!}{\xy
(1,-5)*{\ldots},
(-13,-7)*{_1},
(-8,-7)*{_2},
(-3,-7)*{_3},
(7,-7)*{_{p-1}},
(13,-7)*{_p},
 (0,0)*{\circ}="a",
(0,5)*{}="0",
(-12,-5)*{}="b_1",
(-8,-5)*{}="b_2",
(-3,-5)*{}="b_3",
(8,-5)*{}="b_4",
(12,-5)*{}="b_5",
\ar @{-} "a";"0" <0pt>
\ar @{-} "a";"b_2" <0pt>
\ar @{-} "a";"b_3" <0pt>
\ar @{-} "a";"b_1" <0pt>
\ar @{-} "a";"b_4" <0pt>
\ar @{-} "a";"b_5" <0pt>
\endxy}\Ea,
\ \ \
\Ba{c}\resizebox{21mm}{!}{\xy
(1,-5)*{\ldots},
(-13,-7)*{_1},
(-8,-7)*{_2},
(-3,-7)*{_3},
(7,-7)*{_{k-1}},
(13,-7)*{_k},
 (0,0)*{\bullet}="a",
(0,5)*{}="0",
(-12,-5)*{}="b_1",
(-8,-5)*{}="b_2",
(-3,-5)*{}="b_3",
(8,-5)*{}="b_4",
(12,-5)*{}="b_5",
\ar @{--} "a";"0" <0pt>
\ar @{-} "a";"b_2" <0pt>
\ar @{-} "a";"b_3" <0pt>
\ar @{-} "a";"b_1" <0pt>
\ar @{-} "a";"b_4" <0pt>
\ar @{-} "a";"b_5" <0pt>
\endxy}\Ea
\ \ \ ,
\Ba{c}\resizebox{21mm}{!}{\xy
(1,-5)*{\ldots},
(-13,-7)*{_1},
(-8,-7)*{_2},
(-3,-7)*{_3},
(7,-7)*{_{q-1}},
(13,-7)*{_q},
 (0,0)*{\circ}="a",
(0,5)*{}="0",
(-12,-5)*{}="b_1",
(-8,-5)*{}="b_2",
(-3,-5)*{}="b_3",
(8,-5)*{}="b_4",
(12,-5)*{}="b_5",
\ar @{--} "a";"0" <0pt>
\ar @{--} "a";"b_2" <0pt>
\ar @{--} "a";"b_3" <0pt>
\ar @{--} "a";"b_1" <0pt>
\ar @{--} "a";"b_4" <0pt>
\ar @{--} "a";"b_5" <0pt>
\endxy}\Ea
\right\rangle_{p,q\geq 2, k\geq 1}
$$
and equipped with a differential which is given on white corollas of both colours by
formula (\ref{3: Lie_infty differential}) and on  the black corollas by the following formula
  \Beqr\label{2: d on MorLie corollas}
\p
\Ba{c}\resizebox{21mm}{!}{\xy
(1,-5)*{\ldots},
(-13,-7)*{_1},
(-8,-7)*{_2},
(-3,-7)*{_3},
(7,-7)*{_{k-1}},
(13,-7)*{_k},
 (0,0)*{\bullet}="a",
(0,5)*{}="0",
(-12,-5)*{}="b_1",
(-8,-5)*{}="b_2",
(-3,-5)*{}="b_3",
(8,-5)*{}="b_4",
(12,-5)*{}="b_5",
\ar @{--} "a";"0" <0pt>
\ar @{-} "a";"b_2" <0pt>
\ar @{-} "a";"b_3" <0pt>
\ar @{-} "a";"b_1" <0pt>
\ar @{-} "a";"b_4" <0pt>
\ar @{-} "a";"b_5" <0pt>
\endxy}\Ea
&=&
- \sum_{A\varsubsetneq [n]\atop
\# A\geq 2}
\Ba{c}\resizebox{25mm}{!}{
\begin{xy}
<10mm,0mm>*{\bu},
<10mm,0.8mm>*{};<10mm,5mm>*{}**@{--},
<0mm,-10mm>*{...},
<12mm,-5mm>*{\ldots},
<13mm,-7mm>*{\underbrace{\ \ \ \ \ \ \ \ \ \ \ \ \  }},
<14mm,-10mm>*{_{[k]\setminus A}};
<10.0mm,0mm>*{};<20mm,-5mm>*{}**@{-},
<10.0mm,-0mm>*{};<5mm,-5mm>*{}**@{-},
<10.0mm,-0mm>*{};<8mm,-5mm>*{}**@{-},
<10.0mm,0mm>*{};<0mm,-4.4mm>*{}**@{-},
<10.0mm,0mm>*{};<16.5mm,-5mm>*{}**@{-},
<0mm,-5mm>*{\circ};
<-5mm,-10mm>*{}**@{-},
<-2.7mm,-10mm>*{}**@{-},
<2.7mm,-10mm>*{}**@{-},
<5mm,-10mm>*{}**@{-},
<0mm,-12mm>*{\underbrace{\ \ \ \ \ \ \ \ \ \ }},
<0mm,-15mm>*{_{A}},
\end{xy}}
\Ea
\nonumber\\
&& +\ \, \sum_{l=2}^n \sum_{[k]=B_1\sqcup\ldots\sqcup B_l
\atop \inf B_1<\ldots< \inf B_l} \pm
\Ba{c}\resizebox{32mm}{!}{
\xy
(-15.5,-7)*{...},
(19,-7)*{...},
(7.5,0)*{\ldots},
(-17.8,-12)*{_{B_1}},
(-3.2,-12)*{_{B_2}},
(17.8,-12)*{_{B_l}},
(-1.8,-7)*{...},
%
%(-0.5,0)*{_{s_{(2)}}},
(-3.2,-9)*{\underbrace{\ \ \ \ \ \  \ \ \ \   }},
%
%(-14.5,3)*{_{s_{(1)}}},
(-17.8,-9)*{\underbrace{\ \ \ \ \ \  \ \ \ \   }},
%
%(15.5,3)*{_{s_{(k)}}},
(16.8,-9)*{\underbrace{\ \ \ \ \ \  \ \ \ \ \  }},
%
%(1,11.5)*{^{s_{(0)}}},
 (0,7)*{\circ}="a",
(-14,0)*{\bullet}="b_0",
%(-10,0)*{\bullet}="b_1",
(-4.5,0)*{\bullet}="b_2",
(14,0)*{\bullet}="b_3",
(0,13)*{}="0",
(1,-7)*{}="c_1",
(-8,-7)*{}="c_2",
(-5,-7)*{}="c_3",
(-22,-7)*{}="d_1",
(-19,-7)*{}="d_2",
(-13,-7)*{}="d_3",
(12,-7)*{}="e_1",
(15,-7)*{}="e_2",
(22,-7)*{}="e_3",
\ar @{--} "a";"0" <0pt>
\ar @{--} "a";"b_0" <0pt>
\ar @{--} "a";"b_2" <0pt>
\ar @{--} "a";"b_3" <0pt>
\ar @{-} "b_2";"c_1" <0pt>
\ar @{-} "b_2";"c_2" <0pt>
\ar @{-} "b_2";"c_3" <0pt>
\ar @{-} "b_0";"d_1" <0pt>
\ar @{-} "b_0";"d_2" <0pt>
\ar @{-} "b_0";"d_3" <0pt>
\ar @{-} "b_3";"e_1" <0pt>
\ar @{-} "b_3";"e_2" <0pt>
\ar @{-} "b_3";"e_3" <0pt>
\endxy}
\Ea.
\Eeqr
}

\subsubsection{\bf Example} As $\overline{C}_2(\R^n)=S^{d-1}$, the space
 $\widehat{\fC}_2(\R^d)$ is the closure
of the embedding
$$
\Ba{ccccccc}
\fC_2(\R^d) & \lon &  S^{d-1} &\times & (0,1) &\hook &  S^{d-1} \times  [0,+\infty]\\
(p_1,p_2) &\lon &   \frac{p_1-p_2}{|p_1-p_2|} &\times & \frac{|p_1-p_2|}{1+ |p_1-p_2|} &&
\Ea
$$
and hence can be identified with the closed $d$-dimensional cylinder
\Beq\label{2: fC_2(C) cylinder}
\overline{\fC}_2(\R^d)=\
\Ba{c}
{S^{d-1}_{out}}\\
\xy
(-8,0)*{}="a",
(8,0)*{}="b",
(-8,-15)*{}="a1",
(8,-15)*{}="b1",
(-8,0)*-{};(8,0)*-{};
**\crv{(0,6)}
**\crv{(0,-6)};
\ar @{-} "a";"a1" <0pt>
\ar @{-} "b";"b1" <0pt>
\endxy\vspace{-3mm}\\
\xy
(-8,0)*{},
(8,0)*{},
(-8,0)*-{};(8,0)*-{};
**\crv{(0,6)}
**\crv{(0,-6)};
\endxy
\vspace{1mm}
\\
{{S^{d-1}_{in}}}
\Ea.
\Eeq
where $S^{d-1}_{in}$ is the boundary corresponding to $|p_1-p_2|\rar 0$, and  $S^{d-1}_{out}$ is the boundary corresponding to $|p_1-p_2|\rar +\infty$. This is in accordance with
 the r.h.s.\ of (\ref{2: d on MorLie corollas}) for $k=2$  which is the sum of two terms, the first term corresponding to the bottom ``in" sphere $S^{d-1}$ (``two points collapsing to each other")
and upper ``out" sphere $S^{d-1}$ (``two points going infinitely far away from each other").

\mip

\bip

{\Large
\section{\bf Transcendental formulae for a class of  $\Holied$ morphisms}
}

\bip

\subsection{De Rham theories on operads of manifolds with corners} Let $X=\{X_k\}$
be a (a possibly coloured) operad on the category of semialgebraic manifolds (or smooth manifolds with corners), and $\fG=\{\fG(k)\}$ some dg  cooperad of graphs with the same set of coloures
(e.g., the dual cooperad of the operad $\caD\cG ra_d$ or $\caD \cG ra_d^{or}$ from \S 2). A {\em de Rham  $\fG$-theory on the operad $X$}\, is by definition a collection
of $\bS_n$-equivariant (and respecting colours) morphisms of complexes,
$$
\Ba{rccc}
\Omega_k: & \fG(k)& \lon & \Omega^\bu(X_k)\\
& \Ga & \lon & \Omega_\Ga
\Ea
$$
where  $\Omega^\bu(X_k)$ stands for the de Rham algebra of piecewise semialgebraic differential forms on
$X_k$, which satisfy the following compatibility condition: for any $k,l\in \N$ and any $i\in [k]$
the associated operad composition
$$
\circ_i: X_k \times X_l \lon X_{k+l-1}
$$
and the cooperad co-composition
$$
\Delta_i: \fG(k+l-1) \lon \fG(k)\ot \fG(l)
$$
makes the following diagram commutative,
\[
 \xymatrix{
  \fG(k+l-1)\ar[r]^{\Omega_{k+l-1}}\ar[d]_{\Delta_i} & \Omega^\bu(X_{k+l-1})\ar[r]^{\circ^*_i} & \Omega^\bu(X_{k}\times X_l)       \\
 \fG(k)\ot_\K \fG(l)\ar[r]_{\Omega_{k}\ot \Omega_l\ \ } &  \ \Omega^\bu(X_k)\ot_\K \Omega^\bu(X_{l})\ar[ur]_i &
}
\]
where
$$
\Ba{rccc}
i: & \Omega^\bu(X_k)\ot_\K \Omega^\bu(X_{l}) & \lon &  \Omega^\bu(X_{k}\times X_l)\\
&  \om_1 \ot \om_2 &\lon & \om_1 \wedge \om_2
\Ea
$$
is the natural inclusion.

\subsubsection{\bf Proposition}\label{3: From DRhamFT to MC element} {\em Let $\fG$ be the cooperad dual to the operad $\caD \cG ra_d$
(resp., to $\caD\cG ra^{or}_d$) equipped with the trivial differential. Then a de Rham $\fG$-theory on the operad of configuration spaces
 $\overline{C}_\bu(\R^d)=\{\overline{C}_k(\R^d)\}_{k\geq 2}$ gives rise to the following
 Maurer-Cartan element
 \Beq\label{3: MC element from DRham theory}
 \Upsilon:=\sum_{k\geq 2} \sum_{\Ga\in \fG(k)} \left(\int_{\overline{C}_k(\R^d)} \Omega_\Ga\right) \Ga
 \Eeq
in the (non-differential) Lie algebra\, $\dfGC_d$
(respectively,  in\, $\fGC_d^{or}$).}

 \bip

The second summation in (\ref{3: MC element from DRham theory}) runs over the set of generators of the vector space $\caD \cG ra_d(k)$ (resp., $\cG ra^{or}_d(k)$), and we assume $\int_{\overline{C}_k(\R^d)} \Omega_\Ga=0$ if $\deg \Omega_\Ga\neq \dim \overline{C}_k(\R^d)$.
This proposition is just a reformulation of Theorem 4.2.1 in \cite{Me0} so that we refer to that paper for its proof. It is worth noting that only connected graphs can give a non-zero contribution into the sum (\ref{3: MC element from DRham theory}).

\subsection{De Rham $\fG$-theories from propagators}
There is a large class of de Rham $\fG$-theories on $\overline{C}_\bu(\R^d)=\{\overline{C}_k(\R^d)\}_{k\geq 2}$ constructed as follows. Let $\om$
be an arbitrary differential top degree differential form on the sphere
$$
{C}_2(\R^d)=\overline{C}_2(\R^d)=S^{d-1}
$$
normalized so that
$$
\int_{S^{d-1}} \om=1.
$$
We call such a differential form a {\em propagator}. For any pair of distinct ordered numbers
$(i,j)$ with $i,j\in [k]$, consider a smooth map
$$
\Ba{rccc}
p_{ij}: & C_{k}(\R^d)  &\lon & C_{2}(\R^d) \\
        & (p_1, \ldots, p_k) & \lon & \frac{p_i-p_j}{|p_i-p_j|},
\Ea
$$
The pullback $\pi_{ij}^*(\om)$ defines a degree $d-1$ differential form on $C_{k}(\R^n)$
which extends smoothly to the compactification $\overline{C}_k(\R^d)$. In particular, for any directed graph $\Ga$ with $k$ labelled vertices   and any edge $e\in E(\Ga)$ there is an associated
 differential form $p_e^*(\om)\in \Omega^{d-1}_{\overline{C}_k(\R^d)}$, where $p_e:=p_{ij}$ if the edge $e$ begins at the vertex labelled by
 $i$ and ends at the vertex labelled by $j$.
Then, for $\fG$ being the cooperad dual to the operad $\caD \cG ra_d$,
consider a collection of maps
$$
\Ba{rccc}
\Omega_k: &  \fG(k) & \lon &  \Omega^\bu_{\overline{C}_k(\R^d)}\\
& \Ga & \lon & \Omega_\Ga:=\displaystyle \bigwedge_{e\in E(\Ga)}\hspace{-2mm}
{\pi}^*_e\left(\om\right).
\Ea
$$
It defines a de Rham $\fG$-theory on the operad $\overline{C}_\bu(\R^d)$
which in turn gives rise to a Maurer-Cartan element (\ref{3: MC element from DRham theory}) in $\fGC_d$ which in turn induces a $\Holied$ structure in $\A_d^{(n)}$,
\Beq\label{3: om Lie infty structure}
\mu^{\om}=\left\{\mu^{\om}_k: \ot^n \A_d^{(n)} \lon \A_d^{(n)}[ d+1-kd]   \right\}
\Eeq
given explicitly by
\Beq\label{3: mu_k^omega Lie infty structure}
\mu^{\om}_k=\sum_{\Ga\in \fG(k)} \left(\int_{\overline{C}_k(\R^d)} \Omega_\Ga\right) \Phi_\Ga.
\Eeq
As $\wedge^N \om=0$ for sufficiently large $N$, graphs with too many edges between any pair of vertices do not contribute into the sum in the r.h.s.\ of (\ref{3: mu_k^omega Lie infty structure}) so that the sum is finite and the formula is well-defined.

\sip

Note that an (oriented) graphs $\Ga$  with $k$ vertices can make a non-zero contribution  into (\ref{3: MC element from DRham theory}) or into $\mu^{\om}_k$ only if $d-1\mid kd-d-1$, i.e.\ if and only if
     $k=(d-1)l+2$ for some $l\in \N$;
     in that case the number of edges
    of $\Ga$ must be equal to $\frac{kd-d-1}{d-1}=dl+1$.
%kd-d-1=m(d-1), where m is the number of edges. Hence kd=d(m+1) -m +1, so m-1=dl for some l
% kd=d(dl+1)+d -dl so k=dl+1+1-l=l(d-1)+2 and the number of edges is dl+1

\mip

  Denote by $\sG_{k,l}$ (respectively, $\sG_{k,l}^{or}$)  the subset of the set $\cG_{k,l}$ (respectively, $\cG_{k,l}^{or}$)
  of directed (oriented) graphs consisting of connected graphs $\Ga$ such every vertex of $\Ga$ has valency $\geq 2$. Then we have the following sharpening of Proposition {\ref{3: From DRhamFT to MC element}}.

\subsubsection{\bf Proposition}\label{3: Prop on Upsilon^om and mu^om} {\em  For any propagator $\om$ on $S^{d-1}$, $d\geq 2$,  there is an associated Maurer-Cartan element
\Beq\label{3: Upsilon^omega}
\Upsilon^{\om}=
\Ba{c}\resizebox{7mm}{!}{\xy
(0,2)*+{_1}*\cir{}="b",
(8,2)*+{_2}*\cir{}="c",
\ar @{->} "b";"c" <0pt>
\endxy}
\Ea  - (-1)^d
\Ba{c}\resizebox{7mm}{!}{\xy
(0,2)*+{_2}*\cir{}="b",
(8,2)*+{_1}*\cir{}="c",
\ar @{->} "b";"c" <0pt>
\endxy}
\Ea
+ \sum_{l\geq 1}  \sum_{\Ga\in \sG_{l(d-1)+2,ld+1}} \left(\int_{\overline{C}_{l(d-1)+2}(\R^d)}
\displaystyle \bigwedge_{e\in E(\Ga)}\hspace{-2mm}
{\pi}^*_e\left(\om\right)\right) \Ga
\Eeq
in $\dfGC_d$, and an associated $\Holied$ algebra structure (\ref{3: om Lie infty structure}) can have  $\mu_k^\om$ non-vanishing only for $k= l(d-1)+2$ for some $l\in \N$, and with $\mu_2^{\om}$ given by the standard Schouten bracket (\ref{2: standard Lie bracket in A_d}).}

\begin{proof}
 It remains to check that (i) disconnected graphs and (ii) connected  directed graphs with univalent vertices do not contribute into the sum over $l\geq 1$. Let us show the second claim, the proof of the first claim being analogous (cf.\ \cite{Ko}).

 Let $\Ga\in \sG_{l(d-1)+2,ld+1}$, $l\geq 1$ be a connected directed graph with a univalent vertex $v\in V(\Ga)$, and  let $v'$ be the unique vertex connected to $v$ by the unique edge $e_{v,v'}$. Note that $v'$ has valency at least $2$ (as the $\Ga$ is connected and has $\geq 3$ vertices) so that
 there is a vertex $v''\in V(\Ga)\setminus v $ which is connected by an edge to $v'$.

  Let a $p_{v'}$ and  $p_{v''}$ be two different  points in $\R^d$ corresponding to the vertices $v'$ and respectively $v''$.
 Using the action of the group $\R^+ \ltimes \R^d$ on $\R^d$ we can put $p_{v'}$ into $0\in \R^d$ and $p_{v''}$ on the unital sphere
$S^{d-1}$ with center at $0$. The integral factorizes as follows
$$
\int_{{C}_{l(d-1)+2}(\R^d)}
\displaystyle \bigwedge_{e\in E(\Ga)}\hspace{-2mm}
{\pi}^*_e\left(\om\right)= \int_{\Conf_1(\R^d)} \pi^*_{e_{v,v'}}(\om) \cdot
\int_{M\subset \Conf_{l(d-1)}(\R^d) }
\displaystyle \bigwedge_{e\in E(\Ga)\setminus e_{v,v'}}\hspace{-2mm}
{\pi}^*_e\left(\om\right)
$$
The form $\bigwedge_{e\in E(\Ga)\setminus e_{v,v'}}\hspace{-2mm}
{\pi}^*_e\left(\om\right)$ has degree $ld(d-1)$ and $M$ is a subspace in $\Conf_{l(d-1)}(\R^d)$ of dimension $ld(d-1)-1$ (as one of the configuration points, $p_{v''}$, is restricted to lie on $S^{d-1}$). Hence the form $\bigwedge_{e\in E(\Ga)\setminus e_{v,v'}}\hspace{-2mm}
{\pi}^*_e\left(\om\right)$ vanishes identically on $M$ and the claim is proven.
\end{proof}

\subsubsection{\bf Example: the standard Schouten type bracket}
If one chooses the propagator
$$
\om_0:=\mathrm{Vol}_{S^{d-1}}
$$
to be the standard homogeneous (normalized to $1$) volume form on $S^{d-1}$
then, thanks to Kontsevich's Vanishing Lemma (proven for $d=2$ case in \cite{Ko} and for $d\geq 3$
in \cite{Ko0}), all integrals in the sum (\ref{3: Upsilon^omega}) over $l\geq 1$ vanish so that
\Beq\label{3: Upsilon_0}
\Upsilon^{\om_0}= \Ba{c}\resizebox{6.3mm}{!}{\xy
(0,1)*+{_1}*\cir{}="b",
(8,1)*+{_2}*\cir{}="c",
\ar @{->} "b";"c" <0pt>
\endxy}
\Ea  - (-1)^d
\Ba{c}\resizebox{7mm}{!}{\xy
(0,1)*+{_2}*\cir{}="b",
(8,1)*+{_1}*\cir{}="c",
\ar @{->} "b";"c" <0pt>
\endxy}=
\Ea=:\xy
 (0,0)*{\bullet}="a",
(5,0)*{\bu}="b",
\ar @{->} "a";"b" <0pt>
\endxy
\Eeq
The associated  $\Holied$ structure $\mu^{\om_0}$ in $\A_d^{(n)}$ is just the standard
Schouten bracket (\ref{2: standard Lie bracket in A_d}).

\subsubsection{\bf Example: a class of $\caL ie_\infty$ structures given by oriented graphs}
Let $g(x)$ be a non-negative function on the sphere
$$
S^{d-1}=\{(x_1,\ldots, x_d)\in \R^d\ \mid x_1^2 + \ldots x_d^2=1\}
$$
with compact support in the the upper ($x_d>0$) half of $S^{d-1}$
and normalized so that
$$
\int_{S^{d-1}} g\, \mathrm{Vol}_{S^{d-1}} =1.
$$
We can and will assume from now on that the function $g(x)$ on $S^{d-1}$ is invariant under the reflection in the $x_d$-axis,
$$
\sigma: \{x_i\rar -x_i\}_{1\leq i \leq d-1}, x_d\rar x_d.
$$
so that the propagator
\Beq\label{3: omega_g propagator}
\om_g:=  g\, \mathrm{Vol}_{S^{d-1}}
\Eeq
satisfies
\Beq\label{3: reflection sigma}
\sigma^*(\om_g)=(-1)^{d-1} \om_g
\Eeq
It is clear that only {\em oriented}\, graphs can give a non-trivial contribution into
the associated Maurer-Cartan element (\ref{3: Upsilon^omega}) (so that $\Upsilon^{\om_g}\in
\dfGC_d^{or}$) and that the associated $\Holied$ structure $\mu_{\om_g}$ on $\A_d^{(n)}$ is well-defined in the limit $n\rar +\infty$.

\sip

The imposed symmetry property (\ref{3: reflection sigma}) leads to vanishing of many terms in the sum (\ref{3: Upsilon^omega}).

\subsubsection{\bf Proposition}\label{3: Prop on Upsilon^om_g} {\em For any propagator $\om_g$ as above the associted MC element in $\dfGC_d^{or}$ has the form
\Beq\label{3: Usilon^omega_g}
\Upsilon^{\om_g}= \xy
 (0,0)*{\bullet}="a",
(5,0)*{\bu}="b",
\ar @{->} "a";"b" <0pt>
\endxy
+ \sum_{p\geq 1}  \sum_{\Ga\in \sG^{or}_{2p(d-1)+2,2pd+1}} \left(\int_{\overline{C}_{2p(d-1)+2}(\R^d)} \bigwedge_{e\in E(\Ga)}\hspace{-2mm}
{\pi}^*_e\left(\om_g\right)\right) \Ga
\Eeq
so that the associated $\Holied$ structure in $\A_d^{(n)}$ can have linear maps
$
\mu_k^{\om_g}\neq 0
$ only for $k=2p(d-1)+2 $,  $p\in N$.
}

\begin{proof} By Proposition {\ref{3: Prop on Upsilon^om and mu^om}}, $\mu_k^{\om_g}$ can be non-zero if and only if $k=(d-1)l+2$ for some $l\in \N$.  Let
$$
C_\Ga:=\int_{\overline{C}_{(d-1)l+2}(\R^d)} \bigwedge_{e\in E(\Ga)}\hspace{-2mm}
{\pi}^*_e\left(\om_g\right)
$$
be the weight of a summand $\Gamma\in G_{(d-1)l+2,dl+1}$ in $\mu_{(d-1)l+2}^{\om_g}$ or in $\Upsilon^g$. Using the translation freedom we can fix one of the vertices of $\Ga$ at $0\in \R^d$. If $\sigma$ stands for the reflection in the $x_d$ axis we have (cf.\ \cite{Ko,Sh}),
$$
\int_{\sigma(\overline{C}_{(d-1)l+2}(\R^d))} \bigwedge_{e\in E(\Ga)}\hspace{-2mm} {\pi}^*_e\left(\om\right)=
\int_{\overline{C}_{(d-1)l+2}(\R^d)}
\sigma^*\left({\pi}^*_e\left(\om\right)\right).
$$
As $\sigma(\overline{\R^d}_{(d-1)l+2}(\C))$ is equal to $\overline{C}_{(d-1)l+2}(\R^d)$ with orientation changed by the factor
$(-1)^{(k-1)(d-1)}$ and as $\sigma^*(\om_g)=(-1)^{d-1}\om_g$, we obtain an equality
$$
(-1)^{((d-1)l+2-1)(d-1)} C_\Ga= (-1)^{(dl+1)(d-1)}C_\Ga
$$
which implies $C_\Ga=0$ unless
$$
(d-1)l+1\equiv dl+1\bmod 2\Z
$$
i.e. unless $l=2p$ for some $p\in \N$.
\end{proof}

\subsubsection{\bf Example: a Kontsevich-Shoikhet $\caL ie_\infty$ structure}
If $d=2$, then only oriented graphs $\Ga$ with even number $2p$ of vertices contribute into $\Upsilon^g$, and the leading terms are given explicitly by \cite{Sh}
\Beq\label{3: Upsilon_g^2}
\Upsilon_{KS}^g:=\Upsilon^{\om_g}=  \xy
 (0,0)*{\bullet}="a",
(5,0)*{\bu}="b",
\ar @{->} "a";"b" <0pt>
\endxy \ + \ %\frac{1}{24}
\la\Ba{c}\resizebox{11mm}{!}{\xy
(0,0)*{\bullet}="1",
(-7,16)*{\bullet}="2",
(7,16)*{\bullet}="3",
(0,10)*{\bullet}="4",
\ar @{<-} "2";"4" <0pt>
\ar @{<-} "3";"4" <0pt>
\ar @{<-} "4";"1" <0pt>
\ar @{<-} "2";"1" <0pt>
\ar @{<-} "3";"1" <0pt>
\endxy}\Ea
+
%\frac{1}{12}
2\la
\Ba{c}\resizebox{11mm}{!}{\xy
(0,0)*{\bullet}="1",
(-6,6)*{\bullet}="2",
(6,10)*{\bullet}="3",
(0,16)*{\bullet}="4",
\ar @{<-} "4";"3" <0pt>
\ar @{<-} "4";"2" <0pt>
\ar @{<-} "3";"2" <0pt>
\ar @{<-} "2";"1" <0pt>
\ar @{<-} "3";"1" <0pt>
\endxy}\Ea
+
%\frac{1}{24}
\la
 \Ba{c}\resizebox{11mm}{!}{\xy
(0,16)*{\bullet}="1",
(-7,0)*{\bullet}="2",
(7,0)*{\bullet}="3",
(0,6)*{\bullet}="4",
\ar @{->} "2";"4" <0pt>
\ar @{->} "3";"4" <0pt>
\ar @{->} "4";"1" <0pt>
\ar @{->} "2";"1" <0pt>
\ar @{->} "3";"1" <0pt>
\endxy}\Ea \ \ \ +  \ldots =:\sum_{p\geq 1} \Upsilon_{2p}
\Eeq
for some $\la\in \R\setminus 0$.
In view of the homotopy uniqueness  of the Kontsevich-Shoikhet
element $\Upsilon_{KS}\in \fGC_3^{or}$, the sum $\Upsilon_{KS}^g$ must be gauge equivalent (with the gauge depending on the choice of a function $g$) to an element $\Upsilon_{KS}$ constructed by induction
in \S {\ref{2: subsec on oriented graph complexes}}.

Thus the propagator $\om_g$ induces a Kontsevich-Shoikhet $\Holie_2$  structure $\mu_{KS}^g$ in $\A_d^{(2)}= \cT_{poly}(\R^{n})$
\Beq\label{3: KS g-structure in T_poly}
[\ ,...,\ ]_{2p}:= \sum_{\Ga\in \sG_{2p, 4p-3}^{or}} \left(\int_{\overline{C}_{2p}(\R^d)} \displaystyle \bigwedge_{e\in E(\Ga)}\hspace{-2mm}
{\pi}^*_e\left(\om_g\right)\right) \Phi_\Ga: \cT_{poly}(\R^{n})^{\ot 2p} \rar  \cT_{poly}(\R^n)[3-4p]
\Eeq
which is isomorphic to the one
introduced in \cite{Sh}.

\sip

 We have two explicit $\Holied$ structures in $\A_2^{(n)}$,
the standard one (\ref{2: standard Lie bracket in A_d}) corresponding to the propagator $\om_0$
and the Kontsevich-Shoikhet one $\om_g$ corresponding to the propagator (\ref{3: omega_g propagator}). Shoikhet conjectured in \cite{Sh} that for $d=2$ these two structures are $\Holie_2$
isomorphic, i.e.\ there is a $\Holie_2$ isomorphism
$$
F:  \left(\cT_{poly}(\R^n), [\ ,\ ]_S\right) \lon  \left(\cT_{poly}(\R^n), [\ ,...,\ ]_{2p}, p\geq 1\right)
$$
 Stated in terms of graphs, this conjecture says that as Maurer-Cartan elements in $\mathsf{dfGC}_2$  the expressions (\ref{3: Upsilon_0}) and (\ref{3: Upsilon_g^2}) are gauge-equivalent to each other,
\Beq\label{2: gauge equiv of d and d_KS}
\Upsilon_S = e^{\mathrm{ad}_\Theta} \Upsilon_{KS}^g= e^{\mathrm{ad}_\Theta}\left(\sum_{p=1}^\infty \Upsilon_{2p}\right)
\Eeq
for some degree zero element $\Theta$  in $\mathsf{fGC}_2$. That this relation holds true is far from obvious. Indeed, let us attempt to construct $\Theta$ by induction on the number of vertices (as we managed to construct  $\Upsilon_{KS}$ above). The first step is easy --- the term $\Upsilon_4$ is
$\delta$-exact
in $\mathsf{dfGC}_2$,
$$
\Upsilon_4=\la \delta \left(\Ba{c}\resizebox{10mm}{!}{\xy
(0,0)*{\bullet}="1",
(-12,6)*{\bullet}="2",
(0,12)*{\bullet}="3",
   {\ar@/^0.6pc/(0,12)*{\bu};(0,0)*{\bu}};
 {\ar@/^0.6pc/(0,0)*{\bu};(0,12)*{\bu}};
\ar @{->} "2";"3" <0pt>
\ar @{<-} "1";"2" <0pt>
\endxy}\Ea
+
\Ba{c}\resizebox{10mm}{!}{\xy
(0,0)*{\bullet}="1",
(-12,6)*{\bullet}="2",
(0,12)*{\bullet}="3",
   {\ar@/^0.6pc/(0,12)*{\bu};(0,0)*{\bu}};
 {\ar@/^0.6pc/(0,0)*{\bu};(0,12)*{\bu}};
\ar @{<-} "2";"3" <0pt>
\ar @{->} "1";"2" <0pt>
\endxy}\Ea
     \right)
$$
and we can use the sum of two graphs of degree zero  inside the brackets to gauge away $\Upsilon_4$. However the next obstruction becomes a {\em wheeled}\, graph $\Upsilon_6'$ from $\mathsf{dfGC}_2$ so that starting with the second step all the obstruction classes land in $H^1(\mathsf{dfGC}_2)=H^1(\GC_2)$ (rather than in $H^1(\GC_2^{or})$), the same cohomology group where, according to Kontsevich \cite{Ko-f}, all the obstructions for the universal deformation  quantization of Poisson structures  lie.
Therefore, the formula for $\Theta$  must be as highly non-trivial as the Kontsevich quantization formula in \cite{Ko}.
One of our main results in this paper is such an explicit formula for $\Theta$ proving thereby  Shoikhet's conjecture (in fact, we show that it holds true for {\em any}\, value of the integer parameter $d$).

\sip

An MC element of the $\Holie_2$ algebra $\mu_{KS}^{\om_g}$
can be defined (to assure convergence) as a degree $2$  formal power series $\pi=\hbar \pi^\diamond(\hbar)$   for some
$\pi^\diamond(\hbar) \in \cT_{poly}(\R^n)[[\hbar]]$  satisfying the equation
$$
\frac{1}{2}[\pi,\pi]_2 + \frac{1}{4!} [\pi,\pi,\pi,\pi]_4
+ \ldots =0,
$$
or, equivalently,
$$
\frac{1}{2}[\pi^\diamond,\pi^\diamond]_2 + \frac{\hbar^2}{4!} [\pi^\diamond,\pi^\diamond,\pi^\diamond,\pi^\diamond]_4 + \frac{\hbar^4}{6!} [\pi^\diamond,\pi^\diamond,\pi^\diamond,\pi^\diamond,\pi^\diamond,\pi^\diamond]_6
+ \ldots =0.
$$
The equation is invariant under $\hbar\rar -\hbar$ so that it makes sense to look for solutions
$\pi^\diamond(\hbar)$ which are also invariant under $\hbar\rar -\hbar$, i.e.\ which are
formal power series in $\hbar^2$. Such solutions are precisely what we call quantizable Poisson structures, and making the change $\hbar^2 \rar \hbar$ we arrive at the defining equations in the Subsection 1.2.

\subsubsection{\bf Example: a class of Kontsevich-Shoikhet type $\caL ie_\infty$ structures in the case $d=3$} In this case one can apply a refined version {\ref{A: propos on hat{sG}} of  Proposition {\ref{3: Prop on Upsilon^om_g}} and write explicitly the associated Maurer-Cartan
\Beqr
\Upsilon^{\om_g}&=&\xy
 (0,0)*{\bullet}="a",
(5,0)*{\bu}="b",
\ar @{->} "a";"b" <0pt>
\endxy
+ \sum_{p\geq 2}  \sum_{\Ga\in \hat{\sG}^{or}_{4p+2,6p+1}} \left(\int_{\overline{C}_{4p+2}(\R^3)} \bigwedge_{e\in E(\Ga)}\hspace{-2mm}
{\pi}^*_e\left(\om_g\right)\right) \Ga \label{3: Upsilon om_g for d=3}
\Eeqr
and the associated  $\caH olie_3$ structure $\mu^{\om_{\bar{g}}}=\{\mu^{\om_{\bar{g}}}_{4p+2}\}_{p\geq 2}$ in $\A_3^{(n)}$ for any $n\in N$
\Beq\label{3: Holied om_bar{g} structure in A_3}
\mu_2=[\ ,\ ]_S \ \  \text{and}\ \  \mu^{\om_{\bar{g}}}_{4p+2}:= \sum_{\Ga\in \hat{\sG}^{or}_{4p+2, 6p+1}} \left(\int_{\overline{C}_{4p+2}(\R^3)} \displaystyle \bigwedge_{e\in E(\Ga)}\hspace{-2mm}
{\pi}^*_e\left(\om_g\right)\right) \Phi_\Ga
%: \ot^{4p+2}\A_3^{(n)} \rar  \A_3^{(n)}[-2-12p]
\ \
\text{for}\ p\geq 2.
\Eeq
using the subset of graphs $\hat{\sG}_{4p+2,6p+1}\subset \sG_{4p+2,6p+1}$ introduced in the Appendix A.
This gives us  a $3$-dimensional analogue of the Kontsevich-Shoikhet structure on polyvector fields.

\sip

Maurer-Cartan elements of the Lie algebra $(\A_3^{(n)}, [\ ,\ ]_S)$ are precisely (strongly homotopy)
Lie bialgebra structures in the vector space $V=\mathrm{span}\langle x_1,\ldots, x_n\rangle$.
Maurer-Cartan elements in the continuous $\Holie_3$ algebra $(\A_3^{(n)}[[\hbar]], \mu^{\om_{g}})$, that is, degree $3$ elements $\pi^\diamond\in \A_3^{(n)}[[\hbar]]$
satisfying the equation
$$
[\pi^\diamond,\ga^\diamond]_S+ \sum_{p\geq 2} \frac{\hbar^{p}}{(4p+2)!} \mu_{4p+2}^{\om_{{g}}}(\pi^\diamond,\pi^\diamond,\ldots, \pi^\diamond) =0,
$$
are called {\em quantizable Lie bialgebras}. We show in Section 7 below that the latter structures can be easily deformation quantized via an explicit formula.
We also show below an explicit formula for a universal (i.e.\ independent of a particular value of $n$) $\Holie_2$ morphism
$$
(\A_3^{(n)}, [\ ,\ ]_S) \lon (\A_3^{(n)}[[\hbar]], \mu^{\om_{{g}}}).
$$
The two formulae provide us with an explicit universal quantization of ordinary Lie bialgebras.
\sip

We shall be interested below in a special class of propagators of type  $\om_{g}$ on $S^2$ constructed as follows.  Consider the upper-half hemisphere,
$$
S^2_+:=\{(x,y,z)\in \R^3 \ \mid \ x^2+y^2+z^2=1,\ z>0\}
$$
and a well-defined smooth map
$$
\Ba{rccc}
\nu_+: & S_+^2 & \lon & S^1\times S^1\\
       & (x,y,z) & \lon & (\Arg(x+iz), \Arg(y+iz))
\Ea
$$
Let $\bar{g}(\theta)d\theta$ be a normalized volume form on the circle $S^1=\{e^{i\theta} \ |\ \theta \in [0,2\pi]\}$ as in (\ref{3: omega_g propagator}), i.e.\ the function $\bar{g}(\theta)$ has a compact support in the open interval $(0, \pi)$ and satisfies the standard conditions for $d=2$ propagator,
$$
\bar{g}(\theta)=\bar{g}(\pi-\theta), \ \ \ \  \int_{0}^{2\pi} \bar{g}(\theta)d\theta =1.
$$
Then
\Beq\label{3: propagator omega_bar{g}}
\om_{\bar{g}}:=\nu_+^*\left(\bar{g}(x+iz)\bar{g}(y+iz) dArg(x+iz)\wedge dArg(y+iz) \right)
\Eeq
is a smooth differential form on $S_+^2$ which extends by zero to a smooth differential form
on $S^2$ and which satisfies the standard conditions for $d=3$ propagator,
 $$
 \int_{S^2}\om_{\bar{g}} =1,
 $$
 and
 \Beqrn
 \sigma^*\left(\omega_{\bar{g}}\right)&=&\sigma^*\left(\nu_+^*\left(\bar{g}(x+iz)\bar{g}(y+iz)  dArg(x+iz)\wedge dArg(y+iz) \right)\right)\\
 &=& \nu_+^*\left(\bar{g}(-x+iz)\bar{g}(-y+iz) dArg(-x+iz)\wedge dArg(-y+iz) \right)\\
 &=& \nu_+^*\left(\bar{g}(x+iz)\bar{g}(y+iz) (-1)^2 dArg(x+iz)\wedge dArg(y+iz) \right)\\
 &=& \omega_{\bar{g}}.
 \Eeqrn
Hence the propagator $\om_{\bar{g}}$
belongs to the family of propagators\footnote{We apologize for some abuse of notations ---
the propagator $\om_{\bar{g}}$ is {\em not}\, equal to $\bar{g}\mathrm{Vol}_{S^2}$; the role of the bar in the notation $\om_{\bar{g}}$ is to emphasize this difference.} (\ref{3: omega_g propagator}) so that all the above claims hold true for  $\om_{\bar{g}}$.

\sip

The $1$-form
$$
\Omega_{\bar{g}}(\theta):= \bar{g}(\theta)d\theta
$$
has support in the open interval $(0,\pi)$ and hence it makes sense to consider its {\em iterated integrals},
\Beq\label{3: numbers Lambda_g}
\Lambda_{\bar{g}}^{(p)}:= \int_{0}^\pi\underbrace{\Omega_{\bar{g}}\Omega_{\bar{g}}\ldots \Omega_{\bar{g}}}_{p\ \text{times}}
\Eeq
which are some positive real numbers with $\Lambda_{\bar{g}}^{(1)}=1$.

\subsection{Transcendental formula for a class of  $\Holied$ morphisms}
Consider
an operad $\caD \cG ra_d$ and its 2-coloured version
\Beq\label{3: 2-colourd operad DGrad}
\underline{\caD \cG ra}_d=\left(\caD \cG ra^{out}_d, \caD \cG ra^{mor}_d, \caD \cG ra^{in}_d
\right)
\Eeq
consisting of three copies of $\caD \cG ra_d$: one copy is denoted by $\caD \cG ra^{out}_d$ and has inputs and outputs
in ``dashed" colour, the second copy is denoted by $\caD \cG ra^{mor}_d$ and has inputs in ``solid" colour and the output in the dashed colour, and the third copy is denoted by $\caD \cG ra^{in}_d$ and has both inputs and outputs in the ``solid" colour (cf.\ Proposition {\ref{2: Propos on the face complex of Mor(Lie_infty)}}). Therefore
for any $n,m\in \N$ and any $i\in [n]$ the only non-trivial operadic compositions are of the form
$$
\circ_i: \caD \cG ra^{out}_d(n)\ot \caD \cG ra^{out}_d(m) \lon \caD \cG ra^{out}_d(n+m-1),
\ \ \ \ \ \
\circ_i: \caD \cG ra^{in}_d(n)\ot \caD \cG ra^{in}_d(m) \lon \caD \cG ra^{in}_d(n+m-1),
$$
$$
\circ_i: \caD \cG ra^{out}_d(n)\ot \caD \cG ra^{mor}_d(m) \lon \caD \cG ra^{mor}_d(n+m-1),
\ \ \ \ \ \
\circ_i: \caD \cG ra^{mor}_d(n)\ot \caD \cG ra^{in}_d(m) \lon \caD \cG ra^{mor}_d(n+m-1).
$$
Similarly one defines a 2-coloured operad $\underline{\cG ra}_d^{or}$ of oriented graphs.
Let $\underline{\fG}$ and $\underline{\fG}^{or}$ be the cooperads dual to the operads
$\underline{\caD Gra}_d$ and $\underline{\cG ra}_d^{or}$ respectively.

\subsubsection{\bf Proposition} {\em Let
$
\underline{\fG}=\left(\fG^{in}, \fG^{mor}, \fG^{out}\right)
$
be the 2-coloured cooperad dual to the operad (\ref{3: 2-colourd operad DGrad}).
Then a de Rham $\underline{\fG}$-theory,
$$
\Omega=\left(\Omega^{in}, \Omega^{mor}, \Omega^{out}\right): \left(\fG^{in}, \fG^{mor}, \fG^{out}\right)
\lon \left(\Omega^\bu_{\overline{C}_\bu(\R^d)}, \Omega^\bu_{\overline{\fC}_{\bu}(\R^d)},
\Omega^{\bu}_{\overline{C}_\bu(\R^d)}\right)
$$
 on the 2-coloured operad of compactified configuration spaces $\underline{\fC}(\R^d)$ (see (\ref{3: Lie_infty config topol operad})) provides us with a a $\Holied$-isomorphism between $\Holied$-algebras,
$$
F: \left(\A_d^{(n)}, \mu_\bu^{\Ga_{in}}\right) \lon
\left(\A_d^{(n)}, \mu_\bu^{\Ga_{out}}\right) \ \ \ \ \ \forall\ n\in \N
$$
associated to  Maurer-Cartan elements
$$
\Upsilon_{in}:=\sum_{k\geq 2} \sum_{\Ga\in \fG^{in}(k)} \left(\int_{\overline{C}_k(\R^d)} \Omega_\Ga^{in}\right) \Ga \ \ \ \ \text{and} \ \ \ \ \ \Upsilon_{out}:=\sum_{k\geq 2} \sum_{\Ga\in \fG^{out}(k)} \left(\int_{\overline{C}_k(\R^d)} \Omega_\Ga^{out}\right) \Ga
$$
in $\dfGC_d$. This isomorphism
is given explicitly by the following formulae,
\Beq\label{3: F=F_k for morphism}
F=\left\{ F_k: \ot^k \A_d^{(n)} \lon \A_d^{(n)}[d-dk]\right\}_{k\geq 1}
\Eeq
where
$$
F_k:=  \sum_{\Ga\in \fG^{mor}(k)} \left(\int_{\overline{\fC}_k(\R^d)} \Omega_\Ga^{mor}\right) \Phi_\Ga
$$
}
\begin{proof} The claim follows from the de Rham theorem applied to the family of the compactified configuration spaces
 $\overline{\fC}_\bu(\R^d)$ and Proposition {\ref{2: Propos on the face complex of Mor(Lie_infty)}} (see \S 10.1 in \cite{Me1} for details).
\end{proof}

\subsection{An example} Consider a smooth degree $d-1$ differential form $\varpi$ on $\fC_2(\R^d)= S^{d-1}\times [0,1]$ such that its restrictions $\om_{in}:= \varpi\mid_{t=0}$ and $\om_{out}:=\varpi\mid_{t=1}$ give us top degree differential forms on
$\overline{C}_2(\R^d)=S^{d-1}$ such that $\int_{S^{d-1}}\om_{in}=1$ and
 $\int_{S^{d-1}}\om_{out}=1$. Then the collections of maps, $k\geq 1$,
$$
\Ba{rccc}
\Omega_k^{in}: &  \fG^{in}(k) & \lon &  \Omega^\bu_{\overline{C}_k(\R^d)}\\
& \Ga & \lon & \Omega_\Ga:=\displaystyle \bigwedge_{e\in E(\Ga)}\hspace{-2mm}
{\pi}^*_e\left(\om_{in}\right)
\Ea
\ \ \ \ \ \ \ \ \ \ \ \ \ \ \ \ \ \ \ \ \ \ \ \ \ \
\Ba{rccc}
\Omega_k^{out}: &  \fG^{out}(k) & \lon &  \Omega^\bu_{\overline{C}_k(\R^d)}\\
& \Ga & \lon & \Omega_\Ga:=\displaystyle \bigwedge_{e\in E(\Ga)}\hspace{-2mm}
{\pi}^*_e\left(\om_{out}\right)
\Ea
$$
$$
\Ba{rccc}
\Omega_k^{mor}: &  \fG^{mor}(k) & \lon &  \Omega^\bu_{\overline{\fC}_k(\R^d)}\\
& \Ga & \lon & \Omega_\Ga:=\displaystyle \bigwedge_{e\in E(\Ga)}\hspace{-2mm}
{\pi}^*_e\left(\varpi\right)
\Ea
$$
define a de Rham $\underline{\fG}$-theory on the 2-coloured operad $\underline{\fC}(\R^d)$ (see Theorem 10.1.1 in \cite{Me1} for a proof), and hence a $\Holied$
isomorphism (\ref{3: F=F_k for morphism}) of the associated $\Holied$ algebra structures in $\A_d^{(n)}$ for any $n$.

\sip

The  propagator (\ref{3: omega_g propagator}) satisfies the following equation
$$
\om_g = \mathrm{Vol}_{S^{d-1}} + d\Psi_g
$$
for some degree $d-2$ differential form  $\Psi_g$ on $S^{d-1}$. As $H^{d-2}(S^{d-1})$ equals zero for $d\geq 3$ and $\R$ for $d=2$, we can (and will) choose $\Psi_g$ in such a way
that (cf.\ (\ref{3: reflection sigma}))
% sigma^*(\Psi_g) - (-1)^{d-1}Psi_g=dA_g
%where sigma^*(dA_g)= - (-1)^{d-1} dA_g i.e. sigma^*(A_g) + (-1)^{d-1}A_g= dB_g
% so that dA_g= 1/2 dA_g - 1/2(-1)^{d-1}sigma^*(dA_g)= (-1)^{d-1}dA'_g - \sigma^*(dA'_g)
% for A'_g=1/2(-1)^{d-1}A_g
% i.e. sigma^*(\Psi_g + dA'_g) - (-1)^{d-1}(\Psi_g + dA'_g)=0
$$
\sigma^*(\Psi_g)=(-1)^{d-1} \Psi_g,
$$
where $\sigma: S^{d-1}\rar S^{d-1}$ is the reflection in the $x_d$-axis.

\sip

 Consider
next a differential form on $\fC_2(\R^d)$,
$$
\varpi_g:= \mathsf{Vol}_{S^{d-1}}\left(\frac{p_1-p_2}{|p_1-p_2|}\right) + \frac{|p_1-p_2|}{1+ |p_1-p_2|} d\Psi_g\left(\frac{p_1-p_2}{|p_1-p_2|}\right) +
(-1)^{d-1}\Psi_g\left(\frac{p_1-p_2}{|p_1-p_2|}\right)\wedge  d \left(\frac{|p_1-p_2|}{1+ |p_1-p_2|}\right)
$$
As it satisfies the conditions
$$
\varpi_g\mid_{S^{d-1}_{in}} = \mathsf{Vol}_{S^{d-1}}, \ \ \ \ \ \ \ \
\varpi_g\mid_{S^{d-1}_{out}} = \mathsf{Vol}_{S^{d-1}} + d\Psi_g= \om_g.
$$
and
\Beq\label{reflection for varpi_g}
\sigma^*(\varpi_g):=(-1)^{d-1}\varpi_g.
\Eeq
the associated $\underline{\fG}$-theory on $\underline{\fC}(\R^d)$ gives us almost immediately
the following result (which for $d=2$ proves the Shoikhet conjecture).

% dim C_k(R^d)= kd -d=d(k-1); a graph with k vertices can have non-zero weight iff
% the number of edges l satisfies the eqn d(k-1)=(d-1)l. Thus l=pd for some p and
% k-1=p(d-1). Thus only graphs Ga from G_{1+ p(d-1), pd} can have non-zero weights C_Ga.
% If we make a reflection we get an equality
% (-1)^{p(d-1)(d-1)} C_\Ga= (-1)^{(d-1)pd} C_\Ga so the weight can be non-zero only for
% p(d-1)=pd + 2N i.e. p=2q for some q.
% Conclusion: only graphs from G_{1+2q(d-1), 2qd} can contribute.

\subsubsection{\bf Theorem}\label{3: Theorem on KS iso for A_d^n} {\em
For any $d\geq 2$ and  any $n\geq 1$ there is a $\Holied$ isomorphism between the $\Holied$ algebras,
$$
F^{\om_g}: \left(\A_d^{(n)},  [\ ,\ ]_S\right) \lon \left(\A_d^{(n)}, \mu^{\om_g}           \right)
$$
which is given by (\ref{3: F=F_k for morphism}) with $F_k^{\om_g}$ possibly non-zero only for $k=1+2q(d-1)$, $q\in \Z^{\geq 0}$,
\Beq\label{3: formula for Holie_2 F_k}
F^{\om_g}_{1+2q(d-1)}=\sum_{\Ga\in \sG_{1+2q(d-1), 2qd}} \left(\int_{\overline{\fC}_k(\R^d)} \bigwedge_{e\in E(\Ga)}\hspace{-2mm}
{\pi}^*_e\left(\varpi_g\right) \right) \Phi_\Ga.
\Eeq
}
\begin{proof}
We have only to check that a connected directed graph $\Ga$ with all vertices of valency $\geq 2$ can give a non-trivial contribution to
the above formulae if and only if it belongs to the set $\sG_{1+2q(d-1), 2qd}$ for the non-negative integer $q$.

\sip

As $\dim \fC_k(\R^d)= kd -d=d(k-1)$ a directed graph $\Ga$ with $k$ vertices can have non-zero weight
$$
c_\Ga:=\int_{\overline{\fC}_k(\R^d)}\bigwedge_{e\in E(\Ga)}\hspace{-2mm}
{\pi}^*_e\left(\varpi_g\right)
$$
if and only if its number of edges, say $l$, satisfies the equation
 $$
 d(k-1)=(d-1)l.
$$ Thus
$l=pd$ for some $p\in \Z^{\geq 0}$ and
hence  $k-1=p(d-1)$. Thus only graphs $\Ga$ from $\sG_{1+ p(d-1), pd}$ can have  $c_\Ga\neq 0$.

\sip
 Using the translation freedom we can fix one of the vertices of $\Ga$ at $0\in \R^d$. Using the reflection $\sigma$ in the $x_d$ as in the proof of Proposition {\ref{3: Prop on Upsilon^om_g}} and formula (\ref{reflection for varpi_g}), we obtain an equality
 $$
  (-1)^{p(d-1)(d-1)} c_\Ga= (-1)^{(d-1)pd}c_\Ga
 $$
which implies $c_\Ga=0$ unless
$
p=2q
$
 for some non-negative integer $q$.
\end{proof}

\sip

This Theorem gives us an explicit gauge equivalence between between the Maurer-Cartan elements
$\Upsilon_S$ and $\Upsilon_{KS}^g$.
We use it below in the case $d=2$ to show that such gauge equivalences (and hence the homotopy classes of the associated universal $\Holied$ morphisms) are classified by the set of Drinfeld associators. In particular, the Grothendieck-Teichm\"uller group $GRT_1$ acts effectively and transitively on such gauge equivalences.

\subsubsection{\bf Corollary} {\em Given a Maurer-Cartan element $\pi\in \A_d^{(n)}$,
$$
[\pi,\pi]_S=0
$$
of the Lie algebra $(\A_d^{(n)}, [\ ,\ ]_S)$, the associated formal power series
\Beq\label{3: pi-diamond formula}
\pi^\diamond=\pi + \sum_{q=1}^\infty \frac{\hbar^q}{(1+2q(d-1))!} F^{\om_g}_{1+2q(d-1)}(\pi,\ldots,\pi)
\Eeq
in $\A_d^{(n)}[[\hbar]]$ satisfies the equation
\Beq\label{3: Mc eqn for quantizable str in general d}
[\pi^\diamond, \pi^\diamond]_S + \sum_{p\geq 1} \frac{\hbar^p}{(2p(d-1)+2)!} \mu_{2p(d-1)+2}^{\om_g}(\pi^\diamond,\ldots,\pi^\diamond)=0
\Eeq
}

\bip

In particular, the transcendental morphism $F^{\om_g}$ sends ordinary Poisson and Lie bialgebra structures
into {\em quantizable}\, ones
%(i.e.\ the ones satisfying (\ref{3: Mc eqn for quantizable str in general d}) for $d=2$ and %respectively $d=3$)
establishing thereby a 1-1 correspondence between their gauge equivalence classes: (i) given an ordinary Poisson/Lie bialgebra structure $\pi$ in $\R^n$, the above formal power series gives us a quantizable Poisson/Lie bialgebra structure $\pi^\diamond$, (ii) given a quantizable Poisson/Lie bialgebra structure $\pi^\diamond$ in $\R^n$, the initial term $\pi:=\pi^\diamond|_{\hbar=0}$ is an ordinary Poisson structure/Lie bialgebra.

\subsection{Remark} The Kontsevich-Shoikhet $\Holie_2$ structure on polyvector fields and the associated
$\Holie_2$ isomorphism  (\ref{3: formula for Holie_2 F_k})          have been defined above on the {\em affine}\, space $\R^n$ (as the formulae are
invariant only under the affine group, not under the group of diffeomorphisms). However both structures can be globalized, i.e.\ can be well-defined an arbitrary manifold $M$ using a torsion-free connection on $M$
as they both do not involve graphs with vertices which are univalent or have precisely one incoming edge and precisely one outgoing edge.

\bip

{\Large
\section{\bf A new  explicit  formula for universal quantizations of Poisson structures}
}

\bip

\subsection{The Kontsevich formula for a formality map}
Let $\Conf_{n,m}(\overline{\bbH})$ be the configuration space of injections
$z:[m+n]\hook \overline{\bbH}$ of the set $[m+n]$ into the closed upper-half plane such that
the following conditions are satisfied
\Bi
\item[(i)] for  $1 \leq i \leq m$ one has $z_i:=z(i)\in \R=\p \overline{\bbH}$\, and $z_1<z_2<\ldots < z_m$;
\item[(ii)] for $m+1 \leq i \leq m+n$ one has $z_i\in \bbH$.
\Ei
The group $\R^{+}\ltimes \R$ acts on this configuration space freely via $z_i\rar \la z_i +a$,
$\la\in \R^+$, $a\in \C$, so that the quotient space
$$
C_{n,m}({\bbH}):= \frac{\Conf_{n,m}(\overline{\bbH})}{ \R^{+}\ltimes \R}, \ \ \ \ \ \ 2n+m\geq 2,
$$
is a $2n+m-2$-dimensional manifold. Maxim Kontsevich constructed in \cite{Ko} its compactification
$\overline{C}_{n,m}({\bbH})$ as a smooth manifold with corners, and used it to construct an explicit $\Holie_2$ quasi-isomorphism of dg Lie algebras (for any $n\in \N)$,
$$
 \cF^K:  \left(\cT_{poly}(\R^n), [\ ,\ ]_S\right) \lon \left(C^\bu(\f_{\R^n}, \f_{\R^n})[1], d_H,\ \ [\ ,\ ]_{\mathrm{G}}\right)
 $$
where $(C^\bu(\f_{\R^n}, \f_{\R^n})[1], d_H)$ is the (degree shifted) Hochschild complex of the graded
commutative algebra $\f_{\R^n}=\K[[x_1,\ldots,x_n]]$ and $[\ ,\ ]_G$ are the Gerstenhaber brackets. This quasi-isomorphism
\Beq\label{4: F_k,l}
\cF^K=\left\{ \cF_{k,l}^K: \ot^k\f_{\R^n} \bigotimes \ot^l \cT_{poly}(\R^n)\lon   \f_{\R^n}\right\}_{2k+l\geq 2}
\Eeq
is given explicitly by
$$
\cF_{k,l}=\sum_{\Ga\in G_{k+l, l+2k-2}} \left(\int_{\overline{C}_{l,k}({\bbH})} \bigwedge_{e\in E(\Ga)}\hspace{-2mm}
{\nu}^*_e\left(\om_H\right) \right) \Phi_\Ga
$$
where
\Bi
\item $G_{k+l, l+2k-2}$ is the set of directed graphs with $k+l$ numbered vertices and
$l+2k-2$ edges such that the vertices with labels in the range from $1$ to $k$ have no outgoing edges, and for any $\Ga\in G_{k+l, l+2k-2}$ the associated operator $\Phi_\Ga:\ot^k\f_{\R^n} \bigotimes \ot^l \cT_{poly}(\R^n)\lon   \f_{\R^n}$ is given explicitly in \cite{Ko};

\item for an edge $e\in \Ga$ connecting a vertex with label $i$ to the vertex labelled $j$
$$
\nu_e: \overline{C}_{l,k}({\bbH}) \rar \overline{C}_{2,0}({\bbH})
$$
is the map forgetting all the points in the configuration space except $z_i$ and $z_j$;
\item $\om_H$ is a smooth 1-form on  $\overline{C}_{2,0}({\bbH})$ given explicitly by
$$
\om_H(z_i,z_j)=\frac{1}{2\pi} dArg\frac{z_i-z_j}{\overline{z}_i - z_j}
$$
\Ei

\subsection{A new explicit formula for the formality map}
Note that the 1-form (cf.\ (\ref{3: omega_g propagator}))
$$
\om_g(z_i,z_j)=g\left(\frac{\overline{z}_j-\overline{z}_i}{|z_i-z_j|}\right) dArg(\bar{z}_j -\bar{z}_i)
$$
is well defined on $\overline{C}_{2,0}({\bbH})$ so that it makes sense to consider
a collection of maps $\bar{\cF}=\{\cF_{k,l}\}_{2k+l\geq 2}$ as in (\ref{4: F_k,l}) with
\Beq\label{4: F_k,l from_g}
\bar{\cF}_{k,l}:=\sum_{\Ga\in G_{k+l, l+2k-2)}} \left(\int_{\overline{C}_{l,k}(\overline{\bbH})} \bigwedge_{e\in E(\Ga)}\hspace{-2mm}
{\nu}^*_e\left(\om_g\right) \right) \Phi_\Ga.
\Eeq
The propagator $\om_g$ does {\em not}\, satisfy Kontsevich's Vanishing Lemma 6.4 in \cite{Ko}  so that  many graphs $\Ga$ have non-trivial weights on the strata corresponding to groups of points collapsing to a point
inside $\bbH$; however all such graphs $\Ga$ are easy to describe --- they are precisely the ones which generate the Kontsevich-Shoikhet
$\Holie_2$ structure  $\{[\ ,...,\ ]_{2p}\}_{p\geq 1}$ in $\cT_{poly}(\R^n)$ so that Kontsevich's arguments lead us to the following

\subsubsection{\bf Proposition \cite{B}} {\em The formulae ({\ref{4: F_k,l from_g}}) provide us with
an explicit $\Holie_2$ quasi-isomorphism of $\Holie_2$ algebras
\Beq\label{4: F quasi-iso from KS to Hochshield}
\bar{\cF}: \left(\cT_{poly}(\R^n), \{[\ ,...,\ ]_{2p}\}_{p\geq 1}\right) \lon \left(C^\bu(\f_{\R^n}, \f_{\R^n})[1], d_H,\ \ [\ ,\ ]_{\mathrm{G}}\right).
\Eeq
Moreover, this quasi-isomorphism holds true in infinite dimensions, i.e.\ in the limit
$n\rar +\infty$.}

\begin{proof}
It remains to show the last claim about the limit $n\rar +\infty$. However it is obvious as
the only graphs $\Ga$ which can give a non-trivial contribution into the formula
(\ref{4: F_k,l from_g}) are {\em oriented}\, graphs, i.e.\ the ones which have no closed paths of directed edges.
\end{proof}

The above formulae are transcendental, i.e.\ involve an integration over configuration spaces.
However this $\Holie_2$ quasi-isomorphism can be constructed by a trivial (in the sense, independent of the choice of an associator) induction \cite{Wi2,B}.

\subsubsection{\bf Theorem} \label{4: Theorem on uniquenes of KS quantizations} {\em For any $n$ (including the limit $n\rar +\infty$) there is, up to homotopy equivalence, a unique
$\Holie_2$ quasi-isomorphism of $\Holie_2$ algebras as in (\ref{4: F quasi-iso from KS to Hochshield}).}

\mip
We refer to \cite{Wi2} and \cite{B} for two different proofs of this Theorem.

\sip

Now we can assemble  the previous results into a new proof of the Kontsevich formality theorem which gives us also an new explicit formula
for such a formality map ({\em not}\, involving the  2-dimensional hyperbolic geometry).

\subsubsection{\bf Kontsevich Formality Theorem} {\em For finite natural number $n$
there is a $\Holie_2$ quasi-isomorphism of dg Lie algebras}
$$
 \cF:  \left(\cT_{poly}(\R^n), [\ ,\ ]_S\right) \lon \left(C^\bu(\f_{\R^n}, \f_{\R^n})[1], d_H,\ \ [\ ,\ ]_{\mathrm{G}}\right)
 $$
\begin{proof} Let $g$ be an arbitrary smooth function on the circle $S^1$ with compact support in the  upper half of $S^1$ and normalized so that $\int_{S^1} g \mathrm{Vol_{S^1}}=1$. Then
there is an associated $\Holie_2$ isomorphism of $\Holie_2$ algebras
$$
F: \left(\cT_{poly}(\R^n),  [\ ,\ ]_S\right) \lon \left(\cT_{poly}(\R^n), [\ ,...,\ ]_{2p}, p\geq 1,              \right)
$$
given explicitly by formulae (\ref{3: F=F_k for morphism}), and a quasi-isomorphism
of $\Holie_2$ algebras (\ref{4: F quasi-iso from KS to Hochshield}) given by explicit formulae
(\ref{4: F_k,l from_g}). Hence we obtain the required $\Holie_2$ quasi-isomorphism
as the composition
\Beq\label{4: composition cF}
\cF: \left(\cT_{poly}(\R^n),  [\ ,\ ]_S\right) \stackrel{F}{\lon} \left(\cT_{poly}(\R^n), [\ ,...,\ ]_{2p}, p\geq 1,              \right)
\stackrel{\bar{\cF}}{\lon} \left(C^\bu(\f_{\R^n}, \f_{\R^n})[1], d_H,\ \ [\ ,\ ]_{\mathrm{G}}\right)
\Eeq
which is also given by explicit formulae with weights obtained from integrations on two different families of configuration spaces.
\end{proof}

It was proven in \cite{Do,Wi3}    that the set of homotopy classes of universal formality maps
$\{\cF\}$ can be identified with the set of Drinfeld associators, i.e.\ it is a torsor over the Grothendieck-Teichm\"uler group $GRT$. It follows from Theorem {\ref{3: Theorem on KS iso for A_d^n}} for $d=2$ that every such a quasi-isomorphism can be split as the composition
(\ref{4: composition cF}) with, by Theorem {\ref{4: Theorem on uniquenes of KS quantizations}},
the map $\bar{\cF}$ being unique (up to homotopy). Hence we obtain the following result.

\subsubsection{\bf Corollary}\label{4: Corollary on Holie_2 iso and Drinfeld ass} {\em The set of homotopy classes of universal $\Holie_2$ isomorphims
$$
F: \left(\cT_{poly}(\R^n),  [\ ,\ ]_S\right) \lon \left(\cT_{poly}(\R^n), [\ ,...,\ ]_{2p}, p\geq 1,              \right), \ \ \ \forall\ n\in \N,
$$
can be identified with the set of Drinfeld associators, i.e.\ it is a torsor over the Grothendieck-Teichm\"uler group $GRT_1$.}

\mip

We conclude that a construction of a non-commutative associative star product in $\f_{\R^n}$ out of an arbitrary ordinary Poisson structure $\pi$ can be split in two steps:
\Bi
\item[\sc Step 1] Associate to $\pi$ a quantizable Poisson structure $\pi^\diamond$. This step is most non-trivial and requires a choice of an associator; it can by given by an explicit formula (\ref{3: pi-diamond formula}).

\item[\sc Step 2] Construct a star product in $\f_{\R^n}$ using the unique (up to homotopy) quantization formulae
(\ref{4: F_k,l from_g}).
\Ei

We shall use a similar procedure below to obtain explicit and relatively simple formulae for the universal deformation  quantization
of arbitrary finite-dimensional Lie bialgebras.

\bip

\bip

{\Large
\section{\bf Props governing associative bialgebras, Lie bialgebras\\ and the formality maps}
}

\bip

\subsection{Prop of associative bialgebras and its minimal resolution.}
A prop  $\Assb$ governing associative bialgebras is the quotient,
$$
\Assb:= {\cF ree\langle A_0 \rangle}/(R)
$$
of the free prop $\cF ree\langle A_0 \rangle$ generated
%\footnote{Later we shall work with 2-coloured props so we reserve from now on the ``dashed colour"  %to $\Assb_\infty$ operations.}
by an $\bS$-bimodule $A_0=\{A_0(m,n)\}$\footnote{Here and everywhere all internal edges and legs in the graphical representation of an element of a prop are assumed to be implicitly oriented from the bottom of a graph to its top.},
\[
A_0(m,n):=\left\{
\Ba{rr}
\K[\bS_2]\ot \id_1\equiv\mbox{span}\left\langle
\begin{xy}
 <0mm,-0.55mm>*{};<0mm,-2.5mm>*{}**@{.},
 <0.5mm,0.5mm>*{};<2.2mm,2.2mm>*{}**@{.},
 <-0.48mm,0.48mm>*{};<-2.2mm,2.2mm>*{}**@{.},
 <0mm,0mm>*{\circ};<0mm,0mm>*{}**@{},
 <0mm,-0.55mm>*{};<0mm,-3.8mm>*{_1}**@{},
 <0.5mm,0.5mm>*{};<2.7mm,2.8mm>*{^2}**@{},
 <-0.48mm,0.48mm>*{};<-2.7mm,2.8mm>*{^1}**@{},
 \end{xy}
\,
,\,
\begin{xy}
 <0mm,-0.55mm>*{};<0mm,-2.5mm>*{}**@{.},
 <0.5mm,0.5mm>*{};<2.2mm,2.2mm>*{}**@{.},
 <-0.48mm,0.48mm>*{};<-2.2mm,2.2mm>*{}**@{.},
 <0mm,0mm>*{\circ};<0mm,0mm>*{}**@{},
 <0mm,-0.55mm>*{};<0mm,-3.8mm>*{_1}**@{},
 <0.5mm,0.5mm>*{};<2.7mm,2.8mm>*{^1}**@{},
 <-0.48mm,0.48mm>*{};<-2.7mm,2.8mm>*{^2}**@{},
 \end{xy}
   \right\rangle  & \mbox{if}\ m=2, n=1,\vspace{3mm}\\
\id_1\ot \K[\bS_2]\equiv
\mbox{span}\left\langle
\begin{xy}
 <0mm,0.66mm>*{};<0mm,3mm>*{}**@{.},
 <0.39mm,-0.39mm>*{};<2.2mm,-2.2mm>*{}**@{.},
 <-0.35mm,-0.35mm>*{};<-2.2mm,-2.2mm>*{}**@{.},
 <0mm,0mm>*{\circ};<0mm,0mm>*{}**@{},
   <0mm,0.66mm>*{};<0mm,3.4mm>*{^1}**@{},
   <0.39mm,-0.39mm>*{};<2.9mm,-4mm>*{^2}**@{},
   <-0.35mm,-0.35mm>*{};<-2.8mm,-4mm>*{^1}**@{},
\end{xy}
\,
,\,
\begin{xy}
 <0mm,0.66mm>*{};<0mm,3mm>*{}**@{.},
 <0.39mm,-0.39mm>*{};<2.2mm,-2.2mm>*{}**@{.},
 <-0.35mm,-0.35mm>*{};<-2.2mm,-2.2mm>*{}**@{.},
 <0mm,0mm>*{\circ};<0mm,0mm>*{}**@{},
   <0mm,0.66mm>*{};<0mm,3.4mm>*{^1}**@{},
   <0.39mm,-0.39mm>*{};<2.9mm,-4mm>*{^1}**@{},
   <-0.35mm,-0.35mm>*{};<-2.8mm,-4mm>*{^2}**@{},
\end{xy}
\right\rangle
\ & \mbox{if}\ m=1, n=2, \vspace{3mm}\\
0 & \mbox{otherwise}
\Ea
\right.
\]
modulo the ideal generated by relations
\Beq\label{2: bialgebra relations}
R:\left\{
\Ba{c}
\begin{xy}
 <0mm,0mm>*{\circ};<0mm,0mm>*{}**@{},
 <0mm,-0.49mm>*{};<0mm,-3.0mm>*{}**@{.},
 <0.49mm,0.49mm>*{};<1.9mm,1.9mm>*{}**@{.},
 <-0.5mm,0.5mm>*{};<-1.9mm,1.9mm>*{}**@{.},
 <-2.3mm,2.3mm>*{\circ};<-2.3mm,2.3mm>*{}**@{},
 <-1.8mm,2.8mm>*{};<0mm,4.9mm>*{}**@{.},
 <-2.8mm,2.9mm>*{};<-4.6mm,4.9mm>*{}**@{.},
   <0.49mm,0.49mm>*{};<2.7mm,2.3mm>*{^3}**@{},
   <-1.8mm,2.8mm>*{};<0.4mm,5.3mm>*{^2}**@{},
   <-2.8mm,2.9mm>*{};<-5.1mm,5.3mm>*{^1}**@{},
 \end{xy}\Ea
\ - \
\Ba{c}
\begin{xy}
 <0mm,0mm>*{\circ};<0mm,0mm>*{}**@{},
 <0mm,-0.49mm>*{};<0mm,-3.0mm>*{}**@{.},
 <0.49mm,0.49mm>*{};<1.9mm,1.9mm>*{}**@{.},
 <-0.5mm,0.5mm>*{};<-1.9mm,1.9mm>*{}**@{.},
 <2.3mm,2.3mm>*{\circ};<-2.3mm,2.3mm>*{}**@{},
 <1.8mm,2.8mm>*{};<0mm,4.9mm>*{}**@{.},
 <2.8mm,2.9mm>*{};<4.6mm,4.9mm>*{}**@{.},
   <0.49mm,0.49mm>*{};<-2.7mm,2.3mm>*{^1}**@{},
   <-1.8mm,2.8mm>*{};<0mm,5.3mm>*{^2}**@{},
   <-2.8mm,2.9mm>*{};<5.1mm,5.3mm>*{^3}**@{},
 \end{xy}\Ea=0, \ \ \ \ \
%%%%%%%%%%%%%%%%%%%%%%%%%%%%%%%%%%%%%%
 \Ba{c}\begin{xy}
 <0mm,0mm>*{\circ};<0mm,0mm>*{}**@{},
 <0mm,0.69mm>*{};<0mm,3.0mm>*{}**@{.},
 <0.39mm,-0.39mm>*{};<2.4mm,-2.4mm>*{}**@{.},
 <-0.35mm,-0.35mm>*{};<-1.9mm,-1.9mm>*{}**@{.},
 <-2.4mm,-2.4mm>*{\circ};<-2.4mm,-2.4mm>*{}**@{},
 <-2.0mm,-2.8mm>*{};<0mm,-4.9mm>*{}**@{.},
 <-2.8mm,-2.9mm>*{};<-4.7mm,-4.9mm>*{}**@{.},
    <0.39mm,-0.39mm>*{};<3.3mm,-4.0mm>*{^3}**@{},
    <-2.0mm,-2.8mm>*{};<0.5mm,-6.7mm>*{^2}**@{},
    <-2.8mm,-2.9mm>*{};<-5.2mm,-6.7mm>*{^1}**@{},
 \end{xy}\Ea
\ - \
 \Ba{c}\begin{xy}
 <0mm,0mm>*{\circ};<0mm,0mm>*{}**@{},
 <0mm,0.69mm>*{};<0mm,3.0mm>*{}**@{.},
 <0.39mm,-0.39mm>*{};<2.4mm,-2.4mm>*{}**@{.},
 <-0.35mm,-0.35mm>*{};<-1.9mm,-1.9mm>*{}**@{.},
 <2.4mm,-2.4mm>*{\circ};<-2.4mm,-2.4mm>*{}**@{},
 <2.0mm,-2.8mm>*{};<0mm,-4.9mm>*{}**@{.},
 <2.8mm,-2.9mm>*{};<4.7mm,-4.9mm>*{}**@{.},
    <0.39mm,-0.39mm>*{};<-3mm,-4.0mm>*{^1}**@{},
    <-2.0mm,-2.8mm>*{};<0mm,-6.7mm>*{^2}**@{},
    <-2.8mm,-2.9mm>*{};<5.2mm,-6.7mm>*{^3}**@{},
 \end{xy}\Ea=0,\ \ \ \ \ \
%%%%%%%%%%%%%%%%%%%%%%%%%%%%%%%%%%%%%%
\Ba{c} \begin{xy}
 <0mm,2.47mm>*{};<0mm,-0.5mm>*{}**@{.},
 <0.5mm,3.5mm>*{};<2.2mm,5.2mm>*{}**@{.},
 <-0.48mm,3.48mm>*{};<-2.2mm,5.2mm>*{}**@{.},
 <0mm,3mm>*{\circ};<0mm,3mm>*{}**@{},
  <0mm,-0.8mm>*{\circ};<0mm,-0.8mm>*{}**@{},
<0mm,-0.8mm>*{};<-2.2mm,-3.5mm>*{}**@{.},
 <0mm,-0.8mm>*{};<2.2mm,-3.5mm>*{}**@{.},
     <0.5mm,3.5mm>*{};<2.8mm,5.7mm>*{^2}**@{},
     <-0.48mm,3.48mm>*{};<-2.8mm,5.7mm>*{^1}**@{},
   <0mm,-0.8mm>*{};<-2.7mm,-5.2mm>*{^1}**@{},
   <0mm,-0.8mm>*{};<2.7mm,-5.2mm>*{^2}**@{},
\end{xy}\Ea
\ - \
\Ba{c}\begin{xy}
 <0mm,0mm>*{\circ};<0mm,0mm>*{}**@{},
 <0mm,-0.49mm>*{};<0mm,-3.0mm>*{}**@{.},
 <-0.5mm,0.5mm>*{};<-3mm,2mm>*{}**@{.},
 <-3mm,2mm>*{};<0mm,4mm>*{}**@{.},
 <0mm,4mm>*{\circ};<-2.3mm,2.3mm>*{}**@{},
 <0mm,4mm>*{};<0mm,7.4mm>*{}**@{.},
<0mm,0mm>*{};<2.2mm,1.5mm>*{}**@{.},
 <6mm,0mm>*{\circ};<0mm,0mm>*{}**@{},
 <6mm,4mm>*{};<3.8mm,2.5mm>*{}**@{.},
 <6mm,4mm>*{};<6mm,7.4mm>*{}**@{.},
 <6mm,4mm>*{\circ};<-2.3mm,2.3mm>*{}**@{},
 <0mm,4mm>*{};<6mm,0mm>*{}**@{.},
<6mm,4mm>*{};<9mm,2mm>*{}**@{.},
<6mm,0mm>*{};<9mm,2mm>*{}**@{.},
<6mm,0mm>*{};<6mm,-3mm>*{}**@{.},
   <-1.8mm,2.8mm>*{};<0mm,7.8mm>*{^1}**@{},
   <-2.8mm,2.9mm>*{};<0mm,-4.3mm>*{_1}**@{},
<-1.8mm,2.8mm>*{};<6mm,7.8mm>*{^2}**@{},
   <-2.8mm,2.9mm>*{};<6mm,-4.3mm>*{_2}**@{},
 \end{xy}
\Ea=0
\right.
\Eeq
 Note that the relations are not quadratic (it is proven, however, in \cite{MV} that $\Assb$ is {\em homotopy Koszul}). A minimal resolution, $(\Assb_\infty,\delta)$ of $\Assb$ exists \cite{Ma1} and is generated by an
$\bS$-bimodule $ A=\{ A(m,n)\}_{m,n\geq 1, m+n\geq 3}$,
\[
 A(m,n):= \K[\bS_m]\ot \K[\bS_n][m+n-3]=\mbox{span}\left\langle
\Ba{c}
\resizebox{19mm}{!}{\begin{xy}
 <0mm,0mm>*{\circ};<0mm,0mm>*{}**@{},
 <0mm,0mm>*{};<-8mm,5mm>*{}**@{.},
 <0mm,0mm>*{};<-4.5mm,5mm>*{}**@{.},
 <0mm,0mm>*{};<-1mm,5mm>*{\ldots}**@{},
 <0mm,0mm>*{};<4.5mm,5mm>*{}**@{.},
 <0mm,0mm>*{};<8mm,5mm>*{}**@{.},
   <0mm,0mm>*{};<-10.5mm,5.9mm>*{^{\tau(1)}}**@{},
   <0mm,0mm>*{};<-4mm,5.9mm>*{^{\tau(2)}}**@{},
   <0mm,0mm>*{};<10.0mm,5.9mm>*{^{\tau(m)}}**@{},
 <0mm,0mm>*{};<-8mm,-5mm>*{}**@{.},
 <0mm,0mm>*{};<-4.5mm,-5mm>*{}**@{.},
 <0mm,0mm>*{};<-1mm,-5mm>*{\ldots}**@{},
 <0mm,0mm>*{};<4.5mm,-5mm>*{}**@{.},
 <0mm,0mm>*{};<8mm,-5mm>*{}**@{.},
   <0mm,0mm>*{};<-10.5mm,-6.9mm>*{^{\sigma(1)}}**@{},
   <0mm,0mm>*{};<-4mm,-6.9mm>*{^{\sigma(2)}}**@{},
   <0mm,0mm>*{};<10.0mm,-6.9mm>*{^{\sigma(n)}}**@{},
 \end{xy}}\Ea
\right\rangle_{\tau\in \bS_n\atop \sigma\in \bS_m},
\]
The differential $\delta$ in $\Assb_\infty$ is not
quadratic, and its explicit value on generic $(m,n)$-corolla is not
known at present, but we can  (and will) assume from now on that $\delta$ preserves
the {\em path grading}\, of  $\Assb_\infty$ (which associates to any decorated graph $G$ from
$\Assb_\infty$ the total number of directed paths connecting input legs of $G$ to the output ones).
\mip

 Let $V$ be a $\Z$-graded
vector space over a field $\K$ of characteristic zero. The associated symmetric tensor algebra
$\f_V:= {\odot^{\bullet}} V= \oplus_{n\geq 0} \odot^n V$
comes equipped  with the standard graded commutative and
co-commutative bialgebra structure, i.e.\ there a non-trivial representation,
\Beq\label{2: rho_0}
\rho_0: \Assb \lon \cE nd_{\f_V}.
\Eeq
According to \cite{MV}, the (extended) deformation complex
$$
C_{GS}^\bu\left(\f_V,\f_V\right):=\Def\left(\Assb \stackrel{\rho_0}{\lon} \cE nd_{\f_V}\right)\simeq \prod_{m,n\geq 1}\Hom(\f_V^{\ot m}, \f_V^{\ot n})[2-m-n]
$$
and its polydifferential subcomplex $C_{poly}^\bu\left(\f_V,\f_V\right)$ come equipped with a $\caL ie_\infty$ algebra structure, $
\left\{\mu_n: \wedge^n C_{GS}^\bu(\f_V,\f_V)\lon C_{GS}^\bu(\f_V,\f_V)[2-n]\right\}_{n\geq 1},
$
such that $\mu_1$ coincides precisely with the Gerstenhaber-Shack differential \cite{GS}.
According to \cite{GS}, the cohomology of the complex $(C_{GS}^\bu(\f_V,\f_V), \mu_1)$ is precisely
the deformation complex
$$
\fg_V:=
\Def(\LB\stackrel{0}{\lon} \cE nd_V)
$$
controlling deformations of the zero morphism $0: \LB\rar \cE nd_V$, where $\LB$ is the prop
of Lie bialgebras which we discuss below.

\subsection{Prop governing Lie bialgebras and its minimal resolution}
The prop $\LB$ is defined \cite{D} as a quotient,
$$
\LB:= {\cF ree \langle E_0\rangle}/(R)
$$
of the free prop generated by an $\bS$-bimodule $E_0=\{E_0(m,n)\}$,
\Beq\label{5: module E_0 generating Lieb}
E_0(m,n):=\left\{
\Ba{rr}
sgn_2\ot \id_1\equiv\mbox{span}\left\langle
\begin{xy}
 <0mm,-0.55mm>*{};<0mm,-2.5mm>*{}**@{-},
 <0.5mm,0.5mm>*{};<2.2mm,2.2mm>*{}**@{-},
 <-0.48mm,0.48mm>*{};<-2.2mm,2.2mm>*{}**@{-},
 <0mm,0mm>*{\bu};<0mm,0mm>*{}**@{},
 <0mm,-0.55mm>*{};<0mm,-3.8mm>*{_1}**@{},
 <0.5mm,0.5mm>*{};<2.7mm,2.8mm>*{^2}**@{},
 <-0.48mm,0.48mm>*{};<-2.7mm,2.8mm>*{^1}**@{},
 \end{xy}
=-
\begin{xy}
 <0mm,-0.55mm>*{};<0mm,-2.5mm>*{}**@{-},
 <0.5mm,0.5mm>*{};<2.2mm,2.2mm>*{}**@{-},
 <-0.48mm,0.48mm>*{};<-2.2mm,2.2mm>*{}**@{-},
 <0mm,0mm>*{\bu};<0mm,0mm>*{}**@{},
 <0mm,-0.55mm>*{};<0mm,-3.8mm>*{_1}**@{},
 <0.5mm,0.5mm>*{};<2.7mm,2.8mm>*{^1}**@{},
 <-0.48mm,0.48mm>*{};<-2.7mm,2.8mm>*{^2}**@{},
 \end{xy}
   \right\rangle  & \mbox{if}\ m=2, n=1,\vspace{3mm}\\
\id_1\ot sgn_2\equiv
\mbox{span}\left\langle
\begin{xy}
 <0mm,0.66mm>*{};<0mm,3mm>*{}**@{-},
 <0.39mm,-0.39mm>*{};<2.2mm,-2.2mm>*{}**@{-},
 <-0.35mm,-0.35mm>*{};<-2.2mm,-2.2mm>*{}**@{-},
 <0mm,0mm>*{\bu};<0mm,0mm>*{}**@{},
   <0mm,0.66mm>*{};<0mm,3.4mm>*{^1}**@{},
   <0.39mm,-0.39mm>*{};<2.9mm,-4mm>*{^2}**@{},
   <-0.35mm,-0.35mm>*{};<-2.8mm,-4mm>*{^1}**@{},
\end{xy}=-
\begin{xy}
 <0mm,0.66mm>*{};<0mm,3mm>*{}**@{-},
 <0.39mm,-0.39mm>*{};<2.2mm,-2.2mm>*{}**@{-},
 <-0.35mm,-0.35mm>*{};<-2.2mm,-2.2mm>*{}**@{-},
 <0mm,0mm>*{\bu};<0mm,0mm>*{}**@{},
   <0mm,0.66mm>*{};<0mm,3.4mm>*{^1}**@{},
   <0.39mm,-0.39mm>*{};<2.9mm,-4mm>*{^1}**@{},
   <-0.35mm,-0.35mm>*{};<-2.8mm,-4mm>*{^2}**@{},
\end{xy}
\right\rangle
\ & \mbox{if}\ m=1, n=2, \vspace{3mm}\\
0 & \mbox{otherwise}
\Ea
\right.
\Eeq
modulo the ideal generated by the following relations
\Beq\label{3: LieB relations}
R:\left\{
\Ba{l}
\Ba{c}\begin{xy}
 <0mm,0mm>*{\bu};<0mm,0mm>*{}**@{},
 <0mm,-0.49mm>*{};<0mm,-3.0mm>*{}**@{-},
 <0.49mm,0.49mm>*{};<1.9mm,1.9mm>*{}**@{-},
 <-0.5mm,0.5mm>*{};<-1.9mm,1.9mm>*{}**@{-},
 <-2.3mm,2.3mm>*{\bu};<-2.3mm,2.3mm>*{}**@{},
 <-1.8mm,2.8mm>*{};<0mm,4.9mm>*{}**@{-},
 <-2.8mm,2.9mm>*{};<-4.6mm,4.9mm>*{}**@{-},
   <0.49mm,0.49mm>*{};<2.7mm,2.3mm>*{^3}**@{},
   <-1.8mm,2.8mm>*{};<0.4mm,5.3mm>*{^2}**@{},
   <-2.8mm,2.9mm>*{};<-5.1mm,5.3mm>*{^1}**@{},
 \end{xy}
\ + \
\begin{xy}
 <0mm,0mm>*{\bu};<0mm,0mm>*{}**@{},
 <0mm,-0.49mm>*{};<0mm,-3.0mm>*{}**@{-},
 <0.49mm,0.49mm>*{};<1.9mm,1.9mm>*{}**@{-},
 <-0.5mm,0.5mm>*{};<-1.9mm,1.9mm>*{}**@{-},
 <-2.3mm,2.3mm>*{\bu};<-2.3mm,2.3mm>*{}**@{},
 <-1.8mm,2.8mm>*{};<0mm,4.9mm>*{}**@{-},
 <-2.8mm,2.9mm>*{};<-4.6mm,4.9mm>*{}**@{-},
   <0.49mm,0.49mm>*{};<2.7mm,2.3mm>*{^2}**@{},
   <-1.8mm,2.8mm>*{};<0.4mm,5.3mm>*{^1}**@{},
   <-2.8mm,2.9mm>*{};<-5.1mm,5.3mm>*{^3}**@{},
 \end{xy}
\ + \
\begin{xy}
 <0mm,0mm>*{\bu};<0mm,0mm>*{}**@{},
 <0mm,-0.49mm>*{};<0mm,-3.0mm>*{}**@{-},
 <0.49mm,0.49mm>*{};<1.9mm,1.9mm>*{}**@{-},
 <-0.5mm,0.5mm>*{};<-1.9mm,1.9mm>*{}**@{-},
 <-2.3mm,2.3mm>*{\bu};<-2.3mm,2.3mm>*{}**@{},
 <-1.8mm,2.8mm>*{};<0mm,4.9mm>*{}**@{-},
 <-2.8mm,2.9mm>*{};<-4.6mm,4.9mm>*{}**@{-},
   <0.49mm,0.49mm>*{};<2.7mm,2.3mm>*{^1}**@{},
   <-1.8mm,2.8mm>*{};<0.4mm,5.3mm>*{^3}**@{},
   <-2.8mm,2.9mm>*{};<-5.1mm,5.3mm>*{^2}**@{},
 \end{xy}\Ea =0
 \ \ \ \ , \ \ \
%%%%%%%%%%%%%% Lie %%%%%%%%%%%%%%%%%%%%%%%%
 \Ba{c}\begin{xy}
 <0mm,0mm>*{\bu};<0mm,0mm>*{}**@{},
 <0mm,0.69mm>*{};<0mm,3.0mm>*{}**@{-},
 <0.39mm,-0.39mm>*{};<2.4mm,-2.4mm>*{}**@{-},
 <-0.35mm,-0.35mm>*{};<-1.9mm,-1.9mm>*{}**@{-},
 <-2.4mm,-2.4mm>*{\bu};<-2.4mm,-2.4mm>*{}**@{},
 <-2.0mm,-2.8mm>*{};<0mm,-4.9mm>*{}**@{-},
 <-2.8mm,-2.9mm>*{};<-4.7mm,-4.9mm>*{}**@{-},
    <0.39mm,-0.39mm>*{};<3.3mm,-4.0mm>*{^3}**@{},
    <-2.0mm,-2.8mm>*{};<0.5mm,-6.7mm>*{^2}**@{},
    <-2.8mm,-2.9mm>*{};<-5.2mm,-6.7mm>*{^1}**@{},
 \end{xy}
\ + \
 \begin{xy}
 <0mm,0mm>*{\bu};<0mm,0mm>*{}**@{},
 <0mm,0.69mm>*{};<0mm,3.0mm>*{}**@{-},
 <0.39mm,-0.39mm>*{};<2.4mm,-2.4mm>*{}**@{-},
 <-0.35mm,-0.35mm>*{};<-1.9mm,-1.9mm>*{}**@{-},
 <-2.4mm,-2.4mm>*{\bu};<-2.4mm,-2.4mm>*{}**@{},
 <-2.0mm,-2.8mm>*{};<0mm,-4.9mm>*{}**@{-},
 <-2.8mm,-2.9mm>*{};<-4.7mm,-4.9mm>*{}**@{-},
    <0.39mm,-0.39mm>*{};<3.3mm,-4.0mm>*{^2}**@{},
    <-2.0mm,-2.8mm>*{};<0.5mm,-6.7mm>*{^1}**@{},
    <-2.8mm,-2.9mm>*{};<-5.2mm,-6.7mm>*{^3}**@{},
 \end{xy}
\ + \
 \begin{xy}
 <0mm,0mm>*{\bu};<0mm,0mm>*{}**@{},
 <0mm,0.69mm>*{};<0mm,3.0mm>*{}**@{-},
 <0.39mm,-0.39mm>*{};<2.4mm,-2.4mm>*{}**@{-},
 <-0.35mm,-0.35mm>*{};<-1.9mm,-1.9mm>*{}**@{-},
 <-2.4mm,-2.4mm>*{\bu};<-2.4mm,-2.4mm>*{}**@{},
 <-2.0mm,-2.8mm>*{};<0mm,-4.9mm>*{}**@{-},
 <-2.8mm,-2.9mm>*{};<-4.7mm,-4.9mm>*{}**@{-},
    <0.39mm,-0.39mm>*{};<3.3mm,-4.0mm>*{^1}**@{},
    <-2.0mm,-2.8mm>*{};<0.5mm,-6.7mm>*{^3}**@{},
    <-2.8mm,-2.9mm>*{};<-5.2mm,-6.7mm>*{^2}**@{},
 \end{xy}\Ea =0
 \\
%%%%%%%%%%%%%%%%%%%%%%% Lie[1]Bi %%%%%%%%%%%%%%%
 \begin{xy}
 <0mm,2.47mm>*{};<0mm,0.12mm>*{}**@{-},
 <0.5mm,3.5mm>*{};<2.2mm,5.2mm>*{}**@{-},
 <-0.48mm,3.48mm>*{};<-2.2mm,5.2mm>*{}**@{-},
 <0mm,3mm>*{\bu};<0mm,3mm>*{}**@{},
  <0mm,-0.8mm>*{\bu};<0mm,-0.8mm>*{}**@{},
<-0.39mm,-1.2mm>*{};<-2.2mm,-3.5mm>*{}**@{-},
 <0.39mm,-1.2mm>*{};<2.2mm,-3.5mm>*{}**@{-},
     <0.5mm,3.5mm>*{};<2.8mm,5.7mm>*{^2}**@{},
     <-0.48mm,3.48mm>*{};<-2.8mm,5.7mm>*{^1}**@{},
   <0mm,-0.8mm>*{};<-2.7mm,-5.2mm>*{^1}**@{},
   <0mm,-0.8mm>*{};<2.7mm,-5.2mm>*{^2}**@{},
\end{xy}
\  - \
\begin{xy}
 <0mm,-1.3mm>*{};<0mm,-3.5mm>*{}**@{-},
 <0.38mm,-0.2mm>*{};<2.0mm,2.0mm>*{}**@{-},
 <-0.38mm,-0.2mm>*{};<-2.2mm,2.2mm>*{}**@{-},
<0mm,-0.8mm>*{\bu};<0mm,0.8mm>*{}**@{},
 <2.4mm,2.4mm>*{\bu};<2.4mm,2.4mm>*{}**@{},
 <2.77mm,2.0mm>*{};<4.4mm,-0.8mm>*{}**@{-},
 <2.4mm,3mm>*{};<2.4mm,5.2mm>*{}**@{-},
     <0mm,-1.3mm>*{};<0mm,-5.3mm>*{^1}**@{},
     <2.5mm,2.3mm>*{};<5.1mm,-2.6mm>*{^2}**@{},
    <2.4mm,2.5mm>*{};<2.4mm,5.7mm>*{^2}**@{},
    <-0.38mm,-0.2mm>*{};<-2.8mm,2.5mm>*{^1}**@{},
    \end{xy}
\  + \
\begin{xy}
 <0mm,-1.3mm>*{};<0mm,-3.5mm>*{}**@{-},
 <0.38mm,-0.2mm>*{};<2.0mm,2.0mm>*{}**@{-},
 <-0.38mm,-0.2mm>*{};<-2.2mm,2.2mm>*{}**@{-},
<0mm,-0.8mm>*{\bu};<0mm,0.8mm>*{}**@{},
 <2.4mm,2.4mm>*{\bu};<2.4mm,2.4mm>*{}**@{},
 <2.77mm,2.0mm>*{};<4.4mm,-0.8mm>*{}**@{-},
 <2.4mm,3mm>*{};<2.4mm,5.2mm>*{}**@{-},
     <0mm,-1.3mm>*{};<0mm,-5.3mm>*{^2}**@{},
     <2.5mm,2.3mm>*{};<5.1mm,-2.6mm>*{^1}**@{},
    <2.4mm,2.5mm>*{};<2.4mm,5.7mm>*{^2}**@{},
    <-0.38mm,-0.2mm>*{};<-2.8mm,2.5mm>*{^1}**@{},
    \end{xy}
\  - \
\begin{xy}
 <0mm,-1.3mm>*{};<0mm,-3.5mm>*{}**@{-},
 <0.38mm,-0.2mm>*{};<2.0mm,2.0mm>*{}**@{-},
 <-0.38mm,-0.2mm>*{};<-2.2mm,2.2mm>*{}**@{-},
<0mm,-0.8mm>*{\bu};<0mm,0.8mm>*{}**@{},
 <2.4mm,2.4mm>*{\bu};<2.4mm,2.4mm>*{}**@{},
 <2.77mm,2.0mm>*{};<4.4mm,-0.8mm>*{}**@{-},
 <2.4mm,3mm>*{};<2.4mm,5.2mm>*{}**@{-},
     <0mm,-1.3mm>*{};<0mm,-5.3mm>*{^2}**@{},
     <2.5mm,2.3mm>*{};<5.1mm,-2.6mm>*{^1}**@{},
    <2.4mm,2.5mm>*{};<2.4mm,5.7mm>*{^1}**@{},
    <-0.38mm,-0.2mm>*{};<-2.8mm,2.5mm>*{^2}**@{},
    \end{xy}
\ + \
\begin{xy}
 <0mm,-1.3mm>*{};<0mm,-3.5mm>*{}**@{-},
 <0.38mm,-0.2mm>*{};<2.0mm,2.0mm>*{}**@{-},
 <-0.38mm,-0.2mm>*{};<-2.2mm,2.2mm>*{}**@{-},
<0mm,-0.8mm>*{\bu};<0mm,0.8mm>*{}**@{},
 <2.4mm,2.4mm>*{\bu};<2.4mm,2.4mm>*{}**@{},
 <2.77mm,2.0mm>*{};<4.4mm,-0.8mm>*{}**@{-},
 <2.4mm,3mm>*{};<2.4mm,5.2mm>*{}**@{-},
     <0mm,-1.3mm>*{};<0mm,-5.3mm>*{^1}**@{},
     <2.5mm,2.3mm>*{};<5.1mm,-2.6mm>*{^2}**@{},
    <2.4mm,2.5mm>*{};<2.4mm,5.7mm>*{^1}**@{},
    <-0.38mm,-0.2mm>*{};<-2.8mm,2.5mm>*{^2}**@{},
    \end{xy}=0
\Ea
\right.
\Eeq
Its minimal resolution,
 $\LBm$,  is a dg free prop,
$$
\LBm=\cF ree \langle E\rangle,
$$
generated by the $\bS$--bimodule $ E=\{ E(m,n)\}_{m,n\geq 1, m+n\geq 3}$,
\Beq\label{5: generators of Lieb_infty}
 E(m,n):= sgn_m\ot sgn_n[m+n-3]=\mbox{span}\left\langle
\Ba{c}\resizebox{14mm}{!}{\begin{xy}
 <0mm,0mm>*{\bu};<0mm,0mm>*{}**@{},
 <-0.6mm,0.44mm>*{};<-8mm,5mm>*{}**@{-},
 <-0.4mm,0.7mm>*{};<-4.5mm,5mm>*{}**@{-},
 <0mm,0mm>*{};<-1mm,5mm>*{\ldots}**@{},
 <0.4mm,0.7mm>*{};<4.5mm,5mm>*{}**@{-},
 <0.6mm,0.44mm>*{};<8mm,5mm>*{}**@{-},
   <0mm,0mm>*{};<-8.5mm,5.5mm>*{^1}**@{},
   <0mm,0mm>*{};<-5mm,5.5mm>*{^2}**@{},
   <0mm,0mm>*{};<4.5mm,5.5mm>*{^{m\hspace{-0.5mm}-\hspace{-0.5mm}1}}**@{},
   <0mm,0mm>*{};<9.0mm,5.5mm>*{^m}**@{},
 <-0.6mm,-0.44mm>*{};<-8mm,-5mm>*{}**@{-},
 <-0.4mm,-0.7mm>*{};<-4.5mm,-5mm>*{}**@{-},
 <0mm,0mm>*{};<-1mm,-5mm>*{\ldots}**@{},
 <0.4mm,-0.7mm>*{};<4.5mm,-5mm>*{}**@{-},
 <0.6mm,-0.44mm>*{};<8mm,-5mm>*{}**@{-},
   <0mm,0mm>*{};<-8.5mm,-6.9mm>*{^1}**@{},
   <0mm,0mm>*{};<-5mm,-6.9mm>*{^2}**@{},
   <0mm,0mm>*{};<4.5mm,-6.9mm>*{^{n\hspace{-0.5mm}-\hspace{-0.5mm}1}}**@{},
   <0mm,0mm>*{};<9.0mm,-6.9mm>*{^n}**@{},
 \end{xy}}\Ea
\right\rangle,
\Eeq
and  with the differential given on generating corollas by \cite{MaVo,Va}
\Beq\label{3: differential in LieBinfty}
\delta
\Ba{c}\resizebox{14mm}{!}{\begin{xy}
 <0mm,0mm>*{\bu};<0mm,0mm>*{}**@{},
 <-0.6mm,0.44mm>*{};<-8mm,5mm>*{}**@{-},
 <-0.4mm,0.7mm>*{};<-4.5mm,5mm>*{}**@{-},
 <0mm,0mm>*{};<-1mm,5mm>*{\ldots}**@{},
 <0.4mm,0.7mm>*{};<4.5mm,5mm>*{}**@{-},
 <0.6mm,0.44mm>*{};<8mm,5mm>*{}**@{-},
   <0mm,0mm>*{};<-8.5mm,5.5mm>*{^1}**@{},
   <0mm,0mm>*{};<-5mm,5.5mm>*{^2}**@{},
   <0mm,0mm>*{};<4.5mm,5.5mm>*{^{m\hspace{-0.5mm}-\hspace{-0.5mm}1}}**@{},
   <0mm,0mm>*{};<9.0mm,5.5mm>*{^m}**@{},
 <-0.6mm,-0.44mm>*{};<-8mm,-5mm>*{}**@{-},
 <-0.4mm,-0.7mm>*{};<-4.5mm,-5mm>*{}**@{-},
 <0mm,0mm>*{};<-1mm,-5mm>*{\ldots}**@{},
 <0.4mm,-0.7mm>*{};<4.5mm,-5mm>*{}**@{-},
 <0.6mm,-0.44mm>*{};<8mm,-5mm>*{}**@{-},
   <0mm,0mm>*{};<-8.5mm,-6.9mm>*{^1}**@{},
   <0mm,0mm>*{};<-5mm,-6.9mm>*{^2}**@{},
   <0mm,0mm>*{};<4.5mm,-6.9mm>*{^{n\hspace{-0.5mm}-\hspace{-0.5mm}1}}**@{},
   <0mm,0mm>*{};<9.0mm,-6.9mm>*{^n}**@{},
 \end{xy}}\Ea
\ \ = \ \
 \sum_{[1,\ldots,m]=I_1\sqcup I_2\atop
 {|I_1|\geq 0, |I_2|\geq 1}}
 \sum_{[1,\ldots,n]=J_1\sqcup J_2\atop
 {|J_1|\geq 1, |J_2|\geq 1}
}\hspace{0mm}
(-1)^{\sigma(I_1\sqcup I_2)+ |I_1||I_2|+|J_1||J_2|}
\Ba{c}\resizebox{20mm}{!}{ \begin{xy}
 <0mm,0mm>*{\bu};<0mm,0mm>*{}**@{},
 <-0.6mm,0.44mm>*{};<-8mm,5mm>*{}**@{-},
 <-0.4mm,0.7mm>*{};<-4.5mm,5mm>*{}**@{-},
 <0mm,0mm>*{};<0mm,5mm>*{\ldots}**@{},
 <0.4mm,0.7mm>*{};<4.5mm,5mm>*{}**@{-},
 <0.6mm,0.44mm>*{};<12.4mm,4.8mm>*{}**@{-},
     <0mm,0mm>*{};<-2mm,7mm>*{\overbrace{\ \ \ \ \ \ \ \ \ \ \ \ }}**@{},
     <0mm,0mm>*{};<-2mm,9mm>*{^{I_1}}**@{},
 <-0.6mm,-0.44mm>*{};<-8mm,-5mm>*{}**@{-},
 <-0.4mm,-0.7mm>*{};<-4.5mm,-5mm>*{}**@{-},
 <0mm,0mm>*{};<-1mm,-5mm>*{\ldots}**@{},
 <0.4mm,-0.7mm>*{};<4.5mm,-5mm>*{}**@{-},
 <0.6mm,-0.44mm>*{};<8mm,-5mm>*{}**@{-},
      <0mm,0mm>*{};<0mm,-7mm>*{\underbrace{\ \ \ \ \ \ \ \ \ \ \ \ \ \ \
      }}**@{},
      <0mm,0mm>*{};<0mm,-10.6mm>*{_{J_1}}**@{},
 <13mm,5mm>*{};<13mm,5mm>*{\bu}**@{},
 <12.6mm,5.44mm>*{};<5mm,10mm>*{}**@{-},
 <12.6mm,5.7mm>*{};<8.5mm,10mm>*{}**@{-},
 <13mm,5mm>*{};<13mm,10mm>*{\ldots}**@{},
 <13.4mm,5.7mm>*{};<16.5mm,10mm>*{}**@{-},
 <13.6mm,5.44mm>*{};<20mm,10mm>*{}**@{-},
      <13mm,5mm>*{};<13mm,12mm>*{\overbrace{\ \ \ \ \ \ \ \ \ \ \ \ \ \ }}**@{},
      <13mm,5mm>*{};<13mm,14mm>*{^{I_2}}**@{},
 <12.4mm,4.3mm>*{};<8mm,0mm>*{}**@{-},
 <12.6mm,4.3mm>*{};<12mm,0mm>*{\ldots}**@{},
 <13.4mm,4.5mm>*{};<16.5mm,0mm>*{}**@{-},
 <13.6mm,4.8mm>*{};<20mm,0mm>*{}**@{-},
     <13mm,5mm>*{};<14.3mm,-2mm>*{\underbrace{\ \ \ \ \ \ \ \ \ \ \ }}**@{},
     <13mm,5mm>*{};<14.3mm,-4.5mm>*{_{J_2}}**@{},
 \end{xy}}\Ea
\Eeq
where $\sigma(I_1\sqcup I_2)$ and $\sigma(J_1\sqcup J_2)$ are the signs of the shuffles
$[1,\ldots,m]\rar I_1\sqcup I_2$ and, respectively, $[1,\ldots,n]\rar J_1\sqcup J_2$. 

Let
$V$ be a dg vector space. According to the general theory \cite{MV}, there is a one-to-one correspondence between the set of representations, $\{
\rho: \LBm \rar \cE nd_V\}$, and the set of Maurer-Cartan elements in the dg Lie algebra
\Beq\label{2: fl_V'}
\Def(\LBm \stackrel{0}{\rar} \cE nd_V)\simeq \prod_{m,n\geq 1}
\wedge^mV^*\ot \wedge^n V[2-m-n]= \prod_{m,n\geq 1} \odot^m(V^*[-1])\ot \odot^n(V[-1]) [2] =: \fg_V
\Eeq
controlling deformations of the zero map $\LB_\infty \stackrel{0}{\rar} \cE nd_V$.
The differential in $\fg_V$ is induced by the differential in $V$ while the Lie bracket
can be described explicitly as follows. First one notices that the completed graded vector space
$$
\fg_V[-2]= \prod_{m,n\geq 1} \odot^m(V^*[-1])\ot \odot^n(V[-1])=\widehat{\odot^{\bu \geq 1}}\left( V^*[-1])\oplus V[-1]\right)
$$
is naturally a 3-algebra with  degree $-2$ Lie brackets, $\{\ ,\ \}$,
given on generators by
\[
\{sv, sw\}=0,\ \ \{s\al, s\be\}=0, \ \ \{s\al, sv\}=<\al,v>, \ \ \forall v,w\in V, \al,\be\in V^*.
\]
where $s: V\rar V[-1]$ and $s: V^*\rar V^*[-1]$ are natural isomorphisms.
Maurer-Cartan elements in $\fg_V$, that is degree 3 elements $\nu$ satisfying the equation
$$
\{\nu,\nu\}=0,
$$
 are in 1-1 correspondence with
representations $\nu: \LB_\infty \rar \cE nd_V$. Such elements satisfying the condition
$$
\nu \in \odot^2(V^*[-1)\ot V[-1] \ \oplus V^*[-1]\ot \odot^2(V[-1])
$$
are precisely Lie bialgebra structures in $V$.

\sip

The properads $\LB$ and $\LBm$ admit filtrations by the number of vertices and we denote by
$\wLB$ and $\wLBm$ their completions with respect to these filtrations.

\subsection{Formality maps as morphisms of props} We introduced in \cite{MW2} an endofunctor $\caD$ in the category of augmented props with the property
that for any representation of a prop $\cP$ in a vector space $V$ the associated prop $\caD\cP$
admits an induced representation on the graded commutative algebra $\odot^\bu V$ given in terms of polydifferential operators. More, we proved that

\Bi
 \item[(i)]  For any choice of a Drinfeld associator $\fA$ there is an associated highly non-trivial (in the sense that it is is non-zero on every generator of $\Assb_\infty$, see formula (\ref{5: Boundary cond for formality map}) below)  morphism of dg props,
\Beq\label{1: formality map F_A}
F_\fA: \Assb_\infty \lon \caD\wLBm.
\Eeq
where $\Assb_\infty$ stands for a minimal resolution of the prop of associative bialgebras, and the construction of the polydifferential prop $\caD\wLBm$ out of $\wLBm$ is explained below.

\item[(ii)] For any  graded vector space $V$, each morphism
$F_\fA$ induces a $\caL ie_\infty$ quasi-isomorphism (called a {\em formality map})
between the dg $\caL ie_\infty$ algebra
$$
C_{GS}^\bu(\f_V,\f_V)=
\Def(\cA ss\cB\stackrel{\rho_0}{\lon} \cE nd_{\f_V})
$$
controlling deformations of the standard graded commutative an co-commutative bialgebra structure $\rho_0$ in $\f_V$, and the Lie algebra
$$
\fg_V=\Def(\LB\stackrel{0}{\lon} \cE nd_V)
$$
controlling deformations of the zero morphism $0: \LB\rar \cE nd_V$.

\item[(iii)]  For any formality morphism $F_\fA$ there is a canonical morphism of complexes
 $$
 \fGC_3^{or} \lon \Def\left(\Assb_\infty \stackrel{F_\fA}{\lon} \caD \wLBm\right)
 $$
 which is a quasi-isomorphism up to one class corresponding to the standard rescaling automorphism
 of the prop of Lie bialgebras $\LB$. 

\item[(iv)]
 The set of homotopy classes of universal formality maps as in (\ref{1: formality map F_A}) can be identified with the set of Drinfeld associators. In particular,
the Grothendieck-Teichm\"uller group $GRT=GRT_1\rtimes \K^*$ acts faithfully and transitively on such
universal formality maps.
\Ei

In the proof of item (i) in \cite{MW2} we used the Etingof-Kazhdan theorem \cite{EK} which says that any Lie bialgebra
can deformation quantized in the sense explained by Drinfeld in \cite{D}, and which can be reformulated in our language
as a morphism of props
$$
f_\fA: \Assb \lon \caD \wLB
$$
satisfying certain non-triviality condition (see below). This morphism gives us universal quantizations of arbitrary, possibly infinitely dimensional,  Lie bialgebras. If one is interested in universal quantization of {\em finite-dimensional}\, Lie bialgebras only, then the above morphism should be replaced by a map
$$
f^\circlearrowright:\Assb \lon \caD \wLB^\circlearrowright
$$
to the polydifferential extension of the {\em wheeled}\, closure $\wLB^\circlearrowright$ (see \cite{MMS}) of
the prop $\wLB$. The morphism $f_\fA$ implies the morphism $f^\circlearrowright$ due to the canonical injection $\caD \wLB \rar \caD \wLB^\circlearrowright$, but not vice versa.
In this paper we show a new proof of the Etingof-Kazhdan theorem for finite-dimensional Lie bialgebras by giving an explicit formula for the morphism $f^\circlearrowright$ above. We also show that the morphism
$f^\circlearrowright$ can be lifted by a trivial induction to a morphism
of dg props
\Beq \label{5: F from Assb-infty to FLieb-wheeld}
F^\circlearrowright: \Assb_\infty \lon \caD \wLB^{\mathrm{min},\circlearrowright}_\infty
\Eeq
satisfying  the  conditions
 \Beq\label{5: Boundary cond for formality map}
\pi_1\circ F^\circlearrowright\left(\Ba{c}\resizebox{13mm}{!}{ \xy
 (0,7)*{\overbrace{\ \ \ \  \ \ \ \ \ \ \ \ \ \ }},
 (0,9)*{^m},
 (0,3)*{^{...}},
 (0,-3)*{_{...}},
 (0,-7)*{\underbrace{  \ \ \ \  \ \ \ \ \ \ \ \ \ \ }},
 (0,-9)*{_n},
 (0,0)*{\circ}="0",
(-7,5)*{}="u_1",
(-4,5)*{}="u_2",
(4,5)*{}="u_3",
(7,5)*{}="u_4",
(-7,-5)*{}="d_1",
(-4,-5)*{}="d_2",
(4,-5)*{}="d_3",
(7,-5)*{}="d_4",
\ar @{.} "0";"u_1" <0pt>
\ar @{.} "0";"u_2" <0pt>
\ar @{.} "0";"u_3" <0pt>
\ar @{.} "0";"u_4" <0pt>
\ar @{.} "0";"d_1" <0pt>
\ar @{.} "0";"d_2" <0pt>
\ar @{.} "0";"d_3" <0pt>
\ar @{.} "0";"d_4" <0pt>
\endxy}\Ea\right)=
\la \Ba{c}\resizebox{16mm}{!}{\xy
(0,7.5)*{\overbrace{\ \ \ \ \ \ \ \ \  \ \ \ \ \ \ \ \ \ \ }},
 (0,9.5)*{^m},
 (0,-7.5)*{\underbrace{\ \ \ \ \ \ \ \ \ \ \ \  \ \ \ \ \ \ }},
 (0,-9.9)*{_n},
  (-6,5)*{...},
  (-6,-5)*{...},
 (-3,5)*{\circ}="u1",
  (-3,-5)*{\circ}="d1",
  (-6,5)*{...},
  (-6,-5)*{...},
  (-9,5)*{\circ}="u2",
  (-9,-5)*{\circ}="d2",
 (3,5)*{\circ}="u3",
  (3,-5)*{\circ}="d3",
  (6,5)*{...},
  (6,-5)*{...},
  (9,5)*{\circ}="u4",
  (9,-5)*{\circ}="d4",
 (0,0)*{\bullet}="a",
\ar @{-} "d1";"a" <0pt>
\ar @{-} "a";"u1" <0pt>
\ar @{-} "d2";"a" <0pt>
\ar @{-} "a";"u2" <0pt>
\ar @{-} "d3";"a" <0pt>
\ar @{-} "a";"u3" <0pt>
\ar @{-} "d4";"a" <0pt>
\ar @{-} "a";"u4" <0pt>
\endxy}\Ea\ \ \text{for some non-zero} \la\in \R,
\Eeq
for all $m+n\geq 3$, $m,n\geq 1$; here $\pi_1$ is the projection to the
 vector subspace in $\wLB^{\mathrm{min},\circlearrowright}_\infty$ spanned by graphs with precisely one black vertex.
Moreover, we conjecture an explicit formula for such an extension $F^\circlearrowright$.

\sip

Morphisms of dg props ({\ref{5: F from Assb-infty to FLieb-wheeld}) satisfying the condition (\ref{5: Boundary cond for formality map}) can be called {\em formality morphisms in finite dimensions}\, as every such a morphism gives rise to a quasi-isomorphism of $\caL ie_\infty$-algebras introduced in the item (ii) above, but only  for {\em finite-dimensional}\, graded vector spaces $V$ (cf.
 \cite{MW2}).

\subsection{Polydifferential functor}\label{5: functor D}
We refer to \cite{MW2} for a detailed definition of the endofunctor $\caD$. In this paper we apply this functor to the props $\wLB$ and $\wLBm$, their wheeled closures $\wLB^\circlearrowright$, $\wLB_\infty^{\mathrm{min}, \circlearrowright}$, and their quantized versions
$\wLB^{\mathrm{quant}}$ and $\wLB_\infty^{\mathrm{quant}}$. It is enough to explain the action of $\caD$ on the prop $\wLBm$, the other cases being completely analogous.

\sip

Roughly speaking, $\caD\wLBm$ is spanned as a vector space by graphs from $\wLBm$
whose input and output legs are labeled by {\em not necessarily different integers}; input legs labelled by the same integer $i$ we show as attached to a new white {\em in-vertex}\, to which we assign label $i$;
the same procedure applies to output legs giving us new white {\em out-vertices}. Moreover, we allow these new white in-vertices and out-vertices with no legs attached. For example,
$$
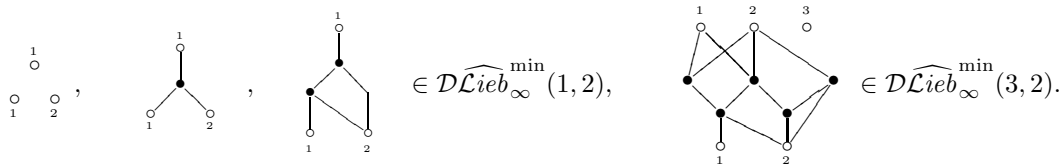
\Ba{c}\resizebox{8mm}{!}{ \xy
(-3,0)*{_{_1}},
(3,0)*{_{_2}},
(0,8)*{^{^1}},
 (0,7)*{\circ}="a",
(-3,2)*{\circ}="b_1",
(3,2)*{\circ}="b_2",
 \endxy}\Ea, \ \ \  \Ba{c}\resizebox{12mm}{!}{
\xy
(-5,0)*{_{_1}},
(5,0)*{_{_2}},
(0,14)*{^{^1}},
(0,13)*{\circ}="0",
 (0,7)*{\bu}="a",
(-5,2)*{\circ}="b_1",
(5,2)*{\circ}="b_2",
(-8,-2)*{}="c_1",
(-2,-2)*{}="c_2",
(2,-2)*{}="c_3",
\ar @{-} "a";"0" <0pt>
\ar @{-} "a";"b_1" <0pt>
\ar @{-} "a";"b_2" <0pt>
\endxy}
\Ea\ \ , \ \ \
\Ba{c}\resizebox{10mm}{!}{ \xy
(-5,-8)*{^{_1}},
(5,-8)*{^{_2}},
(0,14)*{^{^1}},
(0,13)*{\circ}="0",
 (0,7)*{\bu}="a",
(5,-5)*{\circ}="b_1",
(-5,-5)*{\circ}="b_2",
(-5,2)*{\bu}="c",
(5,2)*{}="o",
(-5,2)*{}="c_1",
(2,2)*{}="c_2",
(-2,2)*{}="c_3",
\ar @{-} "a";"0" <0pt>
\ar @{-} "o";"b_1" <0pt>
\ar @{-} "c";"b_2" <0pt>
\ar @{-} "c";"b_1" <0pt>
\ar @{-} "a";"c" <0pt>
\ar @{-} "a";"o" <0pt>
\endxy}
\Ea
\ \ \in \caD \wLBm(1,2), \ \ \ \ \
\Ba{c}\resizebox{22mm}{!}{
\xy
(-5,-5)*{_{_1}},
(5,-5)*{_{_2}},
(-8,16.5)*{^{^1}},
(0,16.5)*{^{^2}},
(8,16.5)*{^{^3}},
(0,15)*{\circ}="u",
(-8,15)*{\circ}="uL",
(8,15)*{\circ}="uR",
%(-5,12)*{\bu}="0",
 (0,7)*{\bu}="a",
(-10,7)*{\bu}="L",
(12,7)*{\bu}="R",
(-5,2)*{\bu}="b_1",
(5,2)*{\bu}="b_2",
(-5,-3)*{\circ}="c_1",
(5,-3)*{\circ}="c_3",
\ar @{-} "a";"b_1" <0pt>
\ar @{-} "a";"b_2" <0pt>
\ar @{-} "b_1";"c_1" <0pt>
\ar @{-} "b_1";"c_3" <0pt>
\ar @{-} "b_2";"c_3" <0pt>
\ar @{-} "R";"c_3" <0pt>
\ar @{-} "b_1";"L" <0pt>
\ar @{-} "u";"L" <0pt>
\ar @{-} "b_2";"R" <0pt>
\ar @{-} "R";"u" <0pt>
\ar @{-} "a";"u" <0pt>
\ar @{-} "a";"uL" <0pt>
\ar @{-} "a";"uL" <0pt>
\ar @{-} "L";"uL" <0pt>
\endxy}\Ea\in \caD \wLBm(3,2).
$$
The linear span of graphs obtained in this way from elements of $\wLBm$ with $n$ in-vertices and $m$-out vertices is denoted by $\caD \wLBm(m,n)$; it is clearly an $\bS_m^{op}\times \bS_n$ module (with elements of the permutation groups acting by relabelling of the
in- and out vertices). The $\bS$-bimodule $\caD \wLBm(m,n)$ has a natural basis $\{\cG_{k;m,n}\}$ where $\cG_{k;m,n}$ is the set of oriented graphs with
$n$ labelled white in-vertices, $m$ labelled white out-vertices and $k$ unlabeled internal (black) vertices and with no edges connecting in-vertices directly to out-vertices.  Any graph $\Ga\in \cG_{k;m,n}$ has its set of edges $E(\Ga)$ decomposed canonically into the disjoint union
$$
E(\Ga)=E_{int}(\Ga) \coprod E_{in}(\Ga)\coprod E_{out}(\Ga)
$$
where $E_{int}(\Ga)$ is the subset of edges connecting two internal vertices,
$E_{in}(\Ga)$ is the subset of edges connecting in-vertices to internal ones,
and $E_{out}(\Ga)$ is the subset of edges connecting internal vertices to out-vertices. As a $\Z$-graded vector space $\caD \wLB_\infty(m,n)$ is defined by
$$
\caD \wLB_\infty(m,n)=\prod_{k\geq 0}\K\langle\cG_{k;m,n}^{or}\rangle
$$
where a graph $\Ga\in \cG_{k;m,n}$ is assigned the following homological degree
$$
|\Ga|=3|V_{int}(\Ga)| -2|E_{int}(\Ga)| -|E_{in}(\Ga)|-|E_{out}(\Ga)|.
$$

\sip

The horizontal composition in $\caD \wLBm$
$$
\Ba{rccc}
\boxtimes: & \caD \wLBm(m,n) \ot \caD \wLBm(m',n') &\lon & \caD \wLBm(m+m',n+n')\\
 & \Ga\ot \Ga' & \lon & \Ga\boxtimes \Ga'
\Ea
$$
is given just by taking the disjoint union of the graphs $\Ga$ and $\Ga'$ and relabelling
in- and  out-vertices of the graph $\Ga'$ accordingly. The vertical composition,
$$
\Ba{rccc}
\circ: & \caD \wLBm(m,n) \ot \caD \wLBm(n,l) &\lon & \caD \wLBm(m,l)\\
 & \Ga\ot \Ga' & \lon & \Ga\circ\Ga',
\Ea
$$
is given by the following two step procedure: (a) erase all $n$ in-vertices of $\Ga$ and all $n$ out-vertices of $\Ga'$, (b) take a sum over all possible ways of attaching the hanging out-legs of $\Ga$ to hanging in-legs of $\Ga'$ (with the same numerical label) as well as to the out-vertices of $\Ga'$, and also attaching the remaining in-legs of $\Ga'$
to in-vertices of  $\Ga$ (see \S 2.2.2 in \cite{MW2} for more details).
For example, a vertical composition of the following two graphs,
$$
\Ba{rccc}
\circ: & \caD \wLBm(2,1) \ot \caD \wLBm(1,2) &\lon &  \caD \wLBm(2,2)\\
&
\Ba{c}\resizebox{7mm}{!}{ \xy
(-5,0)*{_1},
(5,0)*{_2},
(0,13)*{\circ}="0",
 (0,7)*{\bu}="a",
(-5,2)*{\circ}="b_1",
(5,2)*{\circ}="b_2",
\ar @{-} "a";"0" <0pt>
\ar @{-} "a";"b_1" <0pt>
\ar @{-} "a";"b_2" <0pt>
\endxy}\Ea
\ot
\Ba{c}
\resizebox{7mm}{!}{\xy
(-5,15)*{_1},
(5,15)*{_2},
(0,2)*{\circ}="0",
 (0,8)*{\bu}="a",
(-5,13)*{\circ}="b_1",
(5,13)*{\circ}="b_2",
\ar @{-} "a";"0" <0pt>
\ar @{-} "a";"b_1" <0pt>
\ar @{-} "a";"b_2" <0pt>
\endxy} \Ea
 &\lon &  \Ga
\Ea
$$
is given by the following sum
$$
\Ga=
\Ba{c}\resizebox{7mm}{!}{\xy
(-5,0)*{_1},
(5,0)*{_2},
(-5,20)*{^1},
(5,20)*{^2},
%(-2,13)*{_2},
%(-2,7)*{_1},
%
(0,13)*{\bu}="0",
 (0,7)*{\bu}="a",
(-5,2)*{\circ}="b_1",
(5,2)*{\circ}="b_2",
(-5,18)*{\circ}="u_1",
(5,18)*{\circ}="u_2",
\ar @{-} "a";"0" <0pt>
\ar @{-} "a";"b_1" <0pt>
\ar @{-} "a";"b_2" <0pt>
\ar @{-} "0";"u_1" <0pt>
\ar @{-} "0";"u_2" <0pt>
\endxy}
\Ea\ \ \ +\ \ \
\Ba{c}
\resizebox{8mm}{!}{\xy
(-5,0)*{_1},
(5,0)*{_2},
(-5,20)*{^1},
(5,20)*{^2},
%(6,10)*{_2},
%(-6,9)*{_1},
%
(4,10)*{\bu}="0",
 (-4,8)*{\bu}="a",
(-5,2)*{\circ}="b_1",
(5,2)*{\circ}="b_2",
(-5,18)*{\circ}="u_1",
(5,18)*{\circ}="u_2",
\ar @{-} "b_1";"0" <0pt>
\ar @{-} "u_1";"a" <0pt>
\ar @{-} "a";"b_1" <0pt>
\ar @{-} "a";"b_2" <0pt>
\ar @{-} "0";"u_1" <0pt>
\ar @{-} "0";"u_2" <0pt>
\endxy}
\Ea\ \ \ +\ \ \
\Ba{c}
\resizebox{8mm}{!}{\xy
(-5,0)*{_1},
(5,0)*{_2},
(-5,20)*{^1},
(5,20)*{^2},
%(6,10)*{_2},
%(-6,8)*{_1},
%
(4,10)*{\bu}="0",
 (-4,8)*{\bu}="a",
(-5,2)*{\circ}="b_1",
(5,2)*{\circ}="b_2",
(-5,18)*{\circ}="u_1",
(5,18)*{\circ}="u_2",
\ar @{-} "b_2";"0" <0pt>
\ar @{-} "u_1";"a" <0pt>
\ar @{-} "a";"b_1" <0pt>
\ar @{-} "a";"b_2" <0pt>
\ar @{-} "0";"u_1" <0pt>
\ar @{-} "0";"u_2" <0pt>
\endxy}
\Ea\ \ \ +\ \ \
\Ba{c}
\resizebox{8mm}{!}{\xy
(-5,0)*{_1},
(5,0)*{_2},
(-5,20)*{^1},
(5,20)*{^2},
%(6,10)*{_2},
%(-6,8)*{_1},
%
(4,10)*{\bu}="0",
 (-4,8)*{\bu}="a",
(-5,2)*{\circ}="b_1",
(5,2)*{\circ}="b_2",
(-5,18)*{\circ}="u_1",
(5,18)*{\circ}="u_2",
\ar @{-} "b_2";"0" <0pt>
\ar @{-} "u_2";"a" <0pt>
\ar @{-} "a";"b_1" <0pt>
\ar @{-} "a";"b_2" <0pt>
\ar @{-} "0";"u_1" <0pt>
\ar @{-} "0";"u_2" <0pt>
\endxy}
\Ea\ \ \ +\ \ \
\Ba{c}
\resizebox{8mm}{!}{\xy
(-5,0)*{_1},
(5,0)*{_2},
(-5,20)*{^1},
(5,20)*{^2},
%(6,10)*{_2},
%(-6,8)*{_1},
%
(4,10)*{\bu}="0",
 (-4,8)*{\bu}="a",
(-5,2)*{\circ}="b_1",
(5,2)*{\circ}="b_2",
(-5,18)*{\circ}="u_1",
(5,18)*{\circ}="u_2",
\ar @{-} "b_1";"0" <0pt>
\ar @{-} "u_2";"a" <0pt>
\ar @{-} "a";"b_1" <0pt>
\ar @{-} "a";"b_2" <0pt>
\ar @{-} "0";"u_1" <0pt>
\ar @{-} "0";"u_2" <0pt>
\endxy}\Ea
$$
The differential $\delta$ in $\caD \LB_\infty$ acts only on black vertices and splits them as shown in (\ref{3: differential in LieBinfty}).

\mip

For any given representation $\nu: \LB_\infty^{\mathrm{min}}\rar \cE nd_V$, i.e.\ for any Maurer-Cartan
element $\nu$ in the Lie algebra $\fg_V$, there is an associated representation
$\rho_\nu: \caD \LB_\infty \rar \cE nd_{\f_V}$ in $\f_V=\odot^\bu V$ given in terms of polydifferential operators as explained in full details in \S 5.4 of \cite{MW2}. If, for example, $V=\R^n$
with the standard basis denoted by $(x_1, \ldots, x_n)$ (so that $\f_V=\K[x_1, \ldots, x_n]$),
and $\nu$ is a Lie bialgebra structure in $V$ with the structure constants for the Lie bracket and,
respectively, Lie cobracket given by
$$
[x_i,x_j]=:\sum_{k=1}^n C_{ij}^k x_k, \ \ \ \ \ \triangle(x_k)=\sum_{i,j=1}^n \Phi^{ij}_k x_i\wedge x_j
$$
then one has,
$$
\Ba{rccc}
\rho^\nu\left( \Ba{c}\resizebox{7mm}{!}{ \xy
(-5,0)*{_1},
(5,0)*{_2},
(0,13)*{\circ}="0",
 (0,7)*{\bu}="a",
(-5,2)*{\circ}="b_1",
(5,2)*{\circ}="b_2",
\ar @{-} "a";"0" <0pt>
\ar @{-} "a";"b_1" <0pt>
\ar @{-} "a";"b_2" <0pt>
\endxy}\Ea\right): & \f_V\ot \f_V & \lon \f_V\\
& f_1 \ot f_2 &  \lon & \displaystyle\sum_{i,j,k=1}^n  x_k C^k_{ij} \frac{\p f_1}{\p x_i}
\frac{\p f_2}{\p x_j}
\Ea
$$
$$
\Ba{rccc}
\rho^\nu\left( \Ba{c}
\resizebox{7mm}{!}{
\xy
(-5,0)*{_1},
(5,0)*{_2},
(-5,20)*{^1},
(5,20)*{^2},
%(-2,13)*{_2},
%(-2,7)*{_1},
%
(0,13)*{\bu}="0",
 (0,7)*{\bu}="a",
(-5,2)*{\circ}="b_1",
(5,2)*{\circ}="b_2",
(-5,18)*{\circ}="u_1",
(5,18)*{\circ}="u_2",
\ar @{-} "a";"0" <0pt>
\ar @{-} "a";"b_1" <0pt>
\ar @{-} "a";"b_2" <0pt>
\ar @{-} "0";"u_1" <0pt>
\ar @{-} "0";"u_2" <0pt>
\endxy}
\Ea
\right): & \f_V & \lon & \f_V\ot \f_V\\
& f & \lon & \displaystyle\sum_{i,j,k,m,n=1}^n (x_m\ot x_n)\cdot \Phi_k^{mn} C^k_{ij} \Delta( \frac{\p f_1}{\p x_i}
\frac{\p f_2}{\p x_j})
\Ea
$$
while
$
\rho^\nu\left( \Ba{c}\resizebox{6mm}{!}{ \xy
(-3,0)*{_1},
(3,0)*{_2},
 (0,7)*{\circ}="a",
(-3,2)*{\circ}="b_1",
(3,2)*{\circ}="b_2",
 \endxy}\Ea \right): \f_V^{\ot 2}\rar \f_V$ and  $\Delta:=\rho^\nu\left(   \Ba{c}\resizebox{6mm}{!}{ \xy
(-3,9)*{^1},
(3,9)*{^2},
 (0,2)*{\circ}="a",
(-3,7)*{\circ}="b_1",
(3,7)*{\circ}="b_2",
 \endxy}\Ea  \right): \f_V\rar \f_V^{\ot 2}$
are the standard commutative multiplication and, respectively,  co-commutative   comultiplication
 in $\f_V$. Representations of the completed props $\wLB$ and $\wLBm$ (and hence of $\caD\wLB$ and
 $\caD\wLBm$) are considered
in the subsection \S {\ref{5: subsec on repr of wLBq}} below --- they require an introduction of a formal parameter to insure convergence.

\subsection{Properad of quantizable Lie bialgebras} Let us denote by $\widehat{\LB}_\infty$  the {\em non-differential}\, properad $\wLBm$, i.e.\ the completed free properad generated by the same $\bS$-bimodule $E$ but with the differential set to zero. Let $\widehat{\LB}_\infty^+$ be the free extension of $\widehat{\LB}_\infty$  by one extra $(1,1)$ generator $\Ba{c}
\xy
%%%
 (0,0)*{\bullet}="a",
(0,3)*{}="b",
(0,-3)*{}="c",
%%%
\ar @{-} "a";"b" <0pt>
\ar @{-} "a";"c" <0pt>
\endxy\Ea$ of homological degree one. In \cite{MW} (see formula (11) there) we constructed a map of Lie algebras
$$
\Ba{rccc}
f: & \dfGC_3^{or} & \lon & \Der(\widehat{\LB}_\infty^+) \\
& \Ga & \lon & \displaystyle\sum_{m,n\geq 1} \sum_{s:[n]\rar V(\Ga)\atop \hat{s}:[m]\rar V(\Ga)}  \Ba{c}\resizebox{11mm}{!}  {\xy
 (-6,7)*{^1},
(-3,7)*{^2},
(2.5,7)*{},
(7,7)*{^m},
(-3,-8)*{_2},
(3,-6)*{},
(7,-8)*{_n},
(-6,-8)*{_1},
(0,4.5)*+{...},
(0,-4.5)*+{...},
(0,0)*+{\Ga}="o",
(-6,6)*{}="1",
(-3,6)*{}="2",
(3,6)*{}="3",
(6,6)*{}="4",
(-3,-6)*{}="5",
(3,-6)*{}="6",
(6,-6)*{}="7",
(-6,-6)*{}="8",
\ar @{-} "o";"1" <0pt>
\ar @{-} "o";"2" <0pt>
\ar @{-} "o";"3" <0pt>
\ar @{-} "o";"4" <0pt>
\ar @{-} "o";"5" <0pt>
\ar @{-} "o";"6" <0pt>
\ar @{-} "o";"7" <0pt>
\ar @{-} "o";"8" <0pt>
\endxy}\Ea
\Ea
$$
 where the second sum in the r.h.s.\  is taken over all ways of attaching the incoming and outgoing legs to the graph $\Ga$, and then setting to zero every resulting graph if it contains a vertex with valency $\leq 2$ or with no input legs or no output legs
    (there is an implicit rule of signs in-built into this formula and its version in Proposition {\ref{5: Prop on isomorphisms of lieb properads}} below  which is completely analogous to the one explained in \S 7 of \cite{MaVo}).
Here $\Der(\widehat{\LB}_\infty^+)$ is the Lie algebra of continuous derivations of the topological properad $\widehat{\LB}_\infty^+$.  Note that for many graphs $\Ga\in \dfGC_3^{or}$ the associated  $(m=1,n=1)$ summand
$
\Ba{c}\resizebox{3mm}{!}  {\xy
 (0,5)*{^1},
(0,-5)*{_1},
(0,0)*+{\Ga}="o",
(0,4)*{}="1",
(0,-4)*{}="8",
\ar @{-} "o";"1" <0pt>
\ar @{-} "o";"8" <0pt>
\endxy}\Ea
$
 in $f(\Ga)$ can be highly non-trivial, and this phenomenon explains  the need for the  extension $\widehat{\LB}_\infty\rar \widehat{\LB}_\infty^+$ above.

\sip

If $\Upsilon$ is a Maurer-Cartan element in $\dfGC_3^{or}$, then $f(\Upsilon)$ is a differential
in  $\widehat{\LB}_\infty^+$ which acts on the generating $(m,n)$-corolla as follows
$$
f(\Upsilon)\left( \Ba{c}\resizebox{14mm}{!}{\begin{xy}
 <0mm,0mm>*{\bu};<0mm,0mm>*{}**@{},
 <-0.6mm,0.44mm>*{};<-8mm,5mm>*{}**@{-},
 <-0.4mm,0.7mm>*{};<-4.5mm,5mm>*{}**@{-},
 <0mm,0mm>*{};<-1mm,5mm>*{\ldots}**@{},
 <0.4mm,0.7mm>*{};<4.5mm,5mm>*{}**@{-},
 <0.6mm,0.44mm>*{};<8mm,5mm>*{}**@{-},
   <0mm,0mm>*{};<-8.5mm,5.5mm>*{^1}**@{},
   <0mm,0mm>*{};<-5mm,5.5mm>*{^2}**@{},
   <0mm,0mm>*{};<4.5mm,5.5mm>*{^{m\hspace{-0.5mm}-\hspace{-0.5mm}1}}**@{},
   <0mm,0mm>*{};<9.0mm,5.5mm>*{^m}**@{},
 <-0.6mm,-0.44mm>*{};<-8mm,-5mm>*{}**@{-},
 <-0.4mm,-0.7mm>*{};<-4.5mm,-5mm>*{}**@{-},
 <0mm,0mm>*{};<-1mm,-5mm>*{\ldots}**@{},
 <0.4mm,-0.7mm>*{};<4.5mm,-5mm>*{}**@{-},
 <0.6mm,-0.44mm>*{};<8mm,-5mm>*{}**@{-},
   <0mm,0mm>*{};<-8.5mm,-6.9mm>*{^1}**@{},
   <0mm,0mm>*{};<-5mm,-6.9mm>*{^2}**@{},
   <0mm,0mm>*{};<4.5mm,-6.9mm>*{^{n\hspace{-0.5mm}-\hspace{-0.5mm}1}}**@{},
   <0mm,0mm>*{};<9.0mm,-6.9mm>*{^n}**@{},
 \end{xy}}\Ea  \right)= \sum_{s:[n]\rar V(\Ga)\atop \hat{s}:[m]\rar V(\Ga)}  \Ba{c}\resizebox{11mm}{!}  {\xy
 (-6,7)*{^1},
(-3,7)*{^2},
(2.5,7)*{},
(7,7)*{^m},
(-3,-8)*{_2},
(3,-6)*{},
(7,-8)*{_n},
(-6,-8)*{_1},
(0,4.5)*+{...},
(0,-4.5)*+{...},
(0,0)*+{\Upsilon}="o",
(-6,6)*{}="1",
(-3,6)*{}="2",
(3,6)*{}="3",
(6,6)*{}="4",
(-3,-6)*{}="5",
(3,-6)*{}="6",
(6,-6)*{}="7",
(-6,-6)*{}="8",
\ar @{-} "o";"1" <0pt>
\ar @{-} "o";"2" <0pt>
\ar @{-} "o";"3" <0pt>
\ar @{-} "o";"4" <0pt>
\ar @{-} "o";"5" <0pt>
\ar @{-} "o";"6" <0pt>
\ar @{-} "o";"7" <0pt>
\ar @{-} "o";"8" <0pt>
\endxy}\Ea
$$

If $\Upsilon$ is such that the summand $
\Ba{c}\resizebox{3mm}{!}  {\xy
 (0,5)*{^1},
(0,-5)*{_1},
(0,0)*+{\Upsilon}="o",
(0,4)*{}="1",
(0,-4)*{}="8",
\ar @{-} "o";"1" <0pt>
\ar @{-} "o";"8" <0pt>
\endxy}\Ea
$
 in $f(\Upsilon)$ contains at least one vertex of the form
 $\Ba{c}
\xy
%%%
 (0,0)*{\bullet}="a",
(0,3)*{}="b",
(0,-3)*{}="c",
%%%
\ar @{-} "a";"b" <0pt>
\ar @{-} "a";"c" <0pt>
\endxy\Ea$, then the ideal $I^+\subset \widehat{\LB}_\infty^+$ generated by this extra generator
 $\Ba{c}
\xy
%%%
 (0,0)*{\bullet}="a",
(0,3)*{}="b",
(0,-3)*{}="c",
%%%
\ar @{-} "a";"b" <0pt>
\ar @{-} "a";"c" <0pt>
\endxy\Ea$ is respected by the differential $f(\Upsilon)$ so that the latter induces a differential
in the quotient properad
$$
 \widehat{\LB}_\infty=  \widehat{\LB}_\infty^+/I^+.
$$
For example, the standard Maurer-Cartan element
$$
\Upsilon_S:=\xy
 (0,0)*{\bullet}="a",
(5,0)*{\bu}="b",
\ar @{->} "a";"b" <0pt>
\endxy
$$
in $\dfGC_3^{or}$ does have this property as
$$
\Ba{c}\resizebox{7mm}{!}  {\xy
 (0,5)*{^1},
(0,-5)*{_1},
(0,0)*+{\ \Upsilon_S}="o",
(0,4)*{}="1",
(0,-4)*{}="8",
\ar @{-} "o";"1" <0pt>
\ar @{-} "o";"8" <0pt>
\endxy}\Ea = \begin{xy}
<0mm,7mm>*{^1};
 <0mm,-4.4mm>*{_1};
 <0mm,0mm>*{};<0mm,-3mm>*{}**@{-},
 <0mm,0mm>*{};<0mm,6mm>*{}**@{-},
 <0mm,0mm>*{\bullet};
 <0mm,3mm>*{\bullet};
 \end{xy}
$$
and induces  the standard differential $\delta$
in $\widehat{\LB}_\infty$ given by the formula (\ref{3: differential in LieBinfty}).

\sip

Another interesting for us Maurer-Cartan element is given explicitly by (\ref{3: Upsilon om_g for d=3}). It was proven in Lemma {\ref{A: lemma on 4 binary vertices}} that every graph $\Ga\in \hat{\sG}_{4p+2,6p+1}$, $p\geq 2$, contributing into $\Upsilon^{\om_g}$ has at least 4 binary vertices so that again
$$
\Ba{c}\resizebox{9mm}{!}  {\xy
 (0,5)*{^1},
(0,-5)*{_1},
(0,0)*+{\ \ \Upsilon^{\om_g}}="o",
(0,4)*{}="1",
(0,-4)*{}="8",
\ar @{-} "o";"1" <0pt>
\ar @{-} "o";"8" <0pt>
\endxy}\Ea = \begin{xy}
<0mm,7mm>*{^1};
 <0mm,-4.4mm>*{_1};
 <0mm,0mm>*{};<0mm,-3mm>*{}**@{-},
 <0mm,0mm>*{};<0mm,6mm>*{}**@{-},
 <0mm,0mm>*{\bullet};
 <0mm,3mm>*{\bullet};
 \end{xy}
$$
implying that $\Upsilon^{\om_g}$ induces the following differential in $\widehat{\LB}_\infty$
$$
\delta^{\om_g}=f(\Upsilon^{\om_g})\bmod I^+ = \delta + \sum_{p\geq 2}  \sum_{\Ga\in \hat{\sG}^{or}_{4p+2,6p+1}}\displaystyle\sum_{m,n\geq 1\atop m+n\geq 4}  \left(\int_{\overline{C}_{4p+2}(\R^3)}\bigwedge_{e\in E(\Ga)}\hspace{-2mm}
{\pi}^*_e\left(\om_g\right)\right) \sum_{s:[n]\rar V(\Ga)\atop \hat{s}:[m]\rar V(\Ga)}  \Ba{c}\resizebox{9mm}{!}  {\xy
 (-6,7)*{^1},
(-3,7)*{^2},
(2.5,7)*{},
(7,7)*{^m},
(-3,-8)*{_2},
(3,-6)*{},
(7,-8)*{_n},
(-6,-8)*{_1},
(0,4.5)*+{...},
(0,-4.5)*+{...},
(0,0)*+{\Ga}="o",
(-6,6)*{}="1",
(-3,6)*{}="2",
(3,6)*{}="3",
(6,6)*{}="4",
(-3,-6)*{}="5",
(3,-6)*{}="6",
(6,-6)*{}="7",
(-6,-6)*{}="8",
\ar @{-} "o";"1" <0pt>
\ar @{-} "o";"2" <0pt>
\ar @{-} "o";"3" <0pt>
\ar @{-} "o";"4" <0pt>
\ar @{-} "o";"5" <0pt>
\ar @{-} "o";"6" <0pt>
\ar @{-} "o";"7" <0pt>
\ar @{-} "o";"8" <0pt>
\endxy}\Ea.
$$
As every graph in the sum over $p\geq 2$ has at least  $4$ bivalent vertices (see Appendix A), we have, in particular,
$$
\delta^{\om_g}\left(\Ba{c}
\begin{xy}
 <0mm,-0.55mm>*{};<0mm,-2.5mm>*{}**@{-},
 <0.5mm,0.5mm>*{};<2.2mm,2.2mm>*{}**@{-},
 <-0.48mm,0.48mm>*{};<-2.2mm,2.2mm>*{}**@{-},
 <0mm,0mm>*{\bu};<0mm,0mm>*{}**@{},
 <0mm,-0.55mm>*{};<0mm,-3.8mm>*{_1}**@{},
 <0.5mm,0.5mm>*{};<2.7mm,2.8mm>*{^2}**@{},
 <-0.48mm,0.48mm>*{};<-2.7mm,2.8mm>*{^1}**@{},
 \end{xy}
 \Ea\right)=\delta\left(\Ba{c}
\begin{xy}
 <0mm,-0.55mm>*{};<0mm,-2.5mm>*{}**@{-},
 <0.5mm,0.5mm>*{};<2.2mm,2.2mm>*{}**@{-},
 <-0.48mm,0.48mm>*{};<-2.2mm,2.2mm>*{}**@{-},
 <0mm,0mm>*{\bu};<0mm,0mm>*{}**@{},
 <0mm,-0.55mm>*{};<0mm,-3.8mm>*{_1}**@{},
 <0.5mm,0.5mm>*{};<2.7mm,2.8mm>*{^2}**@{},
 <-0.48mm,0.48mm>*{};<-2.7mm,2.8mm>*{^1}**@{},
 \end{xy}
 \Ea\right)=0
\ \ \ \ , \ \ \ \
 \delta^{\om_g}  \left(\begin{xy}
 <0mm,0.66mm>*{};<0mm,3mm>*{}**@{-},
 <0.39mm,-0.39mm>*{};<2.2mm,-2.2mm>*{}**@{-},
 <-0.35mm,-0.35mm>*{};<-2.2mm,-2.2mm>*{}**@{-},
 <0mm,0mm>*{\bu};<0mm,0mm>*{}**@{},
   <0mm,0.66mm>*{};<0mm,3.4mm>*{^1}**@{},
   <0.39mm,-0.39mm>*{};<2.9mm,-4mm>*{^2}**@{},
   <-0.35mm,-0.35mm>*{};<-2.8mm,-4mm>*{^1}**@{},
\end{xy}\right)= \delta \left(\begin{xy}
 <0mm,0.66mm>*{};<0mm,3mm>*{}**@{-},
 <0.39mm,-0.39mm>*{};<2.2mm,-2.2mm>*{}**@{-},
 <-0.35mm,-0.35mm>*{};<-2.2mm,-2.2mm>*{}**@{-},
 <0mm,0mm>*{\bu};<0mm,0mm>*{}**@{},
   <0mm,0.66mm>*{};<0mm,3.4mm>*{^1}**@{},
   <0.39mm,-0.39mm>*{};<2.9mm,-4mm>*{^2}**@{},
   <-0.35mm,-0.35mm>*{};<-2.8mm,-4mm>*{^1}**@{},
\end{xy}\right)=0
$$

The first differential $\delta$  makes $\widehat{\LB}_\infty$ into the standard minimal resolution of the completed properad $\widehat{\LB}$
of Lie bialgebras. The second differential $\delta^{\om_g}$ makes $\widehat{\LB}_\infty$  into a resolution  of a properad $\widehat{\LB}^{\vspace{-1mm}\mathrm{quant}}$ which we call the {\em properad of quantizable Lie bialgebras}\, and which can be defined as follows.

\sip

By contrast to $\widehat{\LB}_\infty^{\mathrm{min}}:= (\widehat{\LB}_\infty, \delta)$, let us abbreviate from now on
$$
\widehat{\LB}_\infty^{\mathrm{quant}}:= (\widehat{\LB}_\infty, \delta^{\om_g})
$$
Let $J$ be the differential closure of the ideal in $\widehat{\LB}_\infty^{\mathrm{quant}}$ generated
by $(m,n)$-corollas with $m+n\geq 4$. The quotient
$$
 \widehat{\LB}^{\mathrm{quant}}:= \widehat{\LB}_\infty^{\mathrm{quant}}/J
$$
is a properad which is concentrated in homological degree zero, and which is generated by the
$\bS$-bimodule (\ref{5: module E_0 generating Lieb}) modulo the following three relations,
$$
0= \Ba{c}\resizebox{6mm}{!}{\begin{xy}
 <0mm,2.47mm>*{};<0mm,0.12mm>*{}**@{-},
 <0.5mm,3.5mm>*{};<2.2mm,5.2mm>*{}**@{-},
 <-0.48mm,3.48mm>*{};<-2.2mm,5.2mm>*{}**@{-},
 <0mm,3mm>*{\bu};<0mm,3mm>*{}**@{},
  <0mm,-0.8mm>*{\bu};<0mm,-0.8mm>*{}**@{},
<-0.39mm,-1.2mm>*{};<-2.2mm,-3.5mm>*{}**@{-},
 <0.39mm,-1.2mm>*{};<2.2mm,-3.5mm>*{}**@{-},
     <0.5mm,3.5mm>*{};<2.8mm,5.7mm>*{^2}**@{},
     <-0.48mm,3.48mm>*{};<-2.8mm,5.7mm>*{^1}**@{},
   <0mm,-0.8mm>*{};<-2.7mm,-5.2mm>*{^1}**@{},
   <0mm,-0.8mm>*{};<2.7mm,-5.2mm>*{^2}**@{},
\end{xy}}\Ea
 -
\Ba{c}\resizebox{8mm}{!}{\begin{xy}
 <0mm,-1.3mm>*{};<0mm,-3.5mm>*{}**@{-},
 <0.38mm,-0.2mm>*{};<2.0mm,2.0mm>*{}**@{-},
 <-0.38mm,-0.2mm>*{};<-2.2mm,2.2mm>*{}**@{-},
<0mm,-0.8mm>*{\bu};<0mm,0.8mm>*{}**@{},
 <2.4mm,2.4mm>*{\bu};<2.4mm,2.4mm>*{}**@{},
 <2.77mm,2.0mm>*{};<4.4mm,-0.8mm>*{}**@{-},
 <2.4mm,3mm>*{};<2.4mm,5.2mm>*{}**@{-},
     <0mm,-1.3mm>*{};<0mm,-5.3mm>*{^1}**@{},
     <2.5mm,2.3mm>*{};<5.1mm,-2.6mm>*{^2}**@{},
    <2.4mm,2.5mm>*{};<2.4mm,5.7mm>*{^2}**@{},
    <-0.38mm,-0.2mm>*{};<-2.8mm,2.5mm>*{^1}**@{},
    \end{xy}}\Ea
  +
\Ba{c}\resizebox{8mm}{!}{\begin{xy}
 <0mm,-1.3mm>*{};<0mm,-3.5mm>*{}**@{-},
 <0.38mm,-0.2mm>*{};<2.0mm,2.0mm>*{}**@{-},
 <-0.38mm,-0.2mm>*{};<-2.2mm,2.2mm>*{}**@{-},
<0mm,-0.8mm>*{\bu};<0mm,0.8mm>*{}**@{},
 <2.4mm,2.4mm>*{\bu};<2.4mm,2.4mm>*{}**@{},
 <2.77mm,2.0mm>*{};<4.4mm,-0.8mm>*{}**@{-},
 <2.4mm,3mm>*{};<2.4mm,5.2mm>*{}**@{-},
     <0mm,-1.3mm>*{};<0mm,-5.3mm>*{^2}**@{},
     <2.5mm,2.3mm>*{};<5.1mm,-2.6mm>*{^1}**@{},
    <2.4mm,2.5mm>*{};<2.4mm,5.7mm>*{^2}**@{},
    <-0.38mm,-0.2mm>*{};<-2.8mm,2.5mm>*{^1}**@{},
    \end{xy}}\Ea
 -
\Ba{c}\resizebox{8mm}{!}{\begin{xy}
 <0mm,-1.3mm>*{};<0mm,-3.5mm>*{}**@{-},
 <0.38mm,-0.2mm>*{};<2.0mm,2.0mm>*{}**@{-},
 <-0.38mm,-0.2mm>*{};<-2.2mm,2.2mm>*{}**@{-},
<0mm,-0.8mm>*{\bu};<0mm,0.8mm>*{}**@{},
 <2.4mm,2.4mm>*{\bu};<2.4mm,2.4mm>*{}**@{},
 <2.77mm,2.0mm>*{};<4.4mm,-0.8mm>*{}**@{-},
 <2.4mm,3mm>*{};<2.4mm,5.2mm>*{}**@{-},
     <0mm,-1.3mm>*{};<0mm,-5.3mm>*{^2}**@{},
     <2.5mm,2.3mm>*{};<5.1mm,-2.6mm>*{^1}**@{},
    <2.4mm,2.5mm>*{};<2.4mm,5.7mm>*{^1}**@{},
    <-0.38mm,-0.2mm>*{};<-2.8mm,2.5mm>*{^2}**@{},
    \end{xy}}\Ea
+
\Ba{c}\resizebox{8mm}{!}{\begin{xy}
 <0mm,-1.3mm>*{};<0mm,-3.5mm>*{}**@{-},
 <0.38mm,-0.2mm>*{};<2.0mm,2.0mm>*{}**@{-},
 <-0.38mm,-0.2mm>*{};<-2.2mm,2.2mm>*{}**@{-},
<0mm,-0.8mm>*{\bu};<0mm,0.8mm>*{}**@{},
 <2.4mm,2.4mm>*{\bu};<2.4mm,2.4mm>*{}**@{},
 <2.77mm,2.0mm>*{};<4.4mm,-0.8mm>*{}**@{-},
 <2.4mm,3mm>*{};<2.4mm,5.2mm>*{}**@{-},
     <0mm,-1.3mm>*{};<0mm,-5.3mm>*{^1}**@{},
     <2.5mm,2.3mm>*{};<5.1mm,-2.6mm>*{^2}**@{},
    <2.4mm,2.5mm>*{};<2.4mm,5.7mm>*{^1}**@{},
    <-0.38mm,-0.2mm>*{};<-2.8mm,2.5mm>*{^2}**@{},
    \end{xy}}\Ea
    + \ \sum_{p\geq 2}  \sum_{\Ga\in \hat{\sG}^{\leq 3}_{4p+2,6p+1}}\displaystyle \hspace{-2mm} \left(\int_{\overline{C}_{4p+2}(\R^3)} \bigwedge_{e\in E(\Ga)}\hspace{-2mm}
{\pi}^*_e\left(\om_g\right)\right)\hspace{-1mm} \sum_{s:[2]\rar V(\Ga)\atop \hat{s}:[2]\rar V(\Ga)}  \Ba{c}\resizebox{10mm}{!}  {\xy
 (-6,7)*{^1},
%(-3,7)*{^2},
(2.5,7)*{},
(7,7)*{^2},
%(-3,-8)*{_2},
(3,-6)*{},
(7,-8)*{_2},
(-6,-8)*{_1},
%
%(0,4.5)*+{...},
%(0,-4.5)*+{...},
%
(0,0)*+{\Ga}="o",
(-5,5)*{}="1",
%(-3,6)*{}="2",
%(3,6)*{}="3",
(5,5)*{}="4",
%(-3,-6)*{}="5",
%(3,-6)*{}="6",
(5,-5)*{}="7",
(-5,-5)*{}="8",
\ar @{-} "o";"1" <0pt>
%\ar @{-} "o";"2" <0pt>
%\ar @{-} "o";"3" <0pt>
\ar @{-} "o";"4" <0pt>
%\ar @{-} "o";"5" <0pt>
%\ar @{-} "o";"6" <0pt>
\ar @{-} "o";"7" <0pt>
\ar @{-} "o";"8" <0pt>
\endxy}\Ea\ ,
$$

$$
0=\Ba{c}\resizebox{9mm}{!}{\begin{xy}
 <0mm,0mm>*{\bu};<0mm,0mm>*{}**@{},
 <0mm,0.69mm>*{};<0mm,3.0mm>*{}**@{-},
 <0.39mm,-0.39mm>*{};<2.4mm,-2.4mm>*{}**@{-},
 <-0.35mm,-0.35mm>*{};<-1.9mm,-1.9mm>*{}**@{-},
 <-2.4mm,-2.4mm>*{\bu};<-2.4mm,-2.4mm>*{}**@{},
 <-2.0mm,-2.8mm>*{};<0mm,-4.9mm>*{}**@{-},
 <-2.8mm,-2.9mm>*{};<-4.7mm,-4.9mm>*{}**@{-},
    <0.39mm,-0.39mm>*{};<3.3mm,-4.0mm>*{^3}**@{},
    <-2.0mm,-2.8mm>*{};<0.5mm,-6.7mm>*{^2}**@{},
    <-2.8mm,-2.9mm>*{};<-5.2mm,-6.7mm>*{^1}**@{},
 \end{xy}}\Ea
\ + \
\Ba{c}\resizebox{9mm}{!}{ \begin{xy}
 <0mm,0mm>*{\bu};<0mm,0mm>*{}**@{},
 <0mm,0.69mm>*{};<0mm,3.0mm>*{}**@{-},
 <0.39mm,-0.39mm>*{};<2.4mm,-2.4mm>*{}**@{-},
 <-0.35mm,-0.35mm>*{};<-1.9mm,-1.9mm>*{}**@{-},
 <-2.4mm,-2.4mm>*{\bu};<-2.4mm,-2.4mm>*{}**@{},
 <-2.0mm,-2.8mm>*{};<0mm,-4.9mm>*{}**@{-},
 <-2.8mm,-2.9mm>*{};<-4.7mm,-4.9mm>*{}**@{-},
    <0.39mm,-0.39mm>*{};<3.3mm,-4.0mm>*{^2}**@{},
    <-2.0mm,-2.8mm>*{};<0.5mm,-6.7mm>*{^1}**@{},
    <-2.8mm,-2.9mm>*{};<-5.2mm,-6.7mm>*{^3}**@{},
 \end{xy}}\Ea
\ + \
\Ba{c}\resizebox{9mm}{!}{ \begin{xy}
 <0mm,0mm>*{\bu};<0mm,0mm>*{}**@{},
 <0mm,0.69mm>*{};<0mm,3.0mm>*{}**@{-},
 <0.39mm,-0.39mm>*{};<2.4mm,-2.4mm>*{}**@{-},
 <-0.35mm,-0.35mm>*{};<-1.9mm,-1.9mm>*{}**@{-},
 <-2.4mm,-2.4mm>*{\bu};<-2.4mm,-2.4mm>*{}**@{},
 <-2.0mm,-2.8mm>*{};<0mm,-4.9mm>*{}**@{-},
 <-2.8mm,-2.9mm>*{};<-4.7mm,-4.9mm>*{}**@{-},
    <0.39mm,-0.39mm>*{};<3.3mm,-4.0mm>*{^1}**@{},
    <-2.0mm,-2.8mm>*{};<0.5mm,-6.7mm>*{^3}**@{},
    <-2.8mm,-2.9mm>*{};<-5.2mm,-6.7mm>*{^2}**@{},
 \end{xy}}\Ea \
    + \ \sum_{p\geq 2}  \sum_{\Ga\in \hat{\sG}^{\leq 3}_{4p+2,6p+1}}\displaystyle  \left(\int_{\overline{C}_{4p+2}(\R^3)} \bigwedge_{e\in E(\Ga)}\hspace{-2mm}
{\pi}^*_e\left(\om_g\right)\right) \sum_{s:[2]\rar V(\Ga)\atop \hat{s}:[2]\rar V(\Ga)}  \Ba{c}\resizebox{10mm}{!}  {\xy
 (0,7)*{^1},
%(-3,7)*{^2},
(2.5,7)*{},
(0,-8)*{_2},
%(-3,-8)*{_2},
(3,-6)*{},
(7,-8)*{_3},
(-6,-8)*{_1},
%
%(0,4.5)*+{...},
%(0,-4.5)*+{...},
%
(0,0)*+{\Ga}="o",
(0,5)*{}="1",
%(-3,6)*{}="2",
%(3,6)*{}="3",
(0,-5)*{}="4",
%(-3,-6)*{}="5",
%(3,-6)*{}="6",
(5,-5)*{}="7",
(-5,-5)*{}="8",
\ar @{-} "o";"1" <0pt>
%\ar @{-} "o";"2" <0pt>
%\ar @{-} "o";"3" <0pt>
\ar @{-} "o";"4" <0pt>
%\ar @{-} "o";"5" <0pt>
%\ar @{-} "o";"6" <0pt>
\ar @{-} "o";"7" <0pt>
\ar @{-} "o";"8" <0pt>
\endxy}\Ea\ ,
$$

$$
0= \Ba{c}\resizebox{9mm}{!}{\begin{xy}
 <0mm,0mm>*{\bu};<0mm,0mm>*{}**@{},
 <0mm,-0.49mm>*{};<0mm,-3.0mm>*{}**@{-},
 <0.49mm,0.49mm>*{};<1.9mm,1.9mm>*{}**@{-},
 <-0.5mm,0.5mm>*{};<-1.9mm,1.9mm>*{}**@{-},
 <-2.3mm,2.3mm>*{\bu};<-2.3mm,2.3mm>*{}**@{},
 <-1.8mm,2.8mm>*{};<0mm,4.9mm>*{}**@{-},
 <-2.8mm,2.9mm>*{};<-4.6mm,4.9mm>*{}**@{-},
   <0.49mm,0.49mm>*{};<2.7mm,2.3mm>*{^3}**@{},
   <-1.8mm,2.8mm>*{};<0.4mm,5.3mm>*{^2}**@{},
   <-2.8mm,2.9mm>*{};<-5.1mm,5.3mm>*{^1}**@{},
 \end{xy}}\Ea
\ + \
\Ba{c}\resizebox{9mm}{!}{\begin{xy}
 <0mm,0mm>*{\bu};<0mm,0mm>*{}**@{},
 <0mm,-0.49mm>*{};<0mm,-3.0mm>*{}**@{-},
 <0.49mm,0.49mm>*{};<1.9mm,1.9mm>*{}**@{-},
 <-0.5mm,0.5mm>*{};<-1.9mm,1.9mm>*{}**@{-},
 <-2.3mm,2.3mm>*{\bu};<-2.3mm,2.3mm>*{}**@{},
 <-1.8mm,2.8mm>*{};<0mm,4.9mm>*{}**@{-},
 <-2.8mm,2.9mm>*{};<-4.6mm,4.9mm>*{}**@{-},
   <0.49mm,0.49mm>*{};<2.7mm,2.3mm>*{^2}**@{},
   <-1.8mm,2.8mm>*{};<0.4mm,5.3mm>*{^1}**@{},
   <-2.8mm,2.9mm>*{};<-5.1mm,5.3mm>*{^3}**@{},
 \end{xy}}\Ea
\ + \
\Ba{c}\resizebox{9mm}{!}{\begin{xy}
 <0mm,0mm>*{\bu};<0mm,0mm>*{}**@{},
 <0mm,-0.49mm>*{};<0mm,-3.0mm>*{}**@{-},
 <0.49mm,0.49mm>*{};<1.9mm,1.9mm>*{}**@{-},
 <-0.5mm,0.5mm>*{};<-1.9mm,1.9mm>*{}**@{-},
 <-2.3mm,2.3mm>*{\bu};<-2.3mm,2.3mm>*{}**@{},
 <-1.8mm,2.8mm>*{};<0mm,4.9mm>*{}**@{-},
 <-2.8mm,2.9mm>*{};<-4.6mm,4.9mm>*{}**@{-},
   <0.49mm,0.49mm>*{};<2.7mm,2.3mm>*{^1}**@{},
   <-1.8mm,2.8mm>*{};<0.4mm,5.3mm>*{^3}**@{},
   <-2.8mm,2.9mm>*{};<-5.1mm,5.3mm>*{^2}**@{},
 \end{xy}}\Ea\
   + \ \sum_{p\geq 2}  \sum_{\Ga\in \hat{\sG}^{\leq 3}_{4p+2,6p+1}}\displaystyle  \left(\int_{\overline{C}_{4p+2}(\R^3)} \bigwedge_{e\in E(\Ga)}\hspace{-2mm}
{\pi}^*_e\left(\om_g\right)\right) \sum_{s:[2]\rar V(\Ga)\atop \hat{s}:[2]\rar V(\Ga)}  \Ba{c}\resizebox{10mm}{!}  {\xy
 (-7,7)*{^1},
%(-3,7)*{^2},
(2.5,7)*{},
(0,7)*{^2},
%(-3,-8)*{_2},
%(3,-6)*{},
(7,7)*{_3},
(0,-8)*{_1},
%
%(0,4.5)*+{...},
%(0,-4.5)*+{...},
%
(0,0)*+{\Ga}="o",
(0,5)*{}="1",
%(-3,6)*{}="2",
%(3,6)*{}="3",
(-5,5)*{}="4",
%(-3,-6)*{}="5",
%(3,-6)*{}="6",
(5,5)*{}="7",
(0,-5)*{}="8",
\ar @{-} "o";"1" <0pt>
%\ar @{-} "o";"2" <0pt>
%\ar @{-} "o";"3" <0pt>
\ar @{-} "o";"4" <0pt>
%\ar @{-} "o";"5" <0pt>
%\ar @{-} "o";"6" <0pt>
\ar @{-} "o";"7" <0pt>
\ar @{-} "o";"8" <0pt>
\endxy}\Ea\ ,
$$
where $\hat{\sG}^{\leq 3}_{4p+2,6p+1}$ is the subset of $\hat{\sG}^{or}_{4p+2,6p+1}$ consisting
of graphs with vertices of valency $\leq 3$ (it was shown in Appendix A that such graphs have precisely $4$ bivalent vertices which explains why there no other relations than the ones shown above).

\subsubsection{\bf Theorem}\label{5: Prop cohomology of LB infty quant}
{\em The natural epimorphism of props
$$
\nu: \widehat{\LB}_\infty^{\mathrm{quant}} \lon \widehat{\LB}^{\mathrm{quant}}
$$
is a quasi-isomorphism.}

\begin{proof}  The morphism $\nu$ respects  complete and exhaustive  filtrations of both  sides by the number of vertices, hence it induces a morphism of the associated spectral sequences,
$$
\nu^r: \cE^r\widehat{\LB}_\infty^{\mathrm{quant}} \lon \cE^r\widehat{\LB}^{\mathrm{quant}}
$$
The term $\cE^0\widehat{\LB}_\infty^{\mathrm{quant}}$ has trivial differential, while
the term $\cE^1\widehat{\LB}_\infty^{\mathrm{quant}}$ can be identified with the dg prop
$\widehat{\LB}_\infty^{\mathrm{min}}$ so that $\cE^2\widehat{\LB}_\infty^{\mathrm{quant}}=\widehat{\LB}$ as an $\bS$-bimodule. On the other hand $\cE^2\wLB^{\mathrm{quant}}$ can also be identified with ${\wLB}$ as an $\bS$-bimodule. Hence the morphism $\nu^2$ is an isomorphism so that, by the Eilenberg-Moore Comparison Theorem 5.5.11 (see \S 5.5 in \cite{Wei}), the morphism $\nu$ is a quasi-isomorphism.
\end{proof}

\sip

Note that graphs in  (\ref{3: formula for Holie_2 F_k}) may contain closed paths of directed edges in general and hence belong to the graph complex $\mathsf{dfGC}_d$ rather than to  $\mathsf{dfGC}_d^{or}$. Therefore
in order to see the meaning of Theorem {{\ref{3: Theorem on KS iso for A_d^n}}} in terms of props one has
to consider
the wheeled closure \cite{MMS} of the prop $\widehat{\LB}_\infty^{\mathrm{min}}$ which we denote by $\widehat{\LB}_\infty^{\mathrm{min}, \circlearrowright}$; by definition, it is generated by the same $\bS$-bimodule (\ref{5: generators of Lieb_infty}) but now using directed graphs with possibly {\em closed}\, directed paths of internal edges.

\sip

Theorem {\ref{3: Theorem on KS iso for A_d^n}} implies almost immediately the following

\subsubsection{\bf Proposition}\label{5: Prop on isomorphisms of lieb properads} {\em There is a morphism of dg props  $\widehat{\LB}_\infty^{\mathrm{quant}}$
$$
F: \widehat{\LB}_\infty^{\mathrm{quant}} \lon \widehat{\LB}_\infty^{\mathrm{min}, \circlearrowright}
$$
given by the following transcendental formula (cf.\ (\ref{3: formula for Holie_2 F_k}))
\Beq\label{5: explicit map F from LB^q_infty to LB_infty wheeled}
F\left(
\Ba{c}\resizebox{14mm}{!}{\begin{xy}
 <0mm,0mm>*{\bu};<0mm,0mm>*{}**@{},
 <-0.6mm,0.44mm>*{};<-8mm,5mm>*{}**@{-},
 <-0.4mm,0.7mm>*{};<-4.5mm,5mm>*{}**@{-},
 <0mm,0mm>*{};<-1mm,5mm>*{\ldots}**@{},
 <0.4mm,0.7mm>*{};<4.5mm,5mm>*{}**@{-},
 <0.6mm,0.44mm>*{};<8mm,5mm>*{}**@{-},
   <0mm,0mm>*{};<-8.5mm,5.5mm>*{^1}**@{},
   <0mm,0mm>*{};<-5mm,5.5mm>*{^2}**@{},
   <0mm,0mm>*{};<4.5mm,5.5mm>*{^{m\hspace{-0.5mm}-\hspace{-0.5mm}1}}**@{},
   <0mm,0mm>*{};<9.0mm,5.5mm>*{^m}**@{},
 <-0.6mm,-0.44mm>*{};<-8mm,-5mm>*{}**@{-},
 <-0.4mm,-0.7mm>*{};<-4.5mm,-5mm>*{}**@{-},
 <0mm,0mm>*{};<-1mm,-5mm>*{\ldots}**@{},
 <0.4mm,-0.7mm>*{};<4.5mm,-5mm>*{}**@{-},
 <0.6mm,-0.44mm>*{};<8mm,-5mm>*{}**@{-},
   <0mm,0mm>*{};<-8.5mm,-6.9mm>*{^1}**@{},
   <0mm,0mm>*{};<-5mm,-6.9mm>*{^2}**@{},
   <0mm,0mm>*{};<4.5mm,-6.9mm>*{^{n\hspace{-0.5mm}-\hspace{-0.5mm}1}}**@{},
   <0mm,0mm>*{};<9.0mm,-6.9mm>*{^n}**@{},
 \end{xy}}\Ea
 \right)=
\sum_{q\geq 0} \sum_{\Ga\in \sG_{1+4q, 6q}} \left(\int_{\overline{\fC}_{1+4q}(\R^d)} \bigwedge_{e\in E(\Ga)}\hspace{-2mm}
{\pi}^*_e\left(\varpi_g\right) \right)  \sum_{s:[n]\rar V(\Ga)\atop \hat{s}:[m]\rar V(\Ga)}  \Ba{c}\resizebox{9mm}{!}  {\xy
 (-6,7)*{^1},
(-3,7)*{^2},
(2.5,7)*{},
(7,7)*{^m},
(-3,-8)*{_2},
(3,-6)*{},
(7,-8)*{_n},
(-6,-8)*{_1},
(0,4.5)*+{...},
(0,-4.5)*+{...},
(0,0)*+{\Ga}="o",
(-6,6)*{}="1",
(-3,6)*{}="2",
(3,6)*{}="3",
(6,6)*{}="4",
(-3,-6)*{}="5",
(3,-6)*{}="6",
(6,-6)*{}="7",
(-6,-6)*{}="8",
\ar @{-} "o";"1" <0pt>
\ar @{-} "o";"2" <0pt>
\ar @{-} "o";"3" <0pt>
\ar @{-} "o";"4" <0pt>
\ar @{-} "o";"5" <0pt>
\ar @{-} "o";"6" <0pt>
\ar @{-} "o";"7" <0pt>
\ar @{-} "o";"8" <0pt>
\endxy}\Ea
\Eeq
 where the third sum in the r.h.s.\  is taken over all ways of attaching the incoming and outgoing legs to the graph $\Ga$, and we set to zero every resulting graph if it contains a vertex with valency $<3$ or
   with no at least one incoming  and at least one outgoing edge.}

\subsubsection{\bf Corollary}\label{5: coroll on f from LB^q to LB^min} {\em  The explicit morphism $F$ in Proposition {\ref{5: Prop on isomorphisms of lieb properads}} induces an explicit morphism $f: \widehat{\LB}^{\mathrm{quant}} \rar \widehat{\LB}^\circlearrowright$. }

\subsubsection{\bf Representations of $\widehat{\LB}_\infty^{\mathrm{quant}}$ and quantizable Lie bialgebras}\label{5: subsec on repr of wLBq}  As properads  $\widehat{\LB}_\infty^{\mathrm{quant}}$ and  $\widehat{\LB}^{\mathrm{quant}}$ are vertex completed one must be careful when defining their representations in a dg space $V$.

Let $F_p\widehat{\LB}_\infty^{\mathrm{quant}}$, $F_p\widehat{\LB}_\infty^{\mathrm{min}}$, $F_p\widehat{\LB}^{\mathrm{quant}}$ and $F_p\widehat{\LB}$ be the sub-properads generated by graphs with $\geq p$ vertices, and let $\la$ be a formal parameter of homological degree zero. By a representation of, say, $\widehat{\LB}_\infty^{\mathrm{quant}}$
in a dg vector space $V$ we mean a morphism of properads
$$
\rho: \widehat{\LB}_\infty^{\mathrm{quant}} \lon \cE nd_V[[\la]]
$$
such that $\rho(F_p\widehat{\LB}_\infty^{\mathrm{quant}})\subset \la^p  \cE nd_V[[\la]]$ where
$\cE nd_V[[\la]]$ is the properad of formal power series in $\la$ with coefficients in $\cE nd_V$,
and  $\la^p  \cE nd_V[[\la]]\subset  \cE nd_V[[\la]]$ is a sub-properad  generated by formal power series which are divisible by $\la^p$. Representations of  $\widehat{\LB}^{\mathrm{quant}}$, $\wLBm$, $\widehat{\LB}$ and of their wheeled versions are defined similarly.

 \sip

 It is clear that there is a 1-1 correspondence between representations of $\widehat{\LB}_\infty^{\mathrm{quant}}$ in $V$ and elements $\pi^\diamond\in \fg_V[-2][[\la]]\simeq \A_3^{(n)}[[\la]]$ (for some $n$ including the case $n=+\infty$) such that the equation holds
$$
[\pi^\diamond, \pi^\diamond]_S + \sum_{p\geq 2} \frac{\la^{4p}}{(4p+2)!} \mu_{4p+2}^{\om_g}(\pi^\diamond,\ldots,\pi^\diamond)=0.
$$
As this equation involves only powers of $\la^4$, it makes sense to introduce $\hbar:=\la^4$
and consider a subclass of solutions $\pi^\diamond$ which belong to $\A_3^{(n)}[[\hbar]]$; in the case $V=\R^n$ these are precisely {\em quantizable Lie bialgebra}\, structures introduced above.

\mip

In the next subsection we construct an explicit morphism of props
$$
f^q: \cA ssB \lon \caD\widehat{\LB}^{\mathrm{quant}}
$$
and show that it lifts by a naive  induction to a morphism of dg props
$$
\cF^q: \cA ssb_\infty \lon \caD\widehat{\LB}_\infty^{\mathrm{quant}}
$$
satisfying the boundary condition (\ref{5: Boundary cond for formality map}).
Such a morphism composed with the explicit  isomorphism $\caD F:  \caD\widehat{\LB}_\infty^{\mathrm{quant}}\lon
\caD\widehat{\LB}_\infty^{\mathrm{min, \circlearrowright}}$ from Proposition {\ref{5: Prop on isomorphisms of lieb properads}}, gives us the required formality map,
$$
\caD F \circ \cF^q: \cA ssb_\infty \lon \caD\widehat{\LB}_\infty^{\mathrm{min},\circlearrowright}
$$
for finite-dimensional Lie bialgebras.

\subsection{Open problems} The prop $\widehat{\LB}^{\mathrm{quant}}$
and the dg prop $\widehat{\LB}_\infty^{\mathrm{quant}}$ have been defined with the help of explicit transcendental formulae. However it is very hard to compute the integrals
given in that formulae. For example the weights of the graphs $\ga_{10}^{2,2}$,
 $\ga_{10}^{1,3}$ and  $\ga_{10}^{3,1}$  (the first possibly non-trivial contributions) introduced in the Appendix A involve
 integrals of top-degree differential forms over $24$-dimensional configuration spaces.
 In principle all these weights might be zero so that $\widehat{\LB}_\infty^{\mathrm{quant}}$ might be identical to $\widehat{\LB}_\infty^{\mathrm{min}}$. If this is the case, then our explicit formulae
 for universal quantization of Lie bialgebras become even much simpler --- the quantization job would be done solely by the map $f^q$ given by the explicit formulae
 (\ref{7: explicit morphism f^q}).
 We conjecture, however, that the situation is quite the opposite:

 \subsubsection{\bf Conjectures} {\em (i) The set of homotopy classes of morphisms of dg
 props $
F: \widehat{\LB}_\infty^{\mathrm{quant}} \lon \widehat{\LB}_\infty^{\mathrm{min}, \circlearrowright}
$ is a torsor over the Grothendieck-Teichm\:uller group $GRT$.

\sip

(ii) The set of homotopy classes of morphisms of dg props $\cF^q: \cA ssb_\infty \lon \caD\widehat{\LB}_\infty^{\mathrm{quant}}$ consists of a single point.
}

\mip

These are open problems which we hope to address in the future. Another open problem is to construct an {\em explicit}\, isomorphism of dg  props
$$
\widehat{\LB}_\infty^{\mathrm{quant}} \lon \widehat{\LB}_\infty^{\mathrm{min}},
$$
i.e.\ to construct an analogue of our explicit morphism $F$ in (\ref{5: explicit map F from LB^q_infty to LB_infty wheeled}) which does {\em not}\, involve graphs with wheels.

\bip

%%%%%%%%%%%%%%%%%%%%%%%%%%%%%%%%%%%%%%%%%%%%%%%%%%%%%%%%%%%%%%%%%%%%%%%%%%%%%%%%%%%%%%%%%%%%%%%%%%%%%%%
%%%%%%%%%%%%%%%%%%%%%%%%%%%%%%%%%%%%%%%%%%%%%%%%%%%%%%%%%%%%%%%%%%%%%%%%%%%%%%%%%%%%%%%%%%%%%%%%%%%%%%%

\bip

{\Large
\section{\bf An explicit formula for universal quantizations of Lie bialgebras}
}

\bip

\subsection{Kontsevich compactified configuration spaces}
 Let $\overline{\bbH}=\{z=x+it\in \C | t\geq 0\}$ be  the closed  upper-half plane. Its open subset $\{z=x+it\in \C | t> 0\}$ is denoted by $\bbH$; we also consider $\p\overline{\bbH}:=\overline{\bbH}\setminus \bbH\simeq \R$. The group
 $G_2:=\R^+\rtimes \R$  acts on $\overline{\bbH}$
$$
\Ba{ccc}
G_2 \times \overline{\bbH} &\lon & \overline{\bbH}\\
(\la\in \R^+,h\in \R)\times z &\lon & \ \la z + h.
\Ea
$$
 Let $A$ and $I$ be some finite sets, and let
$$
\Conf_{A,I}(\overline{\bbH}):=\{f: A\hook \bbH, \ i: I\hook \p\overline{\bbH}\}
$$
be the configuration space of injections of $A$ into the upper half-plane, and of $I$ into the real line $\R\simeq \p\overline{\bbH}$.
 %To any such a configuration $p=\{z_i, x^0_\al\}_{i\in [k], \al\in [m]}$ with $m\geq 1$ we associate a number,
%$$
%||p||=\sum_{i=1}^k\frac{1}{k} |z_i-x_c^0| + \sum_{\al=1}^m\frac{1}{m} |x^0_i-x_c^0|, \ \ \ \ \ x_c^0:= \sum_{al=1}^m %\frac{1}{m} x_\al^0.
%$$
This is a smooth manifold of dimension $2\# A+ \#I$.
The group $G_2$ acts naturally on it, $(f(A), i(I)\rar (\la f(A)+ h, \la i(I) +h)$, and this action is free provided $2|A|+ |I|\geq 2$. The quotient space
$$
C_{A,I}(\bbH):= \Conf_{A,I}(\overline{\bbH}) /G_2
$$
is a smooth manifold of dimension $2|A|+|I|-2$.
Kontsevich constructed in \cite{Ko}
 its compactification, $\overline{C}_{A,I}(\bbH)$, which is a smooth manifold with corners, and which we use below for a construction of a new family of compactified configuration spaces. If $A=[k]$
 and $I=[n]$, we abbreviate $C_{A,I}(\bbH)$ to $C_{k,n}(\bbH)$.

\subsection{Configuration spaces of points in $\R\times \R$} Let $C_n(\R)$, $n\geq 2$,
be the configuration space of injections $\{p: [n]\rar \R\}$ modulo the action of the group
$\R^+\ltimes \R$ sending an injection $p$ into an injection $\la p + \nu$, $\la\in \R^+$,
$\nu\in \R$. We remind in Appendix B its compactification $\overline{C}_n(\R)$ which gives
us a geometric realization (in the category of semialgebraic manifolds) of Jim Stasheff's
{\em associahedra}.

\sip

Boris Shoikhet introduced  in \cite{Sh1} (with a reference to Maxim Kontsevich's informal suggestion) the configuration space $C_{m,n}(\R\times \R)$ of pairs of injections $\{p': [n]\rar \R, [m]\rar \R\}$, $m,n\geq 1$, $m+n\geq 3$, modulo the action of the group
$\R^+\ltimes \R^2$ sending a pair of injections $(p',p'')$ into $(\la p' + \nu', \la^{-1} p''+ \nu'' )$ for any  $\la\in \R^+$,
$\nu',\nu''\in \R$. We remind its compactification $\overline{C}_{m,n}(\R\times \R)$
in Appendix B, and also prove that the family of compactifications $\{\overline{C}_{m,n}(\R\times \R)\}$ gives us a geometric realization (in the category of semialgebraic manifolds) of the (pre)biassociahedra posets introduced by Martin Markl in \cite{Ma} following an earlier work by  Samson Saneblidze and Ron Umble \cite{SU}. This result gives us a nice combinatorial tool to control the boundary strata of the semialgebraic manifolds $\overline{C}_{m,n}(\R\times \R)$.

\subsection{Configuration space $C_{A;I,J}(\caH)$ and its compactification}\label{7: subsection on C_A;I,J} Let $\bbH'=\{(x,t)\in \R\times \R^{> 0}\}$ and $\bbH''=\{(y,\hat{t})\in \R\times \R^{> 0}\}$ be two copies of the  upper-half plane, and let $\overline{\bbH}'=\{(x,t)\in \R\times \R^{\geq 0}\}$ and $\overline{\bbH}''=\{(y,\hat{t})\in \R\times \R^{\geq 0}\}$ be their closures.
Consider a subspace  $\caH \subset \bbH'\times \bbH''$ given by the equation
 $t\widehat{t}=1$,  and denote by $\overline{\caH}$
its closure under the embedding into $\overline{\bbH}'\times \overline{\bbH}''$.
The space
$\overline{\caH}$ has two distinguished lines, $\bfX:=\{(x\in \R,y=0, t=0\}$ and $\bfY:=\{(x=0,y\in \R, \widehat{t}=0\}$;
it also has a natural structure of a smooth manifold with boundary.
$$
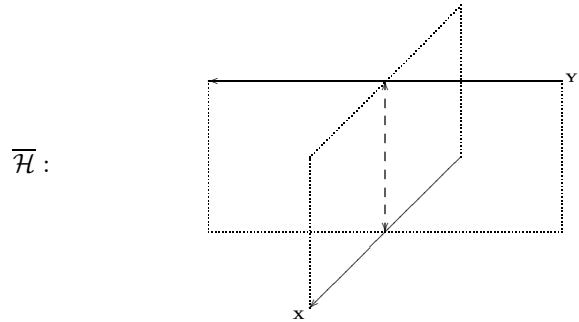
\overline{\caH}:\ \ \ \ \ \ \ \ \ \ \ \ \ \ \
\Ba{c}\resizebox{50mm}{!}{
\xy
(37,30)*{^{\mathbf Y}},
(-17,-17)*{^{\mathbf X}},
(-35,30)*{}="a",
(35,30)*{}="b",
(-35,0)*{}="a1",
(35,0)*{}="b1",
(0,30)*{}="c",
(0,0)*{}="d",
(-15,-15)*{}="a'",
(15,15)*{}="b'",
(0,0)*{}="c'",
(0,30)*{}="d'",
(-15,15)*{}="a''",
(15,45)*{}="b''",
\ar @{<-} "a";"b" <0pt>
\ar @{.} "a1";"b1" <0pt>
\ar @{.} "a1";"a" <0pt>
\ar @{.} "b1";"b" <0pt>
\ar @{-->} "c";"d" <0pt>
\ar @{<-} "a'";"b'" <0pt>
\ar @{.} "a''";"b''" <0pt>
\ar @{.} "a''";"a'" <0pt>
\ar @{.} "b''";"b'" <0pt>
\ar @{-->} "c'";"d'" <0pt>
\endxy}
\Ea
$$
 The group
$G_3:=\R^+ \rtimes \R^2$ acts on  $\overline{\caH}$,
$$
\Ba{ccccc}
\R^+ \rtimes \R^2 &\times& \overline{\caH} &\lon & \overline{\caH}\\
(\la, a,b) &\times& (x,y,t) &\lon & (\la x + a, \la^{-1}y +b , \la t).
\Ea
$$
For finite sets $A$, $I$ and $J$ let us consider a configuration space
$$
\Conf_{A;I,J}(\caH):= \{i: A\hook \caH, i':J\hook \bfX, i'': I\hook \bfY \}
$$
of injections. This is a
$(3\# A + \# I + \# J)$-dimensional smooth manifold. The group $G_3$ acts on it smoothly and,
in the case $3\# A + \# I + \# J\geq 3$ freely. {\em We assume from now on that  conditions\,
$3\# A + \# I + \# J\geq 3$, $\# I\geq 1$ and $\# J\geq 1$ hold true}, and denote by
$$
C_{A;I,J}(\caH)= \Conf_{A;I,J}(\caH)/G_3
$$
the associated smooth manifold of $G_3$-orbits.
 If $A=[k]$, $I=[m]$ and $J=[n]$ for some non-negative integers $k,m,n\in \Z^{\geq 0}$ (with $3k+m+n\geq 3$, $m,n\geq 1$), then we abbreviate $C_{[k];[m],[n]}(\caH)$
to $C_{k;m,n}(\caH)$.

 \sip

 A point $p\in  C_{A;I,J}(\caH)$ can be understood as a collection of numbers
 $$
 p=\left\{(x_a,y_a,t_a=\frac{1}{\hat{t}_a}), x^0_{\al}, y^0_{\be}\right\}_{a\in [A], \al\in [J], \be\in [I]}
 $$
defined modulo the following transformation
$$
\left\{(x_a,y_a,t_a), x^0_{\al}, y^0_{\be}\right\}\lon
\left\{(\la x_a+ h',\la^{-1}y_a +h'',\la t_a), \la x^0_{\al}+h', \la^{-1} y^0_{\be} + h''\right\}
$$
for some $\la\in \R^+$, $h',h''\in \R$.

\sip

The space $C_{0;I,J}(\caH)$ can be identified with $C_{I,J}(\R\times \R)\simeq C_{m,n}(\R\times \R)$ studied in detail in Appendix B,
and we define its compactification $\overline{C}_{0;I,J}(\caH)$ as $\overline{C}_{I,J}(\R\times \R)$.

\sip

The space $C_{A;I,J}(\caH)$ with $\# A\geq 1$ admits a canonical projection
$$
\pi: C_{A;I,J}(\caH) \lon C_{I,J}(\R\times \R)
$$
which forgets {\em internal}\, points in $\caH$ (where we assume  $C_{1,1}(\R\times \R)$ to be the one point set for consistency), and, for any $a\in A$, the following two projections
$$
\Ba{rccc}
\pi_a': & C_{A;I,J}(\caH) &\lon & C_{a,J}(\bbH')\simeq C_{1,n}(\bbH)\vspace{1mm}\\
       &    p            & \lon &  \{z_a':=x_a+it_a, x^0_{\al}\}_{\al\in I}
       \Ea
$$

$$
\Ba{rccc}
\pi_a'': & C_{k;m,n}(\caH) &\lon & C_{1,I}(\bbH'')\simeq C_{1,m}(\bbH)\vspace{1mm}\\
       &    p            & \lon &  \{z_a'':=y_a+i\frac{1}{t_a}, y^0_{\be}\}_{\be\in [m]}.
       \Ea
$$
We use these projections to construct  the following continuous map for $\# A\geq 1$

$$
\Ba{ccccccccccc}
f: C_{A;I,J}(\caH) \hspace{-2mm} &\lon & \hspace{-2mm} \displaystyle\prod_{a\in A} \overline{C}_{a,J}(\bbH') \hspace{-2mm}&\times  &\hspace{-2mm}\displaystyle \prod_{a\in A}
\overline{C}_{a,I}(\bbH'')\hspace{-2mm} &\times& \hspace{-2mm}  \overline{C}_{I,J}(\R\times \R)  &\times &  \hspace{-2mm} (S^2)^{k(k-1)} \hspace{-2mm}& \times & \hspace{-2mm}[0,+\infty]^{ k(k-1)(k-2)}\vspace{3mm}\\
p\hspace{-2mm} &\lon& \hspace{-2mm}\displaystyle \sqcap_{a\in A}\pi'_a(p) && \displaystyle \sqcap_{a\in A}\pi''_a(p) && \hspace{-4mm} \pi(p)
 &&
\hspace{-4mm} {\displaystyle \underset{_{a,b\in A\atop a\neq b}}{\sqcap}} \pi_{ab}(p)
 && \hspace{-7mm}{\displaystyle \underset{_{a,b,c\in A \atop \#\{a,b,c\}|=3}}{\sqcap}}
 \pi_{abc}(p)
\\
\Ea
$$
where
\Beq\label{5: map pi_a,b to the sphere}
 \pi_{ab}(p):= \frac{\left(x_a-x_b,  t_at_b(y_a-y_b), t_a-t_b,\right)}{\sqrt{(x_a-x_b)^2 +
(t_a-t_b)^2 +  t^2_at^2_b(y_a-y_b)^2 }}, \ \ \
\Eeq
$$
 \pi_{abc}(p):= \frac{\sqrt{(x_a-x_b)^2 +
(t_a-t_b)^2 +  t^2_at^2_b(y_a-y_b)^2 }}
{\sqrt{(x_b-x_c)^2 +
(t_b-t_c)^2 +  t^2_bt^2_c(y_b-y_c)^2 }}, \ \ \
$$
Here we assume that the last factor in the r.h.s.\ is omitted for $k<3$, and the last two factors are omitted for $k<2$ (as they have no sense in these cases). It is not hard to check that the above map is an embedding (it is essentially enough to check the cases $C_{1;1,1}(\caH)$ and $C_{2,1,1}(\caH)$)
so that we can define a {\em compactified}\, configuration space $\overline{C}_{A;I,J}(\caH)$ as the closure of the image of ${C}_{A;I,J}(\caH)$ under the map $f$. It clearly has the structure of an oriented smooth manifold with corners and also  of a semi-algebraic manifold.

\subsection{A class of differential forms on $\overline{C}_{A,I,J}(\caH)$} Consider the circle
$$
S^1=\{z\in \C: z=e^{i\theta}, \theta=Arg(z)\in [0,2\pi]\}
$$ and a 1-form on $S^1$ of the form
$\frac{1}{2\pi}\bar{g}(\theta)d\theta$ which satisfies the conditions
$$
\int_0^{2\pi} \frac{1}{2\pi}\bar{g}(\theta)d\theta=1
$$
 and
 $$
 \text{supp} (\bar{g}(\theta))\subset (0,\pi).
 $$
 Thus this $1$-form is concentrated in the upper-half of the circle. We shall use this $1$-form
 to construct a class of closed differential forms $\Omega_\Ga$ on  $\overline{C}_{k;n,m}(\caH)$
 parameterized by a set of graphs $\Ga$ we describe next.

\subsubsection{\bf A family of graphs $\cG_{k;m,n}$} The prop $\caD \wqLB=\{\caD \wqLB(m,n)\}$
introduced in \S {\ref{5: functor D}} is identical as graded vector space 
to the prop $\wLBm$ and hence admits the same set $\{\cG_{k;m,n}\}$ of basis vectors. 
%The graded space $\caD \wqLB(m,n)$
%is spanned by graphs from the sets
% of  graphs with the integer $k$ ranging from $0$ to $+\infty$.  By definition, %$\cG_{k;m,n}$
%consists of  oriented  graphs $\Ga$ with $k$ vertices called {\em internal}, $m$ %vertices called
%{\em out}-vertices and $n$ vertices called {\em in}-vertices and satisfying the %following conditions
%\Bi
%\item[(i)] every internal vertex of $\Ga$ is at least trivalent, and has at %least one incoming edge
%and one outgoing edge,
%\item[(ii)] every in-vertex can have only outgoing edges (called {\em %out-legs}),
%\item[(iii)] every out-vertex can have only ingoing edges (called {\em %in-legs}),
%\item[(iv)]  there are no edges connecting in-vertices to out-vertices,
%\item[(v)] bijections  $V_{in}(\Ga)\rar [n]$, $V_{out}(\Ga)\rar [m]$
%are fixed,
%\item[(vi)] the sets  $E_{in}(\Ga)$ and $E_{out}(\Ga)$ are totally ordered up to %an even permutation,
%\Ei
For example (we omit labellings of white vertices by integers),
$$
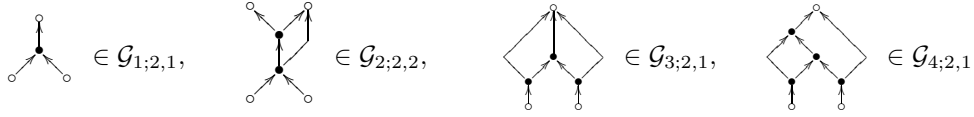
\Ba{c}\resizebox{11mm}{!}{
\xy
(0,13)*{\circ}="0",
 (0,7)*{\bu}="a",
(-5,2)*{\circ}="b_1",
(5,2)*{\circ}="b_2",
(-8,-2)*{}="c_1",
(-2,-2)*{}="c_2",
(2,-2)*{}="c_3",
\ar @{->} "a";"0" <0pt>
\ar @{<-} "a";"b_1" <0pt>
\ar @{<-} "a";"b_2" <0pt>
\endxy}
\Ea\in \cG_{1;2,1}, \ \ \ \
\Ba{c}\resizebox{10mm}{!}{
\xy
(0,13)*{\bu}="0",
 (0,7)*{\bu}="a",
(5,12)*{}="R",
(-5,2)*{\circ}="b_1",
(5,2)*{\circ}="b_2",
(-5,18)*{\circ}="u_1",
(5,18)*{\circ}="u_2",
\ar @{->} "a";"0" <0pt>
\ar @{<-} "a";"b_1" <0pt>
\ar @{<-} "a";"b_2" <0pt>
\ar @{->} "0";"u_1" <0pt>
\ar @{->} "0";"u_2" <0pt>
\ar @{-} "a";"R" <0pt>
\ar @{->} "R";"u_2" <0pt>
\endxy}
\Ea\in \cG_{2;2,2},
\ \ \
\ \ \
\Ba{c}\resizebox{14mm}{!}{
\xy
(0,17)*{\circ}="u",
%(-5,12)*{\bu}="0",
 (0,7)*{\bu}="a",
(-10,7)*{}="L",
(10,7)*{}="R",
(-5,2)*{\bu}="b_1",
(5,2)*{\bu}="b_2",
(-5,-3)*{\circ}="c_1",
(5,-3)*{\circ}="c_3",
\ar @{<-} "a";"b_1" <0pt>
\ar @{<-} "a";"b_2" <0pt>
\ar @{<-} "b_1";"c_1" <0pt>
\ar @{<-} "b_2";"c_3" <0pt>
\ar @{-} "b_1";"L" <0pt>
\ar @{<-} "u";"L" <0pt>
\ar @{-} "b_2";"R" <0pt>
\ar @{->} "R";"u" <0pt>
\ar @{->} "a";"u" <0pt>
\endxy}
\Ea\in \cG_{3;2,1},
\ \ \
\Ba{c}\resizebox{14mm}{!}{
\xy
(0,17)*{\circ}="u",
(-5,12)*{\bu}="0",
 (0,7)*{\bu}="a",
(-10,7)*{}="L",
(10,7)*{}="R",
(-5,2)*{\bu}="b_1",
(5,2)*{\bu}="b_2",
(-5,-3)*{\circ}="c_1",
(5,-3)*{\circ}="c_3",
\ar @{->} "a";"0" <0pt>
\ar @{<-} "a";"b_1" <0pt>
\ar @{<-} "a";"b_2" <0pt>
\ar @{<-} "b_1";"c_1" <0pt>
\ar @{<-} "b_2";"c_3" <0pt>
\ar @{-} "b_1";"L" <0pt>
\ar @{<-} "0";"L" <0pt>
\ar @{-} "b_2";"R" <0pt>
\ar @{->} "R";"u" <0pt>
\ar @{->} "0";"u" <0pt>
\endxy}
\Ea\in \cG_{4;2,1}
 $$

Thus graphs from $\cG_{k;m,n}$ admit a flow which we always assume in our pictures to be directed from the bottom to the top (so that there is no need to show directions of the edges
anymore). As before, $E_{int}(\Ga)$  stands for the set of {\em internal}\, edges, $E_{in}(\Ga)$ for the set of in-legs, 
    $E_{out}(\Ga)$ for the set of out-legs.

\subsubsection{\bf From graphs to differential forms}
Consider a graph $\Ga\in \cG_{k;m,n}$ with $3k+m+n\geq 3$, and an associated configuration space
$$
C(\Ga):=C_{E_{int}(\Ga); E_{out}(\Ga), E_{in}(\Ga)}(\caH)\simeq C_{k;n,m}(\caH).
$$
Let $\sC(\Ga)$ be a subspace of $C(\Ga)$  consisting of points
 $$
 p=\left\{(x_a,y_a,t_a=\frac{1}{\hat{t}_a}), x^0_{\al}, y^0_{\be}\right\}_{a\in E_{int}(\Ga), \al\in E_{in}(\Ga), \be\in E_{out}(\Ga)}
 $$
with
$$
z'_a(p):=x_a+it_a \neq z'_b(p):=x_b+it_b \ \ \text{and}\ \ \ z''_a(p):=y_a+i\frac{1}{t_a} \neq z''_b(p):=y_b+i\frac{1}{t_b}\ \ \ \ \ \ \forall\ a\neq b \in V_{int}(\Ga),
$$
i.e.\ with projections of internal vertices on planes $\bbH'$ and $\bbH''$ being different (so that
differential forms  $dArg(z'_a(p) - z'_b(p))$ and $dArg(z''_a(p) - z''_b(p))$ are well-defined
on  $\sC(\Ga)$).

 \sip

We define a smooth top degree differential form $\Omega_\Ga$ on $\sC(\Ga)$,
\Beq\label{1> Def of Omega_Ga}
\Omega_\Ga:= %(-1)^{\Ga}
\bigwedge_{e\in E_{in}(\Ga)} \om'_e\ \ \wedge\ \ \bigwedge_{e\in E_{int}(\Ga)}\Omega_e\ \ \wedge \ \
 \bigwedge_{e\in E_{out}(\Ga)} \om''_e
\Eeq
where $\om'_e$ and $\om''_e$ are 1-forms and $\Omega_e$ is a 2-form defined as follows.
%{\bf TODO: define the sign $(-1)^\Ga$, or the orderings associated with any particular %graph.}
Identifying vertices of $\Ga$ with their images in $\overline{\caH}$ under injections
$(i,i',i'')$,
%then
%, setting for any pair of distinct vertices $v_1,v_2$ of $\Ga$,
%$$
%\theta'_{v_1,v_2}:=  Arg(z'(v_2) - z'(v_1))\in (0,\pi)\ \ \ %\theta''_{v_1,v_2}:=  Arg(\overline{z''(v_2) - z''(v_1)})\in (0,\pi)
%$$
we define,
\Bi
\item[(i)] for any in-leg $e=\xy
(0,2)*{^{v_1}},
(8,2)*{^{v_2}},
(0,0)*{\circ}="0",
 (8,0)*{\bu}="a"
\ar @{->} "0";"a" <0pt>
\endxy\in E_{in}(\Ga)$,  $\om'_e:=\frac{1}{2\pi}\bar{g}\left(Arg\left(z'_{v_2} - x^0_{v_1}\right)\right)d Arg\left(z'_{v_2} - x^0_{v_1}\right)$,
\item[(ii)] for any out-leg $e=\xy
(0,2)*{^{v_1}},
(8,2)*{^{v_2}},
(0,0)*{\bu}="0",
 (8,0)*{\circ}="a"
\ar @{->} "0";"a" <0pt>
\endxy\in E_{out}(\Ga)$,  $\om''_e:=\frac{1}{2\pi}\bar{g}\left(Arg\left(\overline{y^0_{v_2} - z''_{v_1}}\right)\right)dArg\left(\overline{y^0_{v_2} - z''_{v_1}}\right)$,
\item[(iii)] for any internal edge $e=\xy
(0,2)*{^{v_1}},
(8,2)*{^{v_2}},
(0,0)*{\bu}="0",
 (8,0)*{\bu}="a"
\ar @{->} "0";"a" <0pt>
\endxy\in E_{int}(\Ga)$,
%\Beq\label{1: Omega_e}
$$
\Omega_e:=\frac{1}{(2\pi)^2}\bar{g}\left(Arg\left(z'(v_2) - z'(v_1)\right)\right) \bar{g}\left(Arg\left(\overline{z''(v_2) - z''(v_1)}\right)\right)d Arg\left(z'(v_2) - z'(v_1)\right)\wedge dArg\left(\overline{z''(v_2) - z''(v_1)}\right)
$$
\Ei
As the function $\bar{g}$ has support in the upper-half of the circle, the differential form
$
\Omega_\Ga$ extends smoothly to the configuration space  $C(\Ga)$ and even to its compactification
$\overline{C}(\Ga):=\overline{C}_{E_{int}(\Ga); E_{out}(\Ga), E_{in}(\Ga)}(\caH)$.

\sip

A subset of $\cG_{k;m,n}$ consisting of graphs $\Ga$ satisfying the condition
$$
3\# E_{int}(\Ga) + \# E_{in}(\Ga) + \# E_{out}(\Ga)= 3k+m+n-3
$$
is denoted by $\cG_{k;m,n}^{top}$ as the associated differential forms $\Omega_\Ga$ give us top-degree
forms on the configuration space $\overline{C}(\Ga)$.

\sip

\sip

Notice that if a graph $\Ga\in \cG_{k;m,n}$ satisfies the condition
$$
    3\# E_{int}(\Ga) + \# E_{in}(\Ga) + \# E_{out}(\Ga)= 3k+m+n-4
$$
then the associated differential form $\Omega_\Ga$ has degree $\dim C(\Ga)-1$ and hence one can apply the Stokes theorem to $d\Omega_\Ga$ which is a top degree form. As $\Omega_\Ga$ is closed, we obtain
$$
0=\int_{\overline{C}(\Ga)} d\Omega_\Ga = \int_{\p\overline{C}(\Ga)}
\Omega_\Ga
$$
Let us check in a few concrete examples all the boundary strata in $\p\overline{C}(\Ga)$ on which
the form $\Omega_\Ga$ does not vanish identically.

\subsubsection{\bf Example} Consider
$$
\Ga=
\Ba{c}\resizebox{9mm}{!}{
\xy
(0,13)*{\bu}="0",
 (0,7)*{\bu}="a",
(-5,2)*{\circ}="b_1",
(5,2)*{\circ}="b_2",
(-5,18)*{\circ}="u_1",
(5,18)*{\circ}="u_2",
\ar @{->} "a";"0" <0pt>
\ar @{<-} "a";"b_1" <0pt>
\ar @{<-} "a";"b_2" <0pt>
\ar @{->} "0";"u_1" <0pt>
\ar @{->} "0";"u_2" <0pt>
%\ar @{-} "a";"R" <0pt>
%\ar @{->} "R";"u_2" <0pt>
\endxy}
\Ea\in \cG_{2;2,2}.
$$
The associated $7$-dimensional configuration space $C(\Ga)$ is given by the data,
\Beq\label{7: C(Ga) for the first example}
\left\{ \left[\Ba{c}z_1'=x_1+it_1 \\ z_2'=x_2+ i{t_2}\\  x_1^0, x_2^0\in \R \Ea\right],  \left[\Ba{c}z_1''=y_1+ \frac{i}{t_1} \\ z_2''=y_2+ \frac{i}{t_2}\\
y_1^0, y_2^0\in \R
\Ea\right]\ \text{with}\ x_1^0<x_2^0, y_1^0< y_2^0
\right\}
\end{equation}
modulo the action of the 3-dimensional group $G_3$. The $6$-form $\Omega_\Ga$ is given by
$$
\Omega_\Ga=\Omega_\Ga'\wedge \Omega_\Ga''
$$
where
$$
\Omega_\Ga':=\Omega_{\bar{g}}(z'_2-z_1')\wedge \Omega_{\bar{g}}(z_1'-x_1^0)\wedge \Omega_{\bar{g}}(z_1'-x_2^0), \ \ \
\Omega_{\Ga}'':=\Omega_{\bar{g}}(\overline{z_2''-z_1''})\wedge
\wedge\Omega_{\bar{g}}(\overline{y_1^0-z_2''})\wedge \Omega_{\bar{g}}(\overline{y_2^0-z_2''})
$$
%$$
%\Omega_\Ga=\Omega_{\bar{g}}(z_2'-z_1')\wedge\Omega_{\bar{g}}(\overline{z_2''-z_1''})\wedge
%\Omega_{\bar{g}}(z_1'-x_1^0)\wedge \Omega_{\bar{g}}(z_1'-x_2^0)
%\wedge\Omega_{\bar{g}}(\overline{z_2''-y_1^0})\wedge \Omega_{\bar{g}}(\overline{z_2''-y_2^0})
%$$
and the 1-form $\Omega_{\bar{g}}$ is given by
$$
\Omega_{\bar{g}}(z_1-z_2):= \frac{1}{2\pi}{\bar{g}}(Arg(z_1-z_2))dArg(z_1-z_2).
$$
Let us classify the boundary strata in $\p \overline{C}(\Ga)$ on which the form $\Omega_\Ga$ does not vanish identically.

\mip

{\sc Case I}.
Consider the boundary strata in which two internal vertices  collapse into one internal vertex, that is, the limit $\var\rar 0$ of the configuration in which $(x_1^0, x_2^0, y_1^0, y_2^0)$ stay constant, and
 $$
\left(z_a' = x_* + it_* + \var (\bx_a + i\bt_a),  z_a''= y_* + \var \by_a +\frac{i}{t_* + \var \bt_a}\right)_{a=1,2}
$$
 It is isomorphic to $C_2(\R^3)\times C(\Ga/\Ga_{V_{int}(\Ga)})$ where $\Ga_{V_{int}(\Ga)}= \xy
 (0,0)*{\bullet}="a",
(5,0)*{\bu}="b",
\ar @{->} "a";"b" <0pt>
\endxy $
is the complete subgraph of $\Ga$ spanned by the two internal vertices, and
$$
\Ga/\Ga_{V_{int}(\Ga)}=\Ba{c}\resizebox{11mm}{!}{
\xy
(0,7)*{\bu}="0",
(-5,2)*{\circ}="b_1",
(5,2)*{\circ}="b_2",
(-5,12)*{\circ}="u_1",
(5,12)*{\circ}="u_2",
\ar @{<-} "0";"b_1" <0pt>
\ar @{<-} "0";"b_2" <0pt>
\ar @{->} "0";"u_1" <0pt>
\ar @{->} "0";"u_2" <0pt>
\endxy}
\Ea
$$
is the quotient graph obtained from $\Ga$ by collapsing the subgraph  $\Ga_{V_{int}(\Ga)}$
into a single internal vertex. As we can fix the unique internal vertex at of the latter graph $(x_*=0,y_*=0,t_*=1)$  and
$$
\lim_{\var\rar 0} Arg(z_2'-z_1')=Arg(\bx_2-\bx_1 + i(\bt_2  - \bt_1))
$$
and
\Beqrn
\lim_{\var\rar 0} Arg(\overline{z_2''-z_1''})&=&\lim_{\var\rar 0}  Arg(\by_2   - \by_1  - \frac{i}{t_* + \var \bt_2} + \frac{i}{t_* + \var \bt_1})\\
&=&Arg(\by_2 -\by_1 + i(\bt_2 - \bt_1))
\Eeqrn
we obtain a factorization
$$
\int_{C_2(\R^3)\times C(\Ga/\Ga_{V_{int}(\Ga)})}\Omega_\Ga=\int_{C_2(\R^3)=S^2}\omega_{\bar{g}} \cdot
\int_{C(\Ga/\Ga_{V_{int}(\Ga)})}\Omega_{\Ga/\Ga_{V_{int}(\Ga)}}=(\Lambda_{\bar{g}}^{(2)})^2.
$$

\mip

{\sc Case II}. Using invariance under the group $G_3$ we can always assume  that the point $(x_2,y_2,t_2)$ is fixed at, say, $(0,0,1)$. Thus it remains to consider limit configurations in which the projection $z_1'$ collapses
to a point $x_*$ in the boundary $t=0$ of $\overline{\bbH}'$,
 $$
\left(z_1'= x_*   + \var (\bx_1 + i\bt_1),   z_2'=x_2 + i t_2, z_1''=\by_1(\var) + \frac{i}{\var \bt_1}, z_2''= y_2^*  +\frac{i}{t_2}\right)\  \ \text{with} \ \var\rar 0.
$$
for some function  $\by_1^*(\var)$ of the parameter $\var$.
The limit
$$
\lim_{\var\rar 0} dArg(z_1' - x_1^0)\wedge  dArg(z_1' - x_2^0)
$$
can be non-zero if and only if  the boundary points $x_2^0$ and $x_1^0$  also collapse to $x_*$,
$$
x_1^0=x_* + \var \bx_1^0, \ \ \ \  x_2^0= x_* + \var \bx_2^0,
$$
so that we get in that limit
$$
\Omega_\Ga'
\ \underset{\var\rar 0}{\lon} \ \Omega_{\bar{g}}(z'_2-x_*)\wedge \Omega_{\bar{g}}(\bz_1'-\bx_1^0)\wedge \Omega_{\bar{g}}(\bz_1'-\bx_2^0)
$$
where $\bz_1=\bx_1+i\bt_1$. To make the form
$$
d Arg(\overline{z_2''-z_1})= dArg (y_2-\by_1(\var)-\frac{i}{t_2} + \frac{i}{\var \bt_1})
$$
non-zero in the limit $\var\rar 0$, we have to assume
$$
\by_1(\var) \thicksim \text{const} + \frac{\by_*}{\var} \ \text{for some}\ \by_*\in \R
$$
and then get in the limit
$$
\Omega_\Ga''
\ \underset{\var\rar 0}{\lon} \ \Omega_{\bar{g}}(\overline{\by_1-\bz_1''})\wedge \Omega_{\bar{g}}(\overline{y_1^0-z_2''})\wedge \Omega_{\bar{g}}(\overline{y_2^0- z_2''})
$$
where $\bz_1''=\by_1 + \frac{i}{\bt_1}$. We  conclude that this boundary strata is isomorphic to $C_{1;1,2}(\caH)\times  C_{1;2,1}(\caH)$
and the integral over it factorizes as follows
$$
\int_{C_{1;1,2}(\caH)\times  C_{1;2,1}(\caH)}\Omega_{\Ga_1}=-\int_{C_{1;1,2}(\caH)}\Omega_{\Ga_1} \cdot
\int_{C_{1;2,1}(\caH)}\Omega_{\Ga_2}=-(\Lambda_{\bar{g}}^{(2)})^2,
$$
where
\Beq\label{7: graphs Ga_1 and Ga_2}
\Ga_1=
 \Ba{c}\resizebox{11mm}{!}{
\xy
%(-5,0)*{_1},
%(5,0)*{_2},
%
(0,13)*{\circ}="0",
 (0,7)*{\bu}="a",
(-5,2)*{\circ}="b_1",
(5,2)*{\circ}="b_2",
(-8,-2)*{}="c_1",
(-2,-2)*{}="c_2",
(2,-2)*{}="c_3",
\ar @{->} "a";"0" <0pt>
\ar @{<-} "a";"b_1" <0pt>
\ar @{<-} "a";"b_2" <0pt>
\endxy}
\Ea \ \ \ \ \ \ \ \ \ \
\Ga_2=
\Ba{c}\resizebox{11mm}{!}{
\xy
%(-5,7)*{_1},
%(5,7)*{_2},
%
(0,-6)*{\circ}="0",
 (0,0)*{\bu}="a",
(-5,5)*{\circ}="b_1",
(5,5)*{\circ}="b_2",
(-8,9)*{}="c_1",
(-2,9)*{}="c_2",
(2,9)*{}="c_3",
\ar @{<-} "a";"0" <0pt>
\ar @{->} "a";"b_1" <0pt>
\ar @{->} "a";"b_2" <0pt>
\endxy}
\Ea
\Eeq
As expected,
 $
 \int_{\p\overline{C}(\Ga)}\Omega_\Ga= - (\Lambda_{\bar{g}}^{(2)})^2 + (\Lambda_{\bar{g}}^{(2)})^2=0.$

\sip

\subsubsection{\bf Example}\label{7: example 2} Consider
$$
\Ga=
\Ba{c}\resizebox{9mm}{!}{
\xy
(+3,11)*{\bu}="0",
 (-3,8)*{\bu}="a",
(-5,2)*{\circ}="b_1",
(5,2)*{\circ}="b_2",
(-5,18)*{\circ}="u_1",
(5,18)*{\circ}="u_2",
\ar @{->} "a";"0" <0pt>
\ar @{<-} "a";"b_1" <0pt>
\ar @{<-} "0";"b_2" <0pt>
\ar @{->} "a";"u_1" <0pt>
\ar @{->} "0";"u_2" <0pt>
%\ar @{-} "a";"R" <0pt>
%\ar @{->} "R";"u_2" <0pt>
\endxy}
\Ea\in
 \cG_{2;2,2}.
$$
The associated $7$-dimensional configuration space $C(\Ga)$ is given by the same data
as in (\ref{7: C(Ga) for the first example}),
while the $6$-form $\Omega_\Ga$ is given by
%$$
%\Omega_\Ga=\Omega_\Ga'\wedge \Omega_\Ga''
%$$
$$
\Omega_\Ga:= \Omega_{\bar{g}}(z_1'-x_1^0)\wedge \Omega_{\bar{g}}(z_2'-x_2^0)\wedge \Omega_{\bar{g}}(z'_2-z_1')\wedge\Omega_{\bar{g}}(\overline{z_2''-z_1''})\wedge
\wedge\Omega_{\bar{g}}(\overline{y_1^0-z_1''})\wedge \Omega_{\bar{g}}(\overline{y_1^0-z_2''})
$$

%$$
%\Omega_\Ga=\Omega_{\bar{g}}(z_2'-z_1')\wedge\Omega_{\bar{g}}(\overline{z_2''-z_1''})\wedge
%\Omega_{\bar{g}}(z_1'-x_1^0)\wedge \Omega_{\bar{g}}(z_2'-x_2^0)
%\wedge\Omega_{\bar{g}}(\overline{z_1''-y_1^0})\wedge \Omega_{\bar{g}}(\overline{z_2''-y_2^0})
%$$

Let us classify again the boundary strata in $\p \overline{C}(\Ga)$ which  can contribute non-trivially into the vanishing
integral $\int_{\p\overline{C}(\Ga)} \Omega_\Ga$.

\mip

{\sc Case 0}. Consider the boundary configurations in which the internal points
stay invariant while (i)  $|x_2^0-x_1^0| \rar 0$, or (ii) $|y_2^0-y_1^0|$,
or (iii)  $|x_2^0-x_1^0| \rar +\infty$, or (iv) $|y_2^0-y_1^0|\rar +\infty$. The forms $\Omega_\Ga$ vanishes
identically on boundary strata of types (iii) and (iv), while on strata
of types (i) and, respectively, (ii) ones obtains the integrals
$$
\int_{C_{2;2,1}(\caH)} \Omega_{\Ga_1'}\ \ \text{and}\ \
\int_{C_{2;1,2}(\caH)} \Omega_{\Ga_2'}
\ \ , \ \
\text{where} \ \ \ \ \
 \Ga_1'=
\Ba{c}\resizebox{9mm}{!}{
\xy
(+3,11)*{\bu}="0",
 (-3,8)*{\bu}="a",
(0,2)*{\circ}="b_1",
(0,2)*{\circ}="b_2",
(-5,18)*{\circ}="u_1",
(5,18)*{\circ}="u_2",
\ar @{->} "a";"0" <0pt>
\ar @{<-} "a";"b_1" <0pt>
\ar @{<-} "0";"b_2" <0pt>
\ar @{->} "a";"u_1" <0pt>
\ar @{->} "0";"u_2" <0pt>
%\ar @{-} "a";"R" <0pt>
%\ar @{->} "R";"u_2" <0pt>
\endxy}
\Ea\ \ \   ,  \ \ \
\Ga_2'=
\Ba{c}\resizebox{9mm}{!}{
\xy
(+3,11)*{\bu}="0",
 (-3,8)*{\bu}="a",
(-5,2)*{\circ}="b_1",
(5,2)*{\circ}="b_2",
(0,18)*{\circ}="u_1",
(0,18)*{\circ}="u_2",
\ar @{->} "a";"0" <0pt>
\ar @{<-} "a";"b_1" <0pt>
\ar @{<-} "0";"b_2" <0pt>
\ar @{->} "a";"u_1" <0pt>
\ar @{->} "0";"u_2" <0pt>
%\ar @{-} "a";"R" <0pt>
%\ar @{->} "R";"u_2" <0pt>
\endxy}
\Ea
$$
(which happen to vanish identically --- one can use the standard reflection argument to check this claim  which plays no role below).

\mip

{\sc Case I} is exactly the same as Case I in the previous example. The boundary strata is isomorphic
to $C_2(\R^3)\times C(\Ga/\Ga_{V{int}(\Ga)})$ and one has
$$
\int_{C_2(\R^3)\times C(\Ga/\Ga_{V_{int}(\Ga)})}\Omega_\Ga=\int_{C_2(\R^3)}\omega_{\bar{g}} \cdot
\int_{C(\Ga/\Ga_{V_{int}(\Ga)})}\Omega_{\Ga/\Ga_{V_{int}(\Ga)}}=(\Lambda_{\bar{g}}^{(2)})^2.
$$

\mip

{\sc Case II}. Using invariance under the group $G_3$ we can always assume  that the point $(x_2,y_2,t_2)$ is fixed at, say, $(0,0,1)$. Thus it remains to consider limit configurations in which the projection $z_1'$ collapses
to a point $x_*$ in the boundary $t=0$ of $\overline{\bbH}'$,
 $$
\left(z_1'= x_*   + \var (\bx_1 + i\bt_1),   z_2'=x_2 + i t_2, z_1''=\by_1(\var) + \frac{i}{\var \bt_1}, z_2''= y_2^*  +\frac{i}{t_2}\right)\  \ \text{with} \ \var\rar 0.
$$
for some function  $\by_1^*(\var)$ of the parameter $\var$. Arguing as in the Case II of the previous example, we conclude that for $\Omega_\Ga$ not to vanish identically we have to assume
$$
x_1^0=x_* +\var \bx_1^0 \ , \ \by_1(\var) = \text{const} + \frac{\by_1}{\var}\ , \ y_1^0=\text{const}+
 \frac{\by_1^0}{\var}\ \text{for some}\ \bx_1^0,\by_1, \by_1^0\in \R
$$
so that we get in the limit
$$
\lim_{\var \rar 0} \Omega_\Ga=-\Omega_{\Ga_2}\wedge \Omega_{\Ga_1}
$$
where
$$
\Omega_{\Ga_2}:=\Omega_{\bar{g}}(\bz'_1-\bx_1^0)\wedge \Omega_{\bar{g}}(\overline{\by_1^0-\bz''_1})\wedge \Omega_{\bar{g}}(\overline{0-\bz_1''})\
\text{and}\
\Omega_{\Ga_1}:=\Omega_{\bar{g}}(\overline{z_2''-x_*})\wedge
\wedge\Omega_{\bar{g}}(\overline{z_2''-x_2^0})\wedge \Omega_{\bar{g}}(\overline{y_2^0-z_2''})
$$
are the differential forms associated to the graphs in (\ref{7: graphs Ga_1 and Ga_2}). This boundary stratum is isomorphic to $C_{1;2,1}(\caH)\times C_{1;1;2}(\caH)=C(\Ga_2)\times C(\Ga_1)$ and we get
$$
\int_{C_{1;2,1}(\caH)\times  C_{1;1,2}(\caH)}\Omega_{\Ga}=-\int_{C(\Ga_2)}\Omega_{\Ga_2}\  \cdot \
\int_{C(\Ga_1)}\Omega_{\Ga_1}=-(\Lambda_{\bar{g}}^{(2)})^2.
$$
\subsubsection{\bf A useful observation}\label{6: useful observation}
Notice that   the only boundary strata in the above two examples which lie in the fibre of the surjection
$$
\pi: \overline{C}(\Ga)) \lon  C_{2,2}(\R\times\R)
$$
over a generic point in the base and contributes non-trivially into the integral is the boundary strata of type $I$.

\sip

Analyzing similarly the graphs
$$
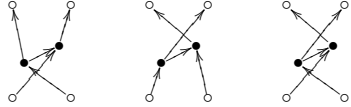
\Ba{c}\resizebox{10mm}{!}{
\xy
(+3,11)*{\bu}="0",
 (-3,8)*{\bu}="a",
(-5,2)*{\circ}="b_2",
(5,2)*{\circ}="b_1",
(-5,18)*{\circ}="u_1",
(5,18)*{\circ}="u_2",
\ar @{->} "a";"0" <0pt>
\ar @{<-} "a";"b_1" <0pt>
\ar @{<-} "0";"b_2" <0pt>
\ar @{->} "a";"u_1" <0pt>
\ar @{->} "0";"u_2" <0pt>
%\ar @{-} "a";"R" <0pt>
%\ar @{->} "R";"u_2" <0pt>
\endxy}
\Ea
\ \ \ \
\Ba{c}\resizebox{10mm}{!}{
\xy
(+3,11)*{\bu}="0",
 (-3,8)*{\bu}="a",
(-5,2)*{\circ}="b_1",
(5,2)*{\circ}="b_2",
(-5,18)*{\circ}="u_2",
(5,18)*{\circ}="u_1",
\ar @{->} "a";"0" <0pt>
\ar @{<-} "a";"b_1" <0pt>
\ar @{<-} "0";"b_2" <0pt>
\ar @{->} "a";"u_1" <0pt>
\ar @{->} "0";"u_2" <0pt>
%\ar @{-} "a";"R" <0pt>
%\ar @{->} "R";"u_2" <0pt>
\endxy}
\Ea
\ \ \ \
\Ba{c}\resizebox{10mm}{!}{
\xy
(+3,11)*{\bu}="0",
 (-3,8)*{\bu}="a",
(-5,2)*{\circ}="b_2",
(5,2)*{\circ}="b_1",
(-5,18)*{\circ}="u_2",
(5,18)*{\circ}="u_1",
\ar @{->} "a";"0" <0pt>
\ar @{<-} "a";"b_1" <0pt>
\ar @{<-} "0";"b_2" <0pt>
\ar @{->} "a";"u_1" <0pt>
\ar @{->} "0";"u_2" <0pt>
%\ar @{-} "a";"R" <0pt>
%\ar @{->} "R";"u_2" <0pt>
\endxy}
\Ea
$$
we obtain the following result
for sum of the push-forwards along the map
$
\pi: \overline{C}(\Ga)\rar C_{2,2}(\R\times\R)
$ and its boundary version $\pi_{\p}: \p \overline{C}(\Ga)\rar C_{2,2}(\R\times\R)$,
\Beqrn
\sum_{\Ga\in \cG_{2;2,2}} \pi_{*\p}\left(\Omega_\Ga \right)\Ga&=&
%\left(\Lambda_{\bar{g}}^{(2)}\right)
\pi_*(\Omega_{\Ga_0})\left(
\Ba{c}\resizebox{8mm}{!}{
\xy
(0,13)*{\bu}="0",
 (0,7)*{\bu}="a",
(-5,2)*{\circ}="b_1",
(5,2)*{\circ}="b_2",
(-5,18)*{\circ}="u_1",
(5,18)*{\circ}="u_2",
\ar @{->} "a";"0" <0pt>
\ar @{<-} "a";"b_1" <0pt>
\ar @{<-} "a";"b_2" <0pt>
\ar @{->} "0";"u_1" <0pt>
\ar @{->} "0";"u_2" <0pt>
%\ar @{-} "a";"R" <0pt>
%\ar @{->} "R";"u_2" <0pt>
\endxy}
\Ea
 -
\Ba{c}\resizebox{8mm}{!}{
\xy
(+3,11)*{\bu}="0",
 (-3,8)*{\bu}="a",
(-5,2)*{\circ}="b_1",
(5,2)*{\circ}="b_2",
(-5,18)*{\circ}="u_1",
(5,18)*{\circ}="u_2",
\ar @{->} "a";"0" <0pt>
\ar @{<-} "a";"b_1" <0pt>
\ar @{<-} "0";"b_2" <0pt>
\ar @{->} "a";"u_1" <0pt>
\ar @{->} "0";"u_2" <0pt>
%\ar @{-} "a";"R" <0pt>
%\ar @{->} "R";"u_2" <0pt>
\endxy}
\Ea
+
\Ba{c}\resizebox{10mm}{!}{
\xy
(+3,11)*{\bu}="0",
 (-3,8)*{\bu}="a",
(-5,2)*{\circ}="b_2",
(5,2)*{\circ}="b_1",
(-5,18)*{\circ}="u_1",
(5,18)*{\circ}="u_2",
\ar @{->} "a";"0" <0pt>
\ar @{<-} "a";"b_1" <0pt>
\ar @{<-} "0";"b_2" <0pt>
\ar @{->} "a";"u_1" <0pt>
\ar @{->} "0";"u_2" <0pt>
%\ar @{-} "a";"R" <0pt>
%\ar @{->} "R";"u_2" <0pt>
\endxy}
\Ea
+
\Ba{c}\resizebox{8mm}{!}{
\xy
(+3,11)*{\bu}="0",
 (-3,8)*{\bu}="a",
(-5,2)*{\circ}="b_1",
(5,2)*{\circ}="b_2",
(-5,18)*{\circ}="u_2",
(5,18)*{\circ}="u_1",
\ar @{->} "a";"0" <0pt>
\ar @{<-} "a";"b_1" <0pt>
\ar @{<-} "0";"b_2" <0pt>
\ar @{->} "a";"u_1" <0pt>
\ar @{->} "0";"u_2" <0pt>
%\ar @{-} "a";"R" <0pt>
%\ar @{->} "R";"u_2" <0pt>
\endxy}
\Ea
-
\Ba{c}\resizebox{8mm}{!}{
\xy
(+3,11)*{\bu}="0",
 (-3,8)*{\bu}="a",
(-5,2)*{\circ}="b_2",
(5,2)*{\circ}="b_1",
(-5,18)*{\circ}="u_2",
(5,18)*{\circ}="u_1",
\ar @{->} "a";"0" <0pt>
\ar @{<-} "a";"b_1" <0pt>
\ar @{<-} "0";"b_2" <0pt>
\ar @{->} "a";"u_1" <0pt>
\ar @{->} "0";"u_2" <0pt>
%\ar @{-} "a";"R" <0pt>
%\ar @{->} "R";"u_2" <0pt>
\endxy}
\Ea
\right)
\\
&=& %\left(\Lambda_{\bar{g}}^{(2)}\right)
\pi_*(\Omega_{\Ga_0})\delta \Ga_0
\Eeqrn
where $\delta$ is the differential in $\caD \LB_\infty$ and
$
\Ga_0=\Ba{c}\resizebox{10mm}{!}{
\xy
(0,7)*{\bu}="0",
(-5,2)*{\circ}="b_1",
(5,2)*{\circ}="b_2",
(-5,12)*{\circ}="u_1",
(5,12)*{\circ}="u_2",
\ar @{<-} "0";"b_1" <0pt>
\ar @{<-} "0";"b_2" <0pt>
\ar @{->} "0";"u_1" <0pt>
\ar @{->} "0";"u_2" <0pt>
\endxy}
\Ea
$.

\subsection{An explicit formula for quantization of Lie bialgebras}
Let $\cG_{k;m,n}^{(3)}$ be a subset of $\cG_{k;m,n}$ consisting of  graphs  forming a basis of
the $\bS$-bimodule $\caD\widehat{\LB}^{\mathrm{quant}}$ (these graphs have, in particular, all their internal vertices $3$-valent).

\subsubsection{\bf Theorem}\label{7: Theorem on f^q} {\em There is a morphism of props
\Beq\label{7: explicit morphism f^q}
f^{q}: \cA ss\cB \lon \caD\widehat{\LB}^{\mathrm{quant}}
\Eeq
given explicitly on the generators of $\cA ss\cB$ as follows,
$$
f^{q}\left(\begin{xy}
 <0mm,0.66mm>*{};<0mm,3mm>*{}**@{.},
 <0.39mm,-0.39mm>*{};<2.2mm,-2.2mm>*{}**@{.},
 <-0.35mm,-0.35mm>*{};<-2.2mm,-2.2mm>*{}**@{.},
 <0mm,0mm>*{\circ};<0mm,0mm>*{}**@{},
   %<0mm,0.66mm>*{};<0mm,3.4mm>*{^1}**@{},
   <0.39mm,-0.39mm>*{};<2.9mm,-4mm>*{^2}**@{},
   <-0.35mm,-0.35mm>*{};<-2.8mm,-4mm>*{^1}**@{},
\end{xy}\right):=
\Ba{c}\resizebox{8mm}{!}{ \xy
(-3,0)*{_1},
(3,0)*{_2},
 (0,7)*{\circ}="a",
(-3,2)*{\circ}="b_1",
(3,2)*{\circ}="b_2",
 \endxy}\Ea
\ +\
\sum_{k\geq 1}
\sum_{\Ga \in \cG_{k;1,2}^{(3)}} \left(\int_{\overline{C}_{k;1,2}(\caH)}\Omega_\Ga\right) \Ga =: \Ba{c}\resizebox{8mm}{!}{ \xy
(-3,0)*{_1},
(3,0)*{_2},
 (0,7)*{\circ}="a",
(-3,2)*{\circ}="b_1",
(3,2)*{\circ}="b_2",
 \endxy}\Ea
\ +\ f^{q}_{\geq 1}\left(\begin{xy}
 <0mm,0.66mm>*{};<0mm,3mm>*{}**@{.},
 <0.39mm,-0.39mm>*{};<2.2mm,-2.2mm>*{}**@{.},
 <-0.35mm,-0.35mm>*{};<-2.2mm,-2.2mm>*{}**@{.},
 <0mm,0mm>*{\circ};<0mm,0mm>*{}**@{},
   %<0mm,0.66mm>*{};<0mm,3.4mm>*{^1}**@{},
   <0.39mm,-0.39mm>*{};<2.9mm,-4mm>*{^2}**@{},
   <-0.35mm,-0.35mm>*{};<-2.8mm,-4mm>*{^1}**@{},
\end{xy}\right)
$$
$$
f^{q}\left(\begin{xy}
 <0mm,-0.55mm>*{};<0mm,-2.5mm>*{}**@{.},
 <0.5mm,0.5mm>*{};<2.2mm,2.2mm>*{}**@{.},
 <-0.48mm,0.48mm>*{};<-2.2mm,2.2mm>*{}**@{.},
 <0mm,0mm>*{\circ};<0mm,0mm>*{}**@{},
 %<0mm,-0.55mm>*{};<0mm,-3.8mm>*{_1}**@{},
 <0.5mm,0.5mm>*{};<2.7mm,2.8mm>*{^2}**@{},
 <-0.48mm,0.48mm>*{};<-2.7mm,2.8mm>*{^1}**@{},
 \end{xy}\right):=
 =\Ba{c}\resizebox{8mm}{!}{ \xy
(-3,9)*{_1},
(3,9)*{_2},
 (0,0)*{\circ}="a",
(-3,7)*{\circ}="b_1",
(3,7)*{\circ}="b_2",
 \endxy}\Ea
\ +\
\sum_{k\geq 1}
\sum_{\Ga \in \cG_{k;2,1}^{(3)}} \left(\int_{\overline{C}_{k;2,1}(\caH)}\Omega_\Ga\right) \Ga :=
\Ba{c}\resizebox{8mm}{!}{ \xy
(-3,9)*{_1},
(3,9)*{_2},
 (0,0)*{\circ}="a",
(-3,7)*{\circ}="b_1",
(3,7)*{\circ}="b_2",
 \endxy}\Ea
\ +\
f^{q}_{\geq 1}\left(\begin{xy}
 <0mm,-0.55mm>*{};<0mm,-2.5mm>*{}**@{.},
 <0.5mm,0.5mm>*{};<2.2mm,2.2mm>*{}**@{.},
 <-0.48mm,0.48mm>*{};<-2.2mm,2.2mm>*{}**@{.},
 <0mm,0mm>*{\circ};<0mm,0mm>*{}**@{},
 %<0mm,-0.55mm>*{};<0mm,-3.8mm>*{_1}**@{},
 <0.5mm,0.5mm>*{};<2.7mm,2.8mm>*{^2}**@{},
 <-0.48mm,0.48mm>*{};<-2.7mm,2.8mm>*{^1}**@{},
 \end{xy}\right)
$$
where  the differential form $\Omega_\Ga$ is defined in (\ref{1> Def of Omega_Ga}).}

\begin{proof}
If $\Ga\in \cG_{k;m,n}^{(3)}$ with $m+n=4$ then
$\deg \Omega_\Ga=3k=\dim \overline{C}_{k;m,n}(\caH) -1$ so that it makes sense to apply the Stokes theorem to the vanishing differential form $d\Omega_\Ga$,
\Beq\label{6: Stokes for Omega_Ga trivalent}
0=\int_{ \overline{C}_{k;m,n}(\caH)} d\Omega_\Ga=\int_{\p  \overline{C}_{k;m,n}(\caH)} \Omega_\Ga, \ \ \ m+n=4, m,n\geq 1.
\Eeq
We claim that the equation
\Bi
\item[(i)]
$
0=\sum_{k\geq 0}\sum_{\Ga\in \cG_{k;1,3}^{(3)}} \int_{\p \overline{C}_{k;1,3}(\caH)} \Omega_\Ga \Ga$
implies that $f^{q}$ respects the first (associativity) relations in (\ref{2: bialgebra relations}),
\item[(ii)]
$
0=\sum_{k\geq 0}\sum_{\Ga\in \cG_{k;3,1}^{(3)}} \int_{\p \overline{C}_{k;3,1}(\caH)} \Omega_\Ga \Ga$
implies that $f^{q}$ respects the second (co-associativity) relations in (\ref{2: bialgebra relations}).
\item[(iii)]
$
0=\sum_{k\geq 0}\sum_{\Ga\in \cG_{k;2,2}^{(3)}} \int_{\p \overline{C}_{k;2,2}(\caH)} \Omega_\Ga \Ga$
implies that $f^{q}$ respects the third (compatibility) relations in (\ref{2: bialgebra relations}).
\Ei

We show the proof of the most difficult step (iii) --- the proofs of the first two steps (i) and (ii) are analogous.

%\sip

%Let $\Ga\in \cG_{k;2,2}$ with $k\geq 2$. For any two vertices %$v_1,v_2\in  V_{int}(\Ga)$ we write
% $v_2\geq v_1$ if there is an internal edge directed from $v_1$ to %$v_2$, or there is no internal edges between $v_1$ and $v_2$ at all %(we also write $v_2>v_1$ if there is an edge connecting $v_1$ to %$v_2$).
%For two disjoint non-empty subsets $I,J\subset V_{int}(\Ga)$ we write %$J\geq I$ if $v_2\geq v_1$ for any $v_2\in J$ and $v_1\in I$.

\sip

Let us classify all the boundary strata
on which the differential forms $\Omega_\Ga$ do not vanish identically. Let us notice that the product the function $|x_2^0-x_1^0||y_2^0-y_1^0|$
can take the following values on the codimension 1 boundary configurations:
\Bi
\item[I:] the value $|x_2^0-x_1^0||y_2^0-y_1^0|$ stays finite;
\item[II:]
$|x_2^0-x_1^0|\rar 0$ while $|y_2^0-y_1^0|$ stays finite, or
 $|y_2^0-y_1^0|\rar 0$ while $|x_2^0-x_1^0|$ stays finite;
\item[III:]
$|y_2^0-y_1^0|\rar +\infty$ while $|x_2^0-x_1^0|$ stays finite, or
 $|x_2^0-x_1^0|\rar +\infty$ while $|y_2^0-y_1^0|$ stays finite.
\Ei
Let us consider each case separately.
\sip

{\bf Case I} corresponds to the boundary strata --- which we denote by  $\p_I\overline{C}_{k;2,2}(\caH)\subset \p\overline{C}_{k;2,2}(\caH)$ --- in which  several internal points collapse into an internal point (see examples in \S {\ref{6: useful observation}}). By Proposition {\ref{3: Prop on Upsilon^om_g}}
for the case $d=3$ the following sum
$$
\sum_{k\geq 0} \sum_{\Ga\in \cG^{(3)}_{k;2,2}}
\left(\int_{\p_I \overline{C}_{k;m,n}(\caH)}\Omega_\Ga)\right) \Ga =
 \sum_{s:[2]\rar V(\Ga)\atop \hat{s}:[2]\rar V(\Ga)}
  \Ba{c}\resizebox{9mm}{!}  {\xy
 (-6,7)*{^1},
(7,7)*{^2},
(7,-8)*{_2},
(-6,-8)*{_1},
(0,0)*+{\Ga^{\om_{\bar{g}}}}="o",
(-6,6)*{\circ}="1",
(6,6)*{\circ}="4",
(-6,-6)*{\circ}="5",
(6,-6)*{\circ}="8",
\ar @{-} "o";"1" <0pt>
\ar @{-} "o";"4" <0pt>
\ar @{-} "o";"5" <0pt>
\ar @{-} "o";"8" <0pt>
\endxy}\Ea \equiv  0
$$
gives an identically vanishing element in $ \caD\widehat{\LB}^{\mathrm{quant}}$ (here the sum is taken over all possible ways of attaching four legs to the MC element $\Ga^{\om_{\bar{g}}}$ and setting to zero every graph which has at least one non-trivalent internal vertex or an internal vertex with no at least ingoing half-edge and at least one outgoing half-edge). Hence  we can skip type I boundary strata in equation (\ref{6: Stokes for Omega_Ga trivalent}).

\mip

{\bf Case II}.
Denote the associated boundary strata by $\p_{II}\overline{C}_{k;2,2}(\caH)$.
If, for example, we consider a limit configuration with $|x_2^0-x_1^0|\rar 0$ but $|y_2^0-y_1^0|$ finite, then the boundary points $x_1^0$, $x_2^0$ and, perhaps, some (possibly empty) subset $I\subset V_{int}(\Ga)$ of internal points tend in the limit $\var\rar 0$ to a point $x_*\in {\mathbf X}$,
\Beqrn
z_i'&=& x_*+ \var(\bx_i + i\bt_i), \ \ \  z_i''= y_i(\var) +  \frac{i}{\var\bt_i}   \ i\in I,\\
x_1'&=& x_* + \var \bx_1^0\\
x_2'&=& x_* + \var \bx_2^0
\Eeqrn
for some functions $y_i(\var)$ of the parameter $\var$ (it is easy to see that if $I\neq \emptyset$, then the differential form $\Omega_\Ga$ has a chance not to vanish identically on such a boundary stratum if and only if  $y_i(\var)\simeq  \frac{\by_i}{\var}\ \text{as}\ \var\rar 0$ for some $\by_i\in \R$).

\sip

Consider (as an elementary illustration) a special case $I=\emptyset$ (and denote the associated strata
in $\p_{II}\overline{C}_{k;2,2}(\caH)$ by $\p_{II\emptyset}\overline{C}_{k;2,2}(\caH)$).
It is clear that in this case we have
$$
\sum_{k\geq 1}
\sum_{\Ga \in \cG_{k;2,2}^{(3)}} \left(\int_{\p_{II\emptyset}C(\Ga)}\Omega_\Ga\right)= -\frac{f^{q}_{\geq 1}\left(\begin{xy}
 <0mm,-0.55mm>*{};<0mm,-2.5mm>*{}**@{.},
 <0.5mm,0.5mm>*{};<2.2mm,2.2mm>*{}**@{.},
 <-0.48mm,0.48mm>*{};<-2.2mm,2.2mm>*{}**@{.},
 <0mm,0mm>*{\circ};<0mm,0mm>*{}**@{},
 %<0mm,-0.55mm>*{};<0mm,-3.8mm>*{_1}**@{},
 %<0.5mm,0.5mm>*{};<2.7mm,2.8mm>*{^2}**@{},
 %<-0.48mm,0.48mm>*{};<-2.7mm,2.8mm>*{^1}**@{},
 \end{xy}\right)}{\begin{xy}
 <0mm,0.66mm>*{};<0mm,3mm>*{}**@{.},
 <0.39mm,-0.39mm>*{};<2.2mm,-2.2mm>*{}**@{.},
 <-0.35mm,-0.35mm>*{};<-2.2mm,-2.2mm>*{}**@{.},
 <0mm,0mm>*{\circ};<0mm,0mm>*{}**@{},
   %<0mm,0.66mm>*{};<0mm,3.4mm>*{^1}**@{},
   %<0.39mm,-0.39mm>*{};<2.9mm,-4mm>*{^2}**@{},
   %<-0.35mm,-0.35mm>*{};<-2.8mm,-4mm>*{^1}**@{},
\end{xy}}
$$
An analogue of this formula in the case  $|y_2^0-y_1^0|\rar 0$ while $|x_2^0-x_1^0|$ stays finite and no internal vertices collapse to the line $\mathbf Y$  would be of course the following one
$$
\sum_{k\geq 1}
\sum_{\Ga \in \cG_{k;2,2}^{(3)}} \left(\int_{\p_{IIa}C(\Ga)}\Omega_\Ga\right)= -\frac{\begin{xy}
 <0mm,-0.55mm>*{};<0mm,-2.5mm>*{}**@{.},
 <0.5mm,0.5mm>*{};<2.2mm,2.2mm>*{}**@{.},
 <-0.48mm,0.48mm>*{};<-2.2mm,2.2mm>*{}**@{.},
 <0mm,0mm>*{\circ};<0mm,0mm>*{}**@{},
 %<0mm,-0.55mm>*{};<0mm,-3.8mm>*{_1}**@{},
 %<0.5mm,0.5mm>*{};<2.7mm,2.8mm>*{^2}**@{},
 %<-0.48mm,0.48mm>*{};<-2.7mm,2.8mm>*{^1}**@{},
 \end{xy}}{f^q_{\geq 1}\left(\begin{xy}
 <0mm,0.66mm>*{};<0mm,3mm>*{}**@{.},
 <0.39mm,-0.39mm>*{};<2.2mm,-2.2mm>*{}**@{.},
 <-0.35mm,-0.35mm>*{};<-2.2mm,-2.2mm>*{}**@{.},
 <0mm,0mm>*{\circ};<0mm,0mm>*{}**@{},
   %<0mm,0.66mm>*{};<0mm,3.4mm>*{^1}**@{},
   %<0.39mm,-0.39mm>*{};<2.9mm,-4mm>*{^2}**@{},
   %<-0.35mm,-0.35mm>*{};<-2.8mm,-4mm>*{^1}**@{},
\end{xy}\right)}
$$
where we use fraction type notation for prop compositions introduced in
\cite{Ma1}  e.g.\
$$
\frac{\Ba{c}  \begin{xy}
 <0mm,-0.55mm>*{};<0mm,-2.5mm>*{}**@{.},
 <0.5mm,0.5mm>*{};<2.2mm,2.2mm>*{}**@{.},
 <-0.48mm,0.48mm>*{};<-2.2mm,2.2mm>*{}**@{.},
 <0mm,0mm>*{\circ};<0mm,0mm>*{}**@{},
 \end{xy}\vspace{-3mm} \ \\    \ \Ea}{\Ba{c}\begin{xy}
 <0mm,0.66mm>*{};<0mm,3mm>*{}**@{.},
 <0.39mm,-0.39mm>*{};<2.2mm,-2.2mm>*{}**@{.},
 <-0.35mm,-0.35mm>*{};<-2.2mm,-2.2mm>*{}**@{.},
 <0mm,0mm>*{\circ};<0mm,0mm>*{}**@{},
\end{xy}\Ea}:= \Ba{c}\begin{xy}
 <0mm,2.47mm>*{};<0mm,-0.5mm>*{}**@{.},
 <0.5mm,3.5mm>*{};<2.2mm,5.2mm>*{}**@{.},
 <-0.48mm,3.48mm>*{};<-2.2mm,5.2mm>*{}**@{.},
 <0mm,3mm>*{\circ};<0mm,3mm>*{}**@{},
  <0mm,-0.8mm>*{\circ};<0mm,-0.8mm>*{}**@{},
<0mm,-0.8mm>*{};<-2.2mm,-3.5mm>*{}**@{.},
 <0mm,-0.8mm>*{};<2.2mm,-3.5mm>*{}**@{.},
\end{xy}\Ea \ \ \ \ \  , \ \ \ \ \ \ \
\Ba{c}\resizebox{12mm}{!}{\xy
(-10,0)*{}="1L",
(10,0)*{}="1R",
(4,10)*{}="0",
 (4,6)*{\circ}="a",
(1,2)*{}="u_1",
(7,2)*{}="u_2",
(-4,10)*{}="0'",
 (-4,6)*{\circ}="a'",
(-1,2)*{}="u_1'",
(-7,2)*{}="u_2'",
(-1,-2)*{}="du1",
(-7,-2)*{}="du2",
(-4,-6)*{\circ}="vu",
 (-4,-10)*{}="vd", %left under line
(4,-10)*{}="xd",
 (4,-6)*{\circ}="x",
(1,-2)*{}="x_1",   %right under line
(7,-2)*{}="x_2",
\ar @{.} "a";"0" <0pt>
\ar @{.} "a";"u_1" <0pt>
\ar @{.} "a";"u_2" <0pt>
\ar @{.} "a'";"0'" <0pt>
\ar @{.} "a'";"u_1'" <0pt>
\ar @{.} "a'";"u_2'" <0pt>
\ar @{.} "vd";"vu" <0pt>
\ar @{.} "vu";"du1" <0pt>
\ar @{.} "vu";"du2" <0pt>
\ar @{.} "x";"xd" <0pt>
\ar @{.} "x";"x_1" <0pt>
\ar @{.} "x";"x_2" <0pt>
\ar @{-} "1L";"1R" <0pt>
\endxy}\Ea :=
\Ba{c}\resizebox{10mm}{!}{\begin{xy}
 <0mm,0mm>*{\circ};<0mm,0mm>*{}**@{},
 <0mm,-0.49mm>*{};<0mm,-3.0mm>*{}**@{.},
 <-0.5mm,0.5mm>*{};<-3mm,2mm>*{}**@{.},
 <-3mm,2mm>*{};<0mm,4mm>*{}**@{.},
 <0mm,4mm>*{\circ};<-2.3mm,2.3mm>*{}**@{},
 <0mm,4mm>*{};<0mm,7.4mm>*{}**@{.},
<0mm,0mm>*{};<2.2mm,1.5mm>*{}**@{.},
 <6mm,0mm>*{\circ};<0mm,0mm>*{}**@{},
 <6mm,4mm>*{};<3.8mm,2.5mm>*{}**@{.},
 <6mm,4mm>*{};<6mm,7.4mm>*{}**@{.},
 <6mm,4mm>*{\circ};<-2.3mm,2.3mm>*{}**@{},
 <0mm,4mm>*{};<6mm,0mm>*{}**@{.},
<6mm,4mm>*{};<9mm,2mm>*{}**@{.},
<6mm,0mm>*{};<9mm,2mm>*{}**@{.},
<6mm,0mm>*{};<6mm,-3mm>*{}**@{.},
 \end{xy}}\Ea
$$

The general case is no more difficult. Let $J:=V_{int}(\Ga)\setminus I$ be the complementary subset corresponding to points which have $\bbH'$-projections {\em not}\, tending to $x_*$ as $\var\rar 0$. We can represent each graph $\Ga$ in the sum
$$
\sum_{k\geq 0}\sum_{\Ga\in \cG_{k;2,2}^{(3)}} \int_{\p_{II} \overline{C}_{k;2,2}(\caH)} \Omega_\Ga \Ga
$$
 in the form
 $$
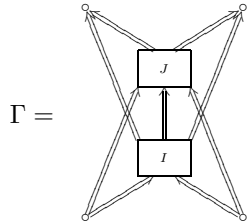
\Ga=
\Ba{c}\resizebox{23mm}{!}{
\xy
(0,28.5)*{_{J}},
(-5,32)*{}="ulUL",
(5,32)*{}="urUL",
(-5,25)*{}="dlUL",
(5,25)*{}="drUL",
(0,11.5)*{_{I}},
(-5,15)*{}="ulDL",
(5,15)*{}="urDL",
(-5,8)*{}="dlDL",
(5,8)*{}="drDL",
(-2,32)*{}="UL1L",
(2,32)*{}="UL1R",
(0,25)*{}="UL2",
(-5,25)*{}="UL2L",
(5,25)*{}="UL2R",
(0,15)*{}="DL1",
(-5,15)*{}="DL1L",
(5,15)*{}="DL1R",
(-2,8)*{}="DL2L",
(2,8)*{}="DL2R",
(-15,0)*{\circ}="b_1",
(15,0)*{\circ}="b_2",
(-15,40)*{\circ}="u_1",
(15,40)*{\circ}="u_2",
\ar @{-} "ulUL";"urUL" <0pt>
\ar @{-} "ulUL";"dlUL" <0pt>
\ar @{-} "dlUL";"drUL" <0pt>
\ar @{-} "drUL";"urUL" <0pt>
\ar @{-} "ulDL";"urDL" <0pt>
\ar @{-} "ulDL";"dlDL" <0pt>
\ar @{-} "dlDL";"drDL" <0pt>
\ar @{-} "drDL";"urDL" <0pt>
\ar @{=>} "UL1L";"u_1" <0pt>
\ar @{=>} "DL1L";"u_1" <0pt>
\ar @{=>} "DL1R";"u_2" <0pt>
\ar @{=>} "DL1";"UL2" <0pt>
\ar @{=>} "b_1";"DL2L" <0pt>
\ar @{=>} "b_1";"UL2L" <0pt>
\ar @{=>} "b_2";"UL2R" <0pt>
\ar @{=>} "UL1R";"u_2" <0pt>
\ar @{=>} "b_2";"DL2R" <0pt>
\endxy}
\Ea
$$
where directed double edges stand for (possibly empty) sets of directed edges.
Let $\Ga'$ (resp., $\Ga''$ )  be the  element of $\caD \widehat{\LB}^{\mathrm{quant}}(1,2)$
(resp., of $\in \caD \widehat{\LB}^{\mathrm{quant}}(2,1)$ defined as the complete
subgraph of $\Ga$ spanned by vertices from the set $I$ (resp., $J$), together with {\em all}\, edges attached to this set,
$$
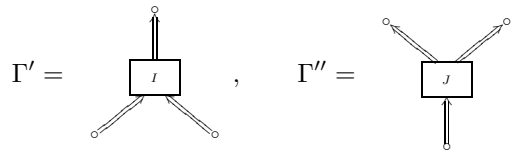
\Ga'=
\Ba{c}\resizebox{18mm}{!}{
\xy
(0,11.5)*{_{I}},
(-5,15)*{}="ulDL",
(5,15)*{}="urDL",
(-5,8)*{}="dlDL",
(5,8)*{}="drDL",
(0,25)*{\circ}="UL2",
(0,15)*{}="DL1",
(-5,15)*{}="DL1L",
(5,15)*{}="DL1R",
(-2,8)*{}="DL2L",
(2,8)*{}="DL2R",
(-12,0)*{\circ}="b_1",
(12,0)*{\circ}="b_2",
\ar @{-} "ulDL";"urDL" <0pt>
\ar @{-} "ulDL";"dlDL" <0pt>
\ar @{-} "dlDL";"drDL" <0pt>
\ar @{-} "drDL";"urDL" <0pt>
\ar @{=>} "DL1";"UL2" <0pt>
\ar @{=>} "b_1";"DL2L" <0pt>
%\ar @{=>} "b_1";"UL2" <0pt>
%\ar @{=>} "b_2";"UL2" <0pt>
%
\ar @{=>} "b_2";"DL2R" <0pt>
\endxy}
\Ea, \ \ \ \ \ \
\Ga''=
\Ba{c}\resizebox{18mm}{!}{
\xy
(0,-11.5)*{_{J}},
(-5,-15)*{}="ulDL",
(5,-15)*{}="urDL",
(-5,-8)*{}="dlDL",
(5,-8)*{}="drDL",
(0,-25)*{\circ}="UL2",
(0,-15)*{}="DL1",
(-5,-15)*{}="DL1L",
(5,-15)*{}="DL1R",
(-2,-8)*{}="DL2L",
(2,-8)*{}="DL2R",
(-12,0)*{\circ}="b_1",
(12,0)*{\circ}="b_2",
\ar @{-} "ulDL";"urDL" <0pt>
\ar @{-} "ulDL";"dlDL" <0pt>
\ar @{-} "dlDL";"drDL" <0pt>
\ar @{-} "drDL";"urDL" <0pt>
\ar @{<=} "DL1";"UL2" <0pt>
\ar @{<=} "b_1";"DL2L" <0pt>
%\ar @{=>} "b_1";"UL2" <0pt>
%\ar @{=>} "b_2";"UL2" <0pt>
%
\ar @{<=} "b_2";"DL2R" <0pt>
\endxy}
\Ea
$$
Note that out-legs  in $\Ga'$ are formed by three types of edges in $\Ga$
(and denoted in $\Ga$ by three different double arrows), the ones which connect vertices of $I$ to the left out-vertex, to the vertices of $J$, and to the right out-vertex. Similarly, the set of in-legs of $\Ga''$ encompasses three different double arrows in $\Ga$. Many different graphs $\Ga$
produce {\em identical}\, associated graphs $\Ga'$ and $\Ga''$ and it is easy to describe this family --- it is precisely the set of non-vanishing summands in the prop composition $\Ga''\ _1\circ_1 \Ga'$! As
 $$
\Omega_\Ga|_{\p_{II}\overline{C}_{k;2,2}(\caH)}=\lim_{\var \rar 0} \Omega_\Ga = \Omega_{\Ga'}\wedge \Omega_{\Ga''}.
$$
we finally get
\Beqrn
-\sum_{k\geq 0} \sum_{\Ga\in \cG^{(3)}_{k;2,2}}\int_{\p_{II}  \overline{C}_{k;2,2}(\caH)} \Omega_\Ga &=&\sum_{k',k''\geq 0} \sum_{\Ga'\in \cG^{(3)}_{k';1,2}} \sum_{\Ga''\in \cG^{(3)}_{k';2,1}}\left(\int_{\overline{C}_{k;1,2}(\caH)} \Omega_{\Ga'}\right) \cdot \left(\int_{\overline{C}_{k;2,1}(\caH)} \Omega_{\Ga'}\right) \Ga''\ _1\circ_1\ \Ga'\\
&=& f^{q}\left(\begin{xy}
 <0mm,-0.55mm>*{};<0mm,-2.5mm>*{}**@{.},
 <0.5mm,0.5mm>*{};<2.2mm,2.2mm>*{}**@{.},
 <-0.48mm,0.48mm>*{};<-2.2mm,2.2mm>*{}**@{.},
 <0mm,0mm>*{\circ};<0mm,0mm>*{}**@{},
 %<0mm,-0.55mm>*{};<0mm,-3.8mm>*{_1}**@{},
 <0.5mm,0.5mm>*{};<2.7mm,2.8mm>*{^2}**@{},
 <-0.48mm,0.48mm>*{};<-2.7mm,2.8mm>*{^1}**@{},
 \end{xy}\right) \ _1\circ_1\ f^{q}\left(\begin{xy}
 <0mm,0.66mm>*{};<0mm,3mm>*{}**@{.},
 <0.39mm,-0.39mm>*{};<2.2mm,-2.2mm>*{}**@{.},
 <-0.35mm,-0.35mm>*{};<-2.2mm,-2.2mm>*{}**@{.},
 <0mm,0mm>*{\circ};<0mm,0mm>*{}**@{},
   %<0mm,0.66mm>*{};<0mm,3.4mm>*{^1}**@{},
   <0.39mm,-0.39mm>*{};<2.9mm,-4mm>*{^2}**@{},
   <-0.35mm,-0.35mm>*{};<-2.8mm,-4mm>*{^1}**@{},
\end{xy}\right)
\\
&=&
f^{q}\left(\Ba{c}\begin{xy}
 <0mm,2.47mm>*{};<0mm,-0.5mm>*{}**@{.},
 <0.5mm,3.5mm>*{};<2.2mm,5.2mm>*{}**@{.},
 <-0.48mm,3.48mm>*{};<-2.2mm,5.2mm>*{}**@{.},
 <0mm,3mm>*{\circ};<0mm,3mm>*{}**@{},
  <0mm,-0.8mm>*{\circ};<0mm,-0.8mm>*{}**@{},
<0mm,-0.8mm>*{};<-2.2mm,-3.5mm>*{}**@{.},
 <0mm,-0.8mm>*{};<2.2mm,-3.5mm>*{}**@{.},
\end{xy}\Ea\right).
\Eeqrn

{\bf Case III}.
 %the value $|y_2^0-y_1^0|\rar +\infty$ while $|x_2^0-x_1^0|$ stays %finite, or
% $|x_2^0-x_1^0|\rar +\infty$ while $|y_2^0-y_1^0|$ stays finite.
 Denote the associated boundary strata by $\p_{III}\overline{C}_{k;2,2}(\caH)$, and
consider for concreteness  limit configuration with $|y_2^0-y_1^0|\rar +\infty$ and $x_2^0$, $x_1^0$ staying constant (the other subcase can be treated similarly). In general a (possibly empty) subset $I_1\subset V_{int}(\Ga)$ (resp,. $I_2$) can collapse to the boundary point $x_1^0$ (resp., $x_2^0$), and a (possibly empty) subset $K_1\subset V_{int}(\Ga)$ (resp., $K_2$) can tend as $\var \rar 0$ to the boundary point $y_1^0$ (resp., $y_2^0$),
\Beqrn
z_{i_1}'&=& x_1^0+ \var(\bx_{i_1} + i\bt_{i_1}), \ \ \  z_{i_1}''= y_{i_1}(\var) +  \frac{i}{\var\bt_{i_1}}   \ , \ \   i_1\in I_1,\\
z_{i_2}'&=& x_2^0+ \var(\bx_{i_2} + i\bt_{i_2}), \ \ \  z_{i_2}''= y_{i_2}(\var) +  \frac{i}{\var\bt_{i_2}}   \ , \ \   i_2\in I_2,\\
y_1^0&=& \frac{\by_1^0}{\var}\ \ , \ \ \ \ y_2^0= \frac{\by_2^0}{\var}\\
z_{k_1}'&=& x_{k_1}+   \frac{i\bt_{k_1}}{\var}, \ \ \  z_{k_1}''= \frac{\by_1^0}{\var} +  \var( \Delta \by_1^0 + \frac{i}{\bt_{k_1}})   \ , \ \   k_1\in K_1,\\
z_{k_2}'&=& x_{k_2}+   \frac{i\bt_{k_2}}{\var}, \ \ \  z_{k_2}''= \frac{\by_1^0}{\var} +  \var( \Delta \by_2^0 + \frac{i}{\bt_{k_2}})   \ , \ \   k_2\in K_2,\\
z'_j&=&x_j + it_j\ , \ \ z''_j=y_j(\var) + \frac{i}{t_j}\ ,\ \  j\in  J:=V_{int}(\Ga)\setminus I_1\sqcup I_2 \sqcup K_1 \sqcup K_2
\Eeqrn
for some functions  $y_\bu(\var)$  of the parameter $\var$ (which we have yet to understand) and some arbitrary constants in bold letters.

\sip

We claim that it is enough to consider the case when the sets $K_1$ and $K_2$ are both empty. Indeed, if at least one of the sets, say $K_1$ is not empty, it has a vertex $k\in K_1$ connected by an edge to a vertex $i$ in the set $J\sqcup I_1\sqcup I_2\sqcup \{x_1^0\}\sqcup \{x_2^0\}$
%(otherwise $\Omega_\Ga|_{\p_{III}\overline{C}_{k;2,2}(\caH)}\equiv 0$) %
which contributes into the form $\Omega_\Ga$  a factor
$$
\lim_{\var\rar 0}dArg\left({z_k'-z_i'}\right)=
\lim_{\var\rar 0}dArg\left({x_{k} + \frac{i\bt_{k}}{\var} - z_i')}\right)\rar 0
$$
which vanishes identically. Hence, for $\Omega_\Ga|_{\p_{II}\overline{C}_{k;2,2}(\caH)}$ not to vanish identically, we can assume assume that $\Ga$ has a form
$$
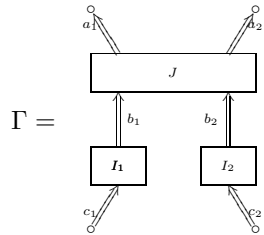
\Ga=
\Ba{c}\resizebox{25mm}{!}{
\xy
(-15,37)*{_{a_1}},
(15,37)*{_{a_2}},
(-15,3)*{_{c_1}},
(15,3)*{_{c_2}},
%
%(-20,27)*{_{b_1}},
%(20,27)*{_{b_2}},
%(-20,13)*{_{d_1}},
%(20,13)*{_{d_2}},
%
(-7,20)*{_{b_1}},
(7,20)*{_{b_2}},
(0,28.5)*{_{J}},
(-15,32)*{}="ulUL",
(-5,32)*{}="urUL",
(-15,25)*{}="dlUL",
(-5,25)*{}="drUL",
(-10,11.5)*{_{I_1}},
(-15,15)*{}="ulDL",
(-5,15)*{}="urDL",
(-15,8)*{}="dlDL",
(-5,8)*{}="drDL",
%
%(10,28.5)*{_{J_2}},
(15,32)*{}="ulUR",
(5,32)*{}="urUR",
(15,25)*{}="dlUR",
(5,25)*{}="drUR",
(-10,11.5)*{_{I_1}},
(-15,15)*{}="ulDL",
(-5,15)*{}="urDL",
(-15,8)*{}="dlDL",
(-5,8)*{}="drDL",
(10,11.5)*{_{I_2}},
(15,15)*{}="ulDR",
(5,15)*{}="urDR",
(15,8)*{}="dlDR",
(5,8)*{}="drDR",
(-10,32)*{}="UL1",
(-10,25)*{}="UL2",
(-10,15)*{}="DL1",
(-10,8)*{}="DL2",
(10,32)*{}="UR1",
(10,25)*{}="UR2",
(10,15)*{}="DR1",
(10,8)*{}="DR2",
(-15,0)*{\circ}="b_1",
(15,0)*{\circ}="b_2",
(-15,40)*{\circ}="u_1",
(15,40)*{\circ}="u_2",
\ar @{-} "ulUL";"urUR" <0pt>
\ar @{-} "ulUL";"dlUL" <0pt>
\ar @{-} "dlUL";"drUR" <0pt>
%\ar @{-} "drUL";"urUL" <0pt>
%
\ar @{-} "ulDL";"urDL" <0pt>
\ar @{-} "ulDL";"dlDL" <0pt>
\ar @{-} "dlDL";"drDL" <0pt>
\ar @{-} "drDL";"urDL" <0pt>
\ar @{-} "ulUR";"urUR" <0pt>
\ar @{-} "ulUR";"dlUR" <0pt>
\ar @{-} "dlUR";"drUR" <0pt>
%\ar @{-} "drUR";"urUR" <0pt>
%
\ar @{-} "ulDR";"urDR" <0pt>
\ar @{-} "ulDR";"dlDR" <0pt>
\ar @{-} "dlDR";"drDR" <0pt>
\ar @{-} "drDR";"urDR" <0pt>
\ar @{=>} "UL1";"u_1" <0pt>
%\ar @{=>} "DR1";"u_2" <0pt>
\ar @{=>} "DL1";"UL2" <0pt>
\ar @{=>} "b_1";"DL2" <0pt>
%\ar @{=>} "b_2";"UR2" <0pt>
%
%\ar @{=>} "DL1";"u_1" <0pt>
\ar @{=>} "UR1";"u_2" <0pt>
\ar @{=>} "DR1";"UR2" <0pt>
%\ar @{=>} "b_1";"UL2" <0pt>
\ar @{=>} "b_2";"DR2" <0pt>
%\ar @{=>} "DR1";"UL2" <0pt>
%\ar @{=>} "DL1";"UR2" <0pt>
\endxy}
\Ea
$$
where some edges ingoing into a box can continue as outgoing edges without ``hitting"
an internal vertex inside the box. Note that no edge can connect a vertex $i_1$ from $I_1$ a to a vertex $i_2$ from $I_2$ as otherwise the differential form $\Omega_\Ga$ vanishes identically in the limit $\var\rar 0$ due to the presence of the factor
$$
\lim_{\var\rar 0} dArg(z_{i_1}' - z_{i_2}')=dArg(x_1^0-x_2^0)=0.
$$

If the set $J$ is empty, then $\Ga$ takes the form
$$
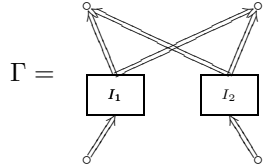
\Ga=
\Ba{c}\resizebox{25mm}{!}{
\xy
%
%
%(0,28.5)*{_{J}},
(-15,32)*{}="ulUL",
(-5,32)*{}="urUL",
(-15,25)*{}="dlUL",
(-5,25)*{}="drUL",
(-10,11.5)*{_{I_1}},
(-15,15)*{}="ulDL",
(-5,15)*{}="urDL",
(-15,8)*{}="dlDL",
(-5,8)*{}="drDL",
(15,32)*{}="ulUR",
(5,32)*{}="urUR",
(15,25)*{}="dlUR",
(5,25)*{}="drUR",
(-10,11.5)*{_{I_1}},
(-15,15)*{}="ulDL",
(-5,15)*{}="urDL",
(-15,8)*{}="dlDL",
(-5,8)*{}="drDL",
(10,11.5)*{_{I_2}},
(15,15)*{}="ulDR",
(5,15)*{}="urDR",
(15,8)*{}="dlDR",
(5,8)*{}="drDR",
(-10,32)*{}="UL1",
(-10,25)*{}="UL2",
(-10,15)*{}="DL1",
(-10,8)*{}="DL2",
(10,32)*{}="UR1",
(10,25)*{}="UR2",
(10,15)*{}="DR1",
(10,8)*{}="DR2",
(-15,0)*{\circ}="b_1",
(15,0)*{\circ}="b_2",
(-15,27)*{\circ}="u_1",
(15,27)*{\circ}="u_2",
\ar @{-} "ulDL";"urDL" <0pt>
\ar @{-} "ulDL";"dlDL" <0pt>
\ar @{-} "dlDL";"drDL" <0pt>
\ar @{-} "drDL";"urDL" <0pt>
\ar @{-} "ulDR";"urDR" <0pt>
\ar @{-} "ulDR";"dlDR" <0pt>
\ar @{-} "dlDR";"drDR" <0pt>
\ar @{-} "drDR";"urDR" <0pt>
\ar @{=>} "DL1";"u_1" <0pt>
\ar @{=>} "DL1";"u_2" <0pt>
\ar @{=>} "b_1";"DL2" <0pt>
\ar @{=>} "DR1";"u_2" <0pt>
\ar @{=>} "DR1";"u_1" <0pt>
\ar @{=>} "b_2";"DR2" <0pt>
\endxy}
\Ea
$$
Let $G_{k;2,2}\subset \cG_{k;2,2}^{(3)}$ be the subset of graphs of this special form with $k$ internal vertices. It is clear that
$$
\sum_{k\geq 0}
\sum_{\Ga \in G_{k;2,2}} \left(\int_{\p_{III}C(\Ga)}\Omega_\Ga\right)=
\frac{\Ba{c}\begin{xy}
 <0mm,0.66mm>*{};<0mm,3mm>*{}**@{.},
 <0.39mm,-0.39mm>*{};<2.2mm,-2.2mm>*{}**@{.},
 <-0.35mm,-0.35mm>*{};<-2.2mm,-2.2mm>*{}**@{.},
 <0mm,0mm>*{\circ};<0mm,0mm>*{}**@{},
\end{xy}\Ea \Ba{c}\begin{xy}
 <0mm,0.66mm>*{};<0mm,3mm>*{}**@{.},
 <0.39mm,-0.39mm>*{};<2.2mm,-2.2mm>*{}**@{.},
 <-0.35mm,-0.35mm>*{};<-2.2mm,-2.2mm>*{}**@{.},
 <0mm,0mm>*{\circ};<0mm,0mm>*{}**@{},
\end{xy}\Ea}
{f^q\left(\Ba{c}  \begin{xy}
 <0mm,-0.55mm>*{};<0mm,-2.5mm>*{}**@{.},
 <0.5mm,0.5mm>*{};<2.2mm,2.2mm>*{}**@{.},
 <-0.48mm,0.48mm>*{};<-2.2mm,2.2mm>*{}**@{.},
 <0mm,0mm>*{\circ};<0mm,0mm>*{}**@{},
 \end{xy}\Ea\right)  f^q\left( \Ba{c}  \begin{xy}
 <0mm,-0.55mm>*{};<0mm,-2.5mm>*{}**@{.},
 <0.5mm,0.5mm>*{};<2.2mm,2.2mm>*{}**@{.},
 <-0.48mm,0.48mm>*{};<-2.2mm,2.2mm>*{}**@{.},
 <0mm,0mm>*{\circ};<0mm,0mm>*{}**@{},
 \end{xy}\Ea\right)}
$$

Consider next a more general case $J\neq \emptyset$. Let $J_1\subset J$
(resp., $J_2\subset J$) be the subset of vertices which can be connected by a directed path of edges to the out-vertex $y_1^0$ (resp., $y_2^0$). At least one of the sets $J_1$ and $J_2$ is non-empty. It is easy to see that for $\Omega_\Ga|_{\p_{II}\overline{C}_{k;2,2}(\caH)}$ not to vanish identically, the functions $y_{j_1}(\var)$ and $y_{j_2}(\var)$ in the formulae above must
be of the form as $\var\rar 0$,
$$
y_{j_1}(\var)=\frac{\by_1^0}{\var} + \by_{j_1}, \ \ \ \ \
y_{j_2}(\var)=\frac{\by_2^0}{\var} + \by_{j_2}, \ \ \ \forall j_1\in J_1,\ \forall j_2\in J_2,
$$
for some constants $\by_{j_1}$ and $\by_{j_2}$. In particular,
$J_1\cap J_2=\emptyset$, so that for $\Omega_\Ga|_{\p_{II}\overline{C}_{k;2,2}(\caH)}$ not to vanish identically, the graph $\Ga$ must be of the form
$$
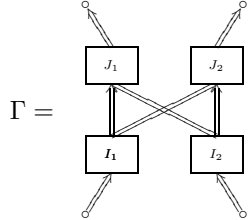
\Ga=
\Ba{c}\resizebox{23mm}{!}{
\xy
(-10,28.5)*{_{J_1}},
(-15,32)*{}="ulUL",
(-5,32)*{}="urUL",
(-15,25)*{}="dlUL",
(-5,25)*{}="drUL",
(-10,11.5)*{_{I_1}},
(-15,15)*{}="ulDL",
(-5,15)*{}="urDL",
(-15,8)*{}="dlDL",
(-5,8)*{}="drDL",
(10,28.5)*{_{J_2}},
(15,32)*{}="ulUR",
(5,32)*{}="urUR",
(15,25)*{}="dlUR",
(5,25)*{}="drUR",
(-10,11.5)*{_{I_1}},
(-15,15)*{}="ulDL",
(-5,15)*{}="urDL",
(-15,8)*{}="dlDL",
(-5,8)*{}="drDL",
(10,11.5)*{_{I_2}},
(15,15)*{}="ulDR",
(5,15)*{}="urDR",
(15,8)*{}="dlDR",
(5,8)*{}="drDR",
(-10,32)*{}="UL1",
(-10,25)*{}="UL2",
(-10,15)*{}="DL1",
(-10,8)*{}="DL2",
(10,32)*{}="UR1",
(10,25)*{}="UR2",
(10,15)*{}="DR1",
(10,8)*{}="DR2",
(-15,0)*{\circ}="b_1",
(15,0)*{\circ}="b_2",
(-15,40)*{\circ}="u_1",
(15,40)*{\circ}="u_2",
\ar @{-} "ulUL";"urUL" <0pt>
\ar @{-} "ulUL";"dlUL" <0pt>
\ar @{-} "dlUL";"drUL" <0pt>
\ar @{-} "drUL";"urUL" <0pt>
\ar @{-} "ulDL";"urDL" <0pt>
\ar @{-} "ulDL";"dlDL" <0pt>
\ar @{-} "dlDL";"drDL" <0pt>
\ar @{-} "drDL";"urDL" <0pt>
\ar @{-} "ulUR";"urUR" <0pt>
\ar @{-} "ulUR";"dlUR" <0pt>
\ar @{-} "dlUR";"drUR" <0pt>
\ar @{-} "drUR";"urUR" <0pt>
\ar @{-} "ulDR";"urDR" <0pt>
\ar @{-} "ulDR";"dlDR" <0pt>
\ar @{-} "dlDR";"drDR" <0pt>
\ar @{-} "drDR";"urDR" <0pt>
\ar @{=>} "UL1";"u_1" <0pt>
%\ar @{=>} "DR1";"u_2" <0pt>
\ar @{=>} "DL1";"UL2" <0pt>
\ar @{=>} "b_1";"DL2" <0pt>
%\ar @{=>} "b_2";"UR2" <0pt>
%
%\ar @{=>} "DL1";"u_1" <0pt>
\ar @{=>} "UR1";"u_2" <0pt>
\ar @{=>} "DR1";"UR2" <0pt>
%\ar @{=>} "b_1";"UL2" <0pt>
\ar @{=>} "b_2";"DR2" <0pt>
\ar @{=>} "DR1";"UL2" <0pt>
\ar @{=>} "DL1";"UR2" <0pt>
\endxy}
\Ea
$$
where some edges ingoing into a box can continue as outgoing edges without ``hitting"
an internal vertex inside the box (note that some of sets $I_1$, $I_2$, $J_1$ and $J_i$ can be empty!).

\mip

If $\Ga$ is a disjoint union of two  graphs, say $\Ga_1$ and $\Ga_2$,
from $\cG^{or}_{n;1,1}$, i.e.\ if it has one of the following two structures,
$$
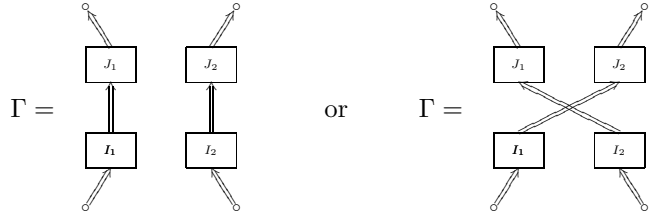
\Ga=
\Ba{c}\resizebox{22mm}{!}{
\xy
(-10,28.5)*{_{J_1}},
(-15,32)*{}="ulUL",
(-5,32)*{}="urUL",
(-15,25)*{}="dlUL",
(-5,25)*{}="drUL",
(-10,11.5)*{_{I_1}},
(-15,15)*{}="ulDL",
(-5,15)*{}="urDL",
(-15,8)*{}="dlDL",
(-5,8)*{}="drDL",
(10,28.5)*{_{J_2}},
(15,32)*{}="ulUR",
(5,32)*{}="urUR",
(15,25)*{}="dlUR",
(5,25)*{}="drUR",
(-10,11.5)*{_{I_1}},
(-15,15)*{}="ulDL",
(-5,15)*{}="urDL",
(-15,8)*{}="dlDL",
(-5,8)*{}="drDL",
(10,11.5)*{_{I_2}},
(15,15)*{}="ulDR",
(5,15)*{}="urDR",
(15,8)*{}="dlDR",
(5,8)*{}="drDR",
(-10,32)*{}="UL1",
(-10,25)*{}="UL2",
(-10,15)*{}="DL1",
(-10,8)*{}="DL2",
(10,32)*{}="UR1",
(10,25)*{}="UR2",
(10,15)*{}="DR1",
(10,8)*{}="DR2",
(-15,0)*{\circ}="b_1",
(15,0)*{\circ}="b_2",
(-15,40)*{\circ}="u_1",
(15,40)*{\circ}="u_2",
\ar @{-} "ulUL";"urUL" <0pt>
\ar @{-} "ulUL";"dlUL" <0pt>
\ar @{-} "dlUL";"drUL" <0pt>
\ar @{-} "drUL";"urUL" <0pt>
\ar @{-} "ulDL";"urDL" <0pt>
\ar @{-} "ulDL";"dlDL" <0pt>
\ar @{-} "dlDL";"drDL" <0pt>
\ar @{-} "drDL";"urDL" <0pt>
\ar @{-} "ulUR";"urUR" <0pt>
\ar @{-} "ulUR";"dlUR" <0pt>
\ar @{-} "dlUR";"drUR" <0pt>
\ar @{-} "drUR";"urUR" <0pt>
\ar @{-} "ulDR";"urDR" <0pt>
\ar @{-} "ulDR";"dlDR" <0pt>
\ar @{-} "dlDR";"drDR" <0pt>
\ar @{-} "drDR";"urDR" <0pt>
\ar @{=>} "UL1";"u_1" <0pt>
%\ar @{=>} "UL1";"u_2" <0pt>
\ar @{=>} "DL1";"UL2" <0pt>
\ar @{=>} "b_1";"DL2" <0pt>
%\ar @{=>} "b_2";"DL2" <0pt>
%
%\ar @{=>} "UR1";"u_1" <0pt>
\ar @{=>} "UR1";"u_2" <0pt>
\ar @{=>} "DR1";"UR2" <0pt>
%\ar @{=>} "b_1";"DR2" <0pt>
\ar @{=>} "b_2";"DR2" <0pt>
%\ar @{=>} "DR1";"UL2" <0pt>
%\ar @{=>} "DL1";"UR2" <0pt>
\endxy}
\Ea
\ \ \ \ \ \ \ \ \mbox{or}\ \ \ \ \ \ \ \
\Ga=
\Ba{c}\resizebox{22mm}{!}{
\xy
(-10,28.5)*{_{J_1}},
(-15,32)*{}="ulUL",
(-5,32)*{}="urUL",
(-15,25)*{}="dlUL",
(-5,25)*{}="drUL",
(-10,11.5)*{_{I_1}},
(-15,15)*{}="ulDL",
(-5,15)*{}="urDL",
(-15,8)*{}="dlDL",
(-5,8)*{}="drDL",
(10,28.5)*{_{J_2}},
(15,32)*{}="ulUR",
(5,32)*{}="urUR",
(15,25)*{}="dlUR",
(5,25)*{}="drUR",
(-10,11.5)*{_{I_1}},
(-15,15)*{}="ulDL",
(-5,15)*{}="urDL",
(-15,8)*{}="dlDL",
(-5,8)*{}="drDL",
(10,11.5)*{_{I_2}},
(15,15)*{}="ulDR",
(5,15)*{}="urDR",
(15,8)*{}="dlDR",
(5,8)*{}="drDR",
(-10,32)*{}="UL1",
(-10,25)*{}="UL2",
(-10,15)*{}="DL1",
(-10,8)*{}="DL2",
(10,32)*{}="UR1",
(10,25)*{}="UR2",
(10,15)*{}="DR1",
(10,8)*{}="DR2",
(-15,0)*{\circ}="b_1",
(15,0)*{\circ}="b_2",
(-15,40)*{\circ}="u_1",
(15,40)*{\circ}="u_2",
\ar @{-} "ulUL";"urUL" <0pt>
\ar @{-} "ulUL";"dlUL" <0pt>
\ar @{-} "dlUL";"drUL" <0pt>
\ar @{-} "drUL";"urUL" <0pt>
\ar @{-} "ulDL";"urDL" <0pt>
\ar @{-} "ulDL";"dlDL" <0pt>
\ar @{-} "dlDL";"drDL" <0pt>
\ar @{-} "drDL";"urDL" <0pt>
\ar @{-} "ulUR";"urUR" <0pt>
\ar @{-} "ulUR";"dlUR" <0pt>
\ar @{-} "dlUR";"drUR" <0pt>
\ar @{-} "drUR";"urUR" <0pt>
\ar @{-} "ulDR";"urDR" <0pt>
\ar @{-} "ulDR";"dlDR" <0pt>
\ar @{-} "dlDR";"drDR" <0pt>
\ar @{-} "drDR";"urDR" <0pt>
\ar @{=>} "UL1";"u_1" <0pt>
%\ar @{=>} "UL1";"u_2" <0pt>
%\ar @{=>} "DL1";"UL2" <0pt>
\ar @{=>} "b_1";"DL2" <0pt>
%\ar @{=>} "b_2";"DL2" <0pt>
%
%\ar @{=>} "UR1";"u_1" <0pt>
\ar @{=>} "UR1";"u_2" <0pt>
%\ar @{=>} "DR1";"UR2" <0pt>
%\ar @{=>} "b_1";"DR2" <0pt>
\ar @{=>} "b_2";"DR2" <0pt>
\ar @{=>} "DR1";"UL2" <0pt>
\ar @{=>} "DL1";"UR2" <0pt>
\endxy}
\Ea
$$
then  $\Omega_{\Ga}|_{\p\overline{C}_{k;2,2}(\caH)}=0$  because of the following

\sip

{\sc Claim.} {\em For any $\Ga\in \cG_{n;1,1}$ the associated integral
$$
\int_{C_{n;1,1}(\caH)}\Omega_{\Ga}
$$
vanishes.} Indeed, let $l'$ be the number of in-legs of $\Ga$,  $l''$  the number of out-legs, and $k$ the number of internal edges. The integral $\int_{C_{n;1,1}(\caH)}\Omega_{\Ga}$ can be non-zero if and only if  $\Omega_\Ga$ has top degree, i.e.\ if and only if
$$
3n-3+2= 2k+ l'+l''
$$
On the other hand, as every internal vertex of $\Ga$ is at least trivalent, one must have
$$
2k+l+l''\geq 3n
$$
These two equations are incompatible which proves the {\sc Claim}.

\mip

Combining all the above observations, we conclude $\Omega_{\Ga}|_{\p_{III}\overline{C}_{k;2,2}(\caH)}$ may not vanish identically only on the boundary strata
of the form
$$
\p_{I_1,I_2,J_1,J_2}\overline{C}_{k;2,2}(\caH):=\overline{C}_{\# I_1;2,1} \times
\overline{C}_{\# I_2;2,1} \times \overline{C}_{\# I_2;1,2} \times \overline{C}_{\# I_2;1,2}
$$
and
$$
\Omega_{\Ga}|_{\p_{I_1,I_2,J_1,J_2} \overline{\fC}(\Ga)}=\Omega_{\Ga_{I_1}}\wedge \Omega_{\Ga_{I_2}}\wedge \Omega_{\Ga_{J_1}}\wedge \Omega_{\Ga_{J_2}}
$$
where the graphs  $\Ga_{I_i}$ and $\Ga_{J_i}$ , $i=1,2$, are given by,
$$
\Ga_{I_i}=\Ba{c}\resizebox{14mm}{!}{
\xy
(0,11.5)*{_{I_i}},
(-5,15)*{}="ulDL",
(5,15)*{}="urDL",
(-5,8)*{}="dlDL",
(5,8)*{}="drDL",
(0,15)*{}="DL1",
(0,8)*{}="DL2",
(0,0)*{\circ}="b_1",
(-8,25)*{\circ}="u_1",
(8,25)*{\circ}="u_2",
\ar @{-} "ulDL";"urDL" <0pt>
\ar @{-} "ulDL";"dlDL" <0pt>
\ar @{-} "dlDL";"drDL" <0pt>
\ar @{-} "drDL";"urDL" <0pt>
\ar @{=>} "DL1";"u_1" <0pt>
\ar @{=>} "DL1";"u_2" <0pt>
\ar @{=>} "b_1";"DL2" <0pt>
\endxy}
\Ea\in \cG^{(3)}_{\# I_i;2,1}, \ \ \ \ \ \ \ \ \ \
\Ga_{J_i}=\Ba{c}\resizebox{14mm}{!}{
\xy
(0,11.5)*{_{J_i}},
(-5,15)*{}="ulDL",
(5,15)*{}="urDL",
(-5,8)*{}="dlDL",
(5,8)*{}="drDL",
(0,15)*{}="DL1",
(0,8)*{}="DL2",
(0,25)*{\circ}="u_1",
(-8,0)*{\circ}="b_1",
(8,0)*{\circ}="b_2",
\ar @{-} "ulDL";"urDL" <0pt>
\ar @{-} "ulDL";"dlDL" <0pt>
\ar @{-} "dlDL";"drDL" <0pt>
\ar @{-} "drDL";"urDL" <0pt>
\ar @{=>} "DL1";"u_1" <0pt>
\ar @{=>} "b_2";"DL2" <0pt>
\ar @{=>} "b_1";"DL2" <0pt>
\endxy}
\Ea\in \cG^{(3)}_{\# J_i;1,2}.
$$
Note that if $I_i$, respectively $J_1$, is empty, then we have to set
$$
\Ga_{I_i}= \Ba{c}\resizebox{9mm}{!}{ \xy
 (0,2)*{\circ}="a",
(-3,7)*{\circ}="b_1",
(3,7)*{\circ}="b_2",
 \endxy}\Ea \ \ \ , \ \ \ \mbox{respectively}\ \
 \Ga_{J_i}= \Ba{c}\resizebox{9mm}{!}{ \xy
 (0,7)*{\circ}="a",
(-3,2)*{\circ}="b_1",
(3,2)*{\circ}="b_2",
 \endxy}\Ea
$$
and $\Omega_{\Ga_{I_i}}=1$, resp.\  $\Omega_{\Ga_{J_i}}=1$.
Therefore we conclude that
$$
\sum_{k\geq 0}\sum_{\Ga\in \cG_{k;2,2}^{(3)}}\left(\sum_{V_{int}(\Ga)=I_1\sqcup I_2\sqcup J_1\sqcup J_2\atop
|I_1|+|I_2|\geq 1, |J_1|+|J_2|\geq 1} \int_{\p_{I_1,I_2,J_1,J_2} \overline{\fC}(\Ga)} \Omega_\Ga\right) \Ga  =\frac{f^{q}\left(\begin{xy}
 <0mm,0.66mm>*{};<0mm,3mm>*{}**@{.},
 <0.39mm,-0.39mm>*{};<2.2mm,-2.2mm>*{}**@{.},
 <-0.35mm,-0.35mm>*{};<-2.2mm,-2.2mm>*{}**@{.},
 <0mm,0mm>*{\circ};<0mm,0mm>*{}**@{},
\end{xy}\right) f^{q}\left(\begin{xy}
 <0mm,0.66mm>*{};<0mm,3mm>*{}**@{.},
 <0.39mm,-0.39mm>*{};<2.2mm,-2.2mm>*{}**@{.},
 <-0.35mm,-0.35mm>*{};<-2.2mm,-2.2mm>*{}**@{.},
 <0mm,0mm>*{\circ};<0mm,0mm>*{}**@{},
\end{xy}\right)}{f^{q}\left(\begin{xy}
 <0mm,-0.55mm>*{};<0mm,-2.5mm>*{}**@{.},
 <0.5mm,0.5mm>*{};<2.2mm,2.2mm>*{}**@{.},
 <-0.48mm,0.48mm>*{};<-2.2mm,2.2mm>*{}**@{.},
 <0mm,0mm>*{\circ};<0mm,0mm>*{}**@{},
 \end{xy}\right) f^{q}\left(\begin{xy}
 <0mm,-0.55mm>*{};<0mm,-2.5mm>*{}**@{.},
 <0.5mm,0.5mm>*{};<2.2mm,2.2mm>*{}**@{.},
 <-0.48mm,0.48mm>*{};<-2.2mm,2.2mm>*{}**@{.},
 <0mm,0mm>*{\circ};<0mm,0mm>*{}**@{},
 \end{xy}\right)}=
 f^{ex}\left(\Ba{c}\resizebox{15mm}{!}{\begin{xy}
 <0mm,0mm>*{\circ};<0mm,0mm>*{}**@{},
 <0mm,-0.49mm>*{};<0mm,-3.0mm>*{}**@{.},
 <-0.5mm,0.5mm>*{};<-3mm,2mm>*{}**@{.},
 <-3mm,2mm>*{};<0mm,4mm>*{}**@{.},
 <0mm,4mm>*{\circ};<-2.3mm,2.3mm>*{}**@{},
 <0mm,4mm>*{};<0mm,7.4mm>*{}**@{.},
<0mm,0mm>*{};<2.2mm,1.5mm>*{}**@{.},
 <6mm,0mm>*{\circ};<0mm,0mm>*{}**@{},
 <6mm,4mm>*{};<3.8mm,2.5mm>*{}**@{.},
 <6mm,4mm>*{};<6mm,7.4mm>*{}**@{.},
 <6mm,4mm>*{\circ};<-2.3mm,2.3mm>*{}**@{},
 <0mm,4mm>*{};<6mm,0mm>*{}**@{.},
<6mm,4mm>*{};<9mm,2mm>*{}**@{.},
<6mm,0mm>*{};<9mm,2mm>*{}**@{.},
<6mm,0mm>*{};<6mm,-3mm>*{}**@{.},
 \end{xy}}\Ea\right)
$$
where the middle expression means the fraction type composition in the prop $\LB^{\mathrm{quant}}$.
Finally, we conclude
\Beqrn
0&=&\sum_{k\geq 0}\sum_{\Ga\in \cG_{k;2,2}^{(3)}}\left( \int_{\p \overline{C}_{k;2,2}(\caH)} \Omega_\Ga\right)  \Ga\\
 &=& \sum_{k\geq 0}\sum_{\Ga\in \cG_{k;2,2}^{(3)}}\left( \int_{\p_{II} \overline{C}_{k;2,2}(\caH)} \Omega_\Ga\right)  \Ga +
 \sum_{k\geq 0}\sum_{\Ga\in \cG_{k;2,2}^{(3)}}\left( \int_{\p_{III} \overline{C}_{k;2,2}(\caH)} \Omega_\Ga\right)  \Ga\\
  &=& \sum_{k\geq 0}\sum_{\Ga\in \cG_{k;2,2}^{(3)}}\left( \int_{\p_{II} \overline{C}_{k;2,2}(\caH)} \Omega_\Ga\right)  +
   \sum_{k\geq 0}\sum_{\Ga\in \cG_{k;2,2}^{(3)}}\sum_{V_{int}(\Ga)=I_1\sqcup I_2\sqcup J_1\sqcup J_2\atop
|I_1|+|I_2|\geq 1, |J_1|+|J_2|\geq 1} \left(\int_{\p_{I_1,I_2,J_1,J_2} \overline{\fC}(\Ga)}\Omega_\Ga\right)
\Ga\\
&=& f^{ex}\left( -\Ba{c} \begin{xy}
 <0mm,2.47mm>*{};<0mm,-0.5mm>*{}**@{.},
 <0.5mm,3.5mm>*{};<2.2mm,5.2mm>*{}**@{.},
 <-0.48mm,3.48mm>*{};<-2.2mm,5.2mm>*{}**@{.},
 <0mm,3mm>*{\circ};<0mm,3mm>*{}**@{},
  <0mm,-0.8mm>*{\circ};<0mm,-0.8mm>*{}**@{},
<0mm,-0.8mm>*{};<-2.2mm,-3.5mm>*{}**@{.},
 <0mm,-0.8mm>*{};<2.2mm,-3.5mm>*{}**@{.},
     <0.5mm,3.5mm>*{};<2.8mm,5.7mm>*{^2}**@{},
     <-0.48mm,3.48mm>*{};<-2.8mm,5.7mm>*{^1}**@{},
   <0mm,-0.8mm>*{};<-2.7mm,-5.2mm>*{^1}**@{},
   <0mm,-0.8mm>*{};<2.7mm,-5.2mm>*{^2}**@{},
\end{xy}\Ea
\ + \
\Ba{c}\begin{xy}
 <0mm,0mm>*{\circ};<0mm,0mm>*{}**@{},
 <0mm,-0.49mm>*{};<0mm,-3.0mm>*{}**@{.},
 <-0.5mm,0.5mm>*{};<-3mm,2mm>*{}**@{.},
 <-3mm,2mm>*{};<0mm,4mm>*{}**@{.},
 <0mm,4mm>*{\circ};<-2.3mm,2.3mm>*{}**@{},
 <0mm,4mm>*{};<0mm,7.4mm>*{}**@{.},
<0mm,0mm>*{};<2.2mm,1.5mm>*{}**@{.},
 <6mm,0mm>*{\circ};<0mm,0mm>*{}**@{},
 <6mm,4mm>*{};<3.8mm,2.5mm>*{}**@{.},
 <6mm,4mm>*{};<6mm,7.4mm>*{}**@{.},
 <6mm,4mm>*{\circ};<-2.3mm,2.3mm>*{}**@{},
 <0mm,4mm>*{};<6mm,0mm>*{}**@{.},
<6mm,4mm>*{};<9mm,2mm>*{}**@{.},
<6mm,0mm>*{};<9mm,2mm>*{}**@{.},
<6mm,0mm>*{};<6mm,-3mm>*{}**@{.},
   <-1.8mm,2.8mm>*{};<0mm,7.8mm>*{^1}**@{},
   <-2.8mm,2.9mm>*{};<0mm,-4.3mm>*{_1}**@{},
<-1.8mm,2.8mm>*{};<6mm,7.8mm>*{^2}**@{},
   <-2.8mm,2.9mm>*{};<6mm,-4.3mm>*{_2}**@{},
 \end{xy}
\Ea   \right)
\Eeqrn
which proves claim (iii).
\end{proof}

\subsubsection{\bf Main Corollary} {\em Composition of the explicit morphism (\ref{7: explicit morphism f^q}) with the explicit morphism $\caD(f)$  (see \S {\ref{5: coroll on f from LB^q to LB^min}}(ii)) gives us an explicit transcendental morphism of props
\Beq\label{7: explicit map from Assb to LB wheeled}
\caD(f) \circ f^q: \Assb   \lon \wLB^\circlearrowright
\Eeq
and hence an explicit universal quantization of finite-dimensional Lie bialgebras.}
 The main purpose of this paper is achieved.

\subsubsection{\bf Other Corollaries}

(i) As the differential 2-forms $\om_g$ and $\varpi_g$ used in the constructions of the maps $f^q$ and $f$ are simple, {\em graphs with multiple edges do not contribute into the map (\ref{7: explicit map from Assb to LB wheeled})}. Essentially this observation  says that {\em our universal quantization formula does not involve graphs which contain a subgraph of the form}\,
$\Ba{c}\xy
(0,0)*{}="0",
 (0,3)*{\bu}="1",
(-3,5)*{}="L",
(3,5)*{}="R",
(0,7)*{\bu}="2",
(0,10)*{}="00",
\ar @{-} "0";"1" <0pt>
\ar @{-} "1";"L" <0pt>
\ar @{-} "1";"R" <0pt>
\ar @{-} "2";"L" <0pt>
\ar @{-} "2";"R" <0pt>
\ar @{-} "2";"00" <0pt>
\endxy\Ea$. It also follows from our explicit formula that all graphs
with at least one black vertex contributing to the universal quantization morphism are {\em connected}.

\sip

(ii) {\em The explicit map (\ref{7: explicit morphism f^q}) lifts by a trivial induction to a morphism of dg props $\cF^q$ which
 fits into a  commutative diagram,
$$
 \xymatrix{
\Assb_\infty\ar[r]^{\cF^q}\ar[d]_p  & \caD\wLB^{\mathrm{quant}}_\infty\ar[d]^\pi\\
 \Assb \ar[r]_{f^q} &
 \caD\wLB^{\mathrm{quant}}}
$$
and which satisfies the condition
$$
\pi_1\circ \cF^q\left(\Ba{c}\resizebox{13mm}{!}{ \xy
 (0,7)*{\overbrace{\ \ \ \  \ \ \ \ \ \ \ \ \ \ }},
 (0,9)*{^m},
 (0,3)*{^{...}},
 (0,-3)*{_{...}},
 (0,-7)*{\underbrace{  \ \ \ \  \ \ \ \ \ \ \ \ \ \ }},
 (0,-9)*{_n},
 (0,0)*{\circ}="0",
(-7,5)*{}="u_1",
(-4,5)*{}="u_2",
(4,5)*{}="u_3",
(7,5)*{}="u_4",
(-7,-5)*{}="d_1",
(-4,-5)*{}="d_2",
(4,-5)*{}="d_3",
(7,-5)*{}="d_4",
\ar @{.} "0";"u_1" <0pt>
\ar @{.} "0";"u_2" <0pt>
\ar @{.} "0";"u_3" <0pt>
\ar @{.} "0";"u_4" <0pt>
\ar @{.} "0";"d_1" <0pt>
\ar @{.} "0";"d_2" <0pt>
\ar @{.} "0";"d_3" <0pt>
\ar @{.} "0";"d_4" <0pt>
\endxy}\Ea\right)=
\la \Ba{c}\resizebox{16mm}{!}{\xy
(0,7.5)*{\overbrace{\ \ \ \ \ \ \ \ \  \ \ \ \ \ \ \ \ \ \ }},
 (0,9.5)*{^m},
 (0,-7.5)*{\underbrace{\ \ \ \ \ \ \ \ \ \ \ \  \ \ \ \ \ \ }},
 (0,-9.9)*{_n},
  (-6,5)*{...},
  (-6,-5)*{...},
 (-3,5)*{\circ}="u1",
  (-3,-5)*{\circ}="d1",
  (-6,5)*{...},
  (-6,-5)*{...},
  (-9,5)*{\circ}="u2",
  (-9,-5)*{\circ}="d2",
 (3,5)*{\circ}="u3",
  (3,-5)*{\circ}="d3",
  (6,5)*{...},
  (6,-5)*{...},
  (9,5)*{\circ}="u4",
  (9,-5)*{\circ}="d4",
 (0,0)*{\bullet}="a",
\ar @{-} "d1";"a" <0pt>
\ar @{-} "a";"u1" <0pt>
\ar @{-} "d2";"a" <0pt>
\ar @{-} "a";"u2" <0pt>
\ar @{-} "d3";"a" <0pt>
\ar @{-} "a";"u3" <0pt>
\ar @{-} "d4";"a" <0pt>
\ar @{-} "a";"u4" <0pt>
\endxy}\Ea\ \ \text{for some non-zero}\ \la\in \R,
$$
for all $m+n\geq 3$, $m,n\geq 1$}. Here $\pi_1$ is the projection to the vector subspace in $\caD\wLB^{\mathrm{quant}}_\infty$ spanned by graphs with precisely one black vertex.

This claim is obvious as surjections $p$ and $\pi^q$ are quasi-isomorphisms.

\sip

(iii) Composition of the maps $\cF^q$ and $\caD(F)$, where $F$ is given by the explicit formula (\ref{5: explicit map F from LB^q_infty to LB_infty wheeled}), gives us a formality map}%\Beq\label{7: explicit formality map from Assb infty to LB wheeled}
$$
\caD(F) \circ \cF^q: \Assb_\infty   \lon \caD\wLB^\circlearrowright_\infty
$$
and hence {\em a universal quantization of finite-dimensional strongly homotopy Lie bialgebras}.

\subsection{An open problem} The above Corollary(ii) gives us an inductive extension of the explicit morphism (\ref{7: explicit morphism f^q}) to some
morphism of dg props $\cF^q: \Assb_\infty\rar \wLB_\infty^{\mathrm{quant}}$.
Can this extension be given by an explicit formula similar to the one for $f^q$?
Here is a conjectural answer.

\subsubsection{\bf Conjecture}\label{7: Conjecture on F^q} {\em There is a morphism of props
\Beq\label{7: explicit morphism F^q}
\cF^{q}: \cA ss\cB_\infty \lon \caD\widehat{\LB}^{\mathrm{quant}}_\infty
\Eeq
given explicitly on the generators of $\cA ss\cB_\infty$ as follows,
$$
\cF^q\left(\Ba{c}\resizebox{13mm}{!}{ \xy
 (0,7)*{\overbrace{\ \ \ \  \ \ \ \ \ \ \ \ \ \ }},
 (0,9)*{^m},
 (0,3)*{^{...}},
 (0,-3)*{_{...}},
 (0,-7)*{\underbrace{  \ \ \ \  \ \ \ \ \ \ \ \ \ \ }},
 (0,-9)*{_n},
 (0,0)*{\circ}="0",
(-7,5)*{}="u_1",
(-4,5)*{}="u_2",
(4,5)*{}="u_3",
(7,5)*{}="u_4",
(-7,-5)*{}="d_1",
(-4,-5)*{}="d_2",
(4,-5)*{}="d_3",
(7,-5)*{}="d_4",
\ar @{.} "0";"u_1" <0pt>
\ar @{.} "0";"u_2" <0pt>
\ar @{.} "0";"u_3" <0pt>
\ar @{.} "0";"u_4" <0pt>
\ar @{.} "0";"d_1" <0pt>
\ar @{.} "0";"d_2" <0pt>
\ar @{.} "0";"d_3" <0pt>
\ar @{.} "0";"d_4" <0pt>
\endxy}\Ea\right)
:=
\sum_{k\geq 1}
\sum_{\Ga \in \cG_{k;m,n}^{or}} \left(\int_{\overline{C}_{k;m,n}(\caH)}\Omega_\Ga\right) \Ga \ \ \ +
\ \ \left\{\Ba{cl} \Ba{c}\resizebox{6mm}{!}{ \xy
(-3,9)*{_1},
(3,9)*{_2},
 (0,0)*{\circ}="a",
(-3,5)*{\circ}="b_1",
(3,5)*{\circ}="b_2",
 \endxy}\Ea& \text{if}\ m=2,n=1\vspace{3mm}\\
  \Ba{c}\resizebox{6mm}{!}{ \xy
(-3,0)*{_1},
(3,0)*{_2},
 (0,6)*{\circ}="a",
(-3,2)*{\circ}="b_1",
(3,2)*{\circ}="b_2",
 \endxy}\Ea & \text{if}\ m=1,n=2\\
 0 & \text{otherwise}
 \Ea\right.
$$
where  the differential form $\Omega_\Ga$ is defined in (\ref{1> Def of Omega_Ga}).}

\mip

Let us provide a strong evidence for this conjecture elucidating  a particular problem which requires a better understanding.

\sip

By construction of the compactified space $\overline{C}_{k;m,n}(\caH)$, we have a natural semialgebraic fibration (see \cite{HLTV})
$$
\pi: \overline{C}_{k;m,n}(\caH) \lon \overline{C}_{m,n}(\R\times \R)
$$
and hence a push-forward map of piecewise semi-algebraic differential forms
$$
\pi_*: \Omega_{ \overline{C}_{k;m,n}(\caH)}^\bu \lon  \Omega^\bu_{\overline{C}_{m,n}(\R\times \R)}
$$
such that for any semialgebraic chain
$$
\phi: M \rar \overline{C}_{m,n}(\R\times \R)
$$
the integral
$$
\int_M \phi^*(\pi_*(\Omega_\Ga))
$$
is well-defined (i.e.\ convergent) for any $\Ga\in \cG_{k;m,n}$. Hence we can consider an $\bS_m^{op}\times \bS_n$ equivariant map
\Beq\label{7: map Phi_n^m from chains to DLieb}
\Ba{rccc}
\Phi_{n}^m: &  Chains(\overline{C}_{m,n}) & \lon & \wqLB_\infty(m,n)\\
   & \phi: M\rar \overline{C}_{m,n}(\R\times \R) & \lon & \displaystyle  \sum_{k\geq 0}\sum_{\Ga\in
   \cG_{k;m,n}} \left(\int_M \phi^*\left(  \pi_*(\Omega_\Ga)\right) \right)\Ga
   \Ea
\Eeq
Note that in our grading conventions  the chain complex $(Chains(\overline{C}_{m,n}),\p)$ is non-positively graded
so that the standard boundary differential $\p$ has degree $+1$. Using arguments almost identical to the ones employed in the proof of Theorem {\ref{7: Theorem on f^q}} one can show the following

\subsubsection{\bf Theorem}\label{7: Theorem on Chains to Dlieb} {\em For any $m,n\geq 1$ with $m+n\geq 3$ the collection of maps $\Phi_{n}^m: Chains(\overline{C}_{m,n})\lon \caD \wqLB_\infty(m,n)$ commutes with the differentials,
$$
\delta^{\om_{\bar{g}}}\circ \Phi_n^m = \Phi_n^m \circ \p
$$
and hence gives us an equivariant morphism of differential $\frac{1}{2}$-props %$\bS$-bimodules
$$
\Phi: Chains(\overline{C}_{\bu,\bu}(\R\times\R))\rar \caD\wqLB_\infty.
$$}

\sip

The restriction of the map $\Phi$ to the Saneblidze-Umble cell complex
$(\cC ell(\sK_\bu^\bu), \p_{cell}) \subset Chains(\overline{C}_{\bu,\bu}(\R\times\R)$ (see Appendix B) gives us precisely the map $\cF^q$ in the Conjecture {\ref{7: Theorem on Chains to Dlieb}}. This map respects the differentials but  at the moment we can not claim it respects {\em all}\, prop compositions as the isomorphism $(\cC ell(\sK_\bu^\bu), \p_{cell})
\simeq \Assb_\infty$ (which is claimed in \cite{SU}) should be  understood better in this context.

\bip

\bip
\appendix
\renewcommand{\thesubsection}{{\bf A.\arabic{subsection}}}
\renewcommand{\thesubsubsection}{{\bf A.\arabic{subsection}.\arabic{subsubsection}}}

\bip

{\Large
\section{\bf Some vanishing Lemmas}\label{App: A}
}

\bip

Let $\om_g$ be a top degree form on $S^2$ given by (\ref{3: omega_g propagator}) for $d=3$.
We shall prove some vanishing results for the weights
$$
C_\Ga=\int_{\overline{C}_{4p+2}(\R^3)} \displaystyle \bigwedge_{e\in E(\Ga)}\hspace{-2mm}
{\pi}^*_e\left(\om_g\right)
$$
of graphs $\Ga\in \sG_{4p+2,6p+1}$ with $p\geq 1$ contributing to the formulae given in Proposition {\ref{3: Prop on Upsilon^om_g}}. 
%In fact, all the results hold true for an arbitrary 2-form
%$\om_g=g \mathrm{Vol}_{S^2}$ on the 2-sphere.

\mip

\subsection{Lemma on binary vertices}\label{A: lemma on 4 binary vertices} {\em Any graph $\Ga\in \cG_{4p+2,6p+1}$ with $p\geq 1$
has at least 4 binary vertices. Moreover, if $\Ga\in \cG_{4p+2,6p+1}$ has precisely $4$ binary vertices, then all other vertices must be trivalent}.

\begin{proof} For a vertex $v\in V(\Ga)$ its valency can be represented as the sum $2+ \Delta v$
for some non-negative integer $\Delta v$. The graph $\Ga$ has $12p+2$ half-edges so we have an equality
$$
\sum_{v\in V(\Ga)} (2+ \Delta v)= 2+12p,
$$
i.e.
$$
\sum_{v\in V(\Ga)} \Delta v=2+12p - 2(2+4p)=4p-2
$$
Therefore at most $4p-2$ vertices can have $\Delta v\geq 1$ which implies that $\Ga$ has at least
$4p+2 -(4p-2)=4$ binary vertices. Moreover, if $\Ga$ has precisely $4$ bivalent vertices, then
the remaining $4p-2$ vertices $v$ must have $\Delta v=1$.
\end{proof}

 Therefore every graph in $\Ga\in \sG_{4p+2,6p+1}$ with $p\geq 1$ has at least four complete\footnote{For a graph $\Ga$ and its pair of vertices $v_1,v_2\in V(\Ga)$ denote by $E_\Ga(v_1,v_2)$ the set of edges connecting $v_1$ to $v_2$. A subgraph $\Ga'$ of graph $\Ga$ is called {\em complete}\, if
  between any pair of its vertices $v_1,v_2 \in V(\Ga')$ we have $E_{\Ga'}(v_1,v_2)=E_{\Ga}(v_1,v_2)$.}   subgraphs of one of the following forms,
$$
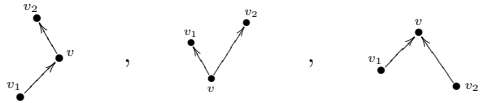
\Ba{c}\resizebox{10mm}{!}{
\xy
(-1,2)*{_{v_1}},
(2,16)*{_{v_2}},
(9,8)*{_{v}},
 (0,0)*{\bullet}="a",
(7,7)*{\bu}="b",
(3,14)*{\bu}="c",
\ar @{->} "a";"b" <0pt>
\ar @{->} "b";"c" <0pt>
\endxy}\Ea \ \ \ \ , \ \ \ \
\Ba{c}\resizebox{10mm}{!}{\xy
(7,5)*{_{v}},
(15,20)*{_{v_2}},
(3,16)*{_{v_1}},
 (14,18)*{\bullet}="a",
(7,7)*{\bu}="b",
(3,14)*{\bu}="c",
\ar @{<-} "a";"b" <0pt>
\ar @{->} "b";"c" <0pt>
\endxy}
\Ea \ \ \ \ , \ \ \ \
\Ba{c}\resizebox{15mm}{!}{\xy
(-1,2)*{_{v_1}},
(17,-3)*{_{v_2}},
(7,9)*{_{v}},
 (0,0)*{\bullet}="a",
(7,7)*{\bu}="b",
(14,-3)*{\bu}="c",
\ar @{->} "a";"b" <0pt>
\ar @{<-} "b";"c" <0pt>
\endxy}\Ea
$$
where the vertex $v$ has no other attached edges except the ones shown in the pictures.

\subsection{Vanishing Lemma}\label{4: vanishing of one bivalent} {\em If\,
$\Ga\in \sG_{4p+2,6p+1}$\, with $p\geq 1$ admits a binary vertex $v$ of the form $\Ba{c}\resizebox{10mm}{!}
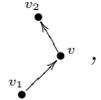
{
\xy
(-1,2)*{_{v_1}},
(2,16)*{_{v_2}},
(9,8)*{_{v}},
%%%
 (0,0)*{\bullet}="a",
(7,7)*{\bu}="b",
(3,14)*{\bu}="c",
%%%
\ar @{->} "a";"b" <0pt>
\ar @{->} "b";"c" <0pt>
\endxy}\Ea$,
then its
 weight $C_\Ga$ vanishes.}

\begin{proof}
We assume here that the propagators are chosen $O(2)$-anti-invariantly, i.e., invariantly for the $SO(2)$ action on the sphere $S^2$, and anti-invariantly for a reflection across a plane containing both poles.
Now, integrating over the position of (the point in a configuration associated to) vertex $v$, the above graph yields a $1$-form on the configuration space of $v_1$ and $v_2$, i.e., on $S^2$.
This 1-form is easily checked to be $O(2)$-anti-invariant, and furthermore closed by Stokes' Theorem. Using standard cylindrical coordinates $(Z,\phi)$ the $O(2)$-anti-invariance implies that the form can be written as
\[
f(Z)d\phi.
\]
for some function $f(Z)$, vanishing at the sphere's poles $Z=\pm 1$ to ensure continuity.
The closedness then implies that in fact $f(Z)\equiv 0$.
 \end{proof}

 \subsection{Vanishing Lemma}\label{Lemma on triangles} {\em If  $\Ga\in \sG_{4p+2,6p+1}$ admits a 3-vertex complete graph (with any possible choice of directions of edges),
$$
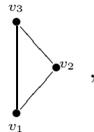
\Ba{c}\resizebox{10mm}{!}{
\xy
(9,8)*{^{v_2}},
(0,18)*{^{v_3}},
(0,-3)*{^{v_1}},
(0,0)*{\bu}="d",
(0,16)*{\bu}="u",
(7,8)*{\bu}="R",
\ar @{-} "d";"u" <0pt>
\ar @{-} "d";"R" <0pt>
\ar @{-} "R";"u" <0pt>
\endxy}
\Ea,
$$
as a subgraph,
then its weight $C_\Ga$ vanishes.}

\begin{proof} The integrand $\Omega_\Ga:=\bigwedge_{e\in E(\Ga)}\hspace{-2mm}
{\pi}^*_e\left(\om_g\right)$ is invariant
under the action of the gauge group $p\rar \R^+p + \R^3$ on points in $\R^3$.
Hence we can place vertex $v_1$ at $0\in \R^3$, and normalized the Euclidean
 distance $|v_2-v_1|$ to be equal to $1$. Then the $6$-form
 $$
 \pi_{v_1,v_2}^*(\om_g)\wedge  \pi_{v_1,v_3}^*(\om_g)\wedge \pi_{v_2,v_3}^*(\om_g)
 $$
depends only on $5$ parameters and hence vanishes identically for degree reasons. Hence the form
$\Omega_\Ga$ is zero.
\end{proof}

  \subsection{Vanishing Lemma}{\em Assume $\Ga\in \cG_{4p+2,6p+1}$ has two  bivalent vertices $v'$ and $v''$ connected by an edge. Then its weight $C_\Ga$ vanishes.}

\begin{proof} It is enough to consider the case when orientations on the subgraph containing
$v'$ and $v''$ and their neighbouring (not necessarily binary) vertices $v_1$ and $v_2$ are as in the following oriented graph,
$$
\Ga_{v_1,v',v'',v_2}:=
\Ba{c}\resizebox{18mm}{!}{
\xy
(-2,1)*{^{v_1}},
(7,9)*{^{v'}},
(15,3)*{^{v''}},
(21,9)*{^{v_2}},
(0,0)*{\bu}="1",
(7,7)*{\bu}="2",
 (14,0)*{\bu}="3",
(21,7)*{\bu}="4",
\ar @{->} "1";"2" <0pt>
\ar @{->} "3";"2" <0pt>
\ar @{->} "3";"4" <0pt>
\endxy}
\Ea\ ,
$$
for all other inequivalent choices the vanishing claim follows from Lemma~{\ref{4: vanishing of one bivalent}}, and, in the case $v_1=v_2$, from Lemma~{\ref{Lemma on triangles}}.
\bip

Let us fix all vertices of the graph except $v'$ and $v''$. We can also fix without loss of generality the vertex $v_1$ at $0\in \R^3$ and the vertex $v_2$ at the unit Euclidean distance from $v_1$. Consider a projection
\Beq\label{map pi}
\pi:  \overline{C}({\Ga_{v_1,v',v'',v_2}}) \lon { C}_{v_1,v_2}(\R^3)
\Eeq
and the function
$$
f:= \pi_*( \underbrace{\pi_{v_1,v'}^*(\om_g)\wedge  \pi_{v',v''}^*(\om_g)\wedge \pi_{v'',v_2}^*(\om_g)}_{\Omega_{\Ga_{v_1,v',v'',v_2}}} )
$$
on  $C_{v_1,v_2}(\R^3)$. By the generalized Stokes Theorem,
$$
d\circ \pi_* =\pm \pi_*\circ d + \pi_{\p *},
$$
so that we have
\Beq\label{df}
df= \pi_{\p *}\left(\Omega_{\Ga_{v_1,v',v'',v_2}}\right)=\al_*(\Omega_{\Ga_{v_1,v'',v_2}}) - \be_*(\Omega_{\Ga_{v_1,v,v_2}})
+ \ga_*(\Omega_{\Ga_{v_1,v',v_2}})
\Eeq
where
$$
\Ga_{v_1,v'',v_2}:=
\Ba{c}\resizebox{13mm}{!}{
\xy
(7,9)*{^{v_1}},
(15,3)*{^{v''}},
(21,9)*{^{v_2}},
%
%(0,0)*{\bu}="1",
(7,7)*{\bu}="2",
 (14,0)*{\bu}="3",
(21,7)*{\bu}="4",
%
%\ar @{->} "1";"2" <0pt>
\ar @{->} "3";"2" <0pt>
\ar @{->} "3";"4" <0pt>
\endxy}
\Ea\ ,\ \
\Ga_{v_1,v,v_2}:=
\Ba{c}\resizebox{10mm}{!}{
\xy
(-1,2)*{_{v_1}},
(2,16)*{_{v_2}},
(9,8)*{_{v}},
 (0,0)*{\bullet}="a",
(7,7)*{\bu}="b",
(3,14)*{\bu}="c",
\ar @{->} "a";"b" <0pt>
\ar @{->} "b";"c" <0pt>
\endxy}\Ea
\ , \ \
\Ga_{v_1,v',v_2}:=
\Ba{c}\resizebox{16mm}{!}{
\xy
(-2,1)*{^{v_1}},
(7,9)*{^{v'}},
(16,1)*{^{v_2}},
(0,0)*{\bu}="1",
(7,7)*{\bu}="2",
 (14,0)*{\bu}="3",
\ar @{->} "1";"2" <0pt>
\ar @{->} "3";"2" <0pt>
\endxy}
\Ea
$$
and
$$
\al: C(\Ga_{v_1,v'',v_2}) \rar C_{v_1,v_2}(\R^3),\ \ \be: C(\Ga_{v_1,v,v_2}) \rar C_{v_1,v_2}(\R^3), \ \ \
\ga: C(\Ga_{v_1,v',v_2}) \rar C_{v_1,v_2}(\R^3)
$$
are the natural forgetful maps. By Lemma~{\ref{4: vanishing of one bivalent}}, the middle term
$\be_*(\Omega_{\Ga_{v_1,v,v_2}})$ vanishes. On the other hand the sum,
$$
\al_*(\Omega_{\Ga_{v_1,v',v_2}})
+ \ga_*(\Omega_{\Ga_{v_1,v'',v_2}})
$$
equals the push down,
$$
p_*\left(\pi_{v_1,v}^*(\om_g)\wedge  \pi_{v,v_2}^*(\om_g)\right)
   $$
   of the $4$-form $\pi_{v_1,v}^*(\om_g)\wedge  \pi_{v,v_2}^*(\om_g)$
along the 3-dimensional fiber of the natural projection,
$$
p: C_{v_1,v,v_2}(\R^3) \lon C_{v_1,v_2}(\R^3).
$$
The latter vanishes by the standard argument using the reflection in the line through vertices $v_1$ and $v_2$ (cf.\ \cite{Ko0}).

\bip

Therefore we conclude that
$$
df=0,
$$
i.e.\ the function $f$ is a constant independent of a  particular position of the vertex $v_2$ (on the sphere). Let us choose $v_2$ to lie in the $(x,t)$-plane. Then the reflection in this plane preserves
the orientation of the fiber of the map (\ref{map pi}) but changes the differential form
$$
\Omega_{\Ga_{v_1,v',v'',v_2}} \lon - \Omega_{\Ga_{v_1,v',v'',v_2}}.
$$
Hence $f=0$ and the proof is completed.
\end{proof}

\mip

Let $\hat{\sG}_{4p+1,6p+1}^{or}$ be the subset of the set of oriented graphs $\hat{\sG}_{4p+1,6p+1}^{or}$ consisting of graphs $\Ga$ which have no
\Bi
\item  binary vertices of arity $(1,1)$, i.e.\ of the form $\Ba{c}\resizebox{2mm}{!}{
\xy
%%%
 (0,0)*{\bullet}="a",
(0,5)*{}="b",
(0,-5)*{}="c",
%%%
\ar @{->} "a";"b" <0pt>
\ar @{<-} "a";"c" <0pt>
\endxy}\Ea$
\item no complete subgraphs of the form $
\Ba{c}\resizebox{10mm}{!}{
\xy
(9,8)*{^{v_2}},
(0,18)*{^{v_3}},
(0,-3)*{^{v_1}},
(0,0)*{\bu}="d",
(0,16)*{\bu}="u",
(7,8)*{\bu}="R",
\ar @{-} "d";"u" <0pt>
\ar @{-} "d";"R" <0pt>
\ar @{-} "R";"u" <0pt>
\endxy}
\Ea,
$
\item no two binary vertices connected by an edge.
\Ei

We proved in this Appendix the following

\subsection{Proposition}\label{A: propos on hat{sG}} {\em In the case $d=3$ Proposition {\ref{3: Prop on Upsilon^om_g}} holds true with the set of graphs
$\sG_{4p+2,6p+1}^{or}$  replaced by its subset  $\hat{\sG}_{4p+1,6p+1}^{or}$.}

\mip

A quick inspection of the case $p=1$ shows that there are no graphs in $\hat{\sG}_{6,7}^{or}$ which satisfy the above three properties so that one gets the following

\subsection{Lemma}\label{A: lemma on p=1} {\em The set $\hat{\sG}_{6,7}^{or}$\, is empty}.

\mip

In the case $p=2$ one has non-trivial examples, e.g.

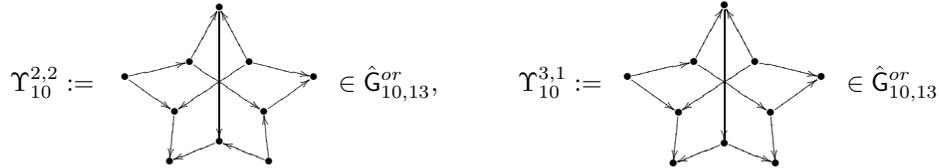
$$
\Upsilon_{10}^{2,2}:=\Ba{c}\resizebox{27mm}{!}{
\xy
(-6,0)*{\bu}="1l",
(6,0)*{\bu}="1r",
(0,11)*{\bu}="u",
(-19,-3)*{\bu}="2l",
(19,-3)*{\bu}="2r",
(-9,-10)*{\bu}="3l",
(9,-10)*{\bu}="3r",
(-10,-20)*{\bu}="4l",
(10,-20)*{\bu}="4r",
(0,-16)*{\bu}="d",
\ar @{->} "1l";"u" <0pt>
\ar @{->} "1r";"u" <0pt>
\ar @{<-} "1l";"2l" <0pt>
\ar @{->} "1r";"2r" <0pt>
\ar @{<-} "3l";"2l" <0pt>
\ar @{->} "3r";"2r" <0pt>
\ar @{->} "3l";"4l" <0pt>
\ar @{<-} "3r";"4r" <0pt>
\ar @{->} "d";"4l" <0pt>
\ar @{<-} "d";"4r" <0pt>
\ar @{->} "u";"d" <0pt>
\ar @{->} "1l";"3r" <0pt>
\ar @{->} "1r";"3l" <0pt>
\endxy}
\Ea \in  \hat{\sG}_{10,13}^{or}, \hspace{10mm}
\Upsilon_{10}^{3,1}:=\Ba{c}
\resizebox{27mm}{!}{\xy
(-6,0)*{\bu}="1l",
(6,0)*{\bu}="1r",
(0,11)*{\bu}="u",
(-19,-3)*{\bu}="2l",
(19,-3)*{\bu}="2r",
(-9,-10)*{\bu}="3l",
(9,-10)*{\bu}="3r",
(-10,-20)*{\bu}="4l",
(10,-20)*{\bu}="4r",
(0,-16)*{\bu}="d",
\ar @{->} "1l";"u" <0pt>
\ar @{->} "1r";"u" <0pt>
\ar @{<-} "1l";"2l" <0pt>
\ar @{->} "1r";"2r" <0pt>
\ar @{<-} "3l";"2l" <0pt>
\ar @{->} "3r";"2r" <0pt>
\ar @{->} "3l";"4l" <0pt>
\ar @{->} "3r";"4r" <0pt>
\ar @{->} "d";"4l" <0pt>
\ar @{->} "d";"4r" <0pt>
\ar @{->} "u";"d" <0pt>
\ar @{->} "1l";"3r" <0pt>
\ar @{->} "1r";"3l" <0pt>
\endxy}
\Ea \in  \hat{\sG}_{10,13}^{or}
$$

The first graph $\Upsilon_{10}^{2,2}$ has two binary vertices have type $(2,0)$ and two binary vertices of type $(0,2)$. The second graph $\Upsilon_{10}^{3,1}$ has three vertices of type $(2,0)$ and one vertex of type $(0,2)$. Reversing all arrows in $\Upsilon_{10}^{3,1}$ one obtains a graph
$$
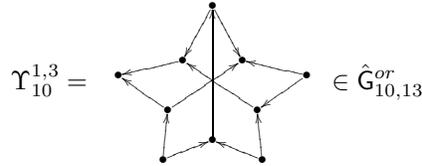
\Upsilon_{10}^{1,3}=\Ba{c}\resizebox{27mm}{!}{
\xy
(-6,0)*{\bu}="1l",
(6,0)*{\bu}="1r",
(0,11)*{\bu}="u",
(-19,-3)*{\bu}="2l",
(19,-3)*{\bu}="2r",
(-9,-10)*{\bu}="3l",
(9,-10)*{\bu}="3r",
(-10,-20)*{\bu}="4l",
(10,-20)*{\bu}="4r",
(0,-16)*{\bu}="d",
\ar @{<-} "1l";"u" <0pt>
\ar @{<-} "1r";"u" <0pt>
\ar @{->} "1l";"2l" <0pt>
\ar @{<-} "1r";"2r" <0pt>
\ar @{->} "3l";"2l" <0pt>
\ar @{<-} "3r";"2r" <0pt>
\ar @{<-} "3l";"4l" <0pt>
\ar @{<-} "3r";"4r" <0pt>
\ar @{<-} "d";"4l" <0pt>
\ar @{<-} "d";"4r" <0pt>
\ar @{<-} "u";"d" <0pt>
\ar @{<-} "1l";"3r" <0pt>
\ar @{<-} "1r";"3l" <0pt>
\endxy}
\Ea \in  \hat{\sG}_{10,13}^{or}
$$
 with three vertices of type $(0,2)$ and one vertex of type $(2,0)$.

\bip

\bip

{\Large
\section{\bf Configuration space models for bipermutahedra\\ and biassociahedra}
}\label{App B}

\bip

\subsection{Associahedron, permutahedron and configuration spaces} Here we remind two well-known constructions \cite{St,Ko, LTV} (see also lecture notes \cite{Me1}) which will be used later.
 Let
$$
\Conf_n(\R):=\{[n]\hook \R\},
$$
be the space of all possible injections of the set $[n]:=\{1,2,\ldots, n\}$ into the real line $\R$.
This space is a disjoint union of $n!$ connected components each of which is isomorphic
 to the space $$
 \Conf_n^{o}(\R)=\{x_{1}< x_{2} <\ldots < x_{n}\}.
 $$
  The set $\Conf_n(\R)$ has a natural structure of an oriented $n$-dimensional manifold
 with orientation on $\Conf_n^{0}(\R)$ given by the volume form $dx_1\wedge dx_2\wedge\ldots\wedge dx_n$;
  orientations of all other connected components are then fixed once we assume that the natural smooth action of $\bS_n$  on $\Conf_n(\R)$ is orientation preserving.
 In fact, we can (and often do) label points by an arbitrary finite set $I$, that is, consider the space
 of injections of sets,
$$
\Conf_I(\R):=\{I\hook \R\}.
$$
 A $2$-dimensional Lie group $G_{2}=\R^+ \ltimes \R$ acts freely on $\Conf_n(\R)$ by the law,
 $$
\Ba{ccccc}
\Conf_n(\R) & \times & \R^+ \ltimes \R & \lon & \Conf_n(\R)\\
p=\{x_1,\ldots,x_n\}&& (\la,\nu) &\lon & \la p+\nu:= \{\la x_1+\nu, \ldots, \la x_n+\nu\}.
\Ea
$$
The action is free
so that the quotient space,
$$
C_n(\R):= \Conf_n(\R)/G_{2},\ \ \  n\geq 2,
$$
is naturally an $(n-2)$-dimensional real oriented manifold equipped with a smooth orientation preserving
action of the group $\bS_n$. In fact,
$$
C_n(\R)=C_n^o(\R)\times \bS_n
$$ with orientation, $\Omega_n$, defined
on $C_n^o(\R):=\Conf^o_n(\R)/G_{2}$ as follows: identify $C_n^o(\R)$ with the subspace of
$\Conf^o_n(\R)$ consisting of points $\{0=x_{1}< x_{2} <\ldots < x_{n}=1\}$ and then set
$\Omega_n:= dx_2\wedge\ldots\wedge dx_{n-1}$.

The space $C_2(\R)$ is closed as it is the disjoint union, $C_2(\R)\simeq \bS_2$,  of two points.
The topological compactification,  $\overline{C}_n(\R)$, of  $C_n(\R)$ for higher $n$
 can be defined as $\overline{C}_n^o(\R)\times \bS_n$ where $\overline{C}_n^o(\R)$ is, by definition, the closure of an embedding,
$$
\Ba{ccc}
C_n^o(\R) & \lon & (\R\P^2)^{n(n-1)(n-2)}\\
(x_{1}, \ldots, x_{n}) & \lon &  \displaystyle\prod_{\#\{i,j,k\}=3}\left[|x_{i}-x_{j}| :
|x_{i}-x_{k}|: |x_{j}-x_{k}|\right].
\Ea
$$
Its codimension one strata are given by
$$
\p \overline{C}_n^o(\R) = \bigsqcup_{A} \overline{C}^o_{n - \# A + 1}(\R)\times
 \overline{C}_{\# A}^o(\R),
$$
where the union runs over {\em connected}\, proper  subsets, $A$, of the ordered set $\{1,2,\ldots,n\}$. The associated collection $\overline{C}(\R)=\{\overline{C}_n(\R)\}$
is a free operad in the category with the set of  generators,
$$
\left\{ {C}_n^o(\R) \simeq \Ba{c}\resizebox{21mm}{!}{\xy
(1,-5)*{\ldots},
(-13,-7)*{_{1}},
(-8,-7)*{_{2}},
(-3,-7)*{_{3}},
(8,-7)*{_{{n-1}}},
(14,-7)*{_{n}},
 (0,0)*{\circ}="a",
(0,5)*{}="0",
(-12,-5)*{}="b_1",
(-8,-5)*{}="b_2",
(-3,-5)*{}="b_3",
(8,-5)*{}="b_4",
(12,-5)*{}="b_5",
\ar @{-} "a";"0" <0pt>
\ar @{-} "a";"b_2" <0pt>
\ar @{-} "a";"b_3" <0pt>
\ar @{-} "a";"b_1" <0pt>
\ar @{-} "a";"b_4" <0pt>
\ar @{-} "a";"b_5" <0pt>
\endxy}\Ea\right\}_{n\geq 2}
$$
With the above graphical notations for the generators, the compactified configuration space
is the disjoint union of sets parameterized by planar rooted (equivalently, directed) trees
$$
 \overline{C}_n^o(\R)=\coprod_{T\in  \cT ree_n} T(\R)
$$
where $ \cT ree_n$ is the set of all planar trees with $n$ input legs whose vertices are at least
trivalent (i.e.\ have at least two input half-edges)\footnote{The set of internal edges of a rooted tree is denoted by $E(T)$, its set of legs by $Leg(G)$, and the set of vertices by $V(T)$; for example, picture (\ref{2: metric tree T}) below shows a rooted tree (with directions of edges tacitly chosen to run from bottom to the top) with $\# E(T)=3$, $\# Leg(T)=7$ and $\# V(T)=4$.
There is a natural partial order on the set $V(T)$: $v_1>  v_2$ if and only if there is a directed path of internal edges starting at $v_2$ and ending at   $v_1$. The set  $\cT ree_n$ also admits a partial order:
$T_1> T_2$ if and only if $T_2$ can be obtained from $T_1$ by contraction of at least one internal edge.} and
$$
T(\R):=\prod_{v\in V(T)}{C}_{\# v} ^o(\R)
$$
 is a set, better to say, a tree ``decorated" by sets. In this decomposition the one-vertex tree
corresponds to the big open cell ${C}_n^o(\R)\subset \overline{C}_n^0(\R)$, while trees with larger number of vertices to the
boundary components of the closed topological space  $\overline{C}_n(\R)$. Therefore the compactified space  $\overline{C}_n^o(\R)$ is homeomorphic, as a stratified topological space, to the $n$-th Stasheff associahedron $\cK_n$, and associated to $\overline{C}_n(\R)$ the operad of fundamental chains gives the minimal resolution, $\cA ss_\infty$, of the operad of associative algebras.

\sip

The trees parameterizing the boundary strata of   $\overline{C}_n^o(\R)$
 can also be used to define a structure of a smooth manifold with corners on $\overline{C}_n^o(\R)$  \cite{Ko}. In particular, a decoration of internal edges of such a tree $T$ with ``small" real parameters defines an smooth open coordinate chart, $\cU_T$, of the boundary strata corresponding to $T$ in $\overline{C}_n^o(\R)$ as follows  (see \cite{Ko,Ga} and lecture notes \cite{Me1} for details)
    $$
    \al_T: [0,\var)^{\# E(T)}\times \prod_{v\in V(T)} {C}^{st}_{\# In(v)}(\R)\simeq  \cU_T\subset
    \overline{C}_n(\R)
    $$
where $E(T)$ is the set of internal edges of $T$, $V(T)$ the set of vertices, $\var\in \R$ is a sufficiently small number (which is in fact depends on coordinates in the factors
  ${C}^{st}_{\# In(v)}(\R)$, i.e.\ strictly speaking the left hand side is a subset of a  smooth bundle over
  $\prod_{v\in V(T)} {C}^{st}_{\# In(v)}(\R)$ but we ignore these unimportant subtleties here),
  and ${C}^{st}_{k}(\R)$ is an $\bS_n$-equivariant section, $\tau: C_n(\R)\rar \Conf_n(\R)$, of the natural projection
$\Conf_n(\R)\rar  C_n(\R)$ defined, for example, by equations $\sum_{i=1}^n x_i=0$ and $\sum_{i}|x_i|^2=1$; clearly, such a section is a smooth manifold so that the l.h.s.\ of the isomorphism $\al_T$ is a smooth manifold with corners and can serve as a coordinate chart indeed.
 For example, a tree \cite{Me1}
 \Beq\label{2: metric tree T}
T=
\Ba{c}\resizebox{23mm}{!}{\xy
(2.0,3.0)*{_{\var_{1}}},
(11.2,3.7)*{_{\var_{2}}},
(-6.8,-4)*{_{\var_{3}}},
(-10.5,-2)*{_1},
(-11,-17)*{_3},
(-2,-17)*{_5},
(3,-10)*{_6},
(8,-10)*{_2},
(14,-10)*{_4},
(21,-10)*{_7},
(0,14)*{}="0",
 (0,8)*{\circ}="a",
(-10,0)*{}="b_1",
(-2,0)*{\circ}="b_2",
(12,0)*{\circ}="b_3",
(2,-8)*{}="c_1",
(-7,-8)*{\circ}="c_2",
(8,-8)*{}="c_3",
(14,-8)*{}="c_4",
(20,-8)*{}="c_5",
(-11,-15)*{}="d_1",
(-3,-15)*{}="d_2",
\ar @{-} "a";"0" <0pt>
\ar @{-} "a";"b_1" <0pt>
\ar @{-} "a";"b_2" <0pt>
\ar @{-} "a";"b_3" <0pt>
\ar @{-} "b_2";"c_1" <0pt>
\ar @{-} "b_2";"c_2" <0pt>
\ar @{-} "b_3";"c_3" <0pt>
\ar @{-} "b_3";"c_4" <0pt>
\ar @{-} "b_3";"c_5" <0pt>
\ar @{-} "c_2";"d_1" <0pt>
\ar @{-} "c_2";"d_2" <0pt>
\endxy}\Ea \ \ \ \ \ \ \var_{1}, \var_{2}, \var_{3}\in [0,\var) \ \mbox{for some}\ 0\leq \var \ll +\infty;
\Eeq
 gives a coordinate chart,
$$
\Ba{ccccccccccc}
 [0,\var)^3 & \hspace{-2mm} \times \hspace{-2mm}& C^{st}_3(\R) & \hspace{-2mm} \times\hspace{-2mm} & C^{st}_2(\R) & \hspace{-2mm} \times \hspace{-2mm} &  C^{st}_3(\R)  &\hspace{-2mm} \times\hspace{-2mm} & C^{st}_2(\R)
& {\lon} & \overline{C}_7(\R) \\
 (\var_1,\var_2,\var_3) & \hspace{-2mm} \times \hspace{-2mm} & (x_1, x',x'') & \hspace{-2mm}\times
 \hspace{-2mm}& (x''', x_6) &\hspace{-2mm} \times \hspace{-2mm}&
(x_2,x_4,x_7) &\hspace{-2mm} \times\hspace{-2mm} & (x_3,x_7) &\lon& (y_1, y_3, y_5,y_6, y_2,y_4, y_7)
\Ea
$$
given explicitly as follows,
$$
\Ba{llllllllrrrr}
y_1 &=& x_1 & &  y_3 &=& x'+   \var_1(x'''+ \var_3 x_3)  &&  y_2 &=&  x''+  \var_2 x_2 \\
&&&&  y_5 &=& x'+  \var_1(x''' + \var_3 x_5)   &&  y_4 &=&  x''+ \var_2 x_4 \\
&&&&    y_6 &=& x'+  \var_1 x_6,            &&  y_7 &=&  x'' +  \var_2 x_7 \\
\Ea
$$
The boundary stratum corresponding to $T$ is given in $\cU_T$ by the equations
$\var_1=\var_2=\var_3=0$. In this atlas the boundary strata gets interpreted as the limit configurations of {\em collapsing}\,
points. However, our configurations are considered only up to an action of the group $G_2$, so that above 3-parameter family of configurations can be equivalently rewritten as
$$
\Ba{llllllllrrrr}
y_1 &=& \frac{1}{\var_1\var_2\var_3} x_1 & &  y_3 &=& \frac{1}{\var_1\var_2\var_3} x'+    \frac{1}{\var_2\var_3}x'''+  \frac{1}{\var_2} x_3  &&  y_2 &=&   \frac{1}{\var_1\var_2\var_3}x''+  \frac{1}{\var_1\var_3} x_2 \\
&&&&  y_5 &=&  \frac{1}{\var_1\var_2\var_3} x'+   \frac{1}{\var_2\var_3}x''' +  \frac{1}{\var_2} x_5   &&  y_4 &=&  \frac{1}{\var_1\var_2\var_3} x''+  \frac{1}{\var_1\var_3} x_4 \\
&&&&    y_6 &=&  \frac{1}{\var_1\var_2\var_3} x'+   \frac{1}{\var_2\var_3} x_6,            &&  y_7 &=&  \frac{1}{\var_1\var_2\var_3} x'' +   \frac{1}{\var_1\var_3} x_7 \\
\Ea
$$
and hence in the corresponding coordinate chart the limit configurations corresponds to points
going in groups {\em infinitely far away from each other}\, (with different relative speeds), i.e.\ as ``exploded" configurations. We shall work below
with configuration spaces of points on a {\em pair of lines}, $\R\times \R$, whose boundary strata are parameterized by pairs of trees (with some extra structure); then it will sometimes be useful to interpret the limit configurations as collapsing ones for one tree (i.e.\ on one copy of the real line), and as exploded ones for another tree (i.e.\ on another copy of $\R$).

\subsubsection{\bf Permutahedron}\label{2: subsection permutahedron} The n-dimensional permutahedron $\cP_n$ is defined as a convex hull in $\R^{n+1}$ of the set  $\{\sigma(1), \sigma(2), \ldots, \sigma(n+1)\}_{\sigma_{\bS_n}}$ of $(n+1)!$ points. The faces of $\cP_n$ are encoded by the ordered partitions of the set
$\{1, 2, \ldots , n + 1\}$, or equivalently, by the set of {\em leveled}\, planar trees with $n+1$
legs (see, e.g., \cite{LTV} or \cite{Ma} for examples and explanations). We recall that
a {\em leveled planar $n$-tree}\,  is  a rooted $n$-tree $T$ together with a surjective map,
$L: V(T) \rar [l]$, from the set of its vertices   to some finite ordinal $[l]=\{1,2,\ldots, l\}$ that respects the standard  partial order on $V(T)$. The set, $\caL \cT ree_{n}$, of  leveled planar trees is partially ordered:  $(T,L) > (T',L')$ if
$(T',L')$ is obtained from $(T,L)$ by a contraction of levels. In particular $(T,L) > (T',L')$ implies
$T \geq T'$. For a level tree $(T,L: V(T)\rar [l])$ we set
$$
|L|:=-l+ \sum_{i=1}^l \#L^{-1}(i).
$$

The configuration space model for the permutahedron was
given in \cite{LTV}. In our context (when we want to keep freedom of interpreting the limit
configurations either as collapsing or as exploded) it is useful to consider
the closure, $\widehat{C}_n^o(\R)$,  of
 ${C}_n^o(\R)$ under the following  embedding (cf.\ \cite{LTV}),
$$
\Ba{ccccc}
C_n^o(\R) & \lon & (\R\P^2)^{n(n-1)(n-2)} & \times & [0,\infty]^{n(n-1)(n-2)(n-3)}\\
(x_{1}, \ldots, x_{n}) & \lon & \displaystyle \prod_{\#\{i,j,k\}=3}\left[|x_{i}-x_{j}| :
|x_{i}-x_{k}|: |x_{j}-x_{k}|\right] &\times&\displaystyle \prod_{\#\{i,j,k,l\}=4}{\frac{|x_i-x_j|}{|x_k-x_l|}}
\Ea
$$
where $[0,\infty]$ is a 1-dimensional compact smooth manifold with corners
with a defining coordinate chart given by
$$
\Ba{ccc}
[0,\infty] &\lon & [0,1]\\
   t       & \lon & \frac{t}{t+1}
\Ea
$$
The set  $\widehat{C}_n^o(\R)$ is
is the disjoint union of sets parameterized by planar rooted level trees
$$
 \widehat{C}_n^o(\R)=\coprod_{T\in  \caL \cT ree_n} T(\R),
$$
and, as a smooth manifold with corners, can be identified with the permutahedron
$\cP_{n-1}$.
For example, the following level trees,
$$
T_1=\Ba{c}\resizebox{16mm}{!}
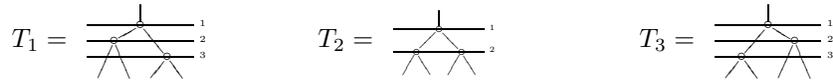
{\xy
(12,0)*{^{_{1}}},
(12,-3)*{^{_{2}}},
(12,-6)*{^{_{3}}},
(0,4)*{}="0",
 (0,0)*{\circ}="a",
(-5,-3)*{\circ}="b_1",
(5,-6)*{\circ}="b_2",
(-8,-10)*{}="c_1",
(-2,-10)*{}="c_2",
(8,-10)*{}="c_3",
(2,-10)*{}="c_4",
 (-10,0)*{}="1L",
(10,0)*{}="1R",
 (-10,-3)*{}="2L",
(10,-3)*{}="2R",
 (-10,-6)*{}="3L",
(10,-6)*{}="3R",
\ar @{-} "a";"0" <0pt>
\ar @{-} "a";"b_1" <0pt>
\ar @{-} "a";"b_2" <0pt>
\ar @{-} "b_1";"c_1" <0pt>
\ar @{-} "b_1";"c_2" <0pt>
\ar @{-} "b_2";"c_3" <0pt>
\ar @{-} "b_2";"c_4" <0pt>
\ar @{-} "1L";"1R" <0pt>
\ar @{-} "2L";"2R" <0pt>
\ar @{-} "3L";"3R" <0pt>
\endxy}\Ea
\ \ \ \ \ \ \ \ \ \ \
T_2=\Ba{c}\resizebox{18mm}{!}{\xy
(12,0)*{^{_{1}}},
(12,-5)*{^{_{2}}},
 (-10,0)*{}="1L",
(10,0)*{}="1R",
 (-10,-5)*{}="2L",
(10,-5)*{}="2R",
(0,4)*{}="0",
 (0,0)*{\circ}="a",
(-5,-5)*{\circ}="b_1",
(5,-5)*{\circ}="b_2",
(-8,-10)*{}="c_1",
(-2,-10)*{}="c_2",
(8,-10)*{}="c_3",
(2,-10)*{}="c_4",
(20,-8)*{}="c_5",
\ar @{-} "a";"0" <0pt>
\ar @{-} "a";"b_1" <0pt>
\ar @{-} "a";"b_2" <0pt>
\ar @{-} "b_1";"c_1" <0pt>
\ar @{-} "b_1";"c_2" <0pt>
\ar @{-} "b_2";"c_3" <0pt>
\ar @{-} "b_2";"c_4" <0pt>
\ar @{-} "1L";"1R" <0pt>
\ar @{-} "2L";"2R" <0pt>
\endxy}\Ea
\ \ \ \ \ \ \ \ \ \ \
T_3=\Ba{c}\resizebox{16mm}{!}{\xy
(12,0)*{^{_{1}}},
(12,-3)*{^{_{2}}},
(12,-6)*{^{_{3}}},
(0,4)*{}="0",
 (0,0)*{\circ}="a",
(-5,-6)*{\circ}="b_1",
(5,-3)*{\circ}="b_2",
(-8,-10)*{}="c_1",
(-2,-10)*{}="c_2",
(8,-10)*{}="c_3",
(2,-10)*{}="c_4",
 (-10,0)*{}="1L",
(10,0)*{}="1R",
 (-10,-3)*{}="2L",
(10,-3)*{}="2R",
 (-10,-6)*{}="3L",
(10,-6)*{}="3R",
\ar @{-} "a";"0" <0pt>
\ar @{-} "a";"b_1" <0pt>
\ar @{-} "a";"b_2" <0pt>
\ar @{-} "b_1";"c_1" <0pt>
\ar @{-} "b_1";"c_2" <0pt>
\ar @{-} "b_2";"c_3" <0pt>
\ar @{-} "b_2";"c_4" <0pt>
\ar @{-} "1L";"1R" <0pt>
\ar @{-} "2L";"2R" <0pt>
\ar @{-} "3L";"3R" <0pt>
\endxy}\Ea
$$
encode, respectively, the following limit configurations (as well as coordinate charts near the limit configurations) in $\cP_3= \widehat{C}_4^o(\R)$:
\Bi
\item[(i)] $T_1$ corresponds to the point in $\cP_3$ obtained in the limit
$\var_1,\var_2\rar +0$  from the configurations,
$$
x_1=-1 - \var_1,\ \ \ \
x_2=-1  + \var_1,\ \ \
x_3=1  - \var_1\var_2,\ \ \
x_4=1  + \var_1\var_2,\ \ \
$$
\item[(ii)] $T_2$ corresponds to the 1-dimensional strata in $\cP_3$ obtained in the limit
$\var\rar +0$  from the configurations,
$$
x_1=-1 - \var x,\ \ \ \
x_2=-1  + \var x,\ \ \
x_3=1  - \var x,\ \ \
x_4=1  + \var x,\ \ \ \ \  x=\frac{x_4-x_3}{x_2-x_1}\in (0,+\infty).
$$
\item[(iii)] $T_3$ corresponds to the point in $\cP_3$ obtained in the limit
$\var_1,\var_2\rar +0$  from the configurations,
$$
x_1=-1 - \var_1\var_2,\ \ \ \
x_2=-1  + \var_1\var_2,\ \ \
x_3=1  - \var_2,\ \ \
x_4=1  + \var_2,\ \ \
$$
\Ei

 For future reference we outline
 a general pattern which associates to a limit configuration, $p=\lim \{x_1,\ldots,x_n\}$, in $\widehat{C}_n^o(\R)$ a {\em levelled}\, tree:
\Bi
\item[(a)] there is a natural projection $\pi:\widehat{C}_n^o(\R)\rar
\overline{C}_n^o(\R)$ which associates to $p$ its image $\pi(p)$ in the associahedron and hence a unique maximal (with respect to the standard partial order in the poset $\cT ree_n$) unlevelled $n$-tree $T\in \cT ree_n$ such that $p\in T(\R)\subset \widehat{C}_n^o(\R)$; the legs of $T$ are naturally labelled by the set $[n]$.

\item[(b)] every vertex $v$ of the unlevelled tree $T$ from (a) stands for a collection of points $\{x_{i_v}\in \R\}_{i_v\in H(v)}$ parameterized by the set $H(v)$ of input half edges at $v\in T_p$ which collapse to a single point $x_v$ in $\R$; we introduce an equivalence relation
    in the set $V(T_p)$ of vertices of the tree $T_p$: $v'\sim v''$ if and only if $\lim \frac{|x_{i_{v'}}-x_{j_{v'}}|}{|x_{k_{v''}}-x_{l_{v''}}|}$ is a non-zero finite number
    for some  (and hence all) $i_{v'}\neq j_{v'}\in H(v')$ and $k_{v''}\neq l_{v''}\in H(v'')$  ; the associated equivalence classes $[v']$ are called {\em levels}; we say that equivalent vertices {\em lie on the same level};

\item[(c)] the natural partial ordering in the set of vertices, $V(T_p)$, induces a well-defined {\em total}\, ordering on the set of its levels.
    Indeed, if $v'$ and $v''$ belong to different levels, then either
     $\lim \frac{|x_{i_{v'}}-x_{j_{v'}}|}{|x_{k_{v''}}-x_{l_{v''}}|}=+\infty$ (in which case the level $[v']$ lies above the level $[v'']$ in the standard pictorial representation of a tree) or $\lim \frac{|x_{i_{v'}}-x_{j_{v'}}|}{|x_{k_{v''}}-x_{l_{v''}}|}=0$ (in which case the level $[v']$ lies below the level $[v'']$).
\Ei
As a result we get  a natural partition of the permutahedron,
$$
 \widehat{C}_n^o(\R)=\coprod_{(T,L)\in  \caL\cT ree_n} T(\R) \times (\R^+)^{|L|},
$$
parameterized by  leveled trees; by analogy to the case of the associahedron,
one can use this partition to introduce a smooth (with corners) atlas on
$\widehat{C}_n^o(\R)$ in which each leveled tree $(T,L)$ (with edges decorated by sufficiently small  parameters and with levels decorated by arbitrary non-negative parameters) gives us  a coordinate chart
near the boundary strata   $T(\R) \times (\R^+)^{|L|}\subset  \widehat{C}_n^o(\R)$. Thus  $\widehat{C}_n^o(\R)=\cP_{n-1}$
can be given a structure of smooth manifold with corners (we do not use in this paper a finer fact that
$\cP_{n-1}$ can be identified with a polytope).

\subsection{Bipermutahedron}\label{2: subsection bipermutohedron} In this and the next subsections we give a configuration space interpretation of the
bipermutahedron and biassociahedron
posets, $\cP^m_n$ and, respectively, $\cK_n^m$, which were introduced and studied by Martin Markl in \cite{Ma}. We show that these posets can be identified
with the boundary posets of certain smooth manifolds with corners (which come equipped with a natural structure of semialgebraic manifolds).

\sip

Consider a configuration space
$$
\Conf^o_{m,n}(\R\times \R):=\Conf^o_{m}(\R) \times \Conf^o_{n}(\R).
$$
A point $p\in \Conf^o_{m,n}(\R\times \R)$ is a pair $(p',p'')$  of collections of real numbers,
$$
p'=\{x_1<\ldots< x_m\}, \ \ \  p''= \{y_1<\ldots< y_n\}.
$$
  The group
$G_3:=R^+\rtimes \R^2$ acts freely on $\Conf^o_{m,n}(\R\times \R)$ for all
$m+n\geq 3$ by rescalings and translations,
$$
\Ba{ccccc}
G_3 &\times &  \Conf^o_{m,n}(\R\times \R) & \lon & \Conf^o_{m,n}(\R\times \R)\\
(\la,a,b) &  & (p',p'') & \lon &
(\la p'+ a; \la^{-1}p'' +b)
\Ea.
$$
The space of orbits,
$$
C_{m,n}(\R\times \R):=\frac{\Conf^o_{m,n}(\R\times \R)}{G_3}
$$
is a $(m+n-3)$-dimensional oriented manifold. It is clear that
$$
C_{1,n}(\R\times \R)=C_{n,1}(\R\times \R)=C_n^o(\R)
$$
and we define their compactifications $\widehat{C_{1,n}}(\R\times \R)$ and $\widehat{C_{n,1}}(\R\times \R)$ as the permutahedron $\widehat{C}_{n}^o(\R)$. For $m,n\geq 2$, there are canonical projections
$$
\pi': C_{m,n}(\R\times \R)\rar C_{m}(\R), \ \ \ \ \pi'': C_{m,n}(\R\times \R)\rar C_{n}(\R)
$$
which can be used to construct the following  embedding
\[
\Ba{ccccccc}
 C_{m,n}(\R\times \R)\hspace{0mm} & \lon &\hspace{0mm} \widehat{C}_m(\R) & \times &  \hspace{0mm}
 \widehat{C}_n(\R) &\times& \hspace{0mm} [0,\infty]^{\frac{nm(n-1)(m-1)}{4}}
 \vspace{1mm}\\
 (p',p'')
 \hspace{0mm} & \lon & \hspace{0mm} p' &\times&
 \hspace{0mm} p''
\hspace{0mm} &\times &\hspace{0mm}
%\Ba{c}\underset{_{\#\{i,j,\al,\be\}=4}}{\sqcap}\Ea\hspace{-7mm}
\displaystyle \prod_{i>j, \al>\be}
{{|x_{i}-x_{j}||y_{\al}-y_{\be}|}}
\hspace{0mm}
\Ea
\]
and define the compactified configuration space $\widehat{C_{m,n}}(\R\times \R)$ as the closure
of the image of $C_{m,n}(\R\times \R)$ under this embedding.
By analogy to the case of permutohedra, the compact space $\widehat{C_{m,n}}(\R\times \R)$ can be given naturally a structure of a smooth manifold with corners; in particular, this space comes with a stratification,
 $$
 \widehat{C_{m,n}}(\R\times \R)\supset \p \widehat{C_{m,n}}(\R\times \R) \supset \p^2 \widehat{C_{m,n}}(\R\times \R)\supset \ldots,
 $$
 and it is not hard to check that  the associated to this stratification poset is precisely the bipermutohedron poset $\cP^m_n$ from \cite{Ma}.
Let us first recall from \cite{Ma} the definition of the poset $\cP_m^n$, $m\geq 1$, $n\geq 1$, $m+n\geq 3$. For $m,n\geq 2$ the set  $\cP_m^n$ is defined as the set of all triples, $(T^\uparrow, T_\downarrow, \ell)$, consisting of an up rooted tree
$T^\uparrow\in \cT ree_n$, of a down-rooted tree $T_\downarrow\in \cT ree_m$,
and a strictly order preserving\footnote{i.e.\ if $v>u$ then $\ell(v)>\ell(u)$.} surjective {\em level  function}\,  $\ell: V(T^\uparrow) \cup V(T_\downarrow)\rar [l]$. For example
$$
\Ba{c}\resizebox{20mm}{!}{\xy
 (-10,0)*{}="1L",
(19,0)*{}="1R",
 (-10,5)*{}="2L",
(19,5)*{}="2R",
 (-10,-5)*{}="3L",
(19,-5)*{}="3R",
(0,4)*{}="0",
 (0,0)*{\circ}="a",
(-5,-5)*{\circ}="b_1",
(5,-5)*{\circ}="b_2",
(-8,-10)*{}="c_1",
(-2,-10)*{}="c_2",
(8,-10)*{}="c_3",
(2,-10)*{}="c_4",
(20,-8)*{}="c_5",
(10,-9)*{}="0'",
 (10,-5)*{\circ}="a'",
(7,0)*{\circ}="b_1'", %2nd graph
(15.5,5)*{}="b_2'",
(4,5)*{}="c_1'",
(9.5,5)*{}="c_2'",
\ar @{-} "a";"0" <0pt>
\ar @{-} "a";"b_1" <0pt>
\ar @{-} "a";"b_2" <0pt>
\ar @{-} "b_1";"c_1" <0pt>
\ar @{-} "b_1";"c_2" <0pt>
\ar @{-} "b_2";"c_3" <0pt>
\ar @{-} "b_2";"c_4" <0pt>
\ar @{-} "a'";"0'" <0pt>
\ar @{-} "a'";"b_1'" <0pt>
\ar @{-} "a'";"b_2'" <0pt>   %2nd graph
\ar @{-} "b_1'";"c_1'" <0pt>
\ar @{-} "b_1'";"c_2'" <0pt>
\ar @{-} "1L";"1R" <0pt>
\ar @{-} "3L";"3R" <0pt>
\endxy}\Ea \in \cP_4^3
$$
We define
$$
|\ell|:=-l+ \sum_{i=1}^l \ell^{-1}(i).
$$
The set $\cP_n^m$ is partially ordered: $(T^\uparrow, T_\downarrow, \ell)>
(\tilde{T}^\uparrow, \tilde{T}_\downarrow, \tilde{\ell})$ if the latter
 can be obtained from  the former by contraction of levels.
 The posets $\cP_n^1$ and $\cP^n_1$ are identified with $\caL \cT ree_n$
 but their elements are still represented as {\em pairs}\, of trees with the help
 of the singular tree $|$ which has no vertices, for example
 $$
 \Ba{c}
\resizebox{6mm}{!}{\xy
 (-4,0)*{}="1L",
(7,0)*{}="1R",
 (-4,-4)*{}="2L",
(7,-4)*{}="2R",
(0,3)*{}="0",
 (0,0)*{\circ}="a",
(-2,-4)*{\circ}="b_1",
(3,-7)*{}="b_2",
(5,3)*{}="r0",
(5,-7)*{}="r00",
(-3.7,-7)*{}="c_1",
(-0.5,-7)*{}="c_2",
\ar @{-} "a";"0" <0pt>
\ar @{-} "a";"b_1" <0pt>
\ar @{-} "a";"b_2" <0pt>
\ar @{-} "b_1";"c_1" <0pt>
\ar @{-} "b_1";"c_2" <0pt>
\ar @{-} "r0";"r00" <0pt>
\ar @{-} "1L";"1R" <0pt>
\ar @{-} "2L";"2R" <0pt>
\endxy}\Ea\in \cP_3^1.
 $$
 To each (limit) configuration, $p=\lim \{x_1,\ldots,x_n\}$, in $\widehat{C}_{m,n}(\R\times \R)$ we associate a uniquely defined leveled bi-tree from
 $\cP_m^n$ by a procedure which is completely analogous to the one described at the end of \S {\ref{2: subsection permutahedron}} and get, therefore, a decomposition,
\Beq\label{2: partition of bipermutohedron}
 \widehat{C}_{m,n}(\R\times \R)=\coprod_{(T^\uparrow,T_\downarrow,\ell)\in  \cP_m^n} T^\uparrow(\R)\times T_\downarrow(\R) \times (\R^+)^{|\ell|}.
\Eeq
 This decomposition can be used to define a smooth (with corners) atlas on
 the bipermutohedron $\widehat{C}_{m,n}(\R\times \R)$.

\subsection{Biassociahedron}\label{6: subsec on biassociahedron} Compactifications
$\overline{C_{1,n}}(\R\times \R)$ and $\overline{C_{n,1}}(\R\times \R)$ of the configuration spaces
${C_{1,n}}(\R\times \R)$ and respectively ${C_{n,1}}(\R\times \R)$ are defined as the associahedron $\overline{C}^o_n(\R)$. For $m,n\geq 2$ we define a  compactification
$\overline{C_{m,n}}(\R\times \R)$ of the configuration space $C_{m,n}(\R\times \R)$ as the closure of  the image of $C_{m,n}(\R\times \R)$ under the following embedding (cf.\ \cite{Sh1}),
\[
\Ba{ccccccc}
 C_{m,n}(\R\times \R)\hspace{0mm} & \lon &\hspace{0mm} \overline{C}_m(\R) & \times &  \hspace{0mm}
 \overline{C}_n(\R) &\times& \hspace{0mm} [0,\infty]^{\frac{nm(n-1)(m-1)}{4}}
 \vspace{1mm}\\
 (p',p'')
 \hspace{0mm} & \lon & \hspace{0mm} p' &\times&
 \hspace{0mm} p''
\hspace{0mm} &\times &\hspace{0mm}
%\Ba{c}\underset{_{\#\{i,j,\al,\be\}=4}}{\sqcap}\Ea\hspace{-7mm}
\displaystyle \prod_{i>j, \al>\be}
{{|x_{i}-x_{j}||y_{\al}-y_{\be}|}}
\hspace{0mm}
\Ea
\]

There is a natural surjection
$$
P: \widehat{C_{m,n}}(\R\times \R) \lon \overline{C_{m,n}}(\R\times \R)
$$
so that the partition (\ref{2: partition of bipermutohedron}) induces a partition of  $\overline{C_{m,n}}(\R\times \R)$. The induced partition is again parameterized by pairs of trees with an extra structure. The difference of the compactification formula for  $\overline{C_{m,n}}(\R\times \R)$ from the one for
$\widehat{C_{m,n}}(\R\times \R)$ is that we have no factors $\frac{|x_{i}-x_{j}|}{|x_k-x_l|}$ and $\frac{|y_{\al}-y_{\be}|}{|y_\ga-y_\delta|}$
which measure relatives speeds of collapsing/exploding groups of points belonging solely to one of the factors in $\R\times \R$. Hence the projection $P$
applied to the stratum  $T^\uparrow(\R)\times T_\downarrow(\R) \times (\R^+)^{|\ell|}$
contracts to single points those factors of $\R^+$ which correspond to the  levels $i\in [l]$ which have the property that either
$\ell^{-1}(i)\cap V(T_\uparrow)=\emptyset$
or   $\ell^{-1}(i)\cap V(T_\downarrow)=\emptyset$.  However such levels do not disappear completely from the induced stratification formula as it still
makes sense to compare  $\ell^{-1}(i)$ with  $\ell^{-1}(j)$ in the cases when $\ell^{-1}(i)\cap V(T_\uparrow)=\emptyset$ and $\ell^{-1}(j)\cap V(T^\uparrow)=\emptyset$.
Thus after the projection $P$ the level
function on $V(T^\uparrow)\sqcup V(T_\downarrow)$ gets transformed into
 a so called {\em zone function}\, \cite{Ma} which, by definition, is a  surjection,
 $$
 \zeta: V(T^\uparrow)\sqcup V(T_\downarrow) \lon [l]
 $$
 satisfying the following conditions:
 \Bi
 \item[(i)] if $v<u$, then $\zeta(v)\leq \zeta (u)$,
 \item[(ii)] for any pair of different elements $i,j\in [l]$ with $\zeta^{-1}(i)$ and $\zeta^{-1}(j)$
 containing vertices from {\em both}\, sets $V(T^\uparrow)$ and $V(T_\downarrow)$, then $i<j$ implies
 $v< u$ for every vertex $v\in \zeta^{-1}(i)$ and every vertex $u\in \zeta^{-1}(j)$ such that the relation
 $v<u$;
 \item[(iii)] there is no $i\in [l]$ such that both subsets $\zeta^{-1}(i)$ and $\zeta^{-1}(i+1)$
 belong to  $V(T^\uparrow)$ or both belong to $V(T_\downarrow)$.
 \Ei
Elements $i\in [l]$ with $\zeta^{-1}(i)\cap V(T^\uparrow)\neq \emptyset$ and  $\zeta^{-1}(i)\cap V(T_\downarrow)\neq \emptyset$ are called {\em barriers}\, and are depicted as solid horizontal lines.
Elements $i\in [l]$ with $\zeta^{-1}(i)\cap V(T^\uparrow)= \emptyset$ are called {\em down-zones},
while elements  $\zeta^{-1}(i)\cap V(T_\downarrow)= \emptyset$ are called {\em up-zones}; they are depicted as dashed
horizontal lines.
Thus condition (i)
says that the zone function is order preserving, condition (ii) says that it is strictly order preserving on barriers, and condition  (iii) says that there are no adjacent zones of the same type. Here are examples,
$$
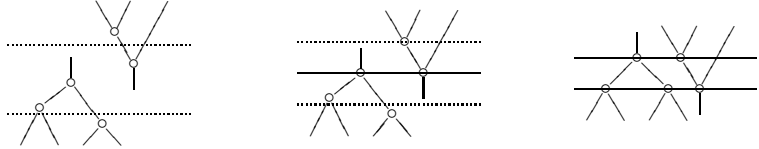
\Ba{c}\resizebox{25mm}{!}{\xy
%(12,0)*{^{_{1}}},
%(12,-5)*{^{_{2}}},
 (-10,6)*{}="1L",
(19,6)*{}="1R",
 (-10,-5)*{}="2L",
(19,-5)*{}="2R",
(0,4)*{}="0",
 (0,0)*{\circ}="a",
(-5,-4)*{\circ}="b_1",
(5,-6.6)*{\circ}="b_2",
(-8,-10)*{}="c_1",
(-2,-10)*{}="c_2",
(8,-10)*{}="c_3",
(2,-10)*{}="c_4",
(20,-8)*{}="c_5",
(10,-1)*{}="0'",
 (10,3)*{\circ}="a'",
(7,8)*{\circ}="b_1'", %2nd graph
(15.5,13)*{}="b_2'",
(4,13)*{}="c_1'",
(9.5,13)*{}="c_2'",
\ar @{-} "a";"0" <0pt>
\ar @{-} "a";"b_1" <0pt>
\ar @{-} "a";"b_2" <0pt>
\ar @{-} "b_1";"c_1" <0pt>
\ar @{-} "b_1";"c_2" <0pt>
\ar @{-} "b_2";"c_3" <0pt>
\ar @{-} "b_2";"c_4" <0pt>
\ar @{-} "a'";"0'" <0pt>
\ar @{-} "a'";"b_1'" <0pt>
\ar @{-} "a'";"b_2'" <0pt>   %2nd graph
\ar @{-} "b_1'";"c_1'" <0pt>
\ar @{-} "b_1'";"c_2'" <0pt>
\ar @{.} "1L";"1R" <0pt>
\ar @{.} "2L";"2R" <0pt>
\endxy}\Ea
\hspace{10mm}
\Ba{c}\resizebox{25mm}{!}{\xy
%(12,0)*{^{_{1}}},
%(12,-5)*{^{_{2}}},
 (-10,0)*{}="1L",
(19,0)*{}="1R",
 (-10,5)*{}="2L",
(19,5)*{}="2R",
 (-10,-5)*{}="3L",
(19,-5)*{}="3R",
(0,4)*{}="0",
 (0,0)*{\circ}="a",
(-5,-4)*{\circ}="b_1",
(5,-6.5)*{\circ}="b_2",
(-8,-10)*{}="c_1",
(-2,-10)*{}="c_2",
(8,-10)*{}="c_3",
(2,-10)*{}="c_4",
(20,-8)*{}="c_5",
(10,-4)*{}="0'",
 (10,0)*{\circ}="a'",
(7,5)*{\circ}="b_1'", %2nd graph
(15.5,10)*{}="b_2'",
(4,10)*{}="c_1'",
(9.5,10)*{}="c_2'",
\ar @{-} "a";"0" <0pt>
\ar @{-} "a";"b_1" <0pt>
\ar @{-} "a";"b_2" <0pt>
\ar @{-} "b_1";"c_1" <0pt>
\ar @{-} "b_1";"c_2" <0pt>
\ar @{-} "b_2";"c_3" <0pt>
\ar @{-} "b_2";"c_4" <0pt>
\ar @{-} "a'";"0'" <0pt>
\ar @{-} "a'";"b_1'" <0pt>
\ar @{-} "a'";"b_2'" <0pt>   %2nd graph
\ar @{-} "b_1'";"c_1'" <0pt>
\ar @{-} "b_1'";"c_2'" <0pt>
\ar @{-} "1L";"1R" <0pt>
\ar @{.} "2L";"2R" <0pt>
\ar @{.} "3L";"3R" <0pt>
\endxy}\Ea
\ \ \ \ \ \ \
\Ba{c}\resizebox{25mm}{!}{\xy
 (-10,0)*{}="1L",
(19,0)*{}="1R",
 (-10,5)*{}="2L",
(19,5)*{}="2R",
 (-10,-5)*{}="3L",
(19,-5)*{}="3R",
(0,4)*{}="0",
 (0,0)*{\circ}="a",
(-5,-5)*{\circ}="b_1",
(5,-5)*{\circ}="b_2",
(-8,-10)*{}="c_1",
(-2,-10)*{}="c_2",
(8,-10)*{}="c_3",
(2,-10)*{}="c_4",
(20,-8)*{}="c_5",
(10,-9)*{}="0'",
 (10,-5)*{\circ}="a'",
(7,0)*{\circ}="b_1'", %2nd graph
(15.5,5)*{}="b_2'",
(4,5)*{}="c_1'",
(9.5,5)*{}="c_2'",
\ar @{-} "a";"0" <0pt>
\ar @{-} "a";"b_1" <0pt>
\ar @{-} "a";"b_2" <0pt>
\ar @{-} "b_1";"c_1" <0pt>
\ar @{-} "b_1";"c_2" <0pt>
\ar @{-} "b_2";"c_3" <0pt>
\ar @{-} "b_2";"c_4" <0pt>
\ar @{-} "a'";"0'" <0pt>
\ar @{-} "a'";"b_1'" <0pt>
\ar @{-} "a'";"b_2'" <0pt>   %2nd graph
\ar @{-} "b_1'";"c_1'" <0pt>
\ar @{-} "b_1'";"c_2'" <0pt>
\ar @{-} "1L";"1R" <0pt>
\ar @{-} "3L";"3R" <0pt>
\endxy}\Ea
$$
of a fixed pair of trees and three different zone functions on the set of their vertices.
For a zone function $\zeta$ on $V(T^\uparrow)\sqcup V(T_\downarrow)$
we denote by $B(\zeta)$ the set of its barriers, and by $|\zeta|$ the non-negative integer,
$$
|\zeta|:=-l + \sum_{i\in B(\zeta)} \# \zeta^{-1}(i).
$$
The compactified configuration space  $\overline{C_{m,n}}(\R\times \R)$, the
{\em  biassociahedron} (cf.\ \cite{Ma}), comes therefore equipped with the induced stratification
\Beq\label{2: C m,n stratification formula}
\overline{C_{m,n}}(\R\times \R)=\bigcup_{(T^\uparrow, T_\downarrow, \zeta)} T^\uparrow(\R)\times T_\downarrow(\R) \times (0,+\infty)^{|\zeta|}
\Eeq
which is parameterized by the poset $\cK_m^n$ consisting of triples
$(T^\uparrow, T_\downarrow, \zeta)$. Therefore we often denote
$\overline{C_{m,n}}(\R\times \R)$ by $\sK_m^n$. This decomposition can be used to define in a standard way a smooth (with corners) atlas on
 the biassociahedron $\sK_m^n=\overline{C_{m,n}}(\R\times \R)$ such that the associated poset
 $$
 \overline{C_{m,n}}(\R\times \R)\supset \p \overline{C_{m,n}}(\R\times \R) \supset \p^2 \overline{C_{m,n}}(\R\times \R)\supset \ldots,
 $$
is precisely the poset  $\cK^m_n$ from \cite{Ma}.

\subsection{Example: $m+n=4$} This is the first non-trivial case. It is clear that
$$
\overline{C_{3,1}}(\R\times\R )\simeq \overline{C_{1,3}}(\R\times \R)\simeq \overline{C_{3}}(\R)\simeq [0,1].
 $$
 Therefore in the cases $(m=1,n=2)$ and $(m=2,n=1)$ the combinatorics of the
natural stratification of the compactified configuration spaces  can be coded by the following pairs of trees (each pair is equipped with the only possible zone
function),
$$
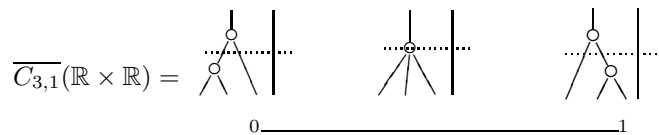
\overline{C_{3,1}}(\R\times\R )=\hspace{3mm} \Ba{c}
\hspace{-3mm}\resizebox{12mm}{!}{\xy
 (-3,-2)*{}="1L",
(7,-2)*{}="1R",
(0,3)*{}="0",
 (0,0)*{\circ}="a",
(-2,-4)*{\circ}="b_1",
(3,-7)*{}="b_2",
(5,3)*{}="r0",
(5,-7)*{}="r00",
(-3.7,-7)*{}="c_1",
(-0.5,-7)*{}="c_2",
\ar @{-} "a";"0" <0pt>
\ar @{-} "a";"b_1" <0pt>
\ar @{-} "a";"b_2" <0pt>
\ar @{-} "b_1";"c_1" <0pt>
\ar @{-} "b_1";"c_2" <0pt>
\ar @{-} "r0";"r00" <0pt>
\ar @{.} "1L";"1R" <0pt>
\endxy}
\ \  \ \ \  \ \ \ \ \
\resizebox{12mm}{!}{\xy
(-3,-1.5)*{}="1L",
(7,-1.5)*{}="1R",
(0,3)*{}="0",
 (0,-1.5)*{\circ}="a",
(3,-7)*{}="b_2",
(5,3)*{}="r0",
(5,-7)*{}="r00",
(-3.7,-7)*{}="c_1",
(-0.5,-7)*{}="c_2",
\ar @{-} "a";"0" <0pt>
\ar @{-} "a";"b_2" <0pt>
\ar @{-} "a";"c_1" <0pt>
\ar @{-} "a";"c_2" <0pt>
\ar @{-} "r0";"r00" <0pt>
\ar @{.} "1L";"1R" <0pt>
\endxy}
\ \  \ \ \  \ \ \ \ \ \
\resizebox{12mm}{!}{\xy
 (-3,-2)*{}="1L",
(7,-2)*{}="1R",
(0,3)*{}="0",
 (0,0)*{\circ}="a",
(2,-4)*{\circ}="b_1",
(-3,-7)*{}="b_2",
(5,3)*{}="r0",
(5,-7)*{}="r00",
(3.7,-7)*{}="c_1",
(0.5,-7)*{}="c_2",
\ar @{-} "a";"0" <0pt>
\ar @{-} "a";"b_1" <0pt>
\ar @{-} "a";"b_2" <0pt>
\ar @{-} "b_1";"c_1" <0pt>
\ar @{-} "b_1";"c_2" <0pt>
\ar @{-} "r0";"r00" <0pt>
\ar @{.} "1L";"1R" <0pt>
\endxy}\\
\vspace{-2mm}
\\
{\xy
 (-5,0)*{^0}="0",
 (44,0)*{^1}="a",
\ar @{-} "a";"0" <0pt>
\endxy}
\Ea
$$
The left pair corresponds to the point $0\in [0,1]$, the middle one to the open interval $(0,1)$, and the right pair of trees to the point
$1\in [0,1]$. Turning the trees above upside down, we get a ``pairs of trees" stratification of $\overline{C_{1,3}}(\R\times\R )$. The trees are not leveled, but it will be useful to understand these trees as trivially {\em zoned} (cf.\ \cite{Ma}), i.e.\ as the ones in which all vertices are assigned one an the same zone value $1$. We shall see below examples of trees with more than one zone.

\sip

The compactification formula says that $\overline{C_{2,2}}(\R\times\R )$ is the closure
of the embedding,
$$
\Ba{cccc}
{C_{2,2}}(\R\times\R ) &\lon & [0,+\infty]\\
(x_1,x_2),(y_1,y_2) & \lon & |x_2-x_1||y_2-y_1|
\Ea
$$
Thus  $\overline{C_{2,2}}(\R\times\R )\simeq [0,1]$, and the stratification $[0,1]=0\sqcup(0.1)\sqcup 1$ can be represented in terms of the pair of trees and three possible zone functions as follows,
$$
\overline{C_{2,2}}(\R\times\R )=\Ba{c}
\resizebox{10mm}{!}{\xy
(-3,-1)*{}="1L",
(7,-1)*{}="1R",
 (-3,1)*{}="2L",
(7,1)*{}="2R",
(0,4)*{}="0",
 (0,-1)*{\circ}="a",
(-2,-4)*{}="b_1",
(2,-4)*{}="b_2",
(4,-3)*{}="00",
 (4,1)*{\circ}="a0",
(2,4)*{}="c_1",
(6,4)*{}="c_2",
\ar @{-} "a";"0" <0pt>
\ar @{-} "a";"b_1" <0pt>
\ar @{-} "a";"b_2" <0pt>
\ar @{-} "a0";"c_1" <0pt>
\ar @{-} "a0";"c_2" <0pt>
\ar @{-} "a0";"00" <0pt>
\ar @{.} "1L";"1R" <0pt>
\ar @{.} "2L";"2R" <0pt>
\endxy}
\ \  \ \ \  \ \ \ \ \
\resizebox{11mm}{!}{
\xy
(-3,0)*{}="1L",
(7,0)*{}="1R",
(0,3)*{}="0",
 (0,0)*{\circ}="a",
(-2,-4)*{}="b_1",
(2,-4)*{}="b_2",
(4,-3)*{}="00",
 (4,0)*{\circ}="a0",
(2,4)*{}="c_1",
(6,4)*{}="c_2",
\ar @{-} "a";"0" <0pt>
\ar @{-} "a";"b_1" <0pt>
\ar @{-} "a";"b_2" <0pt>
\ar @{-} "a0";"c_1" <0pt>
\ar @{-} "a0";"c_2" <0pt>
\ar @{-} "a0";"00" <0pt>
\ar @{-} "1L";"1R" <0pt>
\endxy}
\ \  \ \ \  \ \ \ \ \ \
\resizebox{11mm}{!}{
\xy
(-3,-1)*{}="1L",
(7,-1)*{}="1R",
 (-3,1)*{}="2L",
(7,1)*{}="2R",
(0,4)*{}="0",
 (0,1)*{\circ}="a",
(-2,-4)*{}="b_1",
(2,-4)*{}="b_2",
(4,-4)*{}="00",
 (4,-1)*{\circ}="a0",
(2,4)*{}="c_1",
(6,4)*{}="c_2",
\ar @{-} "a";"0" <0pt>
\ar @{-} "a";"b_1" <0pt>
\ar @{-} "a";"b_2" <0pt>
\ar @{-} "a0";"c_1" <0pt>
\ar @{-} "a0";"c_2" <0pt>
\ar @{-} "a0";"00" <0pt>
\ar @{.} "1L";"1R" <0pt>
\ar @{.} "2L";"2R" <0pt>
\endxy}\\
\vspace{-2mm}
\\
\ \ {\xy
(-12,0)*{^0}="0",
 (46,0)*{^1}="a",
\ar @{-} "a";"0" <0pt>
\endxy}
\Ea
$$
The left pair of trees corresponds to the limit $\var\rar 0$ configuration
$$
(x_1=-\var, x_2=\var),\  (y_1=-1, y_2=1)\ \ \ \sim  \ \ \ (x_1=-1, x_2=1),\  (y_1=-\var, y_2=\var),
$$
with $|x_2-x_1||y_2-y_1|\rar 0$.
 The middle pair of tress  corresponds to the generic configurations,
$$
(x_1=- x, x_2= x),\  (y_1=-y, y_2=y)\ \ \ \sim \ \ \ (x_1=-\var x, x_2=\var x),\  (y_1=-\frac{1}{\var}y, y_2=\frac{1}{\var}y), \ \ \ x,y\in \R^+,
$$
with  $|x_2-x_1||y_2-y_1|$ a positive finite number (so that $|x_2-x_1|\sim |y_2-y_1|$
and the associated vertices are on the same level ). The right pair of trees corresponds to the limit $\var\rar 0$ of the configuration
$$
(x_1=-1, x_2=1), \ \ (y_1=-\frac{1}{\var}, y_2=\frac{1}{\var}) \ \ \ \sim  \ \ \  (x_1=-\frac{1}{\var}, x_2=\frac{1}{\var}), \  (y_1=-1, y_2=1)
$$
with $|x_2-x_1||y_2-y_1|\rar +\infty $.

\subsection{Example: $m+n=5$} The cases $(m=1,n=4)$ and $(m=4,n=1)$ are completely analogous to the example discussed above. The cases $(m=2,n=3)$ and $(m=3,n=2)$ are similar so that we shall study in detail only one of them. The compactification   $\overline{C_{3,2}}(\R\times\R )$ is the closure
of the embedding,
$$
\Ba{ccccc}
{C_{3,2}}(\R\times\R ) &\lon & \R\P^2 &\times & [0,+\infty]^3\\
(x_1,x_2,x_3),(y_1,y_2) & \lon & \left[|x_1-x_2|:|x_1-x_3|:|x_2-x_3|\right] &\times&
\left\{{\Ba{c} |x_2-x_1||y_2-y_1|\\  |x_3-x_1||y_2-y_1| \\  |x_2-x_3||y_2-y_1|\Ea}\right.
\Ea
$$
There are three possible pairs of trees in this case,
$$
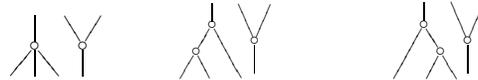
\Ba{c}\resizebox{12mm}{!}{\xy
(0,2)*{}="0",
 (0,-3)*{\circ}="a",
(-4,-8)*{}="b_1",
(0,-8)*{}="b_2",
(4,-8)*{}="b_3",
(8,-8)*{}="0'",
 (8,-3)*{\circ}="a'",
(5,2)*{}="b_1'", %2nd tree
(11,2)*{}="b_2'",
\ar @{-} "a";"0" <0pt>
\ar @{-} "a";"b_1" <0pt>
\ar @{-} "a";"b_2" <0pt>
\ar @{-} "a";"b_3" <0pt>
%
%\ar @{-} "1L";"1R" <0pt>
%\ar @{-} "2L";"2R" <0pt>
\ar @{-} "a'";"0'" <0pt>
\ar @{-} "a'";"b_1'" <0pt>
\ar @{-} "a'";"b_2'" <0pt>
\endxy}\Ea
\ \ \ \ \ \
\Ba{c}\resizebox{12mm}{!}{\xy
(0,4)*{}="0",
 (0,0)*{\circ}="a",
(-3,-5)*{\circ}="b_1",
(5.5,-10)*{}="b_2",
(8,-9)*{}="0'",
 (8,-3)*{\circ}="a'",
(5,4)*{}="b_1'", %2nd tree
(11,4)*{}="b_2'",
(-6,-10)*{}="c_1",
(-0.5,-10)*{}="c_2",
\ar @{-} "a";"0" <0pt>
\ar @{-} "a";"b_1" <0pt>
\ar @{-} "a";"b_2" <0pt>
\ar @{-} "b_1";"c_1" <0pt>
\ar @{-} "b_1";"c_2" <0pt>
%
%\ar @{-} "1L";"1R" <0pt>
%\ar @{-} "2L";"2R" <0pt>
\ar @{-} "a'";"0'" <0pt>
\ar @{-} "a'";"b_1'" <0pt>
\ar @{-} "a'";"b_2'" <0pt>
\endxy}\Ea
\ \ \ \ \ \ \ \ \ \ \
\Ba{c}\resizebox{12mm}{!}{\xy
(0,4)*{}="0",
 (0,0)*{\circ}="a",
(3,-5)*{\circ}="b_2",
(-5.5,-10)*{}="c_1",
(8,-9)*{}="0'",
 (8,-3)*{\circ}="a'",
(5,4)*{}="b_1'", %2nd tree
(11,4)*{}="b_2'",
(6,-10)*{}="c_3",
(0.5,-10)*{}="c_2",
\ar @{-} "a";"0" <0pt>
\ar @{-} "a";"c_1" <0pt>
\ar @{-} "a";"b_2" <0pt>
\ar @{-} "b_2";"c_2" <0pt>
\ar @{-} "b_2";"c_3" <0pt>
\ar @{-} "a'";"0'" <0pt>
\ar @{-} "a'";"b_1'" <0pt>
\ar @{-} "a'";"b_2'" <0pt>
\endxy}\Ea
$$
To check  claim (\ref{2: C m,n stratification formula}) we have to consider the list of all possible zone functions on these pairs, together with the associated boundary strata.

\sip

1) To the zone function
$\Ba{c}\resizebox{12mm}{!}{\xy
(-4,-3)*{}="1L",
(12,-3)*{}="1R",
(0,2)*{}="0",
 (0,-3)*{\circ}="a",
(-4,-8)*{}="b_1",
(0,-8)*{}="b_2",
(4,-8)*{}="b_3",
(8,-8)*{}="0'",
 (8,-3)*{\circ}="a'",
(5,2)*{}="b_1'", %2nd tree
(11,2)*{}="b_2'",
\ar @{-} "a";"0" <0pt>
\ar @{-} "a";"b_1" <0pt>
\ar @{-} "a";"b_2" <0pt>
\ar @{-} "a";"b_3" <0pt>
\ar @{-} "1L";"1R" <0pt>
%\ar @{-} "2L";"2R" <0pt>
\ar @{-} "a'";"0'" <0pt>
\ar @{-} "a'";"b_1'" <0pt>
\ar @{-} "a'";"b_2'" <0pt>
\endxy}\Ea$
 we associate, in accordance with (\ref{2: C m,n stratification formula}), the $2$-dimensional big cell
$$
C_{3,2}(\R\times\R )\simeq \left\{\Ba{c}(x_1=0,x_2=x,x_3=1)\\ (y_1=-y,y_2=y)\Ea\right. \simeq (0,1)\times (0,+\infty)
$$

\sip

2) The zone function $\Ba{c}\resizebox{12mm}{!}{\xy
(-4,-2)*{}="1L",
(12,-2)*{}="1R",
(-4,-5)*{}="2L",
(12,-5)*{}="2R",
(0,0)*{}="0",
 (0,-5)*{\circ}="a",
(-4,-10)*{}="b_1",
(0,-10)*{}="b_2",
(4,-10)*{}="b_3",
(8,-7)*{}="0'",
 (8,-2)*{\circ}="a'",
(5,3)*{}="b_1'", %2nd tree
(11,3)*{}="b_2'",
\ar @{-} "a";"0" <0pt>
\ar @{-} "a";"b_1" <0pt>
\ar @{-} "a";"b_2" <0pt>
\ar @{-} "a";"b_3" <0pt>
\ar @{.} "1L";"1R" <0pt>
\ar @{.} "2L";"2R" <0pt>
\ar @{-} "a'";"0'" <0pt>
\ar @{-} "a'";"b_1'" <0pt>
\ar @{-} "a'";"b_2'" <0pt>
\endxy}\Ea$
corresponds to the $1$-dimensional cell
$$
\lim_{\var\rar 0} \left\{\Ba{c}(x_1=0,x_2=x,x_3=1)\\ (y_1=-\var,y_2=\var)\Ea\right. \simeq (0,1)
$$

3) The zone function $\Ba{c}\resizebox{12mm}{!}{\xy
(-4,-2)*{}="1L",
(12,-2)*{}="1R",
(-4,-5)*{}="2L",
(12,-5)*{}="2R",
(0,3)*{}="0",
 (0,-2)*{\circ}="a",
(-4,-7)*{}="b_1",
(0,-7)*{}="b_2",
(4,-7)*{}="b_3",
(8,-10)*{}="0'",
 (8,-5)*{\circ}="a'",
(5,0)*{}="b_1'",      %2nd tree
(11,0)*{}="b_2'",
\ar @{-} "a";"0" <0pt>
\ar @{-} "a";"b_1" <0pt>
\ar @{-} "a";"b_2" <0pt>
\ar @{-} "a";"b_3" <0pt>
\ar @{.} "1L";"1R" <0pt>
\ar @{.} "2L";"2R" <0pt>
\ar @{-} "a'";"0'" <0pt>
\ar @{-} "a'";"b_1'" <0pt>
\ar @{-} "a'";"b_2'" <0pt>
\endxy}\Ea$
corresponds to the $1$-dimensional cell
$$
\lim_{\var\rar 0} \left\{\Ba{c}(x_1=0,x_2=x,x_3=1)\\ (y_1=-\frac{1}{\var},y_2=\frac{1}{\var})\Ea\right.
 \simeq (0,1)
$$

4) The zone functions $\Ba{c}\resizebox{12mm}{!}{\xy
(-4,6)*{}="1L",
(12,6)*{}="1R",
(-4,-3)*{}="2L",
(12,-3)*{}="2R",
(0,4)*{}="0",
 (0,0)*{\circ}="a",
(-3,-5)*{\circ}="b_1",
(5.5,-10)*{}="b_2",
(8,0)*{}="0'",
 (8,6)*{\circ}="a'",
(5,13)*{}="b_1'",                %2nd tree
(11,13)*{}="b_2'",
(-6,-10)*{}="c_1",
(-0.5,-10)*{}="c_2",
\ar @{-} "a";"0" <0pt>
\ar @{-} "a";"b_1" <0pt>
\ar @{-} "a";"b_2" <0pt>
\ar @{-} "b_1";"c_1" <0pt>
\ar @{-} "b_1";"c_2" <0pt>
\ar @{-} "a'";"0'" <0pt>
\ar @{-} "a'";"b_1'" <0pt>
\ar @{-} "a'";"b_2'" <0pt>
\ar @{.} "1L";"1R" <0pt>
\ar @{.} "2L";"2R" <0pt>
\endxy}\Ea$
and, respectively,
$\Ba{c}\resizebox{12mm}{!}{\xy
(-4,6)*{}="1L",
(12,6)*{}="1R",
(-4,-3)*{}="2L",
(12,-3)*{}="2R",
(0,4)*{}="0",
 (0,0)*{\circ}="a",
(3,-5)*{\circ}="b_1",
(5.5,-10)*{}="c_3",
(8,0)*{}="0'",
 (8,6)*{\circ}="a'",
(5,13)*{}="b_1'",                %2nd tree
(11,13)*{}="b_2'",
(-6,-10)*{}="c_1",
(-0.5,-10)*{}="c_2",
\ar @{-} "a";"0" <0pt>
\ar @{-} "a";"b_1" <0pt>
\ar @{-} "a";"c_3" <0pt>
\ar @{-} "a";"c_1" <0pt>
\ar @{-} "b_1";"c_2" <0pt>
\ar @{-} "a'";"0'" <0pt>
\ar @{-} "a'";"b_1'" <0pt>
\ar @{-} "a'";"b_2'" <0pt>
\ar @{.} "1L";"1R" <0pt>
\ar @{.} "2L";"2R" <0pt>
\endxy}\Ea$
correspond to $2$ points which are boundaries of the  closure of the strata 2)
in $\overline{C_{3,2}}(\R\times \R)$, i.e.\ they correspond, respectively, to the following two limit  configuration
$$
\lim_{\var_1,\var_2\rar 0} \left\{\Ba{c}(x_1=0,x_2=\var_1,x_3=1)\\ (y_1=-\var_2,y_2=\var_2)\Ea\right. \ \  \ \ \ \ \ \ \ \
\lim_{\var_1,\var_2\rar 0} \left\{\Ba{c}(x_1=0,x_2=1-\var_1,x_3=1)\\ (y_1=-\var_2,y_2=\var_2)\Ea\right.
$$

5) The zone functions $\Ba{c}\resizebox{12mm}{!}{\xy
(-4,-7)*{}="1L",
(12,-7)*{}="1R",
(-4,-3)*{}="2L",
(12,-3)*{}="2R",
(0,4)*{}="0",
 (0,0)*{\circ}="a",
(-3,-5)*{\circ}="b_1",
(5.5,-10)*{}="b_2",
(8,-13)*{}="0'",
 (8,-7)*{\circ}="a'",
(5,0)*{}="b_1'",                %2nd tree
(11,0)*{}="b_2'",
(-6,-10)*{}="c_1",
(-0.5,-10)*{}="c_2",
\ar @{-} "a";"0" <0pt>
\ar @{-} "a";"b_1" <0pt>
\ar @{-} "a";"b_2" <0pt>
\ar @{-} "b_1";"c_1" <0pt>
\ar @{-} "b_1";"c_2" <0pt>
\ar @{-} "a'";"0'" <0pt>
\ar @{-} "a'";"b_1'" <0pt>
\ar @{-} "a'";"b_2'" <0pt>
\ar @{.} "1L";"1R" <0pt>
\ar @{.} "2L";"2R" <0pt>
\endxy}\Ea$
and, respectively,
$\Ba{c}\resizebox{12mm}{!}{\xy
(-4,-7)*{}="1L",
(12,-7)*{}="1R",
(-4,-3)*{}="2L",
(12,-3)*{}="2R",
(0,4)*{}="0",
 (0,0)*{\circ}="a",
(3,-5)*{\circ}="b_1",
(5.5,-10)*{}="c_3",
(8,-13)*{}="0'",
 (8,-7)*{\circ}="a'",
(5,0)*{}="b_1'",                %2nd tree
(11,0)*{}="b_2'",
(-6,-10)*{}="c_1",
(-0.5,-10)*{}="c_2",
\ar @{-} "a";"0" <0pt>
\ar @{-} "a";"b_1" <0pt>
\ar @{-} "a";"c_3" <0pt>
\ar @{-} "a";"c_1" <0pt>
\ar @{-} "b_1";"c_2" <0pt>
\ar @{-} "a'";"0'" <0pt>
\ar @{-} "a'";"b_1'" <0pt>
\ar @{-} "a'";"b_2'" <0pt>
\ar @{.} "1L";"1R" <0pt>
\ar @{.} "2L";"2R" <0pt>
\endxy}\Ea$
correspond to $2$ points which are boundaries of the closure of the strata 3)
in $\overline{C_{3,2}}(\R\times \R)$, i.e.\ they correspond, respectively, to the following two limit  configuration
$$
\lim_{\var_1,\var_2\rar 0} \left\{\Ba{c}(x_1=0,x_2=\var_1,x_3=1)\\ (y_1=-\frac{1}{\var_1\var_2},y_2=\frac{1}{\var_1\var_2})\Ea\right. \ \  \ \ \ \ \ \ \ \
\lim_{\var_1,\var_2\rar 0} \left\{\Ba{c}(x_1=0,x_2=1-\var_1\,x_3=1)\\ (y_1=-\frac{1}{\var_1\var_2},y_2=\frac{1}{\var_1\var_2})\Ea\right.
$$

6) The zone functions $\Ba{c}\resizebox{12mm}{!}{\xy
(-4,0)*{}="1L",
(12,0)*{}="1R",
(-4,-5)*{}="2L",
(12,-5)*{}="2R",
(0,4)*{}="0",
 (0,0)*{\circ}="a",
(-3,-5)*{\circ}="b_1",
(5.5,-10)*{}="b_2",
(8,-6)*{}="0'",
 (8,0)*{\circ}="a'",
(5,7)*{}="b_1'",                %2nd tree
(11,7)*{}="b_2'",
(-6,-10)*{}="c_1",
(-0.5,-10)*{}="c_2",
\ar @{-} "a";"0" <0pt>
\ar @{-} "a";"b_1" <0pt>
\ar @{-} "a";"b_2" <0pt>
\ar @{-} "b_1";"c_1" <0pt>
\ar @{-} "b_1";"c_2" <0pt>
\ar @{-} "a'";"0'" <0pt>
\ar @{-} "a'";"b_1'" <0pt>
\ar @{-} "a'";"b_2'" <0pt>
\ar @{-} "1L";"1R" <0pt>
\ar @{.} "2L";"2R" <0pt>
\endxy}\Ea$
and, respectively,
$\Ba{c}\resizebox{12mm}{!}{\xy
(-4,0)*{}="1L",
(12,0)*{}="1R",
(-4,-5)*{}="2L",
(12,-5)*{}="2R",
(0,4)*{}="0",
 (0,0)*{\circ}="a",
(3,-5)*{\circ}="b_1",
(5.5,-10)*{}="c_3",
(8,-6)*{}="0'",
 (8,0)*{\circ}="a'",
(5,7)*{}="b_1'",                %2nd tree
(11,7)*{}="b_2'",
(-6,-10)*{}="c_1",
(-0.5,-10)*{}="c_2",
\ar @{-} "a";"0" <0pt>
\ar @{-} "a";"b_1" <0pt>
\ar @{-} "a";"c_3" <0pt>
\ar @{-} "a";"c_1" <0pt>
\ar @{-} "b_1";"c_2" <0pt>
\ar @{-} "a'";"0'" <0pt>
\ar @{-} "a'";"b_1'" <0pt>
\ar @{-} "a'";"b_2'" <0pt>
\ar @{-} "1L";"1R" <0pt>
\ar @{.} "2L";"2R" <0pt>
\endxy}\Ea$
correspond, respectively,  to the following $1$-dimensional cells,

$$
\lim_{\var\rar 0} \left\{\Ba{c}(x_1=0,x_2=\var,x_3=1)\\ (y_1=-y,y_2=y)
\Ea\right.\simeq (0,+\infty)  \ \  \ \ \ \ \ \ \ \
\lim_{\var\rar 0} \left\{\Ba{c}(x_1=0,x_2=1-\var_1,x_3=1)\\ (y_1=-y,y_2=y)\Ea\right.\simeq (0,+\infty)
$$

7) The zone functions $\Ba{c}\resizebox{12mm}{!}{\xy
(-4,0)*{}="1L",
(12,0)*{}="1R",
(-4,-5)*{}="2L",
(12,-5)*{}="2R",
(-4,-3)*{}="3L",
(12,-3)*{}="3R",
(0,4)*{}="0",
 (0,0)*{\circ}="a",
(-3,-5)*{\circ}="b_1",
(5.5,-10)*{}="b_2",
(8,-9)*{}="0'",
 (8,-3)*{\circ}="a'",
(5,4)*{}="b_1'",                %2nd tree
(11,4)*{}="b_2'",
(-6,-10)*{}="c_1",
(-0.5,-10)*{}="c_2",
\ar @{-} "a";"0" <0pt>
\ar @{-} "a";"b_1" <0pt>
\ar @{-} "a";"b_2" <0pt>
\ar @{-} "b_1";"c_1" <0pt>
\ar @{-} "b_1";"c_2" <0pt>
\ar @{-} "a'";"0'" <0pt>
\ar @{-} "a'";"b_1'" <0pt>
\ar @{-} "a'";"b_2'" <0pt>
\ar @{.} "1L";"1R" <0pt>
\ar @{.} "2L";"2R" <0pt>
\ar @{.} "3L";"3R" <0pt>
\endxy}\Ea$
and, respectively,
$\Ba{c}\resizebox{12mm}{!}{\xy
(-4,0)*{}="1L",
(12,0)*{}="1R",
(-4,-5)*{}="2L",
(12,-5)*{}="2R",
(-4,-3)*{}="3L",
(12,-3)*{}="3R",
(0,4)*{}="0",
 (0,0)*{\circ}="a",
(3,-5)*{\circ}="b_1",
(5.5,-10)*{}="c_3",
(8,-9)*{}="0'",
 (8,-3)*{\circ}="a'",
(5,4)*{}="b_1'",                %2nd tree
(11,4)*{}="b_2'",
(-6,-10)*{}="c_1",
(-0.5,-10)*{}="c_2",
\ar @{-} "a";"0" <0pt>
\ar @{-} "a";"b_1" <0pt>
\ar @{-} "a";"c_3" <0pt>
\ar @{-} "a";"c_1" <0pt>
\ar @{-} "b_1";"c_2" <0pt>
\ar @{-} "a'";"0'" <0pt>
\ar @{-} "a'";"b_1'" <0pt>
\ar @{-} "a'";"b_2'" <0pt>
\ar @{.} "1L";"1R" <0pt>
\ar @{.} "2L";"2R" <0pt>
\ar @{.} "3L";"3R" <0pt>
\endxy}\Ea$
correspond, respectively,  to the following points in $\overline{C_{3,2}}(\R\times \R)$,
$$
\lim_{\var_1,\var_2\rar 0} \left\{\Ba{c}(x_1=0,x_2=\var_1\var_2,x_3=1)\\ (y_1=-\frac{1}{\var_2},y_2=\frac{1}{\var_2})\Ea\right. \ \  \ \ \ \ \ \ \ \
\lim_{\var_2\rar 0} \left\{\Ba{c}(x_1=0,x_2=1-\var_1\var_2\,x_3=1)\\ (y_1=-\frac{1}{\var_2},y_2=\frac{1}{\var_2})\Ea\right.
$$

8) The zone functions $\Ba{c}\resizebox{12mm}{!}{\xy
(-4,0)*{}="1L",
(12,0)*{}="1R",
(-5,-5)*{}="2L",
(12,-5)*{}="2R",
(0,4)*{}="0",
 (0,0)*{\circ}="a",
(-3,-5)*{\circ}="b_1",
(5.5,-10)*{}="b_2",
(8,-11)*{}="0'",
 (8,-5)*{\circ}="a'",
(5,2)*{}="b_1'",                %2nd tree
(11,2)*{}="b_2'",
(-6,-10)*{}="c_1",
(-0.5,-10)*{}="c_2",
\ar @{-} "a";"0" <0pt>
\ar @{-} "a";"b_1" <0pt>
\ar @{-} "a";"b_2" <0pt>
\ar @{-} "b_1";"c_1" <0pt>
\ar @{-} "b_1";"c_2" <0pt>
\ar @{-} "a'";"0'" <0pt>
\ar @{-} "a'";"b_1'" <0pt>
\ar @{-} "a'";"b_2'" <0pt>
\ar @{.} "1L";"1R" <0pt>
\ar @{-} "2L";"2R" <0pt>
\endxy}\Ea$
and, respectively,
$\Ba{c}\resizebox{12mm}{!}{\xy
(-4,0)*{}="1L",
(12,0)*{}="1R",
(-5,-5)*{}="2L",
(12,-5)*{}="2R",
(0,4)*{}="0",
 (0,0)*{\circ}="a",
(3,-5)*{\circ}="b_1",
(5.5,-10)*{}="c_3",
(8,-11)*{}="0'",
 (8,-5)*{\circ}="a'",
(5,2)*{}="b_1'",                %2nd tree
(11,2)*{}="b_2'",
(-6,-10)*{}="c_1",
(-0.5,-10)*{}="c_2",
\ar @{-} "a";"0" <0pt>
\ar @{-} "a";"b_1" <0pt>
\ar @{-} "a";"c_3" <0pt>
\ar @{-} "a";"c_1" <0pt>
\ar @{-} "b_1";"c_2" <0pt>
\ar @{-} "a'";"0'" <0pt>
\ar @{-} "a'";"b_1'" <0pt>
\ar @{-} "a'";"b_2'" <0pt>
\ar @{.} "1L";"1R" <0pt>
\ar @{-} "2L";"2R" <0pt>
\endxy}\Ea$
correspond, respectively,  to the following $1$-dimensional cells,
$$
\lim_{\var\rar 0} \left\{\Ba{c}(x_1=0,x_2=\var ,x_3=1)\\ (y_1=-\frac{y}{\var},y_2=\frac{y}{\var})\Ea\right.\simeq (0,+\infty)  \ \  \ \ \ \ \ \ \ \
\lim_{\var\rar 0} \left\{\Ba{c}(x_1=0,x_2=1-\var ,x_3=1)\\ (y_1=-\frac{y}{\var},y_2=\frac{y}{\var})\Ea\right.\simeq (0,+\infty)
$$

This list exhaust all possible natural strata of and all possible triples $(T^\uparrow\in \cT ree_3, T_\downarrow\in \cT ree_2, \zeta)$. The stratification formula (\ref{2: C m,n stratification formula})
holds true in this case. Not surprisingly,   $\overline{C_{3,2}}(\R\times\R )$ is the hexagon from the multiplihedra family \cite{Ma,SU}

$$
\Ba{c}
\resizebox{90mm}{!}{\xy
%center
(0,16)*{\resizebox{12mm}{!}{\xy
(-4,-3)*{}="1L",
(12,-3)*{}="1R",
(0,2)*{}="0",
 (0,-3)*{\circ}="a",
(-4,-8)*{}="b_1",
(0,-8)*{}="b_2",
(4,-8)*{}="b_3",
(8,-8)*{}="0'",
 (8,-3)*{\circ}="a'",
(5,2)*{}="b_1'", %2nd tree
(11,2)*{}="b_2'",
\ar @{-} "a";"0" <0pt>
\ar @{-} "a";"b_1" <0pt>
\ar @{-} "a";"b_2" <0pt>
\ar @{-} "a";"b_3" <0pt>
\ar @{-} "1L";"1R" <0pt>
%\ar @{-} "2L";"2R" <0pt>
\ar @{-} "a'";"0'" <0pt>
\ar @{-} "a'";"b_1'" <0pt>
\ar @{-} "a'";"b_2'" <0pt>
\endxy}},
%
%down edge
%
(0,-7)*{\resizebox{12mm}{!}{\xy
(-4,-2)*{}="1L",
(12,-2)*{}="1R",
(-4,-5)*{}="2L",
(12,-5)*{}="2R",
(0,0)*{}="0",
 (0,-5)*{\circ}="a",
(-4,-10)*{}="b_1",
(0,-10)*{}="b_2",
(4,-10)*{}="b_3",
(8,-7)*{}="0'",
 (8,-2)*{\circ}="a'",
(5,3)*{}="b_1'", %2nd tree
(11,3)*{}="b_2'",
\ar @{-} "a";"0" <0pt>
\ar @{-} "a";"b_1" <0pt>
\ar @{-} "a";"b_2" <0pt>
\ar @{-} "a";"b_3" <0pt>
\ar @{.} "1L";"1R" <0pt>
\ar @{.} "2L";"2R" <0pt>
\ar @{-} "a'";"0'" <0pt>
\ar @{-} "a'";"b_1'" <0pt>
\ar @{-} "a'";"b_2'" <0pt>
\endxy}},
%
%upper edge
%
(0,38)*{\resizebox{12mm}{!}{\xy
(-4,-2)*{}="1L",
(12,-2)*{}="1R",
(-4,-5)*{}="2L",
(12,-5)*{}="2R",
(0,3)*{}="0",
 (0,-2)*{\circ}="a",
(-4,-7)*{}="b_1",
(0,-7)*{}="b_2",
(4,-7)*{}="b_3",
(8,-10)*{}="0'",
 (8,-5)*{\circ}="a'",
(5,0)*{}="b_1'",      %2nd tree
(11,0)*{}="b_2'",
\ar @{-} "a";"0" <0pt>
\ar @{-} "a";"b_1" <0pt>
\ar @{-} "a";"b_2" <0pt>
\ar @{-} "a";"b_3" <0pt>
\ar @{.} "1L";"1R" <0pt>
\ar @{.} "2L";"2R" <0pt>
\ar @{-} "a'";"0'" <0pt>
\ar @{-} "a'";"b_1'" <0pt>
\ar @{-} "a'";"b_2'" <0pt>
\endxy}},
%
%right down edge
%
(27,3)*{\resizebox{12mm}{!}{\xy
(-4,0)*{}="1L",
(12,0)*{}="1R",
(-5,-5)*{}="2L",
(12,-5)*{}="2R",
(0,4)*{}="0",
 (0,0)*{\circ}="a",
(-3,-5)*{\circ}="b_1",
(5.5,-10)*{}="b_2",
(8,-6)*{}="0'",
 (8,0)*{\circ}="a'",
(5,7)*{}="b_1'",                %2nd tree
(11,7)*{}="b_2'",
(-6,-10)*{}="c_1",
(-0.5,-10)*{}="c_2",
\ar @{-} "a";"0" <0pt>
\ar @{-} "a";"b_1" <0pt>
\ar @{-} "a";"b_2" <0pt>
\ar @{-} "b_1";"c_1" <0pt>
\ar @{-} "b_1";"c_2" <0pt>
\ar @{-} "a'";"0'" <0pt>
\ar @{-} "a'";"b_1'" <0pt>
\ar @{-} "a'";"b_2'" <0pt>
\ar @{-} "1L";"1R" <0pt>
\ar @{.} "2L";"2R" <0pt>
\endxy}},
%
%left down edge
%
(-25,3)*{\resizebox{12mm}{!}{\xy
(-4,0)*{}="1L",
(12,0)*{}="1R",
(-4,-5)*{}="2L",
(12,-5)*{}="2R",
(0,4)*{}="0",
 (0,0)*{\circ}="a",
(3,-5)*{\circ}="b_1",
(5.5,-10)*{}="c_3",
(8,-6)*{}="0'",
 (8,0)*{\circ}="a'",
(5,7)*{}="b_1'",                %2nd tree
(11,7)*{}="b_2'",
(-6,-10)*{}="c_1",
(-0.5,-10)*{}="c_2",
\ar @{-} "a";"0" <0pt>
\ar @{-} "a";"b_1" <0pt>
\ar @{-} "a";"c_3" <0pt>
\ar @{-} "a";"c_1" <0pt>
\ar @{-} "b_1";"c_2" <0pt>
\ar @{-} "a'";"0'" <0pt>
\ar @{-} "a'";"b_1'" <0pt>
\ar @{-} "a'";"b_2'" <0pt>
\ar @{-} "1L";"1R" <0pt>
\ar @{.} "2L";"2R" <0pt>
\endxy}},
%
%right up edge
%
(25,29)*{\resizebox{12mm}{!}{\xy
(-4,0)*{}="1L",
(12,0)*{}="1R",
(-5,-5)*{}="2L",
(12,-5)*{}="2R",
(0,4)*{}="0",
 (0,0)*{\circ}="a",
(-3,-5)*{\circ}="b_1",
(5.5,-10)*{}="b_2",
(8,-11)*{}="0'",
 (8,-5)*{\circ}="a'",
(5,2)*{}="b_1'",                %2nd tree
(11,2)*{}="b_2'",
(-6,-10)*{}="c_1",
(-0.5,-10)*{}="c_2",
\ar @{-} "a";"0" <0pt>
\ar @{-} "a";"b_1" <0pt>
\ar @{-} "a";"b_2" <0pt>
\ar @{-} "b_1";"c_1" <0pt>
\ar @{-} "b_1";"c_2" <0pt>
\ar @{-} "a'";"0'" <0pt>
\ar @{-} "a'";"b_1'" <0pt>
\ar @{-} "a'";"b_2'" <0pt>
\ar @{.} "1L";"1R" <0pt>
\ar @{-} "2L";"2R" <0pt>
\endxy}},
%
%left up edge
%
(-27,28)*{\resizebox{12mm}{!}{\xy
(-4,0)*{}="1L",
(12,0)*{}="1R",
(-5,-5)*{}="2L",
(12,-5)*{}="2R",
(0,4)*{}="0",
 (0,0)*{\circ}="a",
(3,-5)*{\circ}="b_1",
(5.5,-10)*{}="c_3",
(8,-11)*{}="0'",
 (8,-5)*{\circ}="a'",
(5,2)*{}="b_1'",                %2nd tree
(11,2)*{}="b_2'",
(-6,-10)*{}="c_1",
(-0.5,-10)*{}="c_2",
\ar @{-} "a";"0" <0pt>
\ar @{-} "a";"b_1" <0pt>
\ar @{-} "a";"c_3" <0pt>
\ar @{-} "a";"c_1" <0pt>
\ar @{-} "b_1";"c_2" <0pt>
\ar @{-} "a'";"0'" <0pt>
\ar @{-} "a'";"b_1'" <0pt>
\ar @{-} "a'";"b_2'" <0pt>
\ar @{.} "1L";"1R" <0pt>
\ar @{-} "2L";"2R" <0pt>
\endxy}},
%
%bottom left vertex
%
(-31,-26)*{\resizebox{12mm}{!}{
\xy
(-4,6)*{}="1L",
(12,6)*{}="1R",
(-4,-3)*{}="2L",
(12,-3)*{}="2R",
(0,4)*{}="0",
 (0,0)*{\circ}="a",
(-3,-5)*{\circ}="b_1",
(5.5,-10)*{}="b_2",
(8,0)*{}="0'",
 (8,6)*{\circ}="a'",
(5,13)*{}="b_1'",                %2nd tree
(11,13)*{}="b_2'",
(-6,-10)*{}="c_1",
(-0.5,-10)*{}="c_2",
\ar @{-} "a";"0" <0pt>
\ar @{-} "a";"b_1" <0pt>
\ar @{-} "a";"b_2" <0pt>
\ar @{-} "b_1";"c_1" <0pt>
\ar @{-} "b_1";"c_2" <0pt>
\ar @{-} "a'";"0'" <0pt>
\ar @{-} "a'";"b_1'" <0pt>
\ar @{-} "a'";"b_2'" <0pt>
\ar @{.} "1L";"1R" <0pt>
\ar @{.} "2L";"2R" <0pt>
\endxy
}},
%
%bottom right vertex
%
(29,-22)*{\resizebox{12mm}{!}{\xy
(-4,6)*{}="1L",
(12,6)*{}="1R",
(-4,-3)*{}="2L",
(12,-3)*{}="2R",
(0,4)*{}="0",
 (0,0)*{\circ}="a",
(-3,-5)*{\circ}="b_1",
(5.5,-10)*{}="b_2",
(8,0)*{}="0'",
 (8,6)*{\circ}="a'",
(5,13)*{}="b_1'",                %2nd tree
(11,13)*{}="b_2'",
(-6,-10)*{}="c_1",
(-0.5,-10)*{}="c_2",
\ar @{-} "a";"0" <0pt>
\ar @{-} "a";"b_1" <0pt>
\ar @{-} "a";"b_2" <0pt>
\ar @{-} "b_1";"c_1" <0pt>
\ar @{-} "b_1";"c_2" <0pt>
\ar @{-} "a'";"0'" <0pt>
\ar @{-} "a'";"b_1'" <0pt>
\ar @{-} "a'";"b_2'" <0pt>
\ar @{.} "1L";"1R" <0pt>
\ar @{.} "2L";"2R" <0pt>
\endxy}},
%
%right most vertex
(54,16)*{\resizebox{12mm}{!}{\xy
(-4,0)*{}="1L",
(12,0)*{}="1R",
(-4,-5)*{}="2L",
(12,-5)*{}="2R",
(-4,-3)*{}="3L",
(12,-3)*{}="3R",
(0,4)*{}="0",
 (0,0)*{\circ}="a",
(-3,-5)*{\circ}="b_1",
(5.5,-10)*{}="b_2",
(8,-9)*{}="0'",
 (8,-3)*{\circ}="a'",
(5,4)*{}="b_1'",                %2nd tree
(11,4)*{}="b_2'",
(-6,-10)*{}="c_1",
(-0.5,-10)*{}="c_2",
\ar @{-} "a";"0" <0pt>
\ar @{-} "a";"b_1" <0pt>
\ar @{-} "a";"b_2" <0pt>
\ar @{-} "b_1";"c_1" <0pt>
\ar @{-} "b_1";"c_2" <0pt>
\ar @{-} "a'";"0'" <0pt>
\ar @{-} "a'";"b_1'" <0pt>
\ar @{-} "a'";"b_2'" <0pt>
\ar @{.} "1L";"1R" <0pt>
\ar @{.} "2L";"2R" <0pt>
\ar @{.} "3L";"3R" <0pt>
\endxy}},
%
%left most vertex
%
(-56,16)*{\resizebox{12mm}{!}{\xy
(-4,0)*{}="1L",
(12,0)*{}="1R",
(-4,-5)*{}="2L",
(12,-5)*{}="2R",
(-4,-3)*{}="3L",
(12,-3)*{}="3R",
(0,4)*{}="0",
 (0,0)*{\circ}="a",
(3,-5)*{\circ}="b_1",
(5.5,-10)*{}="c_3",
(8,-9)*{}="0'",
 (8,-3)*{\circ}="a'",
(5,4)*{}="b_1'",                %2nd tree
(11,4)*{}="b_2'",
(-6,-10)*{}="c_1",
(-0.5,-10)*{}="c_2",
\ar @{-} "a";"0" <0pt>
\ar @{-} "a";"b_1" <0pt>
\ar @{-} "a";"c_3" <0pt>
\ar @{-} "a";"c_1" <0pt>
\ar @{-} "b_1";"c_2" <0pt>
\ar @{-} "a'";"0'" <0pt>
\ar @{-} "a'";"b_1'" <0pt>
\ar @{-} "a'";"b_2'" <0pt>
\ar @{.} "1L";"1R" <0pt>
\ar @{.} "2L";"2R" <0pt>
\ar @{.} "3L";"3R" <0pt>
\endxy}},
%
%right up vertex
%
(27,54)*{\resizebox{12mm}{!}{\xy
(-4,-7)*{}="1L",
(12,-7)*{}="1R",
(-4,-3)*{}="2L",
(12,-3)*{}="2R",
(0,4)*{}="0",
 (0,0)*{\circ}="a",
(-3,-5)*{\circ}="b_1",
(5.5,-10)*{}="b_2",
(8,-13)*{}="0'",
 (8,-7)*{\circ}="a'",
(5,0)*{}="b_1'",                %2nd tree
(11,0)*{}="b_2'",
(-6,-10)*{}="c_1",
(-0.5,-10)*{}="c_2",
\ar @{-} "a";"0" <0pt>
\ar @{-} "a";"b_1" <0pt>
\ar @{-} "a";"b_2" <0pt>
\ar @{-} "b_1";"c_1" <0pt>
\ar @{-} "b_1";"c_2" <0pt>
\ar @{-} "a'";"0'" <0pt>
\ar @{-} "a'";"b_1'" <0pt>
\ar @{-} "a'";"b_2'" <0pt>
\ar @{.} "1L";"1R" <0pt>
\ar @{.} "2L";"2R" <0pt>
\endxy}},
%
%left up vertex
%
(-24,54)*{\resizebox{12mm}{!}{\xy
(-4,-7)*{}="1L",
(12,-7)*{}="1R",
(-4,-3)*{}="2L",
(12,-3)*{}="2R",
(0,4)*{}="0",
 (0,0)*{\circ}="a",
(3,-5)*{\circ}="b_1",
(5.5,-10)*{}="c_3",
(8,-13)*{}="0'",
 (8,-7)*{\circ}="a'",
(5,0)*{}="b_1'",                %2nd tree
(11,0)*{}="b_2'",
(-6,-10)*{}="c_1",
(-0.5,-10)*{}="c_2",
\ar @{-} "a";"0" <0pt>
\ar @{-} "a";"b_1" <0pt>
\ar @{-} "a";"c_3" <0pt>
\ar @{-} "a";"c_1" <0pt>
\ar @{-} "b_1";"c_2" <0pt>
\ar @{-} "a'";"0'" <0pt>
\ar @{-} "a'";"b_1'" <0pt>
\ar @{-} "a'";"b_2'" <0pt>
\ar @{.} "1L";"1R" <0pt>
\ar @{.} "2L";"2R" <0pt>
\endxy}},
%polytpoe
%
(-25,-15)*{\bu}="1",
 (25,-15)*{\bu}="2",
 (45,15)*{\bu}="3",
 (25,45)*{\bu}="4",
 (-25,45)*{\bu}="6",
(-45,15)*{\bu}="7",
\ar @{-} "1";"2" <0pt>
\ar @{-} "2";"3" <0pt>
\ar @{-} "3";"4" <0pt>
\ar @{-} "4";"6" <0pt>
\ar @{-} "6";"7" <0pt>
\ar @{-} "7";"1" <0pt>
\endxy}\\
\\
\mbox{\sc{Fig.}\ 1:\ Biassociahedron $\sK_3^2$}
\Ea
\ \ \ \ \ \ \ \ \ \ \ \ \ \
\Ba{c}
\resizebox{60mm}{!}{\xy
%center
(0,16)*{\resizebox{9mm}{!}{\xy
(3,2)*{}="u1",
(-3,2)*{}="u2",
 (0,-3)*{\circ}="a",
(-4,-8)*{}="b_1",
(0,-8)*{}="b_2",
(4,-8)*{}="b_3",
\ar @{.} "a";"b_1" <0pt>
\ar @{.} "a";"b_2" <0pt>
\ar @{.} "a";"b_3" <0pt>
\ar @{.} "a";"u1" <0pt>
\ar @{.} "a";"u2" <0pt>
\endxy}},
%
%down edge
%
(0,-7)*{\resizebox{7mm}{!}{\xy
 (0,-5)*{\circ}="a",
(-4,-10)*{}="b_1",
(0,-10)*{}="b_2",
(4,-10)*{}="b_3",
 (0,0)*{\circ}="a'",
(-3,4)*{}="b_1'", %2nd tree
(3,4)*{}="b_2'",
\ar @{.} "a";"a'" <0pt>
\ar @{.} "a";"b_1" <0pt>
\ar @{.} "a";"b_2" <0pt>
\ar @{.} "a";"b_3" <0pt>
\ar @{.} "a'";"b_1'" <0pt>
\ar @{.} "a'";"b_2'" <0pt>
\endxy}},
%
%upper edge
%
(0,40)*{\resizebox{12mm}{!}{\xy
(-12,0)*{}="1L",
(12,0)*{}="1R",
(-4,7)*{}="0",
(5,12)*{}="0",
 (5,7)*{\circ}="a",
(1,2)*{}="b_1",
(5,2)*{}="b_2",
(9,2)*{}="b_3",
(-5,12)*{}="0'",
 (-5,7)*{\circ}="a'",
(-9,2)*{}="b_1'",
(-5,2)*{}="b_2'",
(-1,2)*{}="b_3'",
(0,-12)*{}="01",
 (0,-7)*{\circ}="a1",
(-3,-2)*{}="b_11",      %middle loe tree
(3,-2)*{}="b_21",
(-7,-12)*{}="02",
 (-7,-7)*{\circ}="a2",
(-10,-2)*{}="b_12",      %left loe tree
(-3,-2)*{}="b_22",
(7,-12)*{}="03",
 (7,-7)*{\circ}="a3",
(10,-2)*{}="b_13",      %right loe tree
(3,-2)*{}="b_23",
\ar @{.} "a";"0" <0pt>
\ar @{.} "a";"b_1" <0pt>
\ar @{.} "a";"b_2" <0pt>
\ar @{.} "a";"b_3" <0pt>
\ar @{.} "a'";"0'" <0pt>
\ar @{.} "a'";"b_1'" <0pt>
\ar @{.} "a'";"b_2'" <0pt>
\ar @{.} "a'";"b_3'" <0pt>
\ar @{.} "1L";"1R" <0pt>
\ar @{.} "a1";"01" <0pt>
\ar @{.} "a1";"b_11" <0pt>
\ar @{.} "a1";"b_21" <0pt>
\ar @{.} "a2";"02" <0pt>
\ar @{.} "a2";"b_12" <0pt>
\ar @{.} "a2";"b_22" <0pt>
\ar @{.} "a3";"03" <0pt>
\ar @{.} "a3";"b_13" <0pt>
\ar @{.} "a3";"b_23" <0pt>
\endxy}},
%
%right down edge
%
(29,3)*{\resizebox{10mm}{!}{
\xy
 (0,0)*{\circ}="a",
(-3,-5)*{\circ}="b_1",
(5.5,-10)*{}="b_2",
(-3,5)*{}="b_1'",                %2nd tree
(3,5)*{}="b_2'",
(-6,-10)*{}="c_1",
(-0.5,-10)*{}="c_2",
\ar @{.} "a";"b_1" <0pt>
\ar @{.} "a";"b_2" <0pt>
\ar @{.} "b_1";"c_1" <0pt>
\ar @{.} "b_1";"c_2" <0pt>
\ar @{.} "a";"b_1'" <0pt>
\ar @{.} "a";"b_2'" <0pt>
\endxy}},
%
%left down edge
%
(-27,3)*{\resizebox{10mm}{!}{
\xy
 (0,0)*{\circ}="a",
(3,-5)*{\circ}="b_1",
(-5.5,-10)*{}="b_2",
(-3,5)*{}="b_1'",                %2nd tree
(3,5)*{}="b_2'",
(6,-10)*{}="c_1",
(0.5,-10)*{}="c_2",
\ar @{.} "a";"b_1" <0pt>
\ar @{.} "a";"b_2" <0pt>
\ar @{.} "b_1";"c_1" <0pt>
\ar @{.} "b_1";"c_2" <0pt>
\ar @{.} "a";"b_1'" <0pt>
\ar @{.} "a";"b_2'" <0pt>
\endxy
}},
%
%right up edge
%
(27,29)*{\resizebox{13mm}{!}{
\xy
(-10,0)*{}="1L",
(12,0)*{}="1R",
(4,10)*{}="0",
 (4,6)*{\circ}="a",
(1,2)*{}="u_1",
(7,2)*{}="u_2",
(-4,10)*{}="0'",
 (-4,6)*{\circ}="a'",
(-1,2)*{}="u_1'",
(-7,2)*{}="u_2'",
(-1,-2)*{}="du1",
(-7,-2)*{}="du2",
(-4,-6)*{\circ}="v",
% (-4,-10)*{\circ}="vd", %left under line
(-1,-10)*{}="dd1",
(-7,-10)*{}="dd2",
(4,-10)*{}="xd",
 (4,-6)*{\circ}="x",
(1,-2)*{}="x_1",   %right under line
(7,-2)*{}="x_2",
\ar @{.} "a";"0" <0pt>
\ar @{.} "a";"u_1" <0pt>
\ar @{.} "a";"u_2" <0pt>
\ar @{.} "a'";"0'" <0pt>
\ar @{.} "a'";"u_1'" <0pt>
\ar @{.} "a'";"u_2'" <0pt>
\ar @{.} "v";"du1" <0pt>
\ar @{.} "v";"du2" <0pt>
\ar @{.} "v";"dd1" <0pt>
\ar @{.} "v";"dd2" <0pt>
\ar @{.} "x";"xd" <0pt>
\ar @{.} "x";"x_1" <0pt>
\ar @{.} "x";"x_2" <0pt>
\ar @{.} "1L";"1R" <0pt>
\endxy
}},
%
%left up edge
%
(-29,28)*{\resizebox{13mm}{!}{
\xy
(-10,0)*{}="1L",
(12,0)*{}="1R",
(4,10)*{}="0",
 (4,6)*{\circ}="a",
(1,2)*{}="u_1",
(7,2)*{}="u_2",
(-4,10)*{}="0'",
 (-4,6)*{\circ}="a'",
(-1,2)*{}="u_1'",
(-7,2)*{}="u_2'",
(1,-2)*{}="du1",
(7,-2)*{}="du2",
(4,-6)*{\circ}="v",
% (-4,-10)*{\circ}="vd", %left under line
(1,-10)*{}="dd1",
(7,-10)*{}="dd2",
(-4,-10)*{}="xd",
 (-4,-6)*{\circ}="x",
(-1,-2)*{}="x_1",   %right under line
(-7,-2)*{}="x_2",
\ar @{.} "a";"0" <0pt>
\ar @{.} "a";"u_1" <0pt>
\ar @{.} "a";"u_2" <0pt>
\ar @{.} "a'";"0'" <0pt>
\ar @{.} "a'";"u_1'" <0pt>
\ar @{.} "a'";"u_2'" <0pt>
\ar @{.} "v";"du1" <0pt>
\ar @{.} "v";"du2" <0pt>
\ar @{.} "v";"dd1" <0pt>
\ar @{.} "v";"dd2" <0pt>
\ar @{.} "x";"xd" <0pt>
\ar @{.} "x";"x_1" <0pt>
\ar @{.} "x";"x_2" <0pt>
\ar @{.} "1L";"1R" <0pt>
\endxy
}},
%
%bottom left vertex
%
(-28,-22)*{\resizebox{9mm}{!}{
\xy
(0,4)*{}="0",
 (0,0)*{\circ}="a",
(3,-5)*{\circ}="b_1",
(5.5,-10)*{}="b_2",
 (0,5)*{\circ}="a'",
(-3,10)*{}="b_1'",                %2nd tree
(3,10)*{}="b_2'",
(-6,-10)*{}="c_1",
(-0.5,-10)*{}="c_2",
\ar @{.} "a";"a'" <0pt>
\ar @{.} "a";"b_2" <0pt>
\ar @{.} "a";"c_1" <0pt>
\ar @{.} "b_1";"c_2" <0pt>
\ar @{.} "a'";"b_1'" <0pt>
\ar @{.} "a'";"b_2'" <0pt>
\endxy
}},
%
%bottom right vertex
%
(28,-20)*{\resizebox{9mm}{!}{
\xy
(0,4)*{}="0",
 (0,0)*{\circ}="a",
(-3,-5)*{\circ}="b_1",
(5.5,-10)*{}="b_2",
 (0,5)*{\circ}="a'",
(-3,10)*{}="b_1'",                %2nd tree
(3,10)*{}="b_2'",
(-6,-10)*{}="c_1",
(-0.5,-10)*{}="c_2",
\ar @{.} "a";"a'" <0pt>
\ar @{.} "a";"b_1" <0pt>
\ar @{.} "a";"b_2" <0pt>
\ar @{.} "b_1";"c_1" <0pt>
\ar @{.} "b_1";"c_2" <0pt>
\ar @{.} "a'";"b_1'" <0pt>
\ar @{.} "a'";"b_2'" <0pt>
\endxy
}},
%
%right most vertex
(52,16)*{\resizebox{12mm}{!}{\xy
(-10,0)*{}="1L",
(12,0)*{}="1R",
(4,10)*{}="0",
 (4,6)*{\circ}="a",
(1,2)*{}="u_1",
(7,2)*{}="u_2",
(-4,10)*{}="0'",
 (-4,6)*{\circ}="a'",
(-1,2)*{}="u_1'",
(-7,2)*{}="u_2'",
(-1,-2)*{}="du1",
(-7,-2)*{}="du2",
(-4,-6)*{\circ}="vu",
 (-4,-10)*{\circ}="vd", %left under line
(-1,-14)*{}="dd1",
(-7,-14)*{}="dd2",
(4,-10)*{}="xd",
 (4,-6)*{\circ}="x",
(1,-2)*{}="x_1",   %right under line
(7,-2)*{}="x_2",
\ar @{.} "a";"0" <0pt>
\ar @{.} "a";"u_1" <0pt>
\ar @{.} "a";"u_2" <0pt>
\ar @{.} "a'";"0'" <0pt>
\ar @{.} "a'";"u_1'" <0pt>
\ar @{.} "a'";"u_2'" <0pt>
\ar @{.} "vd";"vu" <0pt>
\ar @{.} "vu";"du1" <0pt>
\ar @{.} "vu";"du2" <0pt>
\ar @{.} "vd";"dd1" <0pt>
\ar @{.} "vd";"dd2" <0pt>
\ar @{.} "x";"xd" <0pt>
\ar @{.} "x";"x_1" <0pt>
\ar @{.} "x";"x_2" <0pt>
\ar @{-} "1L";"1R" <0pt>
\endxy}},
%
%left most vertex
%
(-52,16)*{\resizebox{12mm}{!}{
\xy
(-10,0)*{}="1L",
(12,0)*{}="1R",
(4,10)*{}="0",
 (4,6)*{\circ}="a",
(1,2)*{}="u_1",
(7,2)*{}="u_2",
(-4,10)*{}="0'",
 (-4,6)*{\circ}="a'",
(-1,2)*{}="u_1'",
(-7,2)*{}="u_2'",
(1,-2)*{}="du1",
(7,-2)*{}="du2",
(4,-6)*{\circ}="vu",
 (4,-10)*{\circ}="vd", %left under line
(1,-14)*{}="dd1",
(7,-14)*{}="dd2",
(-4,-10)*{}="xd",
 (-4,-6)*{\circ}="x",
(-1,-2)*{}="x_1",   %right under line
(-7,-2)*{}="x_2",
\ar @{.} "a";"0" <0pt>
\ar @{.} "a";"u_1" <0pt>
\ar @{.} "a";"u_2" <0pt>
\ar @{.} "a'";"0'" <0pt>
\ar @{.} "a'";"u_1'" <0pt>
\ar @{.} "a'";"u_2'" <0pt>
\ar @{.} "vd";"vu" <0pt>
\ar @{.} "vu";"du1" <0pt>
\ar @{.} "vu";"du2" <0pt>
\ar @{.} "vd";"dd1" <0pt>
\ar @{.} "vd";"dd2" <0pt>
\ar @{.} "x";"xd" <0pt>
\ar @{.} "x";"x_1" <0pt>
\ar @{.} "x";"x_2" <0pt>
\ar @{-} "1L";"1R" <0pt>
\endxy}},
%
%
%right up vertex
%
(25,52)*{\resizebox{12mm}{!}{
\xy
(-14,0)*{}="1L",
(14,0)*{}="1R",
(-6.5,16)*{}="l0",
 (-6.5,12)*{\circ}="la",
(-9.5,7)*{\circ}="lb_1",
(-1,2)*{}="lb_2",
(-12.5,2)*{}="lc_1",
(-7,2)*{}="lc_2",
(7,16)*{}="r0",
 (7,12)*{\circ}="ra",
(4,7)*{\circ}="rb_1",
(12.5,2)*{}="rb_2",
(1,2)*{}="rc_1",
(6.5,2)*{}="rc_2",
(0,-12)*{}="01",
 (0,-7)*{\circ}="a1",
(-3,-2)*{}="b_11",      %middle loe tree
(3,-2)*{}="b_21",
(-7,-12)*{}="02",
 (-7,-7)*{\circ}="a2",
(-10,-2)*{}="b_12",      %left loe tree
(-3,-2)*{}="b_22",
(7,-12)*{}="03",
 (7,-7)*{\circ}="a3",
(10,-2)*{}="b_13",      %right loe tree
(3,-2)*{}="b_23",
\ar @{.} "la";"l0" <0pt>
\ar @{.} "la";"lb_1" <0pt>
\ar @{.} "la";"lb_2" <0pt>
\ar @{.} "lb_1";"lc_1" <0pt>
\ar @{.} "lb_1";"lc_2" <0pt>
\ar @{.} "ra";"r0" <0pt>
\ar @{.} "ra";"rb_1" <0pt>
\ar @{.} "ra";"rb_2" <0pt>
\ar @{.} "rb_1";"rc_1" <0pt>
\ar @{.} "rb_1";"rc_2" <0pt>
\ar @{.} "1L";"1R" <0pt>
\ar @{.} "a1";"01" <0pt>
\ar @{.} "a1";"b_11" <0pt>
\ar @{.} "a1";"b_21" <0pt>
\ar @{.} "a2";"02" <0pt>
\ar @{.} "a2";"b_12" <0pt>
\ar @{.} "a2";"b_22" <0pt>
\ar @{.} "a3";"03" <0pt>
\ar @{.} "a3";"b_13" <0pt>
\ar @{.} "a3";"b_23" <0pt>
\endxy
}},
%left up vertex
%
(-25,53)*{\resizebox{12mm}{!}{
\xy
(-14,0)*{}="1L",
(14,0)*{}="1R",
(6.5,16)*{}="l0",
 (6.5,12)*{\circ}="la",
(9.5,7)*{\circ}="lb_1",
(1,2)*{}="lb_2",
(12.5,2)*{}="lc_1",
(7,2)*{}="lc_2",
(-7,16)*{}="r0",
 (-7,12)*{\circ}="ra",
(-4,7)*{\circ}="rb_1",
(-12.5,2)*{}="rb_2",
(-1,2)*{}="rc_1",
(-6.5,2)*{}="rc_2",
(0,-12)*{}="01",
 (0,-7)*{\circ}="a1",
(-3,-2)*{}="b_11",      %middle loe tree
(3,-2)*{}="b_21",
(-7,-12)*{}="02",
 (-7,-7)*{\circ}="a2",
(-10,-2)*{}="b_12",      %left loe tree
(-3,-2)*{}="b_22",
(7,-12)*{}="03",
 (7,-7)*{\circ}="a3",
(10,-2)*{}="b_13",      %right loe tree
(3,-2)*{}="b_23",
\ar @{.} "la";"l0" <0pt>
\ar @{.} "la";"lb_1" <0pt>
\ar @{.} "la";"lb_2" <0pt>
\ar @{.} "lb_1";"lc_1" <0pt>
\ar @{.} "lb_1";"lc_2" <0pt>
\ar @{.} "ra";"r0" <0pt>
\ar @{.} "ra";"rb_1" <0pt>
\ar @{.} "ra";"rb_2" <0pt>
\ar @{.} "rb_1";"rc_1" <0pt>
\ar @{.} "rb_1";"rc_2" <0pt>
\ar @{-} "1L";"1R" <0pt>
\ar @{.} "a1";"01" <0pt>
\ar @{.} "a1";"b_11" <0pt>
\ar @{.} "a1";"b_21" <0pt>
\ar @{.} "a2";"02" <0pt>
\ar @{.} "a2";"b_12" <0pt>
\ar @{.} "a2";"b_22" <0pt>
\ar @{.} "a3";"03" <0pt>
\ar @{.} "a3";"b_13" <0pt>
\ar @{.} "a3";"b_23" <0pt>
\endxy
}},
%
%polytpoe
%
(-25,-15)*{\bu}="1",
 (25,-15)*{\bu}="2",
 (45,15)*{\bu}="3",
 (25,45)*{\bu}="4",
 (-25,45)*{\bu}="6",
(-45,15)*{\bu}="7",
\ar @{-} "1";"2" <0pt>
\ar @{-} "2";"3" <0pt>
\ar @{-} "3";"4" <0pt>
\ar @{-} "4";"6" <0pt>
\ar @{-} "6";"7" <0pt>
\ar @{-} "7";"1" <0pt>
\endxy}\\
\\
\mbox{{\sc Fig}.\ 2:\  $r_3^2\left(\cF Chains(\sK_3^2)\right)$}
\Ea
$$

%%%%%%%%%%%%%%%%%%%%%%%%%%%%%%%%%%%%%%%%%

\subsection{From biassociahedra to  strongly homotopy bialgebras}
  As we saw in the previous subsection, the biassociahedron $\sK_m^n$ is a smooth manifold with corners which
  comes equipped with a boundary stratification parameterized by 
  Markl's poset $\cK_m^n$. In fact, we constructed $\sK_m^n$  as a closed semi-algebraic subset in the product of copies
of 2-spheres $S^2$ and the intervals $[0,1]$. Hence $\sK_m^n$ comes equipped with a structure of a semialgebraic set (which is finer than just the structure of a smooth manifold with corners). Kontsevich and Soibelman introduced in the Appendix 8 of \cite{KS} a suitable theory of singular chains for such semialgebraic spaces $X$ (see \cite{HLTV} for full details); in this theory
$Chains(X)$ is a vector space of a field $\K$ group generated by
(equivalence classes) of  semialgebraic maps $f:Y\rar X$ from oriented compact semialgebraic spaces $Y$. As in \cite{KS} we assume  that the semialgebraic chain complex $(Chains(X),\p)$ is negatively graded  so that the boundary operator has degree $+1$.

  \sip

  This canonical stratification of the  biassociahedron $\sK_m^n$ in terms of zoned trees gives us (i) an obvious $\frac{1}{2}$-structure on the collection of dg $\bS$-bimodules $\{Chains(\sK_m^n)\}_{m,n\in \N}$, and (ii) a $\frac{1}{2}$-subprop
  $\cF Chains(\sK_m^n)\subset Chains(\sK_m^n)$ spanned by fundamental chains which is called  the dg $\frac{1}{2}$-prop of {\em  of fundamental or cellular chains}\, of the biassociahedron. Unfortunately,
  the $\bS$-submodule $\{\cF Chains(\sK_m^n)\}_{m,n\in \N}$ is not a prop.
  Martin Markl constructed by induction  a collection $r=\{r_m^n\}$ of linear {\em monomorphisms} of graded vector spaces in \cite{Ma},
  $$
  r_m^n:  \cF Chains(\sK_m^n)\} \hook \cA ss \cB_\infty,
  $$
  The image under $r_3^2$ of generators of
  $\cF Chains(\sK_3^2)$ is given in \mbox{{\sc Fig.} \hspace{-2mm} 2}.
As we see from this example, the monomorphism $r$ is not even homogeneous:
the upper edge of $\sK_3^2$ (which is a degree $-1$ element in
$\cF Chains(\sK_3^2)$) gets mapped into a degree $-2$ element\footnote{We also use here fraction notations for elements of $\cA ss \cB_\infty$ introduced in \cite{Ma1}.}
$\Ba{c}{\resizebox{11mm}{!}{\xy
(-12,0)*{}="1L",
(12,0)*{}="1R",
(-4,7)*{}="0",
(5,12)*{}="0",
 (5,7)*{\circ}="a",
(1,2)*{}="b_1",
(5,2)*{}="b_2",
(9,2)*{}="b_3",
(-5,12)*{}="0'",
 (-5,7)*{\circ}="a'",
(-9,2)*{}="b_1'",
(-5,2)*{}="b_2'",
(-1,2)*{}="b_3'",
(0,-12)*{}="01",
 (0,-7)*{\circ}="a1",
(-3,-2)*{}="b_11",      %middle loe tree
(3,-2)*{}="b_21",
(-7,-12)*{}="02",
 (-7,-7)*{\circ}="a2",
(-10,-2)*{}="b_12",      %left loe tree
(-3,-2)*{}="b_22",
(7,-12)*{}="03",
 (7,-7)*{\circ}="a3",
(10,-2)*{}="b_13",      %right loe tree
(3,-2)*{}="b_23",
\ar @{.} "a";"0" <0pt>
\ar @{.} "a";"b_1" <0pt>
\ar @{.} "a";"b_2" <0pt>
\ar @{.} "a";"b_3" <0pt>
\ar @{.} "a'";"0'" <0pt>
\ar @{.} "a'";"b_1'" <0pt>
\ar @{.} "a'";"b_2'" <0pt>
\ar @{.} "a'";"b_3'" <0pt>
\ar @{-} "1L";"1R" <0pt>
\ar @{.} "a1";"01" <0pt>
\ar @{.} "a1";"b_11" <0pt>
\ar @{.} "a1";"b_21" <0pt>
\ar @{.} "a2";"02" <0pt>
\ar @{.} "a2";"b_12" <0pt>
\ar @{.} "a2";"b_22" <0pt>
\ar @{.} "a3";"03" <0pt>
\ar @{.} "a3";"b_13" <0pt>
\ar @{.} "a3";"b_23" <0pt>
\endxy}}\Ea$ in $\cA ss \cB_\infty$. Thus we can not use the map $r$
to make $\cF Chains(\sK_\bu^\bu)$ into a prop (however the collection of maps $\{r_{\bu}^\bu\}$ respects
$\frac{1}{2}$-prop compositions in the dg $\bS$-bimodules $\cF Chains(\sK_\bu^\bu)\}$ and
$\cA ss \cB_\infty$).

\mip

It is not hard to see how the complex
$\cF Chains(\sK_3^2)$ should be modified in order to make the map $r_3^2: \cF Chains(\sK_3^2) \rar \cA ss \cB_\infty
$ into a degree zero morphism of {\em complexes}. One has to subdivide the upper edge of $\sK_3^2$ into the union of two edges by adding a new vertex
in the middle. Equivalently, one has to replace
the degree $-2$ element $\Ba{c}{\resizebox{11mm}{!}{\xy
(-12,0)*{}="1L",
(12,0)*{}="1R",
(-4,7)*{}="0",
(5,12)*{}="0",
 (5,7)*{\circ}="a",
(1,2)*{}="b_1",
(5,2)*{}="b_2",
(9,2)*{}="b_3",
(-5,12)*{}="0'",
 (-5,7)*{\circ}="a'",
(-9,2)*{}="b_1'",
(-5,2)*{}="b_2'",
(-1,2)*{}="b_3'",
(0,-12)*{}="01",
 (0,-7)*{\circ}="a1",
(-3,-2)*{}="b_11",      %middle loe tree
(3,-2)*{}="b_21",
(-7,-12)*{}="02",
 (-7,-7)*{\circ}="a2",
(-10,-2)*{}="b_12",      %left loe tree
(-3,-2)*{}="b_22",
(7,-12)*{}="03",
 (7,-7)*{\circ}="a3",
(10,-2)*{}="b_13",      %right loe tree
(3,-2)*{}="b_23",
\ar @{.} "a";"0" <0pt>
\ar @{.} "a";"b_1" <0pt>
\ar @{.} "a";"b_2" <0pt>
\ar @{.} "a";"b_3" <0pt>
\ar @{.} "a'";"0'" <0pt>
\ar @{.} "a'";"b_1'" <0pt>
\ar @{.} "a'";"b_2'" <0pt>
\ar @{.} "a'";"b_3'" <0pt>
\ar @{-} "1L";"1R" <0pt>
\ar @{.} "a1";"01" <0pt>
\ar @{.} "a1";"b_11" <0pt>
\ar @{.} "a1";"b_21" <0pt>
\ar @{.} "a2";"02" <0pt>
\ar @{.} "a2";"b_12" <0pt>
\ar @{.} "a2";"b_22" <0pt>
\ar @{.} "a3";"03" <0pt>
\ar @{.} "a3";"b_13" <0pt>
\ar @{.} "a3";"b_23" <0pt>
\endxy}}\Ea$
with a {\em sum}\, of degree $-1$ elements,
%upper edge
%
$
{\resizebox{11mm}{!}{\xy
(-12,0)*{}="1L",
(12,0)*{}="1R",
(-4,7)*{\mbox{$\Delta$}}="0",
(0,12)*{}="0",
 (0,7)*{\circ}="a",
(-4,2)*{}="b_1",
(0,2)*{}="b_2",
(4,2)*{}="b_3",
(0,-12)*{}="01",
 (0,-7)*{\circ}="a1",
(-3,-2)*{}="b_11",      %middle loe tree
(3,-2)*{}="b_21",
(-7,-12)*{}="02",
 (-7,-7)*{\circ}="a2",
(-10,-2)*{}="b_12",      %left loe tree
(-3,-2)*{}="b_22",
(7,-12)*{}="03",
 (7,-7)*{\circ}="a3",
(10,-2)*{}="b_13",      %right loe tree
(3,-2)*{}="b_23",
\ar @{.} "a";"0" <0pt>
\ar @{.} "a";"b_1" <0pt>
\ar @{.} "a";"b_2" <0pt>
\ar @{.} "a";"b_3" <0pt>
\ar @{-} "1L";"1R" <0pt>
\ar @{.} "a1";"01" <0pt>
\ar @{.} "a1";"b_11" <0pt>
\ar @{.} "a1";"b_21" <0pt>
\ar @{.} "a2";"02" <0pt>
\ar @{.} "a2";"b_12" <0pt>
\ar @{.} "a2";"b_22" <0pt>
\ar @{.} "a3";"03" <0pt>
\ar @{.} "a3";"b_13" <0pt>
\ar @{.} "a3";"b_23" <0pt>
\endxy}}
$,
where $\Delta $ stands for a $A_\infty$ diagonal \cite{SU2,MS},
$$
\Delta\Ba{c}\resizebox{5mm}{!}
{\xy
(0,12)*{}="0",
 (0,7)*{\circ}="a",
(-4,2)*{}="b_1",
(0,2)*{}="b_2",
(4,2)*{}="b_3",
\ar @{.} "a";"0" <0pt>
\ar @{.} "a";"b_1" <0pt>
\ar @{.} "a";"b_2" <0pt>
\ar @{.} "a";"b_3" <0pt>
\endxy}\Ea = \Ba{c}\resizebox{7mm}{!}{\xy
(0,5)*{}="0",
 (0,0)*{\circ}="a",
(-3,-4)*{\circ}="b_1",
(5.5,-8)*{}="b_2",
(-6,-8)*{}="c_1",
(-0.5,-8)*{}="c_2",
\ar @{.} "a";"0" <0pt>
\ar @{.} "a";"b_1" <0pt>
\ar @{.} "a";"b_2" <0pt>
\ar @{.} "b_1";"c_1" <0pt>
\ar @{.} "b_1";"c_2" <0pt>
\endxy}\Ea
\ot
\Ba{c}\resizebox{5mm}{!}
{\xy
(0,12)*{}="0",
 (0,7)*{\circ}="a",
(-4,2)*{}="b_1",
(0,2)*{}="b_2",
(4,2)*{}="b_3",
\ar @{.} "a";"0" <0pt>
\ar @{.} "a";"b_1" <0pt>
\ar @{.} "a";"b_2" <0pt>
\ar @{.} "a";"b_3" <0pt>
\endxy}\Ea
\ \ \ \
+ \ \ \ \
\Ba{c}\resizebox{5mm}{!}
{\xy
(0,12)*{}="0",
 (0,7)*{\circ}="a",
(-4,2)*{}="b_1",
(0,2)*{}="b_2",
(4,2)*{}="b_3",
\ar @{.} "a";"0" <0pt>
\ar @{.} "a";"b_1" <0pt>
\ar @{.} "a";"b_2" <0pt>
\ar @{.} "a";"b_3" <0pt>
\endxy}\Ea
\ot
 \Ba{c}\resizebox{7mm}{!}{\xy
(0,5)*{}="0",
 (0,0)*{\circ}="a",
(3,-4)*{\circ}="b_1",
(-5.5,-8)*{}="b_2",
(6,-8)*{}="c_1",
(0.5,-8)*{}="c_2",
\ar @{.} "a";"0" <0pt>
\ar @{.} "a";"b_1" <0pt>
\ar @{.} "a";"b_2" <0pt>
\ar @{.} "b_1";"c_1" <0pt>
\ar @{.} "b_1";"c_2" <0pt>
\endxy}\Ea
$$

After this subdivision one reads from $\sK_3^2$ the correct formula
for the value of the differential
in $\cA ss \cB_\infty$,
%%%%%%%%%%%%%%%%%%%% delta on 3,2 corolla %%%%%%%%%%%%%
$$
\delta\ \Ba{c}
\resizebox{4.6mm}{!}
{\xy
(3,2)*{}="u1",
(-3,2)*{}="u2",
 (0,-3)*{\circ}="a",
(-4,-8)*{}="b_1",
(0,-8)*{}="b_2",
(4,-8)*{}="b_3",
\ar @{.} "a";"b_1" <0pt>
\ar @{.} "a";"b_2" <0pt>
\ar @{.} "a";"b_3" <0pt>
\ar @{.} "a";"u1" <0pt>
\ar @{.} "a";"u2" <0pt>
\endxy}\Ea \ = \
%
%down edge
\Ba{c}\resizebox{5.6mm}{!}
{\xy
 (0,-5)*{\circ}="a",
(-4,-10)*{}="b_1",
(0,-10)*{}="b_2",
(4,-10)*{}="b_3",
 (0,0)*{\circ}="a'",
(-3,4)*{}="b_1'", %2nd tree
(3,4)*{}="b_2'",
\ar @{.} "a";"a'" <0pt>
\ar @{.} "a";"b_1" <0pt>
\ar @{.} "a";"b_2" <0pt>
\ar @{.} "a";"b_3" <0pt>
\ar @{.} "a'";"b_1'" <0pt>
\ar @{.} "a'";"b_2'" <0pt>
\endxy}\Ea
\ - \
%
%right down edge
%
\Ba{c}\resizebox{8mm}{!}{
\xy
 (0,0)*{\circ}="a",
(-3,-5)*{\circ}="b_1",
(5.5,-10)*{}="b_2",
(-3,5)*{}="b_1'",                %2nd tree
(3,5)*{}="b_2'",
(-6,-10)*{}="c_1",
(-0.5,-10)*{}="c_2",
\ar @{.} "a";"b_1" <0pt>
\ar @{.} "a";"b_2" <0pt>
\ar @{.} "b_1";"c_1" <0pt>
\ar @{.} "b_1";"c_2" <0pt>
\ar @{.} "a";"b_1'" <0pt>
\ar @{.} "a";"b_2'" <0pt>
\endxy}\Ea
\ +\
%left down edge
%
\Ba{c}\resizebox{8mm}{!}{
\xy
 (0,0)*{\circ}="a",
(3,-5)*{\circ}="b_1",
(-5.5,-10)*{}="b_2",
(-3,5)*{}="b_1'",                %2nd tree
(3,5)*{}="b_2'",
(6,-10)*{}="c_1",
(0.5,-10)*{}="c_2",
\ar @{.} "a";"b_1" <0pt>
\ar @{.} "a";"b_2" <0pt>
\ar @{.} "b_1";"c_1" <0pt>
\ar @{.} "b_1";"c_2" <0pt>
\ar @{.} "a";"b_1'" <0pt>
\ar @{.} "a";"b_2'" <0pt>
\endxy}\Ea
\ - \
%right up edge
%
\Ba{c}\resizebox{14mm}{!}{
\xy
(-10,0)*{}="1L",
(12,0)*{}="1R",
(4,10)*{}="0",
 (4,6)*{\circ}="a",
(1,2)*{}="u_1",
(7,2)*{}="u_2",
(-4,10)*{}="0'",
 (-4,6)*{\circ}="a'",
(-1,2)*{}="u_1'",
(-7,2)*{}="u_2'",
(-1,-2)*{}="du1",
(-7,-2)*{}="du2",
(-4,-6)*{\circ}="v",
% (-4,-10)*{\circ}="vd", %left under line
(-1,-10)*{}="dd1",
(-7,-10)*{}="dd2",
(4,-10)*{}="xd",
 (4,-6)*{\circ}="x",
(1,-2)*{}="x_1",   %right under line
(7,-2)*{}="x_2",
\ar @{.} "a";"0" <0pt>
\ar @{.} "a";"u_1" <0pt>
\ar @{.} "a";"u_2" <0pt>
\ar @{.} "a'";"0'" <0pt>
\ar @{.} "a'";"u_1'" <0pt>
\ar @{.} "a'";"u_2'" <0pt>
\ar @{.} "v";"du1" <0pt>
\ar @{.} "v";"du2" <0pt>
\ar @{.} "v";"dd1" <0pt>
\ar @{.} "v";"dd2" <0pt>
\ar @{.} "x";"xd" <0pt>
\ar @{.} "x";"x_1" <0pt>
\ar @{.} "x";"x_2" <0pt>
\ar @{-} "1L";"1R" <0pt>
\endxy}\Ea
\ +\
%left up edge
%
\Ba{c}
\resizebox{14mm}{!}{
\xy
(-10,0)*{}="1L",
(12,0)*{}="1R",
(4,10)*{}="0",
 (4,6)*{\circ}="a",
(1,2)*{}="u_1",
(7,2)*{}="u_2",
(-4,10)*{}="0'",
 (-4,6)*{\circ}="a'",
(-1,2)*{}="u_1'",
(-7,2)*{}="u_2'",
(1,-2)*{}="du1",
(7,-2)*{}="du2",
(4,-6)*{\circ}="v",
% (-4,-10)*{\circ}="vd", %left under line
(1,-10)*{}="dd1",
(7,-10)*{}="dd2",
(-4,-10)*{}="xd",
 (-4,-6)*{\circ}="x",
(-1,-2)*{}="x_1",   %right under line
(-7,-2)*{}="x_2",
\ar @{.} "a";"0" <0pt>
\ar @{.} "a";"u_1" <0pt>
\ar @{.} "a";"u_2" <0pt>
\ar @{.} "a'";"0'" <0pt>
\ar @{.} "a'";"u_1'" <0pt>
\ar @{.} "a'";"u_2'" <0pt>
\ar @{.} "v";"du1" <0pt>
\ar @{.} "v";"du2" <0pt>
\ar @{.} "v";"dd1" <0pt>
\ar @{.} "v";"dd2" <0pt>
\ar @{.} "x";"xd" <0pt>
\ar @{.} "x";"x_1" <0pt>
\ar @{.} "x";"x_2" <0pt>
\ar @{-} "1L";"1R" <0pt>
\endxy}\Ea
\ - \
%upper edge
\Ba{c}\resizebox{14mm}{!}{\xy
(-12,0)*{}="1L",
(12,0)*{}="1R",
(-4,7)*{\mbox{$\Delta$}}="0",
(0,12)*{}="0",
 (0,7)*{\circ}="a",
(-4,2)*{}="b_1",
(0,2)*{}="b_2",
(4,2)*{}="b_3",
(0,-12)*{}="01",
 (0,-7)*{\circ}="a1",
(-3,-2)*{}="b_11",      %middle loe tree
(3,-2)*{}="b_21",
(-7,-12)*{}="02",
 (-7,-7)*{\circ}="a2",
(-10,-2)*{}="b_12",      %left loe tree
(-3,-2)*{}="b_22",
(7,-12)*{}="03",
 (7,-7)*{\circ}="a3",
(10,-2)*{}="b_13",      %right loe tree
(3,-2)*{}="b_23",
\ar @{.} "a";"0" <0pt>
\ar @{.} "a";"b_1" <0pt>
\ar @{.} "a";"b_2" <0pt>
\ar @{.} "a";"b_3" <0pt>
\ar @{.} "1L";"1R" <0pt>
\ar @{.} "a1";"01" <0pt>
\ar @{.} "a1";"b_11" <0pt>
\ar @{.} "a1";"b_21" <0pt>
\ar @{.} "a2";"02" <0pt>
\ar @{.} "a2";"b_12" <0pt>
\ar @{.} "a2";"b_22" <0pt>
\ar @{.} "a3";"03" <0pt>
\ar @{.} "a3";"b_13" <0pt>
\ar @{.} "a3";"b_23" <0pt>
\endxy}\Ea
$$
on the $(2,3)$-corolla.

\sip

Note that the definition of the $\cA ss_\infty$ diagonal $\Delta$ involves choices so that the best one can hope for is to find a (non-uniquely) defined
cellular refinement, $(\cC ell(\sK_\bu^\bu), \p_{cell})$, of the fundamental chain complex of the biassociahedron  together with a monomorphism complexes
$$
r: \cC ell(\sK_\bu^\bu)\lon \cA ss \cB_\infty
$$
such that the free properad generated by ``big" cells $\sK^m_n$ and equipped with the differential
$\p_{cell}$ can be identified via $r$ with some minimal resolution $\cA ss \cB_\infty$ of $\cA ss\cB$.
The existence of such an intermediate complex
$$
\cF\cC hains(\sK_\bu^\bu) \subset  \cC ell(\sK_\bu^\bu) \subset \cC hains(\sK_\bu^\bu)
$$
was claimed by Samson Saneblidze and Ron Umble in \cite{SU}.

\def\cprime{$'$}


\begin{thebibliography}{10}

\bibitem[B]{B}
T.\ Backman, {\em Configuration spaces, props and wheel-free deformation quantization}, PhD thesis (2016), Stockholm University.

\bibitem[Do]{Do} V.\ Dolgushev, {\em Stable Formality Quasi-isomorphisms for Hochschild Cochains I}. preprint (2011).
arXiv:1109.6031.


\bibitem[Dr]{D}
V.\ Drinfeld,
{\em On some unsolved problems in quantum group theory}.
In: Lecture Notes in Math., Springer,  {\bf 1510} (1992), 1-8.

\bibitem[G]{Ga} G.\ Gaiffi, {\em Compactifications of configuration spaces}, In: Algebraic Geometry Seminars, 1998-1999 (Pisa), Scuola Normale Superiore di Pisa, %Pisa, 1999, pp. 87-109.

\bibitem[EK]{EK} P.\ Etingof and D.\ Kazhdan.
\newblock {\em Quantization of Lie bialgebras, I}.
\newblock{ Selecta Math. (N.S.)} {\bf 2} (1996), 1-41.


\bibitem[GJ]{GJ} E.\ Getzler and J.D.S.\ Jones, {\em Operads, homotopy algebra, and
 iterated integrals for double loop spaces},
 preprint  hep-th/9403055.

\bibitem[GS]{GS} M.\ Gerstenhaber and S.D.\ Schack,
\newblock {\em Bialgebra cohomology, deformations, and quantum groups}.
\newblock{Proc.\ Nat.\ Acad.\ Sci.\ USA, } {\bf 87} (1990), 478-481.


\bibitem[GY]{GY} G. Ginot and S. Yalin,
\newblock {\em Deformation theory of bialgebras, higher Hochschild cohomology and formality}.
\newblock Preprint, arXiv:1606.01504.

\bibitem[HLTV]{HLTV} R.\ M. Hardt, P.\ Lambrechts, V. Turchin, and I.\ Voli\'{c},
{\em Real homotopy theory of semi-algebraic sets},
Algebr. Geom. Topol. {\bf 11} (2011), no. 5, 2477-2545.


\bibitem[KM]{KM} M.\ Kapranov and Yu.I.\ Manin, {\em Modules and Morita theorem for operads}. Amer. J. Math. {\bf 123} (2001), no. 5, 811-838.

\bibitem[KMW]{KMW} A. Khoroshkin, S. Merkulov and T. Willwacher,
{\em On Quantizable Odd Lie Bialgebras}. Lett Math Phys {\bf 106} (2016), no 9, pp. 1199--1215.

 \bibitem[Ko1]{Ko0} M.\ Kontsevich, {\em Feynman diagrams and low-dimensional topology}. In  Progr. Math., vol. 120, Birkhäuser, Basel, 1994, pp. 97 - 121.

 \bibitem[Ko2]{Ko-f} M. Kontsevich, {\em Formality Conjecture}, D. Sternheimer et al. (eds.),
Deformation Theory and Symplectic
Geometry, Kluwer 1997, 139-156.

\bibitem[Ko3]{Ko} M.\ Kontsevich, {\em Deformation quantization
 of Poisson manifolds}, Lett.\ Math.\ Phys. {\bf 66} (2003), 157-216.




  \bibitem[KS]{KS} M.\ Kontsevich and Y.\ Soibelman,
 {\em Deformations of algebras over operads and the Deligne conjecture}. In
Conf´erence Mosh´e Flato 1999, Vol. I (Dijon), volume 21 of Math. Phys. Stud., pages 255-307. Kluwer Acad.
Publ., Dordrecht, 2000.


\bibitem[LTV]{LTV} P.\ Lambrechts, V.\ Turchin and I.\ Volic, {\em
Associahedron, cyclohedron, and permutohedron
as compactifications of configuration spaces},
Bull. Belg. Math. Soc. Simon Stevin 17 (2010), no. 2, 303-332.


\bibitem[Ma1]{Ma1} M.\ Markl,
\newblock {\em A resolution (minimal model) of the prop for bialgebras},
J.\ Pure Appl.\ Algebra {\bf 205} (2006), no. 2, 341"1¤7374.

\bibitem[Ma2]{Ma} M.\ Markl,
\newblock {\em Bipermutahedron and biassociahedron},
\newblock arXiv:1209.5193


\bibitem[MMS]{MMS}  M.\ Markl, S.\ Merkulov and S.\ Shadrin, {\em Wheeled props and the master
equation}, preprint math.AG/0610683, J.\ Pure and Appl.\ Algebra {\bf 213} (2009), 496-535.

\bibitem[MaVo]{MaVo}
M.\ Markl and A.A.\ Voronov, {\em PROPped-up graph cohomology}. Algebra, arithmetic, and geometry: in honor of Yu. I. Manin. Vol. II, 249"1¤7281, Progr. Math., 270, Birkhäuser Boston, Inc., Boston, MA, 2009.

\bibitem[MS]{MS} M.\ Markl and S.\ Shnider, {\em Associahedra, cellular $W$-construction and products of $A_\infty$-algebras}, Transactions of the AMS {\bf 358}, Number 6 (2005) 2353-2372.

\bibitem[Me1]{Me0} S.A.\ Merkulov, {\em Exotic automorphisms of the Schouten algebra of polyvector fields }, arXiv:0809.2385 (2008)

\bibitem[Me2]{Me1} S.A.\ Merkulov, {\em Operads, configuration spaces and quantization}.
 In: ``Proceedings of Poisson 2010, Rio de Janeiro", Bull.\ Braz.\ Math.\ Soc., New Series {\bf 42}(4) (2011), 1-99.

\bibitem[Me3]{Me2} S.A.\ Merkulov, {\em Formality Theorem for Quantizations of Lie Bialgebras}, Letters in Mathematical Physics,  Volume 106, Issue 2 (2016) 169-195.


\bibitem[MV]{MV}  S.A.\ Merkulov and  B.\ Vallette,
{\em Deformation theory of representations of prop(erad)s I \& II},
 Journal f\"ur die reine und angewandte Mathematik (Qrelle)  {\bf 634}, 51-106,
 \& {\bf 636}, 123-174 (2009)


\bibitem[MW1]{MW} S.A. Merkulov and T.\ Willwacher, {\em Deformation theory of  Lie bialgebra properads}, preprint  arXiv:1512.05252 (2015)

\bibitem[MW2]{MW2} S.A. Merkulov and T.\ Willwacher, {\em Classification of universal formality maps\\
for quantizations of Lie bialgebras}, preprint arXiv:1605.01282 (2016)

\bibitem[RW]{RW} C.\ Rossi and T.\ Willwacher, {\em P.\ Etingof's conjecture about Drinfeld associators}, preprint  arXiv:1404.2047 (2014)

\bibitem[SU1]{SU} S.\ Saneblidze and R.\ Umble, {\em  Matrads, biassociahedra, and $A_\infty$-bialgebras}, Homology, Homotopy Appl.,
{\bf 13}(1)  (2011), 1-57.

\bibitem[SU2]{SU2} S.\ Saneblidze and R.\ Umble, {\em Diagonals on the permutahedra, multiplihedra and associahedra},
Homology, Homotopy and Applications {\bf 6} (2004) 363-411.

\bibitem[Se]{Se} P. $\check{\mathrm S}$evera, {\em Quantization of Lie bialgebras revisited},
arXiv:1401.6164 (2014)

\bibitem[Sh1]{Sh1} B.\ Shoikhet, {\em An explicit formula for the deformation
quantization of Lie bialgebras}, arXiv:math.QA/0402046.

\bibitem[Sh2]{Sh} B.\ Shoikhet, {\em  An $L_\infty$ algebra structure on polyvector fields
}, preprint arXiv:0805.3363, (2008).

\bibitem[St]{St} J.D.\ Stasheff. {\em On the homotopy
  associativity of $H$-spaces, I \@ II}. {Trans.\ Amer.\
  Math.\ Soc}., {\bf 108} (1963), 272--292 \& 293--312.




\bibitem[T1]{Ta2} D.E.\ Tamarkin, {\em Another proof of M. Kontsevich formality
theorem},  math.QA/9803025, Lett.\ Math.\ Phys. {\bf 66}
(2003) 65-72.

\bibitem[T2]{Ta} D.E.\ Tamarkin, {\em Quantization of lie bialgebras via the
formality of the operad of little disks}, GAFA Geom. funct. anal.
{\bf 17} (2007), 537-604.

\bibitem[V]{Va}
B.\ Vallette,  {\em A Koszul duality for
props}, Trans.\ Amer. Math. Soc., {\bf 359} (2007), 4865--4943.

\bibitem[We]{Wei} C.A.\ Weibel. An introduction to homological algebra,
CUP, 2003.

\bibitem[W1]{Wi} T.\ Willwacher, {\em M.\ Kontsevich's graph complex and the
 Grothendieck-Teichmueller Lie algebra},
Invent. Math. 200 (2015) 671-760.

\bibitem[W2]{Wi3}  T.\ Willwacher, {\em Stable cohomology of polyvector fields}, Math.\ Res.\ Lett.
{\bf 21}, No.\ 6 (2014), 1501-1530.

\bibitem[W3]{Wi2} T.\ Willwacher, {\em Oriented graph complexes},   Comm. Math. Phys. {\bf 334} (2015), no. 3, 1649--1666.

 \bibitem[W4]{Wi4} T.\ Willwacher, {\em Models for the $n$-Swiss Cheese operads}, preprint arXiv:1506.07021 (2015).



 \end{thebibliography}
\end{document}